%% file: 00_main_SISC.tex
\documentclass{article}
\input{ex_shared}

\begin{document}

\maketitle

% REQUIRED
\begin{abstract}
Topology optimization under uncertainty or reliability-based topology optimization is usually numerically very expensive. This is mainly due to the fact that an accurate evaluation of the probabilistic model requires the system to be simulated for a large number of varying parameters. Traditional gradient-based optimization schemes thus face the difficulty that reasonable accuracy and numerical efficiency often seem mutually exclusive. In this work, we propose a stochastic optimization technique to tackle this problem. To be precise, we combine the well-known method of moving asymptotes (MMA) with a stochastic sample-based integration strategy. By adaptively recombining gradient information from previous steps, we obtain a noisy gradient estimator that is asymptotically correct, i.e., the approximation error vanishes over the course of iterations. As a consequence, the resulting stochastic method of moving asymptotes (sMMA) allows us to solve chance constraint topology optimization problems for a fraction of the cost compared to traditional approaches from literature. To demonstrate the efficiency of sMMA, we analyze structural optimization problems in two and three dimensions.
\end{abstract}

% % REQUIRED
% \begin{keywords}
% reliability-based topology optimization, chance constraints, method of moving asymptotes, stochastic optimization
% \end{keywords}

% % REQUIRED
% \begin{MSCcodes}
% 74P05, 90C15, 90C30, 90C59
% \end{MSCcodes}

%################################### Introduction ############################################
%#############################################################################################
\section{Introduction}
In this contribution, we consider structural optimization problems involving uncertainties, i.e., problems of the form
\begin{align}\label{eq:GeneralSetting}
    \begin{split}
    \min_{u\in\U}\quad & \E_{\xi_1}\Big[f_{\xi_1}^{(1)}\Big(\E_{\xi_2}\Big[f_{\xi_2}^{(2)}\Big(\ldots\big(\E_{\xi_N}\big[f_{\xi_N}^{(N)}(u)\big]\big)\ldots\Big)\Big]\Big)\Big], \\
    \text{s.t.}\quad & \E_{\omega_1}\Big[g_{\omega_1}^{(1)}\Big(\E_{\omega_2}\Big[g_{\omega_2}^{(2)}\Big(\ldots\big(\E_{\omega_M}\big[g_{\omega_M}^{(M)}(u)\big]\big)\ldots\Big)\Big]\Big)\Big] \le 0,
    \end{split}
\end{align}
where,
\begin{align*}
\begin{alignedat}{3}
    f_{\xi_i}^{(i)}&:\R^{n_i}\to\R^{n_{i-1}},&\quad n_i\in\N, \quad i=1,\ldots,I, \\
    g_{\omega_j}^{(j)}&:\R^{m_j}\to\R^{m_{j-1}},&\quad m_j\in\N,\quad j=1,\ldots,J
\end{alignedat}
\end{align*}
and $n_0=1$. While such problems arise in countless different contexts, we focus on the field of density based topology optimization~\cite{Bendse1989,BendsoeKikuchi,TopOpt,Mlejnek1992,Sigmund2013}. Here, uncertainty may enter the model as a consequence of, e.g., material/manufacturing imperfections or distributed load cases~\cite{De2020,De2022,Komini2023,Liu2018,Long2018,Maute2014,Nishioka2023,Tootkaboni2012,Torres2021}. 

As motivation, assume that we are tasked to design a pole tent. Our goal is to use as little material as possible, such that the tent is unlikely to break under usual wind loads. To be precise, for a given design $\rho$ and scenario $\omega$, let $g(\rho,\omega)$ describe the behavior of the tent under the wind load described by $\omega$, where $g(\rho,\omega)>0$ corresponds to the tent breaking down under the load.
Given how vague our task is formulated, it is not surprising that several different modelling approaches fall well within the structure of~\eqref{eq:GeneralSetting}.
\begin{enumerate}
    \item \emph{Constraints in expectation.} Directly following the formulation in~\eqref{eq:GeneralSetting}, the optimization problem can be stated as
    \begin{align*}
        \min_{\rho}\quad & \vol(\rho), \\
        \text{s.t.}\quad & \E_\omega\big[g(\rho,\omega)] \le 0.
    \end{align*}
    This approach will yield a design, which is resilient enough \emph{on average}. However, we have no control over how likely the structure is to break, i.e., ${\P_\omega[g(\rho,\omega) > 0]}$ can (theoretically) admit any value in $(0,1)$.
    \item \emph{Robust optimization.} If we want to find a design, which we can guarantee to not break under any random load, we can do so by following a worst case approach
    \begin{align*}
        \min_{\rho}\quad & \vol(\rho), \\
        \text{s.t.}\quad & g(\rho,\omega) \le 0\quad \forall \omega.
    \end{align*}
    Note that the constraint can equivalently be stated as
    \begin{equation*}
        \P_\omega\big[g(\rho,\omega)>0\big] = \E_\omega\big[\chi_{(0,\infty)}\big(g(\rho,\omega)\big)\big] \le 0,
    \end{equation*}
    where $\chi_{\mathcal{M}}:\R\to\R$ denotes the indicator function of a set $\mathcal{M}\subseteq\R$. There are two main reasons, why we do not consider such problems in this contribution. 
    
    First, enforcing the constraint for \emph{all} realizations of $\omega$ will most likely result in a design with enormous cost. For example, if our wind loads are taken from a weather database, the worst case constraint will enforce a tent design capable of withstanding hurricanes and other extremely rare weather conditions. It is debatable whether load cases of such low probabilities should be allowed to have that strong of an impact on the final design. 

    Second, although the robust formulation can be expressed in the form of~\eqref{eq:GeneralSetting}, there are several numerical difficulties in solving the problem in this formulation. For example, by construction, there can not exist a strictly feasible point $\tilde{\rho}$ with ${\P_\omega[g(\tilde{\rho},\omega)]<0}$. As a consequence, we refer the reader to specialized solution techniques that have been developed for this type of problems~\cite{Brittain2011,Elishakoff1994,Greifenstein2020,Guo2009,Lombardi1998}.
    \item \emph{Chance constraints.} As a compromise between the two approaches, we can instead consider designs for which the probability of failure is small enough. That is, we can fix a probability level ${p\in(0,1)}$ and model our problem as
    \begin{align*}
        \min_{\rho}\quad & \vol(\rho), \\
        \text{s.t.}\quad & \E_\omega\big[\chi_{(0,\infty)}\big(g(\rho,\omega)\big)\big] \le p.
    \end{align*}
    This formulation allows us to fine-tune the balance between costs ($p\approx 1$) and security ($p\approx0$). Moreover, we will see in later examples that choosing $p$ small enough can also produce designs that are robust in the above sense.
\end{enumerate}
Chance constraint optimization problems have a rich theory concerning solution techniques, e.g., scenario generation approaches~\cite{Calafiore2004,Campi2010}, sample average approximations~\cite{Luedtke2008,Pagnoncelli2009}, smoothing approximations~\cite{Hong2011,Shan2014} or bundle methods~\cite{deOliveira2014,Kuchlbauer2022}, just to name a few.

Many of the approaches mentioned above rely on extensive sampling of the constraint function $g$, rendering them numerically inefficient in our framework. Further amplifying the difference between these techniques and our setting, structural optimization problems of the above type are usually not called ``chance constrained" and instead are referred to as \emph{reliability-based topology optimization}~\cite{De2021,Gao2021,Kharmanda2004,Li2024}. Standard approaches for such problems typically combine an established optimization method from deterministic topology optimization with validation techniques to estimate the reliability of the structure w.r.t. the considered uncertainties. 

In this contribution, we propose a slightly different approach in which validation and optimization are more closely intertwined. The motivation is straight-forward: If we were able to efficiently evaluate all appearing probabilities, we would simply use a specialized deterministic optimization scheme, e.g., the \emph{method of moving asymptotes} (MMA)~\cite{MMA_original}, to solve the problem. Since this is not the case, we want to at least avoid an expensive recomputation in each optimization step. The \emph{continuous stochastic gradient method} (CSG), proposed in~\cite{CSG_original}, is a stochastic sample-based optimization approach, which aims to minimize the number of function evaluations during the optimization process by recombining information collected in past iterations. In fact, the approximation error of the CSG gradient estimator was shown to almost surely converge to 0 over the course of iterations~\cite{CSG_part1}. Therefore, combining MMA with CSG integration techniques, we obtain the \emph{stochastic method of moving asymptotes} (sMMA), an efficient stochastic optimization scheme for reliability-based topology optimization. 

\subsection*{Structure of the paper}
In~\Cref{sec:Problem}, the general optimization problem for linear elasticity with chance constraints is formulated. Moreover, regularization techniques and regularity assumptions are presented. \Cref{sec:Method} recalls important basic concepts from CSG and MMA theory. Afterwards, we introduce the general sMMA framework and propose algorithmic augmentations in the context of topology optimization. In~\Cref{sec:Appl}, we provide extensive numerical results for sMMA in different setups. The considered examples cover two- and three-dimensional structures and uncertainty entering the state equation through the right hand side as well as the system matrix. All results are compared to MMA, running with a finite discretization of the appearing integrals.

%#############################################################################################
%#################################### Setting ################################################
%#############################################################################################
\section{Problem formulation}\label{sec:Problem}
In our applications, we focus on density based topology optimization for linear elasticity. Here, given a design domain $\mathscr{D}\subset\R^n$, ${n\in\{2,3\}}$, we denote by ${\rho\in\R^d}$ the finite element representation of a pseudo-density function for the material. By convention, ${\rho=1}$ indicates finite elements filled with material, whereas ${\rho=0}$ corresponds to void. To ensure that the state problem is well-posed, material parameters for empty cells are chosen positive, but small. Thus, void is essentially replaced by very weak material.

Given a force $\mathbf{F}$ (load) acting on the design $\rho$, the compliance of the structure is given by ${\mathbf{F}^\top\mathbf{U}}$, where $\mathbf{U}$ is the unique solution to the underlying state equation in algebraic form ${\mathbf{K}(\rho)\mathbf{U}=\mathbf{F}}$, with global stiffness matrix $\mathbf{K}(\rho)$. 

In this formulation, uncertainties can enter the system mainly at two points: the stiffness matrix $\mathbf{K}$, or the load $\mathbf{F}$. Thus, let ${(\Xi,\mathscr{A},\P^{\Xi})}$ and ${(\Omega,\mathscr{B},\P^{\Omega})}$ be probability spaces with ${\Xi\subset\R^{d_\Xi}}$ and ${\Omega\subset\R^{d_\Omega}}$. Let ${\xi:\Xi\to\R^{d_\xi}}$ and ${\omega:\Omega\to\R^{d_\omega}}$ be random vectors defined on these probability spaces with probability distributions $\mu$ and $\nu$. Then, the probabilistic state equation reads
\begin{equation*}
    \mathbf{K}(\rho,\xi)\mathbf{U}=\mathbf{F}(\omega).
\end{equation*}
To keep matters simple, our objective in all applications will be to minimize the volume of material used, while ensuring the compliance to stay below an upper threshold with high probability. Thus, given a probability level $p>0$ and maximum compliance value $\cm>0$, the resulting optimization problem can be stated as follows:
\begin{align}\label{eq:GeneralOptProb}
    \begin{split}
        \min_{\rho\in\R^d}\quad & \vol(\rho),\\
        \text{s.t.}\quad & \P_{\xi,\omega} \left[\mathbf{F}(\omega)^\top\mathbf{U}_{\xi,\omega}>\cm\right]\le p,\\
            &\mathbf{K}(\rho,\xi)\mathbf{U}_{\xi,\omega}=\mathbf{F}(\omega),\\
            &0\le \rho\le 1.
    \end{split}
\end{align}
\subsection{Regularization and assumptions}
As is common in the literature~\cite{Bourdin2001,Sigmund2007,Wang2010}, a linear filter (with filter matrix ${\mathcal{F}\in\R^{d\times d}}$) is applied to the design vector $\rho$. Moreover, the usage of fictitious intermediate material ${0<\rho<1}$ is discouraged by using the solid isotropic material with penalization (SIMP) technique~\cite{Bendse1989,Bendse1999,Mlejnek1993,Rozvany2000} with SIMP-parameter ${s>1}$. Therefore, the material parameters of a design cell ${i\in\{1,\ldots,d\}}$, which are needed to construct $\mathbf{K}$, are modeled as
\begin{equation*}
    [\mathcal{E}(\rho)]_i = \big([\mathcal{F}\rho]_i\big)^s\mathcal{E}^{(1)} + \big(1-\big([\mathcal{F}\rho]_i\big)^s\big)\mathcal{E}^{(0)},
\end{equation*}
where $\mathcal{E}^{(1)}$ and $\mathcal{E}^{(0)}$ denote the parameters for material and void, respectively. Recall that this interpolation is not applied for the volume, as this would only substitute ${\rho\mapsto (\mathcal{F}\rho)^s}$ in the original problem and not yield any penalization of artificial material. As a consequence, we distinguish two different terminologies:
\begin{enumerate}
    \item The volume of a design ${\vol(\mathcal{F}\rho)}$, which is used for the optimization process. The objective function in all of our numerical examples will be given by the relative volume 
    \begin{equation*}
        \rvol(\rho):=\frac{\vol(\mathcal{F}\rho)}{\vol(\mathscr{D})}.
    \end{equation*}
    \item The physical interpretation ${(\mathcal{F}\rho)^s}$ of a design is what would be of interest if we were to manufacture the structure. Thus, it is used for our design visualizations and corresponds to a physical relative volume 
    \begin{equation*}
        \pvol(\rho):=\frac{\vol\big((\mathcal{F}\rho)^s\big)}{\vol(\mathscr{D})}.
    \end{equation*}
\end{enumerate}
Since these values coincide for a truly ``black \& white" design ${\rho\in\{0,1\}^d}$, their difference provides a measure for the greyness of a result, i.e., how much fictitious material is used.

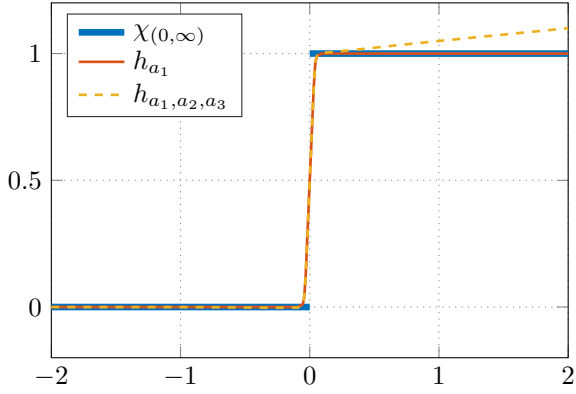
\begin{figure}
  \begin{minipage}[c]{0.55\textwidth}
    \input{cc_smooth}
  \end{minipage}\hfill
  \begin{minipage}[c]{0.4\textwidth}
    \caption{Indicator function $\chi_{(0,\infty)}$ as well as the smooth approximation $h_{a_1}$ and steepened smooth approximation $h_{a_1,a_2,a_3}$. In this illustration, we chose ${a_1=35}$, ${a_2=\tfrac{1}{20}}$ and ${a_3=5}$.}
      \label{fig:cc_smooth}
  \end{minipage}
\end{figure}

As the chance constraint in~\eqref{eq:GeneralOptProb} will later be approximated by sMMA, we can not directly work with the equivalent reformulation
\begin{equation*}
    \E_{\xi,\omega}\left[ \chi_{(0,\infty)}\big( \mathbf{F}(\omega)^\top\mathbf{U}_{\xi,\omega}-\cm \big) \right] \le p,
\end{equation*}
since $\chi_{(0,\infty)}$ is not continuous. Thus, we smoothly approximate $\chi_{(0,\infty)}$ by the family of functions ${h_{a_1}:\R\to\R}$, given by
\begin{equation*}
    h_{a_1}(t):= \frac{1}{2}\big(\tanh(a_1t)+1\big),\quad a_1 > 0.
\end{equation*}
An illustration of the smoothed approximation is given in~\Cref{fig:cc_smooth}. Moreover, since we apply first-order methods to the regularized version of~\eqref{eq:GeneralOptProb}, there arises an additional difficulty from the chance constraint. For ${t\ll0}$ or ${t\gg0}$, the derivative of $h_{a_1}$ is almost zero. As a result, the first-order linearization of the chance constraint in~\eqref{eq:GeneralOptProb} will be approximately constant for designs with ${\mathbf{F}^\top\mathbf{U}\ll\cm}$ or ${\mathbf{F}^\top\mathbf{U}\gg\cm}$. While the first case does not cause a problem, the optimizer will not be able to recover from the second case and might converge to the infeasible point ${\rho\equiv0}$. To circumvent this problem, the constraint gradient can be steepened in the infeasible region. While specialized techniques have been proposed in the literature~\cite{BernhardLiersStingl2024}, we found it sufficient for our numerical experiments to add a smoothmax-term, i.e., a smooth approximation to ${t\mapsto\max\{0,t\}}$, to the definition of $h$:
\begin{equation}\label{eq:smooth_cc_fun}
    h_{a_1,a_2,a_3}(t):=\frac{1}{2}\big(\tanh(a_1t)+1\big)+a_2\left( t-\frac{t}{1+\exp(a_3t)} \right),\quad a_1,a_2,a_3 >0.
\end{equation}
The resulting approximation is shown in~\Cref{fig:cc_smooth} as well. Thus, the final regularized chance constraint can be expressed as
\begin{equation}\label{eq:ConstraintIntegral}
    G(\rho):=\int_{\Xi\times\Omega} h_{a_1,a_2,a_3}\left( \mathbf{F}(\omega)^\top\mathbf{U}_{\xi,\omega} - \cm \right)\dd(\mu\times\nu)(\xi,\omega)\le p.
\end{equation}

Even if the compliance bound ${\mathbf{F}^\top\mathbf{U}\le\cm}$ is satisfied with equality for all ${\xi\in\Xi}$ and ${\omega\in\Omega}$, the smooth approximations provided by $h_{a_1}$ and $h_{a_1,a_2,a_3}$ will yield a chance constraint value of ${G(\rho)=\tfrac{1}{2}}$. Thus, when optimizing the smoothed reformulation 
\begin{align}\label{eq:SmoothOptProb}
    \begin{split}
        \min_{\rho\in\R^d}\quad & \rvol(\rho),\\
        \text{s.t.}\quad & G(\rho)\le p,\\
            &\mathbf{K}(\rho,\xi)\mathbf{U}_{\xi,\omega}=\mathbf{F}(\omega),\\
            &0\le \rho\le 1.
    \end{split}
\end{align}
of~\eqref{eq:GeneralOptProb}, we expect feasible solutions to satisfy the original chance constraint even for a slightly smaller value of $p$. For our numerical experiments, we found the differences between each of the three formulations to be small. In~\Cref{fig:wheel_smoothingeffect}, exemplary comparisons between the values obtained by each method can be found.

We finish this section by collecting regularity assumptions.
\begin{assumption}\label{assu:Regularity}
Both $\mu$ and $\nu$ are Borel probability measures with bounded support. Moreover, for all ${\xi\in\supp(\mu)}$ and ${\omega\in\supp(\nu)}$, the state mapping 
\begin{equation*}
    \mathcal{S}(\rho,\xi,\omega):=\mathbf{K}(\rho,\xi)^{-1}\mathbf{F}(\omega)
\end{equation*}
is differentiable with respect to $\rho$. Lastly, $\mathcal{S}$ and $\nabla_1\mathcal{S}$ are Lipschitz continuous is all arguments.
\end{assumption}

%#############################################################################################
%##################################### Method ################################################
%#############################################################################################
\section{Methods and concepts}\label{sec:Method}
As the proposed method represents a combination of MMA and CSG, we briefly recall key concepts of both optimization schemes. For a more detailed analysis of CSG and MMA, see~\cite{CSG_part2,CSG_part1,CSG_original} and~\cite{MMA_original,MMA_lecture,GCMMA_1,MMA_global_2}, respectively. Afterwards, some possible adaptations for special difficulties arising in topology optimization problems are discussed.
\subsection*{The method of moving asymptotes}
In each MMA iteration, gradient information of the objective and constraint function at the current design is used to build convex first-order approximations. These approximations are then used to construct the current MMA subproblem. To be precise consider the optimization problem
\begin{align*}
        \min_{z\in\R^d}\quad & f(z), \\
        \text{s.t.} \quad & g_i(z)\le 0, \quad i=1,\ldots,m\\
        & \underline{z}_j\le z_j \le \overline{z}_j, \quad j=1,\ldots,d,
\end{align*}
with ${f:\R^d\to\R}$, ${g_i:\R^d\to\R}$ and ${\underline{z},\overline{z}\in\R^d}$. Denote by $z^{(k)}$ the current design vector with entries $z^{(k)}_j$, $j=1,\ldots,d$. Then, the convex first-order approximations are defined as
\begin{equation*}
    \widetilde{g}_i(z) = r^{(k)}_i + \sum_{j=1}^d\left( \frac{p^{(k)}_{ij}}{U^{(k)}_j-z_j} + \frac{q^{(k)}_{ij}}{z_j-L^{(k)}_j} \right),
\end{equation*}
with
\begin{align*}
    p^{(k)}_{ij} &= \begin{cases} (U^{(k)}_j-z^{(k)}_j)^2\frac{\partial}{\partial z_j}g_i(z^{(k)}), & \text{if }\frac{\partial}{\partial z_j}g_i(z^{(k)}) > 0 \\ 0, & \text{else}, \end{cases} \\
    q^{(k)}_{ij} &= \begin{cases} -(z^{(k)}_j-L^{(k)})^2\frac{\partial}{\partial z_j}g_i(z^{(k)}), & \text{if }\frac{\partial}{\partial z_j}g_i(z^{(k)}) < 0 \\ 0, & \text{else}, \end{cases} \\
    r^{(k)}_i &= g_i(z^{(k)}) - \sum_{j=1}^d \left( \frac{p^{(k)}_{ij}}{U^{(k)}_j-z^{(k)}_j} + \frac{q^{(k)}_{ij}}{z^{(k)}_j-L^{(k)}_j} \right).
\end{align*}
The approximation to $f$ is defined analogously. The parameters ${L^{(k)},U^{(k)}}$ are the (heuristically) chosen \emph{moving asymptotes}, which satisfy
\begin{equation*}
    L^{(k)}_j < z^{(k)}_j < U^{(k)}_j, \quad j=1,\ldots,d.
\end{equation*}
By construction, these first-order approximations are convex and separable, i.e., they can be decomposed into a sum of independent single-variable functions. Thus, the resulting MMA subproblem is easy to solve.

To obtain a globally convergent optimization scheme (GCMMA~\cite{GCMMA_1}), it is possible to force a decrease in objective function values between iterations by reconstructing the subproblems with tighter asymptotes, if necessary. However, in practice, it is often sufficient to simply augment the subproblems by a constraint
\begin{equation*}
    z^{(k)}_j-\tau \le z_j \le z^{(k)}_j+\tau,\quad j=1,\ldots,d,
\end{equation*}
where $\tau>0$ is referred to as \emph{move limit}.

If both measures $\mu$ and $\nu$ are known, MMA can directly be applied to~\eqref{eq:SmoothOptProb} by numerically approximating all integrals. However, due to the complex underlying structure, the number of integration points required to achieve a desired accuracy may be extremely high. Coupled with the fact that solving the state equation for fixed $\xi,\omega$ is typically expensive, a direct MMA approach is often computationally infeasible or highly inaccurate, as we will see in our numerical experiments.

\subsection*{The continuous stochastic gradient method}
\begin{figure}
  \begin{minipage}[c]{0.55\textwidth}
    \input{weights_L2}
  \end{minipage}\hfill
  \begin{minipage}[c]{0.4\textwidth}
    \caption{Nearest neighbor approximation for $K=10$, ${\Dim(\U)=\Dim(\X)=1}$ and ${\Vert\cdot\Vert_{\U\times\X}=\Vert\cdot\Vert_2}$. Piecewise constant regions of $\widehat{j}_{10}$ are indicated by the cells in the background. The approximation $\widehat{J}_{10}(u_{10})$ of $J(u_{10})$ is obtained by integrating $\widehat{j}_{10}$ along the solid line. The sets $M_k^{10}$ are given by the colored line segments. Grey cells correspond to cases where $M_k^{10}$ is empty, resulting in an integration weight $\alpha_{k,10}=0$.}
    \label{fig:weights_exact}
  \end{minipage}
\end{figure}
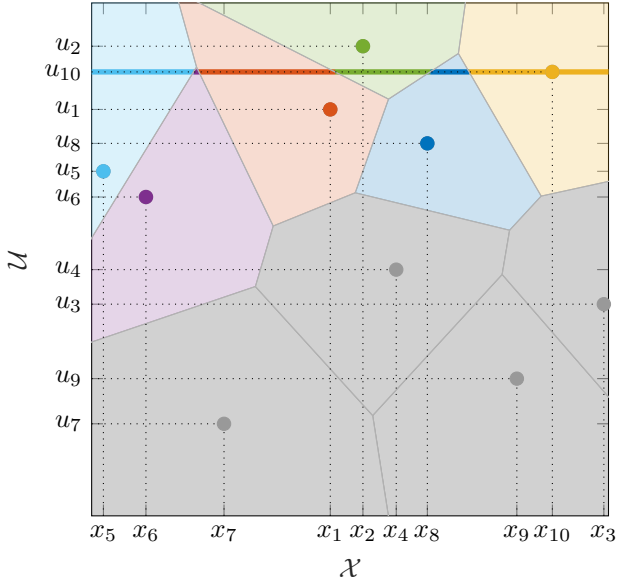

CSG is a stochastic sample-based optimization technique, where old gradient information is used to reconstruct the current gradient via an efficient on-the-fly integration approach. To illustrate the basic idea, let $J:\U\times\X\to\R^d$ be given by
\begin{equation*}
    J(u):= \int_\X j(u,x)\dd\mu(x)
\end{equation*}
and assume that $j$ has been sampled at points ${(u_k,x_k)_{k=1,\ldots,K}}$. Then, the goal is to find design-dependent integration weights $(\alpha_k)_{k=1,\ldots,K}$, such that
\begin{equation*}
    J(u) \approx \sum_{k=1}^K \alpha_k j(u_k,x_k).
\end{equation*}
As a straightforward approach, we might approximate unknown values of $j$ by the value of the nearest available sample point, i.e., we choose a norm $\Vert\cdot\Vert_{\U\times\X}$ and set
\begin{align*}
    &j(u,x) \approx \widehat{j}_K(u,x):= j\left(u_{\tau(u,x;K)},x_{\tau(u,x;K)}\right),\\
    &\tau(u,x;K)\in\argmin_{k=1,\ldots,K} \;\big\Vert (u,x)-(u_k,x_k)\big\Vert_{\U\times\X}.
\end{align*}
Since $\widehat{j}_K$ is piecewise constant, integrating this approximation over $\X$ yields
\begin{equation*}
    J(u)\approx \widehat{J}_K(u):=\int_\X \widehat{j}_K(u,x)\dd\mu(x) = \sum_{k=1}^K j(u_k,x_k)\mu\big(M_k^K(u)\big),
\end{equation*}
where
\begin{equation*}
    M_k^K(u) := \big\{ x\in\X\,:\, \tau(u,x;K)=k \big\}.
\end{equation*}
An illustration for a low-dimensional example is given in~\Cref{fig:weights_exact}. By iterating this procedure, the approximation scheme can be generalized to composite structures, like they appear in~\eqref{eq:GeneralSetting}. In CSG, these approximations are constructed for the objective function as well as its gradient and then used in a standard (projected) gradient scheme.

The key feature of CSG lies within the \emph{approximation property}, which, translated to our setting, states:
\begin{proposition}\label{prop:ApproximationProperty}
Let the sample sequences $(\xi_{n})_{n\in\N}$ and $(\omega_{n})_{n\in\N}$ of the random variables $\xi$ and $\omega$ be independent and identically distributed according to the respective probability measures $\mu$ and $\nu$. Then, for any bounded sequence $(\rho_n)_{n\in\N}$ of designs, under~\Cref{assu:Regularity}, it holds
\begin{equation*}
    \big\Vert G(\rho_n)-\widehat{G}_n(\rho_n)\big\Vert + \big\Vert \nabla G(\rho_n)-\widehat{dG}_n(\rho_n)\big\Vert \xrightarrow{\text{a.s.}}0
\end{equation*}
for ${n\to\infty}$, where $\widehat{G}$ and $\widehat{dG}$ denote the CSG approximations to $G$ and $\nabla G$, as given in~\eqref{eq:ConstraintIntegral}.
\end{proposition}
\begin{proof}
See~\cite[Lemma 4.6 and Remark 4.2]{CSG_part1}.
\end{proof}

It is important to note that, with the exception of $\X$ being one-dimensional, it is usually numerically infeasible to calculate the sets $M_k^K$ explicitly. Thus, a number of more efficient ways to obtain approximations to ${\alpha_{k,K}=\mu(M_k^K)}$ have been proposed in~\cite{CSG_part1}. Moreover, a novel approach, suited for the applications within this study, is proposed later~\Cref{subsec:weights}.
\subsection*{The stochastic method of moving asymptotes}
In material and topology optimization, MMA is usually considered superior to a standard gradient descent scheme. As a result, granted that the integral appearing in~\eqref{eq:ConstraintIntegral} can efficiently be numerically approximated with high enough precision, there is no reason to choose CSG over MMA. However, if the evaluation of~\eqref{eq:ConstraintIntegral} is the limiting factor, MMA might no longer be the better choice. In fact, even for a one-dimensional integral, it was observed in~\cite{CSG_acoustics} that CSG was more efficient than MMA, when counting the number of required state equation solutions. 

In essence, CSG provides a very efficient mechanism of obtaining a gradient approximation, but wastes this information by using it for a simple projected gradient step. In contrast, MMA is highly efficient in utilizing given gradient information, but provides no tools for actually obtaining it.

The stochastic method of moving asymptotes (sMMA) scheme proposed in this contribution is constructed as a direct combination of both approaches, aiming to combine key advantages of each method. In each sMMA iteration, the state equation is solved only for a small batch of realizations of the random variables $\xi$ and $\omega$. The sampled information is then used to build function approximations according to the CSG model. Based on these approximations, the MMA subproblem is constructed and solved as usual. An overview of the method is given in~\Cref{alg:MCMSA}. Note that, by construction, sMMA does not require any changes in the MMA subroutine (lines 7 to 9).
\begin{algorithm}
\caption{sMMA}
\begin{algorithmic}[1]
    \STATE Choose initial design $\rho_1$.
    \FOR{$k=1,2,\ldots$}
        \STATE Draw random samples $\xi_k\sim\mu$ and $\omega_k\sim\nu$.
        \STATE Solve $\mathbf{K}(\rho_k,\xi_k)\mathbf{U}_{\xi_k,\omega_k}=\mathbf{F}(\omega_k)$ and evaluate $\mathbf{F}(\omega_k)^\top\mathbf{U}_{\xi_k,\omega_k}$.
        \STATE Calculate CSG integration weights $(\alpha_{i,k})_{i=1,\ldots,k}$.
        \STATE Calculate CSG approximations $\widehat{G}_k(\rho_k)$ and $\widehat{dG}_k(\rho_k)$ to $G(\rho_k)$ and $\nabla G(\rho_k)$.
        \STATE Construct MMA subproblem based on $\widehat{G}_k(\rho_k)$ and $\widehat{dG}_k(\rho_k)$.
        \STATE Optional: Choose move limits.
        \STATE Calculate solution $\rho^\ast$ of MMA subproblem.
        \STATE Set $\rho_{k+1}=\rho^\ast$.
    \ENDFOR
\end{algorithmic}\label{alg:MCMSA}
\end{algorithm}
\subsection{Pseudoexact weights}\label{subsec:weights}
Consider again the example ${J:\U\times\X\to\R}$ from above. As mentioned before, calculating the full cells $M_k^K$, appearing in the integration weight calculation, is generally very difficult and too time-consuming in practice. However, given any point $x\in\X$, it is very easy to determine, which of the sets $M_k^K$ it belongs to. In fact, all we need to do is calculate the distances ${\Vert (u,x)-(u_k,x_k)\Vert_{\U\times\X}}$ and pick an index minimizing this expression. Therefore, we can obtain approximations to $\mu(M_k^K)$ by discretizing $\X$ into $(\mathbf{x}_t)_{t=1,\ldots,T}$ (according to $\mu$). Then, 
\begin{equation*}
    \mu(M_k^K) \approx \sum_{t=1}^T w_t\chi_{_{M_k^K}}(\mathbf{x}_t), \quad w_t \ge 0, \quad \sum_{t=1}^T w_t = 1.
\end{equation*}
As increasing $T$ only requires additional evaluations of $\Vert\cdot\Vert_{\U\times\X}$ and no additional solutions of the state equation, we can typically choose $T$ rather large in our applications.
Specifically choosing the points $\mathbf{x}_t$ as the sample points $x_k$ and $w_t=T^{-1}$, $\mu$ is approximated by the empirical measure, yielding \emph{empirical integration weights}, as proposed in~\cite{CSG_part1}.

\begin{figure}
  \begin{minipage}[c]{0.55\textwidth}
    \input{weights_pseudoexact}
  \end{minipage}\hfill
  \begin{minipage}[c]{0.4\textwidth}
    \caption{For pseudoexact integration weights, the domain $\X$ is first discretized by the points $(\mathbf{x}_t)_{t=1,\ldots,T}$ (colored diamonds). For each $\mathbf{x}_t$, it is easy to check the closest sample point (indicated by the diamond's color). The measure of a line segment $\mu(M_k^K)$ is then approximated by adding the weights $w_t$ of all discretization points of the same color.}
    \label{fig:weights_pseudoexact}
  \end{minipage}
\end{figure}
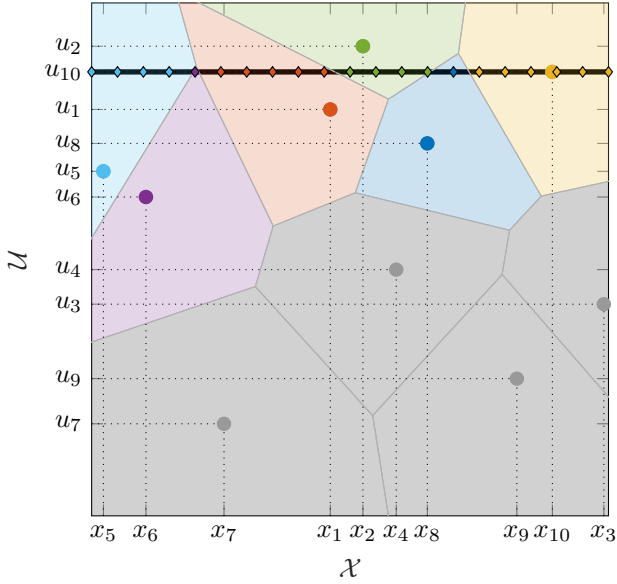

While the empirical measure provides an increasingly better approximation to the underlying true measure if the number of points grows, we might achieve a better practical performance by fixing the amount of integration points $T>0$ beforehand. In essence, this replaces the exact calculation of $\mu(M_k^K)$ by a quadrature rule based on the discretization points $(\mathbf{x}_t)_{t=1,\ldots,T}$. The idea is visualized in~\Cref{fig:weights_pseudoexact}. Note that for uniform probability distributions, we can easily choose $w_t=T^{-1}$.

\subsection{Limited memory}
Since sMMA reuses old gradient information instead of discarding it after each iteration, the memory required to store this information increases with each iteration. Especially considering 3d topology optimization problems, where the number of design elements can easily be well above $10^5$, this may pose problems when aiming for a large number of iterations. Thus, instead of opting to save all available gradient information, we can fix $\mathcal{M}\in\N$ and limit the amount of gradients stored to a maximum of $\mathcal{M}$.

While this does not impact the first iterations of sMMA at all, there inevitably arises the question which gradient information should be kept and which should be discarded, once $\mathcal{M}$ gradients are already stored. Hence, a method to measure the relative importance of a gradient sample in the storage is required, to determine which old sample is erased from memory. This, however, is precisely the information the integration weights $\alpha_k$ encode! Thus, after each iteration, we can simply remove the gradient sample with the smallest associated integration weight, if necessary. The resulting limited memory sMMA is given in~\Cref{alg:LimitedMemory}. If the batch size $\mathcal{B}$ of sMMA is chosen larger than one, we simply remove the gradient samples associated to the $\mathcal{B}$ smallest integration weights.
\begin{algorithm}
\caption{Limited memory sMMA}
\begin{algorithmic}[1]
    \STATE Choose initial design $\rho_1$ and $\mathcal{M}\in\N$.
    \FOR{$k=1,2,\ldots$}
        \STATE Draw random samples $\xi_k\sim\mu$ and $\omega_k\sim\nu$.
        \STATE Solve $\mathbf{K}(\rho_k,\xi_k)\mathbf{U}_{\xi_k,\omega_k}=\mathbf{F}(\omega_k)$ and evaluate $\mathbf{F}(\omega_k)^\top\mathbf{U}_{\xi_k,\omega_k}$.
        \STATE Calculate CSG integration weights $(\alpha_{i,k})_{i=1,\ldots,k}$.
        \STATE Calculate CSG approximations $\widehat{G}_k(\rho_k)$ and $\widehat{dG}_k(\rho_k)$ to $G(\rho_k)$ and $\nabla G(\rho_k)$.
        \STATE Construct MMA subproblem based on $\widehat{G}_k(\rho_k)$ and $\widehat{dG}_k(\rho_k)$.
        \STATE Optional: Choose move limits.
        \STATE Calculate solution $\rho^\ast$ of MMA subproblem.
        \STATE Set $\rho_{k+1}=\rho^\ast$.
        \IF{$k\ge\mathcal{M}$}
        \STATE Find index $j$ of minimal integration weight.
        \STATE Remove from memory gradient and function value information associated to $(\rho_j,\xi_j,\omega_j)$.
        \ENDIF
    \ENDFOR
\end{algorithmic}\label{alg:LimitedMemory}
\end{algorithm}

%#############################################################################################
%################################### Application #############################################
%#############################################################################################
%\FloatBarrier
\section{Numerical examples}\label{sec:Appl}
To analyze the performance of sMMA, we apply it to three examples with increasing complexity. For comparison, each problem is also solved using a standard MMA approach with appropriate discretization of the chance constraint, i.e.,
\begin{align*}
    G(\rho) &= \int_{\Xi\times\Omega} h_{a_1,a_2,a_3}\left( \mathbf{F}(\omega)^\top\mathbf{U}_{\xi,\omega} - \cm \right)\dd(\mu\times\nu)(\xi,\omega) \\
    &\approx \sum_{i=1}^{N_\Xi} \sum_{j=1}^{N_\Omega} \lambda_{ij} h_{a_1,a_2,a_3}\left( \mathbf{F}^\top(\omega_j)\mathbf{U}_{\xi_i,\omega_j}-\cm \right),
\end{align*}
where integration weights ${\lambda_{ij}\in\R}$ and integration points ${\xi_i\in\R^{d_\xi}}$, ${\omega_j\in\R^{d_\omega}}$ are given by a suitable quadrature rule.
\subsection{2d compliance with load uncertainty}
\begin{figure}
  \begin{minipage}[c]{0.6\textwidth}
    \input{wheel_setup_2}
  \end{minipage}\hfill
  \begin{minipage}[c]{0.35\textwidth}
    \caption{Design domain $\mathscr{D}$ (light grey) with fixed inner Dirichlet boundary (green) and material at the outer rim (dark grey). Depending on $\omega$, the structure is loaded by the force $F(\omega)$, acting in normal direction of the Neumann boundary (red). The force intensity $f_\omega$ associated to $F(\omega)$, see~\eqref{eq:force_intensity_wheel}, is indicated by the blue curve.}
    \label{fig:WheelSetup}
  \end{minipage}
\end{figure}
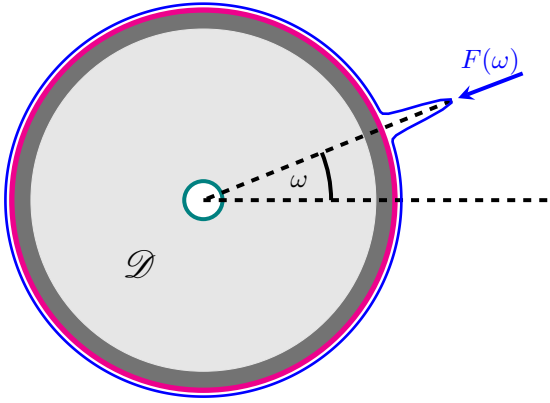

Consider a wheel of radius 1. The wheel is assumed to be fixed at a center section with radius $\tfrac{1}{10}$, and the outer rim ($\tfrac{95}{100}<r\le1$) consists of some given isotropic material. Lastly, for $\omega\in[0,2\pi]$, let $F(\omega)$ be a normal force, acting on the outer boundary of the wheel from direction $\omega$. An illustration of the setup is given in~\Cref{fig:WheelSetup}.

% \begin{figure}
%     \centering
%     \input{wheel_setup_2}
%     \caption{Design domain $\mathscr{D}$ (light grey) with fixed inner Dirichlet boundary (green) and material at the outer rim (dark grey). Depending on $\omega$, the structure is loaded by the force $F(\omega)$, acting in normal direction of the Neumann boundary (red). The force intensity $f_\omega$ associated to $F(\omega)$, see~\eqref{eq:force_intensity_wheel}, is indicated by the blue curve.}
%     \label{fig:WheelSetup}
% \end{figure}

For our calculations, the design domain was discretized into roughly $8\cdot10^4$ unstructured triangular elements. The force at position $x\in\R^2$ on the outer boundary of $\mathscr{D}$ was modeled as ${F(x,\omega) := f_\omega\big(\operatorname{arctan2}(x_1,x_2)\big)n(x)}$, where $n(x)$ denotes the inner normal vector in $x$ and $f_\omega$ is given by
\begin{equation}\label{eq:force_intensity_wheel}
    f_\omega(\beta) := 1 + \tanh\big( 10^3(\cos(\beta-\omega)-1)+10^{-1}\big),
\end{equation}
see~\Cref{fig:WheelSetup}. Moreover, $\omega$ is assumed to follow a uniform distribution on $[0,2\pi]$, while the stiffness of material and void was chosen as 1 and $10^{-4}$, respectively. The parameters for the smoothing function~\eqref{eq:smooth_cc_fun} are ${a_1=50}$, ${a_2=\tfrac{1}{10}}$ as well as ${a_3=5}$ and the probability level was chosen as ${p=0.025}$.

It is important to note that in the resulting optimization problem
\begin{align*}
        \min_{\rho\in\R^d}\quad & \rvol(\rho),\\
        \text{s.t.}\quad & \frac{1}{2\pi}\int_0^{2\pi}h_{a_1,a_2,a_3}\left(\mathbf{F}(\omega)^\top\mathbf{U}_{\omega}-\cm\right)\dd{\omega}\le p,\\
            &\mathbf{K}(\rho)\mathbf{U}_{\omega}=\mathbf{F}(\omega),\\
            &0\le \rho\le 1,
\end{align*}
uncertainties enter only via the right hand side of the state equation. Coupled with the fact that the design dimension is moderate, a direct solver, e.g., sparse Cholesky decomposition, can be used to solve the system of linear equations. Thus, it is possible to solve the state equation for many realizations of $\omega$ simultaneously, without significantly increasing the numerical effort. This allows us to thoroughly analyze the performance of sMMA and MMA, before considering more advanced applications later on. 

Both methods were initialized with ${\rho_1\equiv 0.75}$ and set to perform 400 iterations each. The SIMP parameter was picked as 10 and increased to 15 after 200 iterations. Note that these values are pretty large when compared to typical choices for SIMP parameters (${s\sim3}$) and also much higher than for our other numerical experiments. However, in this setup, fictitious material is very profitable for the optimizer and requires a harsh penalization in order to obtain black \& white designs. Furthermore, note that we did not observe that the optimizers got stuck in poor local minima due to the high penalty parameter. 

In each iteration, each optimizer was allowed to evaluate function and gradient values for a batch size $\mathcal{B}$ of different realizations of $\omega$ (for MMA, they were picked evenly spaced, corresponding to a trapezoidal integration rule on a circle). Moreover, we imposed different move limits $\tau$ on the methods to observe their resulting behavior. All in all, each combination of ${\mathcal{B}\in\{4,8,16,32,64\}}$ and ${\tau\in\{1,\tfrac{3}{4},\tfrac{1}{2},\tfrac{1}{4}\}}$ was tested. 

Since both approaches work with a different internal approximation of the chance constraint (sMMA uses the CSG integration scheme with $\mathcal{B}$ random samples per iteration, whereas MMA uses a trapezoidal quadrature rule with $\mathcal{B}$ equidistantly spaced integration points), we are highly interested in the quality of this approximation. Therefore, each intermediate design is additionally analyzed using 1080 integration points on the circle, which we consider to yield a sufficiently exact value of the chance constraint. 

The relative physical volumes and chance constraint values of all final designs can be found in~\Cref{fig:wheel_pareto}. As it turns out, all designs obtained by MMA with ${\mathcal{B}\in\{4,8,16\}}$ are infeasible, due to the large difference between true and discretized chance constraint. In fact, the discrepancy is so large, that it can directly be seen when looking at the final designs shown in~\Cref{tab:wheel_designs_mcmsa,tab:wheel_designs_mma}. Only for large enough batch sizes ${\mathcal{B}\in\{32,64\}}$ does MMA produce feasible (up to tolerance) designs. In contrast, sMMA does not struggle with small batch sizes and produces feasible designs in each setup. As an example, the chance constraint values during the optimization for $\tau=1$ and ${\mathcal{B}\in\{16,32\}}$ are shown in~\Cref{fig:wheel_cc}. Therein, we can see that the internal sMMA estimator ${\widehat{G}_n(\rho_n)}$ very closely approximates ${G(\rho_n)}$, even in early iterations. 

To analyze the effect of our smoothing techniques on the original problem, a comparison between ${G(\rho_n)}$ and the nonsmooth chance constraint (as appearing in~\eqref{eq:GeneralOptProb}) can be found in~\Cref{fig:wheel_smoothingeffect}. As predicted above, we see that the smooth chance constraint approximation overestimates the actual nonsmooth chance constraint. In fact, all feasible final designs are robust, i.e., they satisfy ${\mathbf{F}(\omega)^\top\mathbf{U}_\omega \le \cm}$ for all ${\omega\in\Omega}$.

Moreover, we see that adding the smoothmax-term to steepen the gradient has only a negligible impact on the resulting value. Note that, for the depicted sMMA design, it actually leads to a smaller chance constraint value. This is due to the fact that 
\begin{equation*}
    t-\frac{t}{1+\exp(a_3t)}
\end{equation*}
admits a local minimum with negative function values at ${t^\ast\approx-\tfrac{1.28}{a_3}}$. As a result, the optimizer is encouraged (very slightly) to undershoot the compliance bound by this margin. This can again be seen as measure of caution, since the obtained designs are more likely to satisfy the nonsmooth chance constraint due to an additional tolerance for errors.

\begin{figure} 
    \centering
    \begin{subfigure}{0.45\textwidth}
        \input{wheel_obj_scatter}
    \end{subfigure}\hfill
    \begin{subfigure}{0.45\textwidth}
        \input{wheel_phy_scatter}%
    \end{subfigure}
    \caption{Relative volume of final design without SIMP interpolation (left) and with SIMP interpolation (right) for sMMA (blue) and MMA (red). Only volumes corresponding to feasible solutions are shown, i.e., all results obtained by sMMA as well as designs found by MMA with ${\mathcal{B}\in\{32,64\}}$. Since an increased batch size corresponds to a numerically more expensive approximation to the true chance constraint, designs obtained with a smaller batch size are expected to have slightly worse objective function values (e.g., ${\mathcal{B}=4}$ in left picture). Differences between $\rvol$ and $\pvol$ are mostly related to the boundary length of ${\{x\in\mathscr{D}\,:\, \rho(x)\approx1\}}$. That is, finer structures in the design lead to larger gaps between $\rvol$ and $\pvol$, see~\Cref{tab:wheel_designs_mcmsa,tab:wheel_designs_mma}.}
    \label{fig:wheel_objphyscatter}
\end{figure}
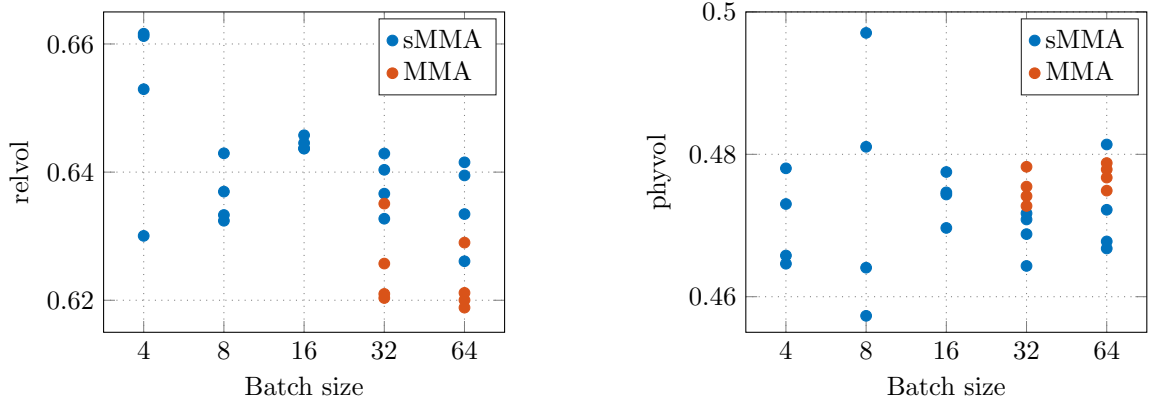

Lastly, we want to analyze the impact of $\mathcal{B}$ on the optimization process. While we have already seen that a large enough value of $\mathcal{B}$ is crucial for MMA to yield feasible designs, this is not the case for sMMA. However, since a larger batch size increases the number of available samples (and numerical cost), we still expect sMMA to benefit from the additional information. This is indeed the case, as the overview of the associated relative volumes and relative physical volumes (\Cref{fig:wheel_objphyscatter}) shows a trend of decreasing objective function values for larger batch sizes.

\begin{figure} 
    \centering
    \begin{subfigure}{0.45\textwidth}
        \input{wheel_pareto}
    \end{subfigure}\hfill
    \begin{subfigure}{0.45\textwidth}
        \input{wheel_pareto_zoom}
    \end{subfigure}
    \caption{Final relative physical volumes and chance constraint values (evaluated using a trapezoidal rule with 1080 load cases) for sMMA (blue circles) and MMA (red triangles). In the right diagram, only feasible designs (all sMMA results and MMA results for ${\mathcal{B}\in\{32,64\}}$) are shown. Of all feasible designs, the lowest relative physical volume corresponds to sMMA with ${\mathcal{B}=8}$ and ${\tau=1}$, while the highest value of $\pvol$ is associated to sMMA with ${\mathcal{B}=8}$ and ${\tau=\tfrac{1}{4}}$. Values of $\rvol$ and $\pvol$, sorted by batch size, can be found in~\Cref{fig:wheel_objphyscatter}.}
    \label{fig:wheel_pareto}
\end{figure}
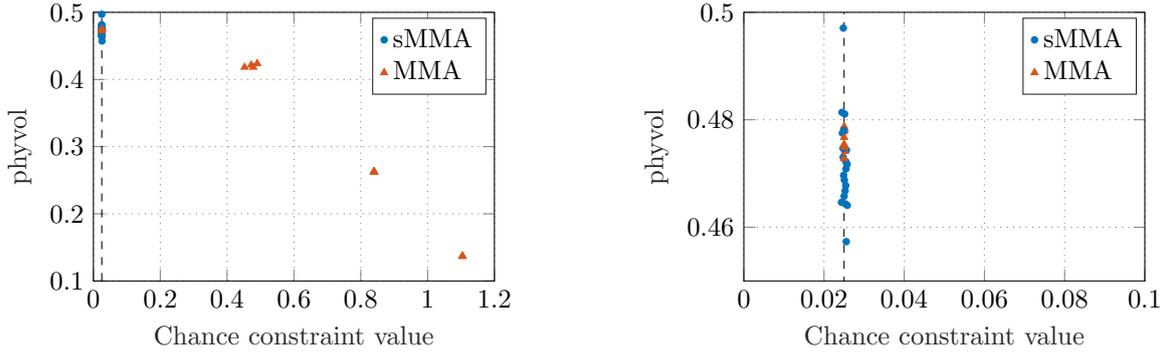

\begin{figure} 
    \centering
    % \begin{subfigure}{0.5\textwidth}
    %     \centering
    %     \input{mcmsa_8_cc}
    % \end{subfigure}%
    % \begin{subfigure}{0.5\textwidth}
    %     \centering
    %     \input{mma_8_cc}%
    % \end{subfigure}
    \begin{subfigure}{0.5\textwidth}
        \centering
        \input{mcmsa_16_cc}
    \end{subfigure}%
    \begin{subfigure}{0.5\textwidth}
        \centering
        \input{mma_16_cc}%
    \end{subfigure}
    \begin{subfigure}{0.5\textwidth}
        \centering
        \input{mcmsa_32_cc}
    \end{subfigure}%
    \begin{subfigure}{0.5\textwidth}
        \centering
        \input{mma_32_cc}%
    \end{subfigure}
    \caption{Chance constraint values during the sMMA (left) and MMA (right) optimization process for ${\tau=1}$ and ${\mathcal{B}\in\{16,32\}}$. The blue curve (internal) corresponds to chance constraint values as approximated by the respective optimizer, i.e., a trapezoidal rule with $\mathcal{B}$ equidistantly spaced integration points for MMA or the sample-based approximation for sMMA. The red curve (exact) is obtained by using a trapezoidal rule with 1080 equidistant integration points. We can see a large error in the MMA discretization for ${\mathcal{B}=16}$. In contrast, sMMA always yields a close approximation to the exact values, even at early iterations.}
    \label{fig:wheel_cc}
\end{figure}
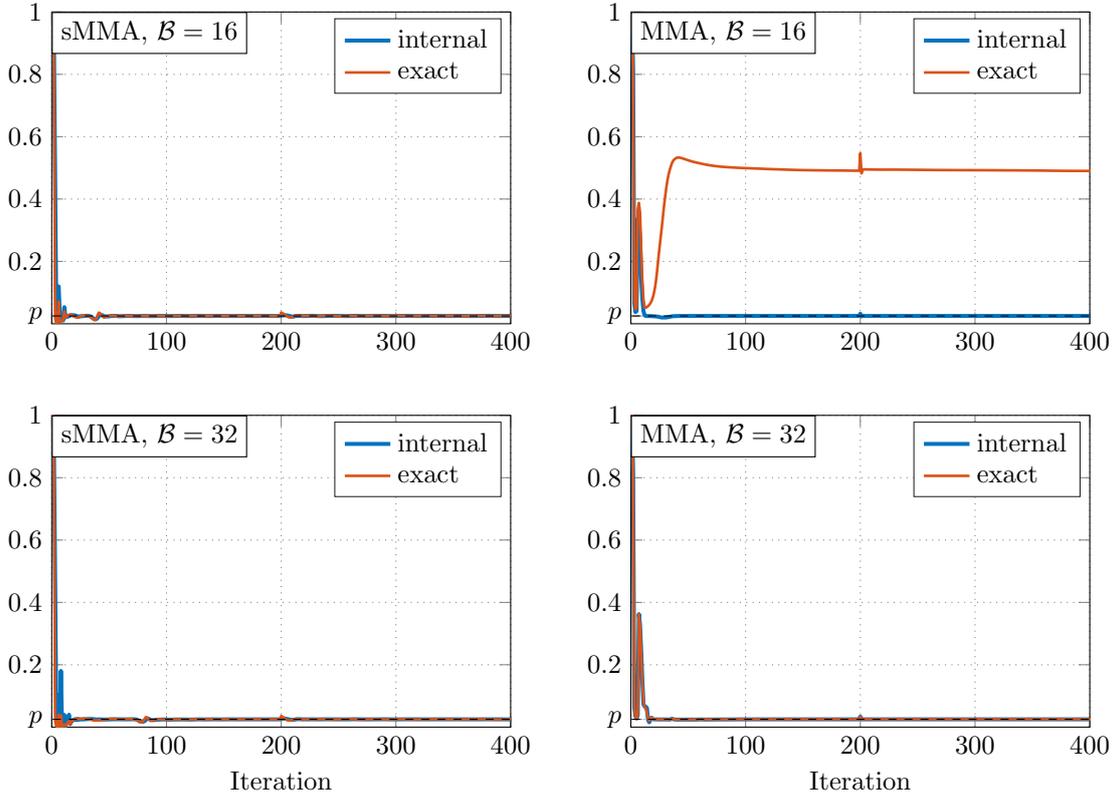

\begin{figure} 
    \centering
    \begin{subfigure}{0.45\textwidth}
        \input{mcmsa_16_1_smoothingeffect}
    \end{subfigure}\hfill
    \begin{subfigure}{0.45\textwidth}
        \input{mma_16_1_smoothingeffect}%
    \end{subfigure}
    \caption{Evolution of chance constraint values for sMMA (left) and MMA (right) with ${\tau=1}$ and ${\mathcal{B}=16}$. The curves are obtained by a trapezoidal rule with 1080 equidistant integration points when modelling the chance constraint with the nonsmooth indicator function ${\chi_{(0,\infty)}}$ (blue), the ${\tanh}$-approximation used in $h_{a_1}$ (red) and the smooth and steep approximation $h_{a_1,a_2,a_3}$ (yellow) used for the optimization, see~\eqref{eq:smooth_cc_fun}. As expected, the difference between the smooth approximations is insignificant and both provide an upper bound to the true chance constraint without regularization.}
    \label{fig:wheel_smoothingeffect}
\end{figure}
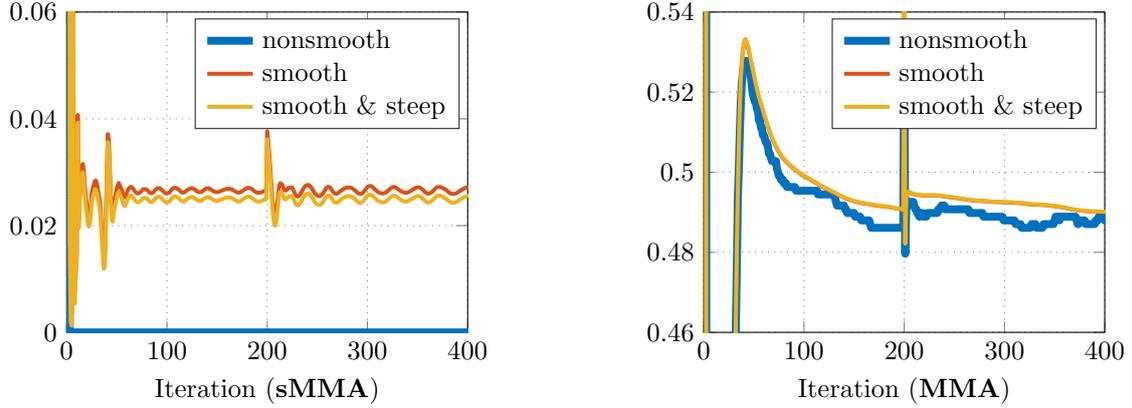

\begin{table}
    \centering
    \begin{tabular}{l|ccccc}
    & $\mathcal{B}=4$ & $\mathcal{B}=8$ & $\mathcal{B}=16$ & $\mathcal{B}=32$ & $\mathcal{B}=64$ \\[1em] \hline
        $\tau=1$ & \includegraphics[width=.135\textwidth,keepaspectratio,valign=c]{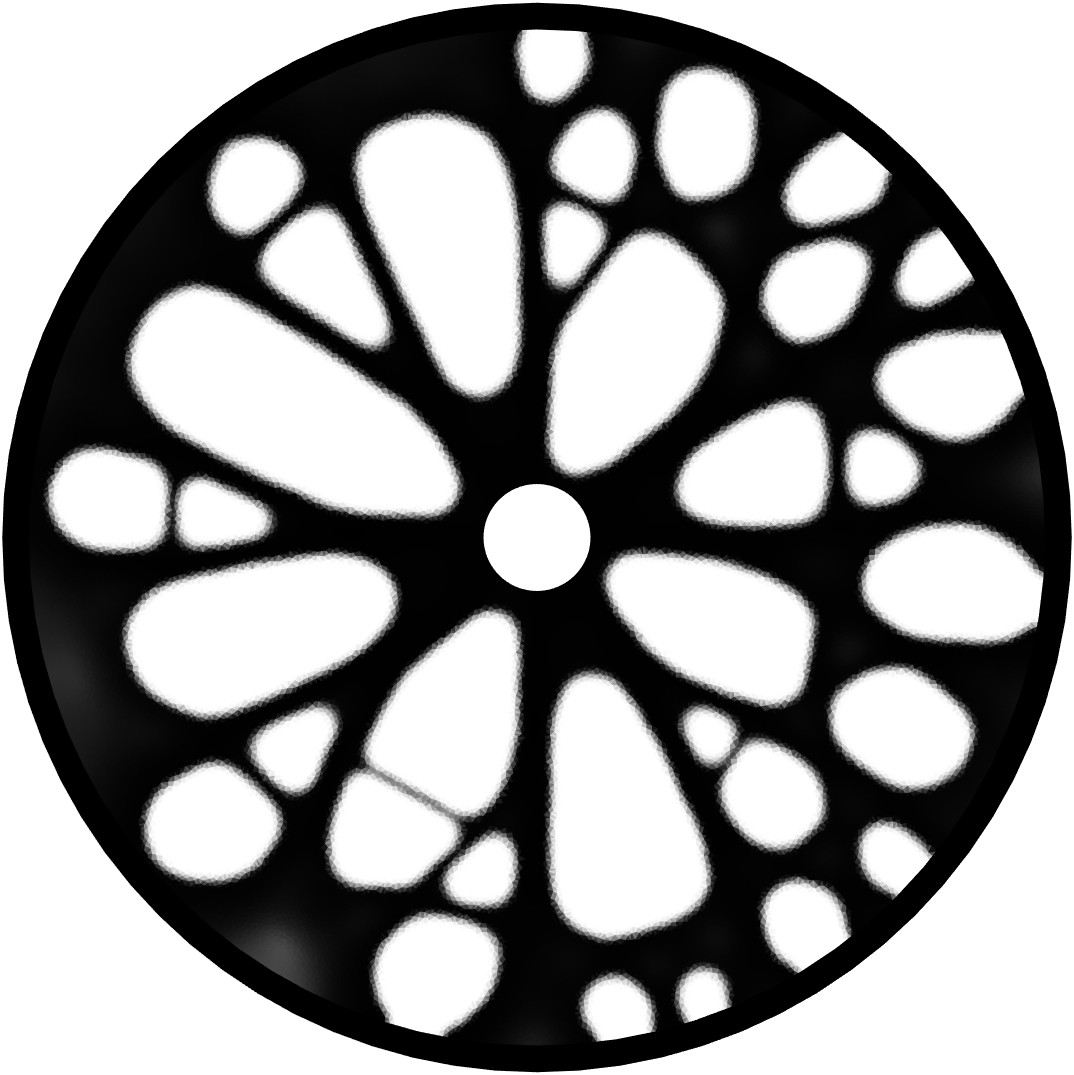} & \includegraphics[width=.135\textwidth,keepaspectratio,valign=c]{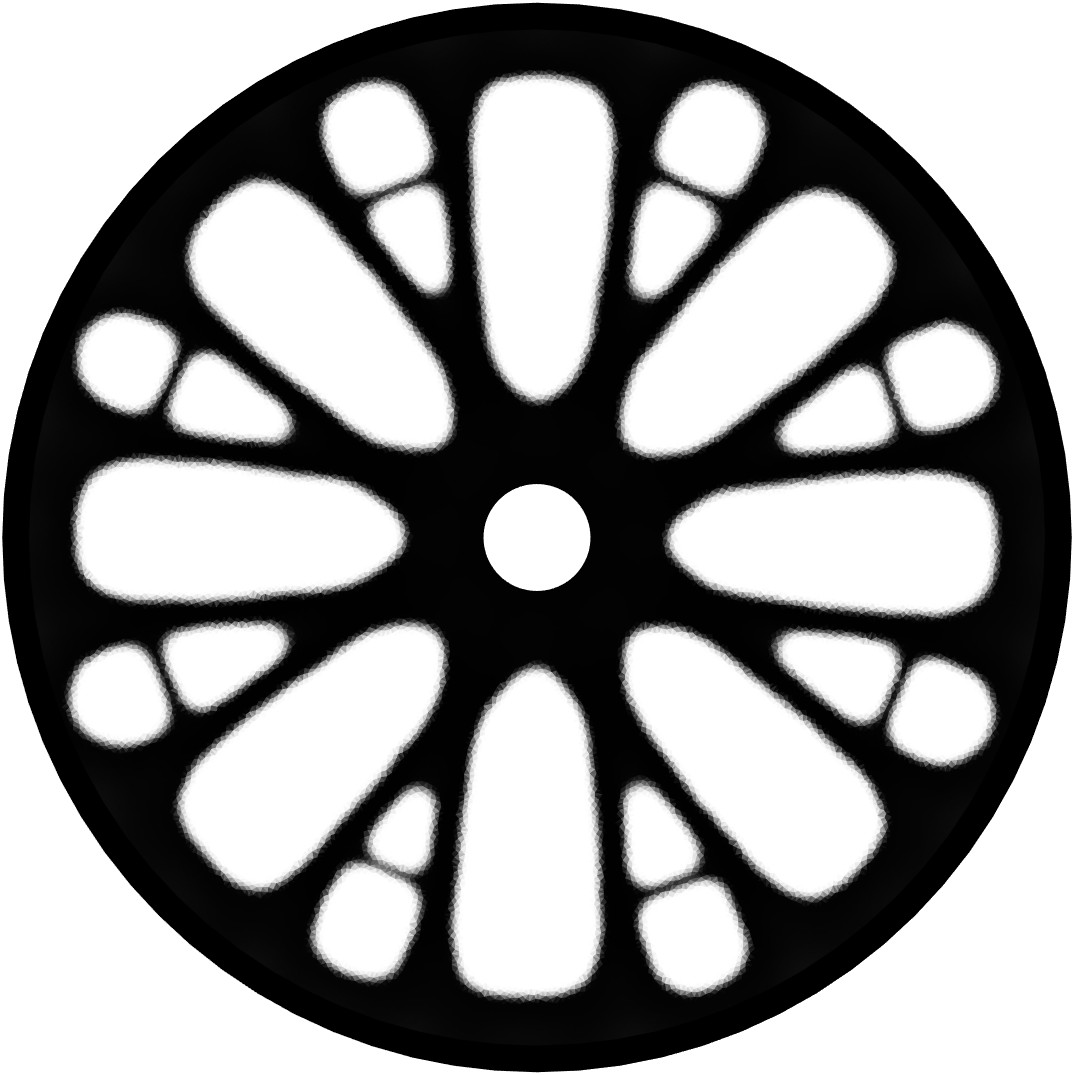} & \includegraphics[width=.135\textwidth,keepaspectratio,valign=c]{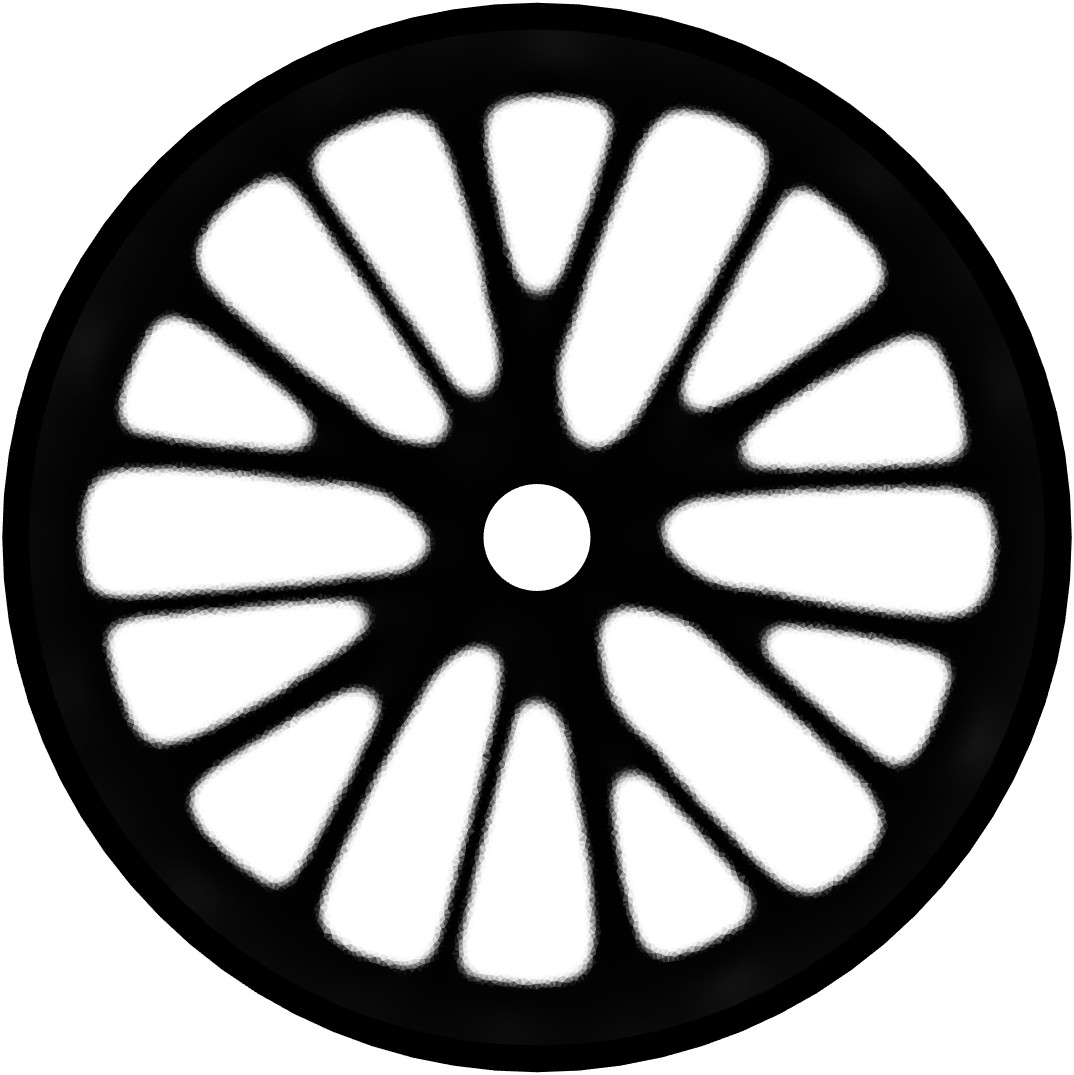} & 
         \includegraphics[width=.135\textwidth,keepaspectratio,valign=c]{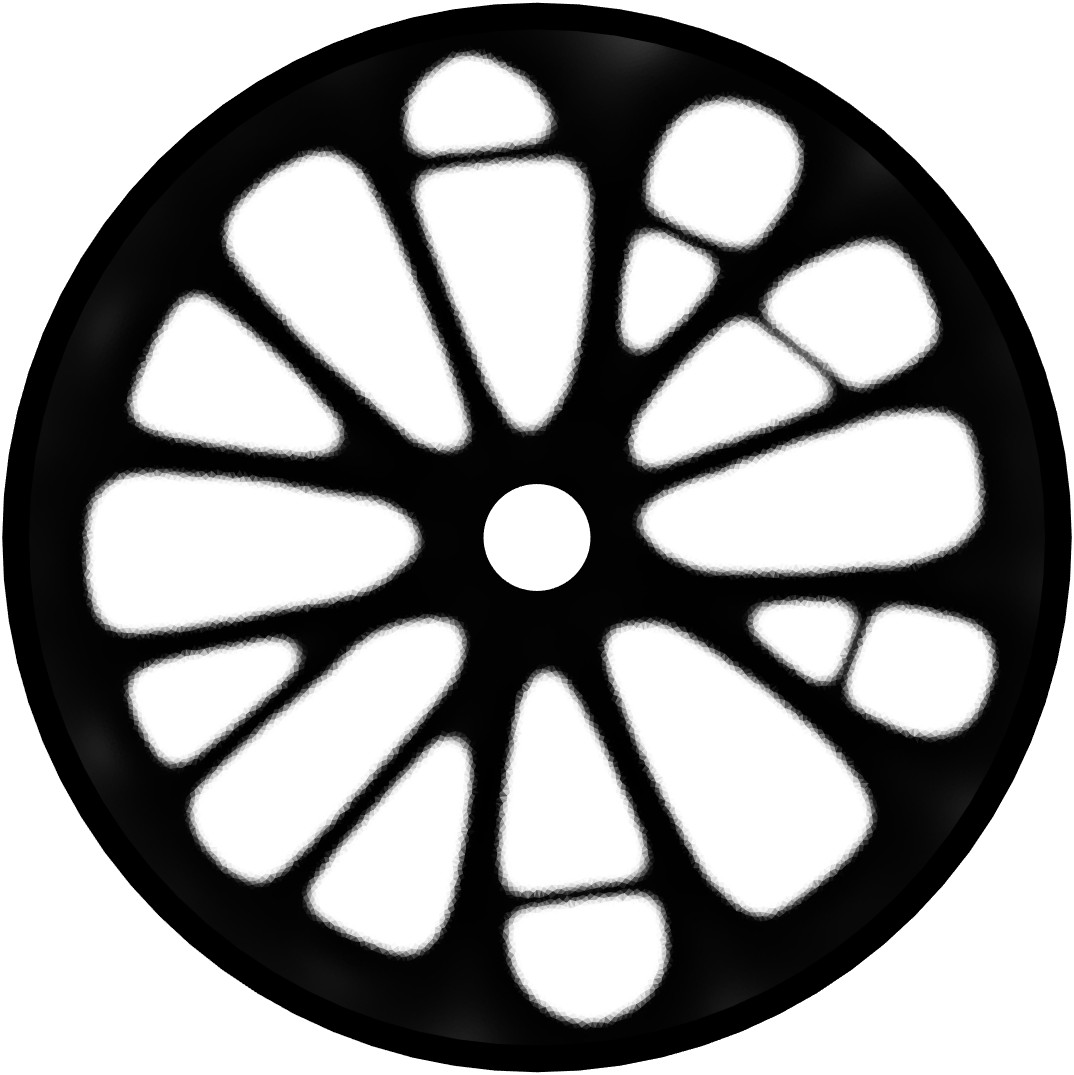} & \includegraphics[width=.135\textwidth,keepaspectratio,valign=c]{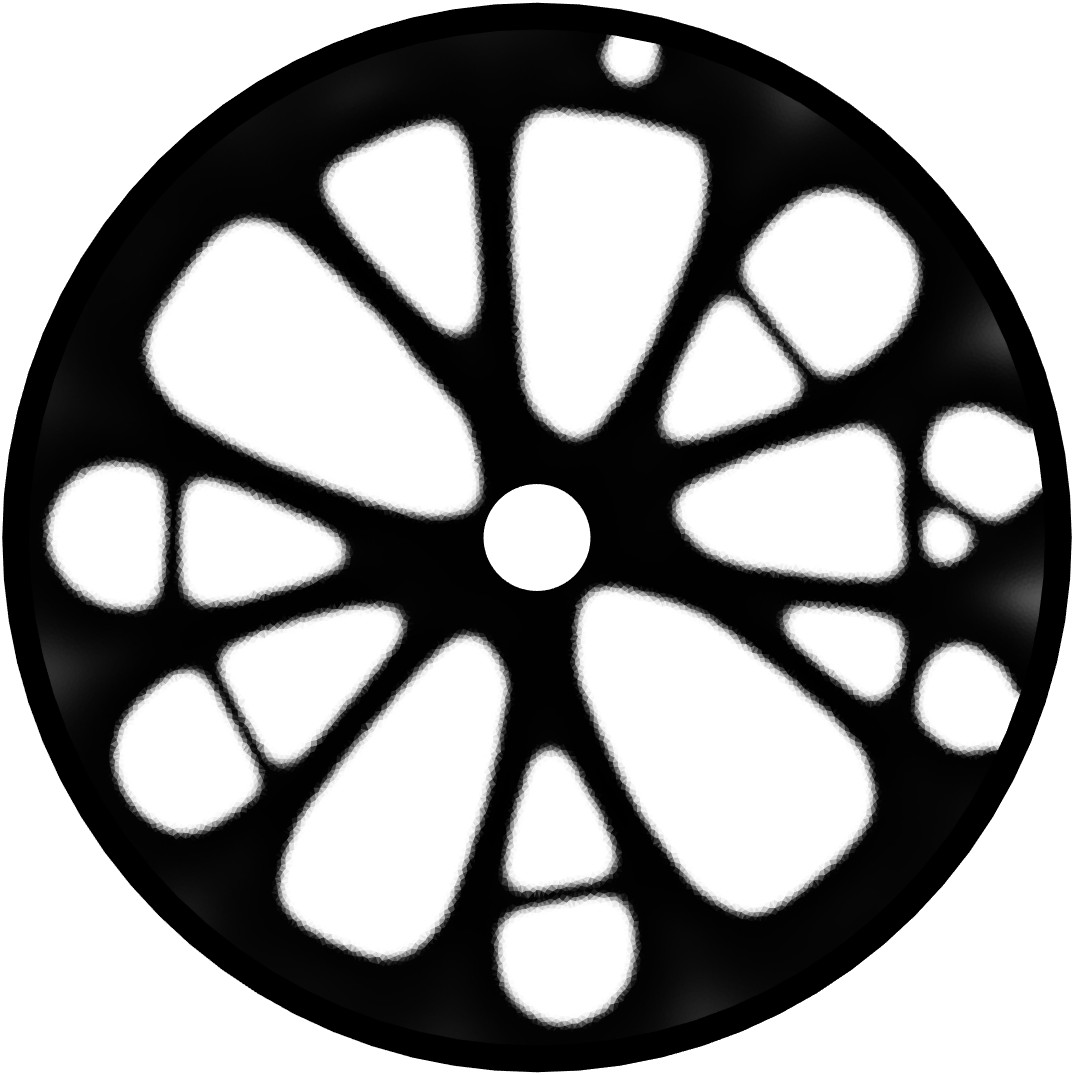}\\
        $\tau=\tfrac{3}{4}$ & \includegraphics[width=.135\textwidth,keepaspectratio,valign=c]{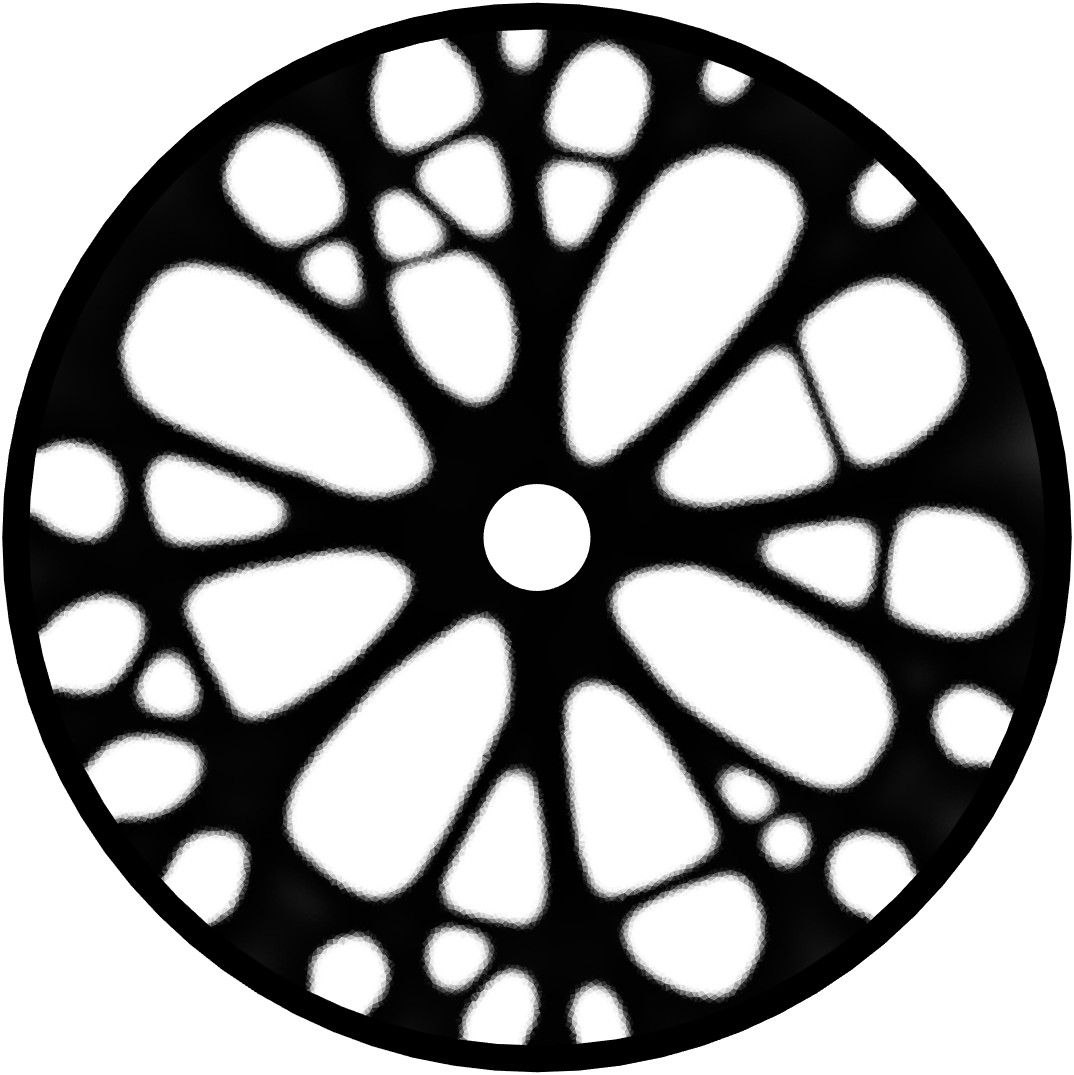} & \includegraphics[width=.135\textwidth,keepaspectratio,valign=c]{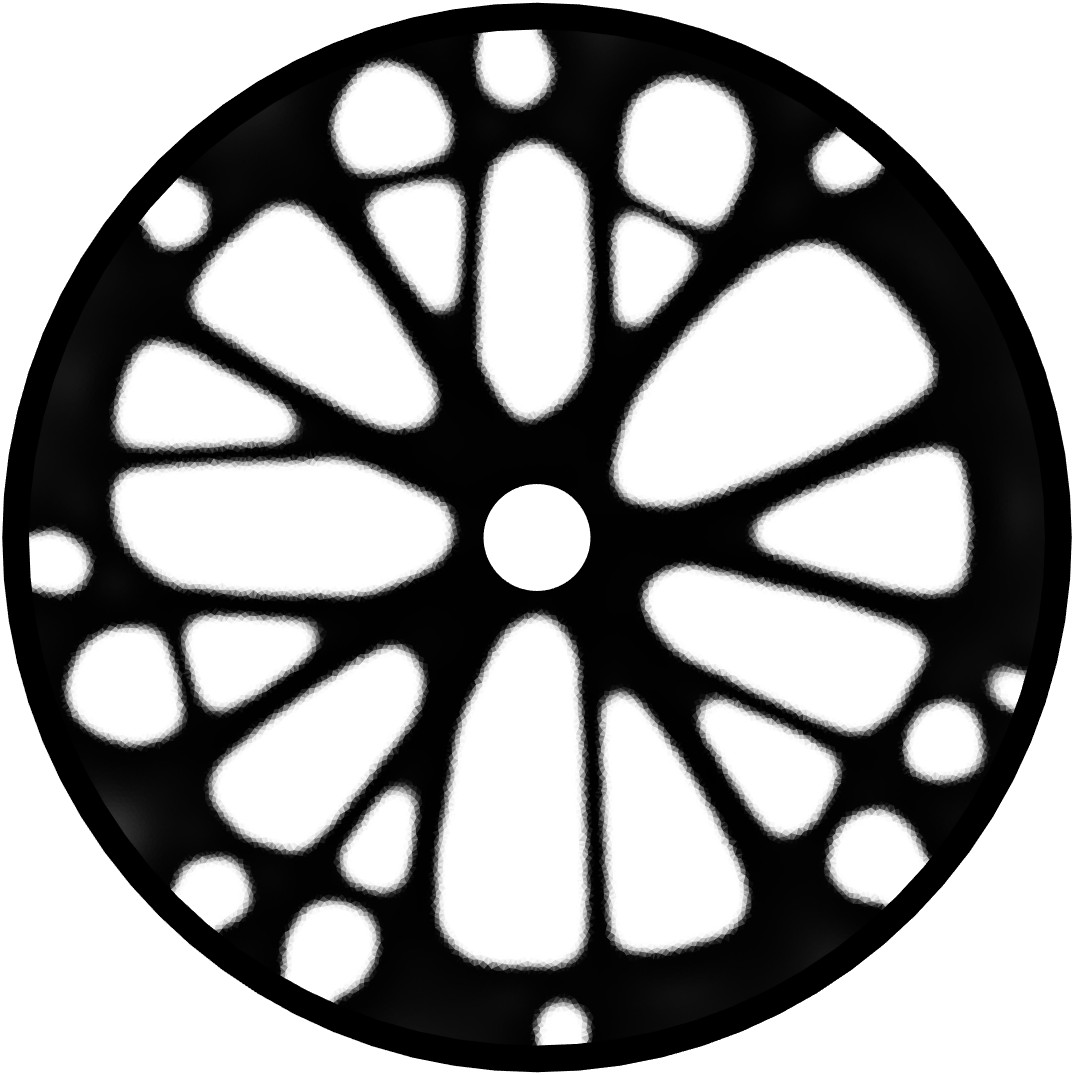} & \includegraphics[width=.135\textwidth,keepaspectratio,valign=c]{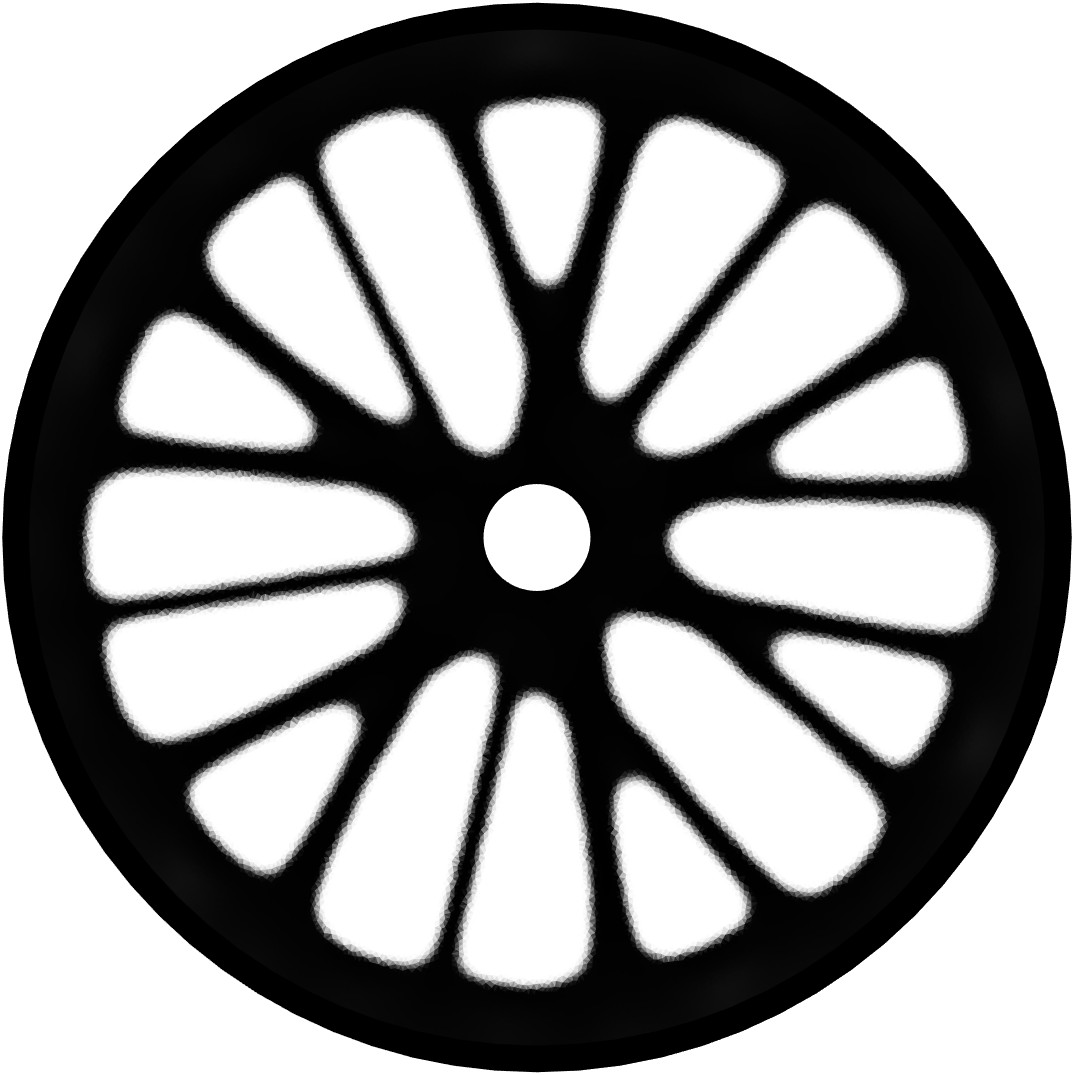} & 
         \includegraphics[width=.135\textwidth,keepaspectratio,valign=c]{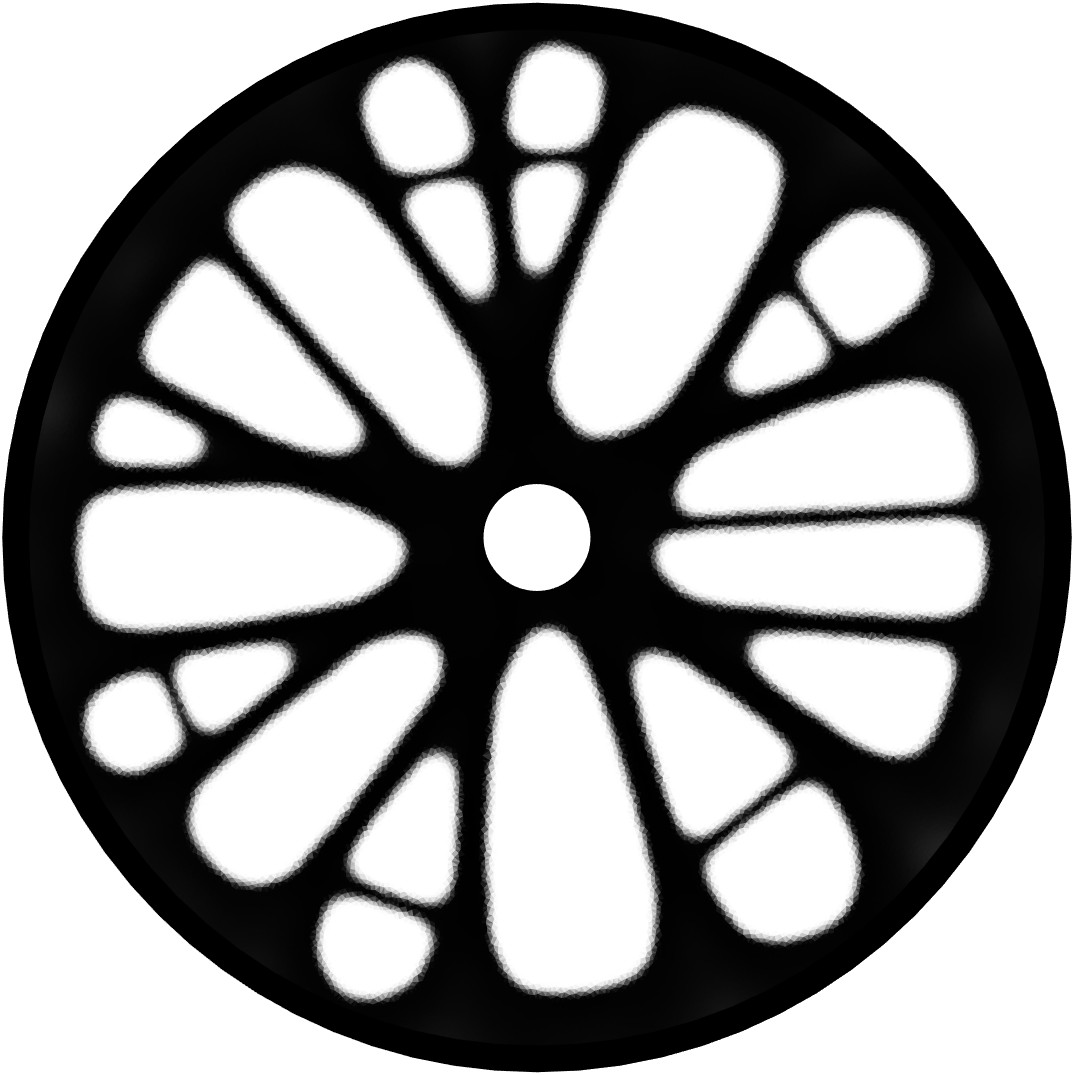} & \includegraphics[width=.135\textwidth,keepaspectratio,valign=c]{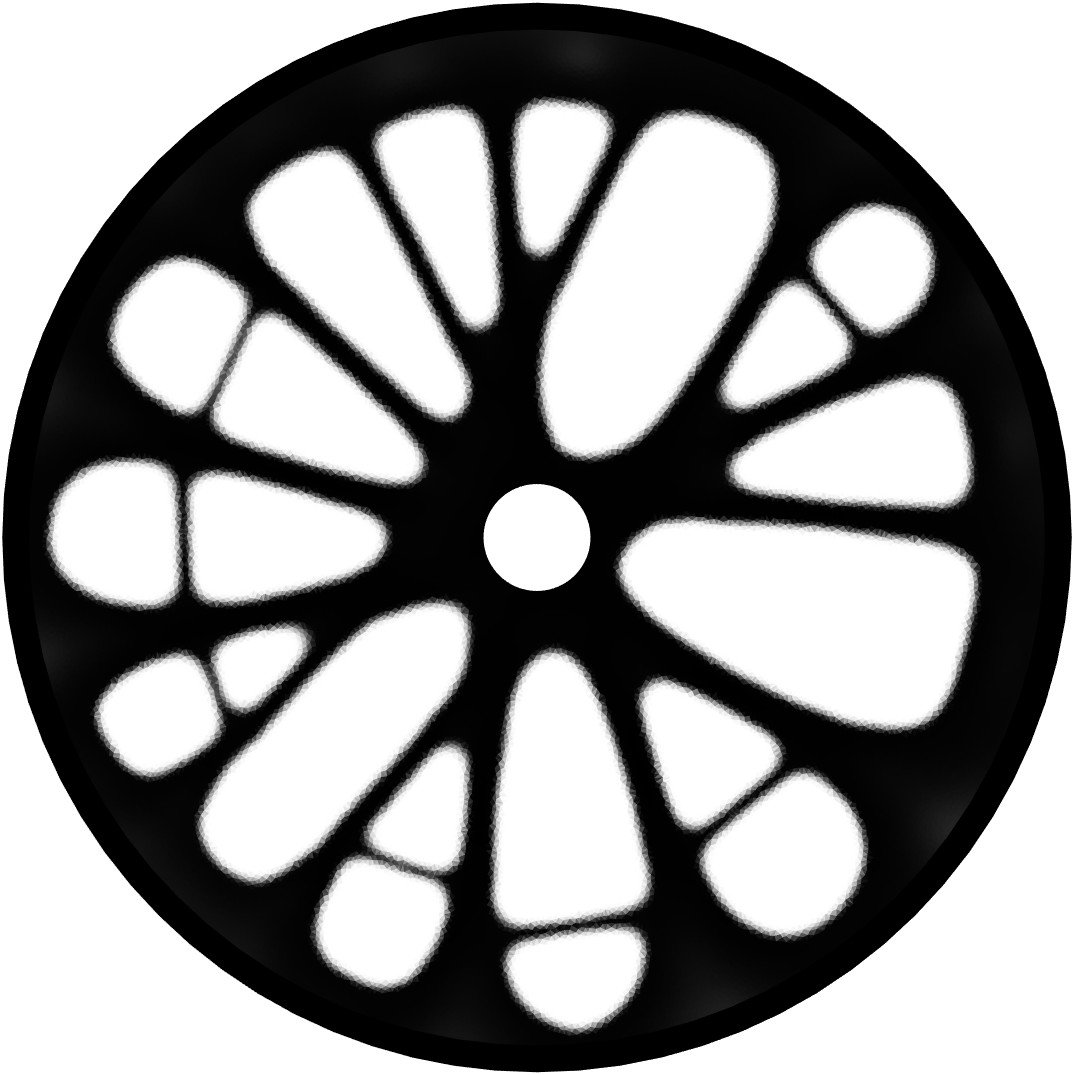} \\
         $\tau=\tfrac{1}{2}$ & \includegraphics[width=.135\textwidth,keepaspectratio,valign=c]{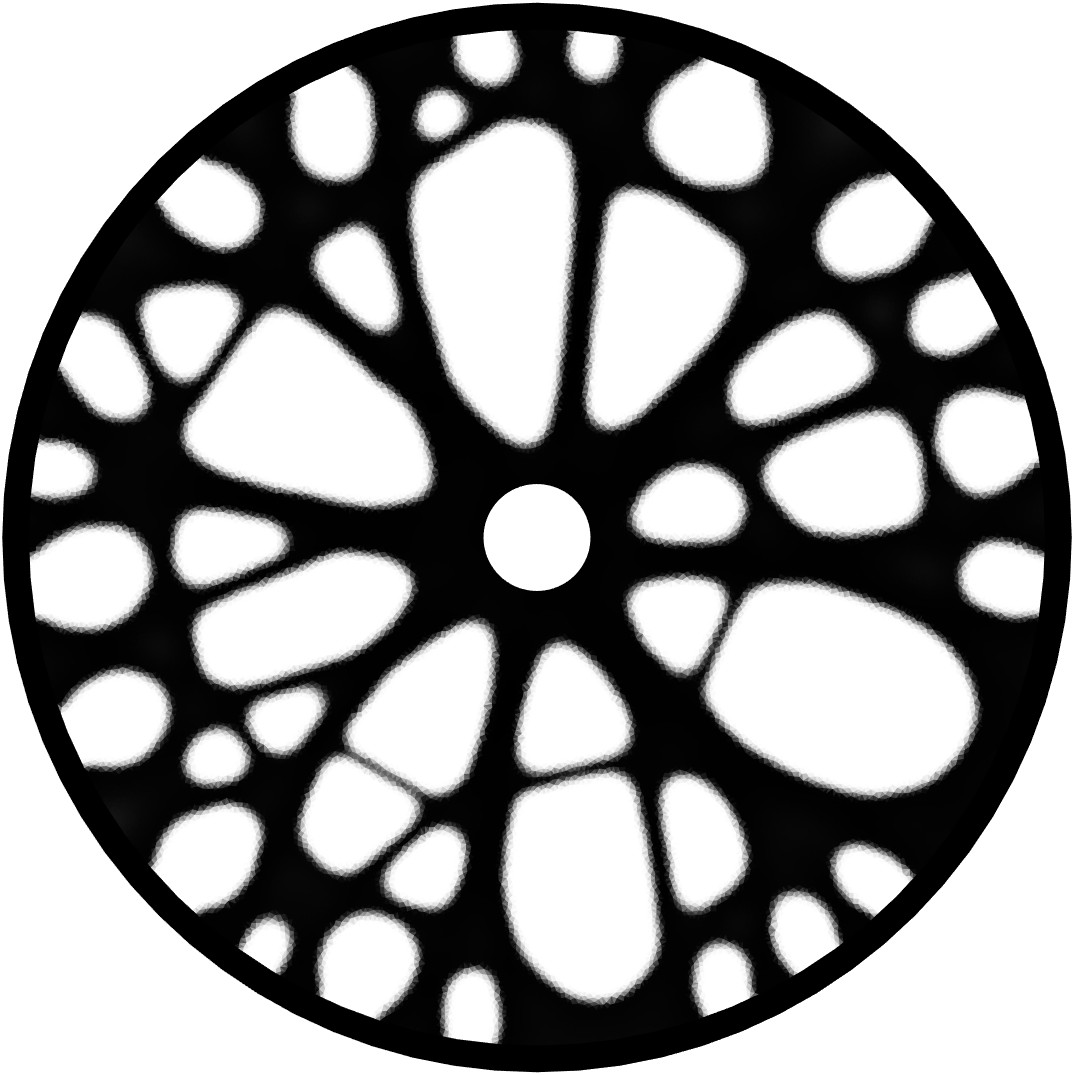} & \includegraphics[width=.135\textwidth,keepaspectratio,valign=c]{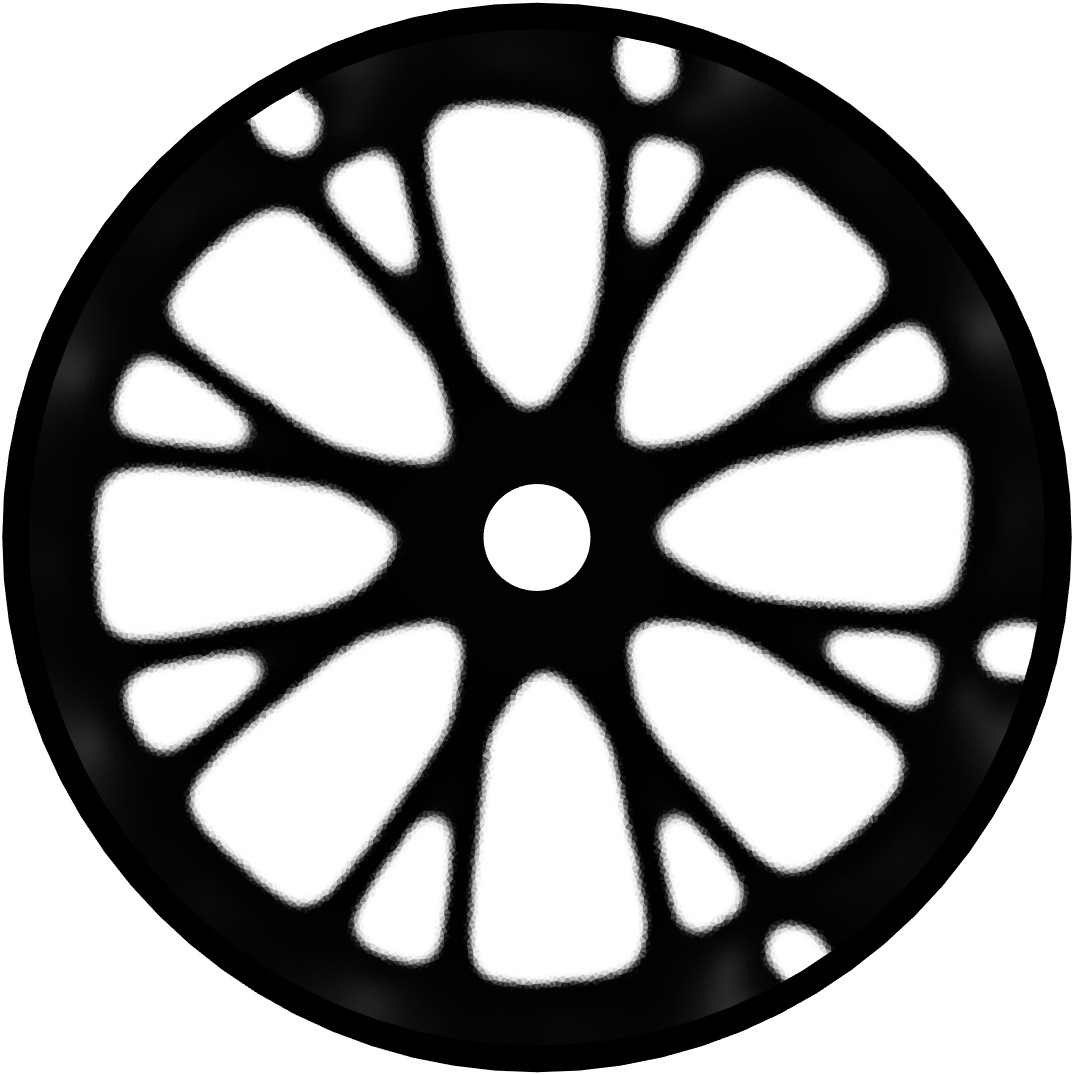} & \includegraphics[width=.135\textwidth,keepaspectratio,valign=c]{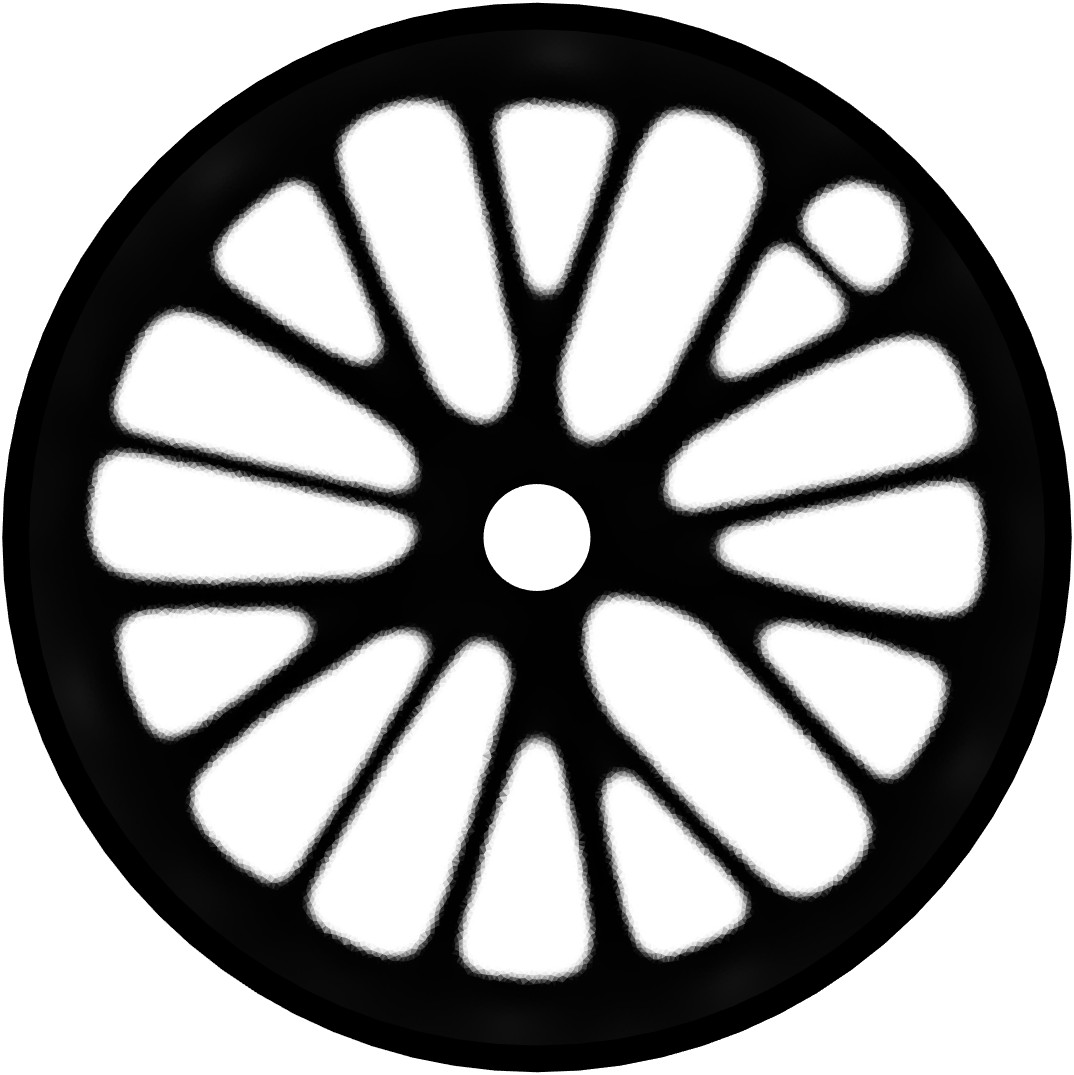} & 
         \includegraphics[width=.135\textwidth,keepaspectratio,valign=c]{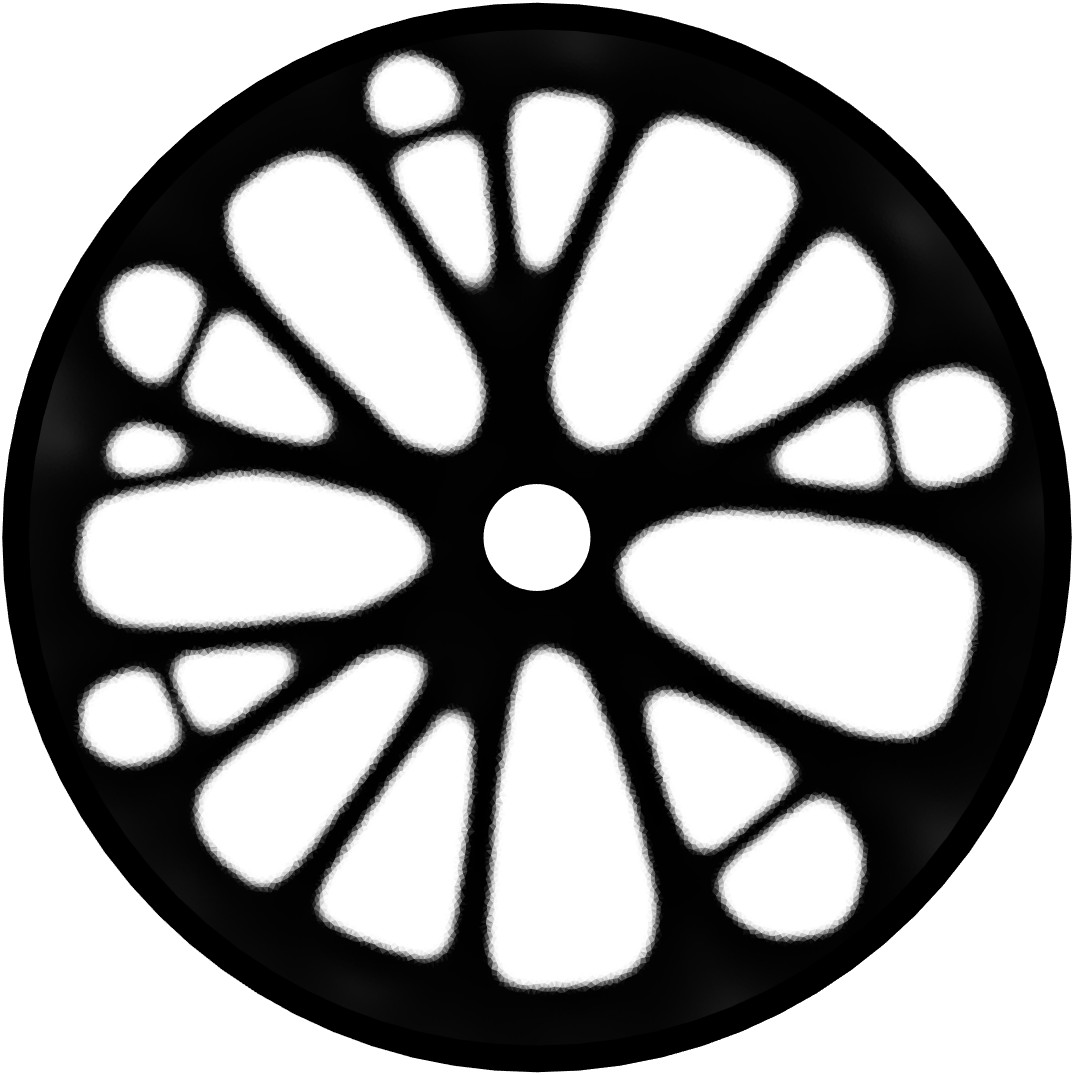} & \includegraphics[width=.135\textwidth,keepaspectratio,valign=c]{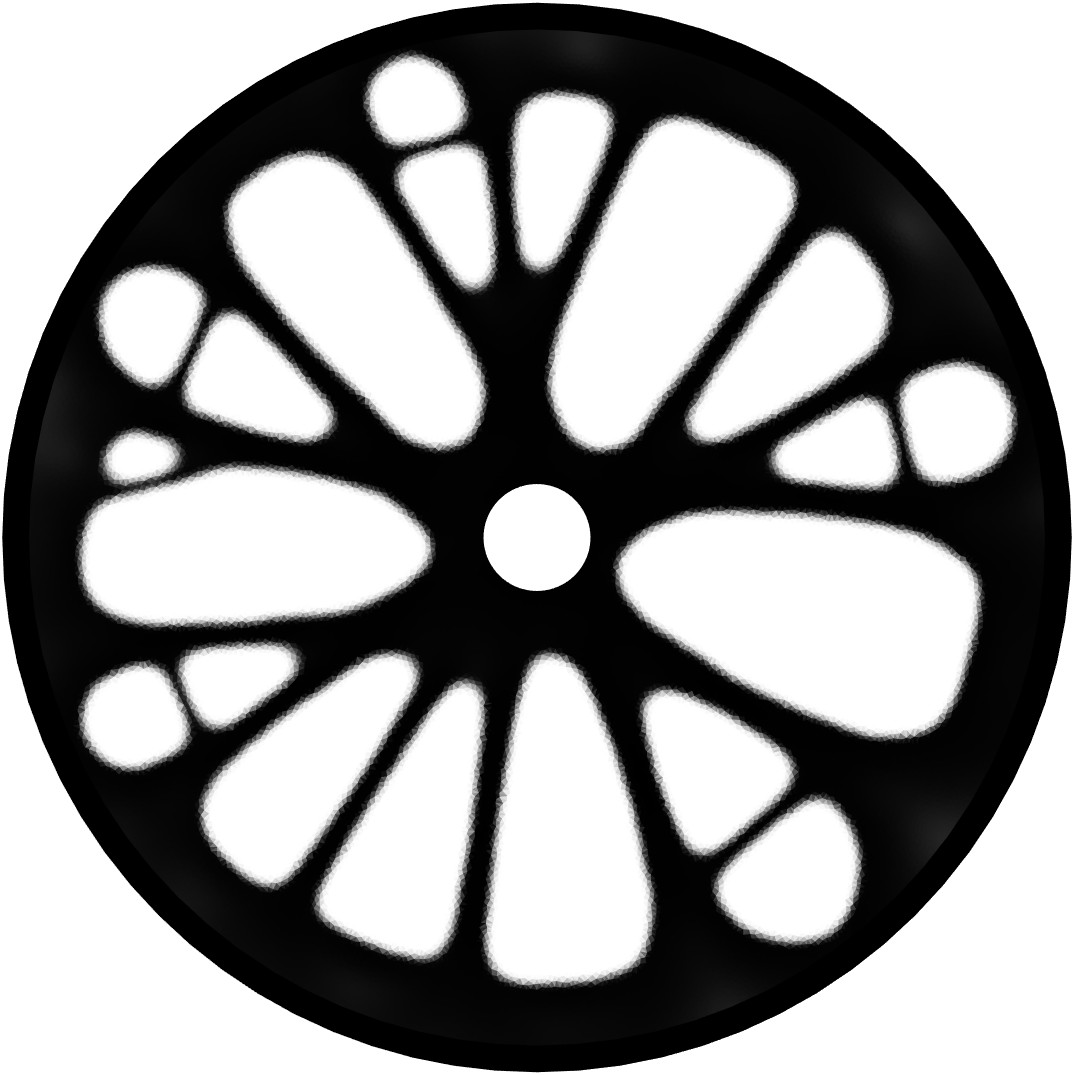}\\
         $\tau=\tfrac{1}{4}$ & \includegraphics[width=.135\textwidth,keepaspectratio,valign=c]{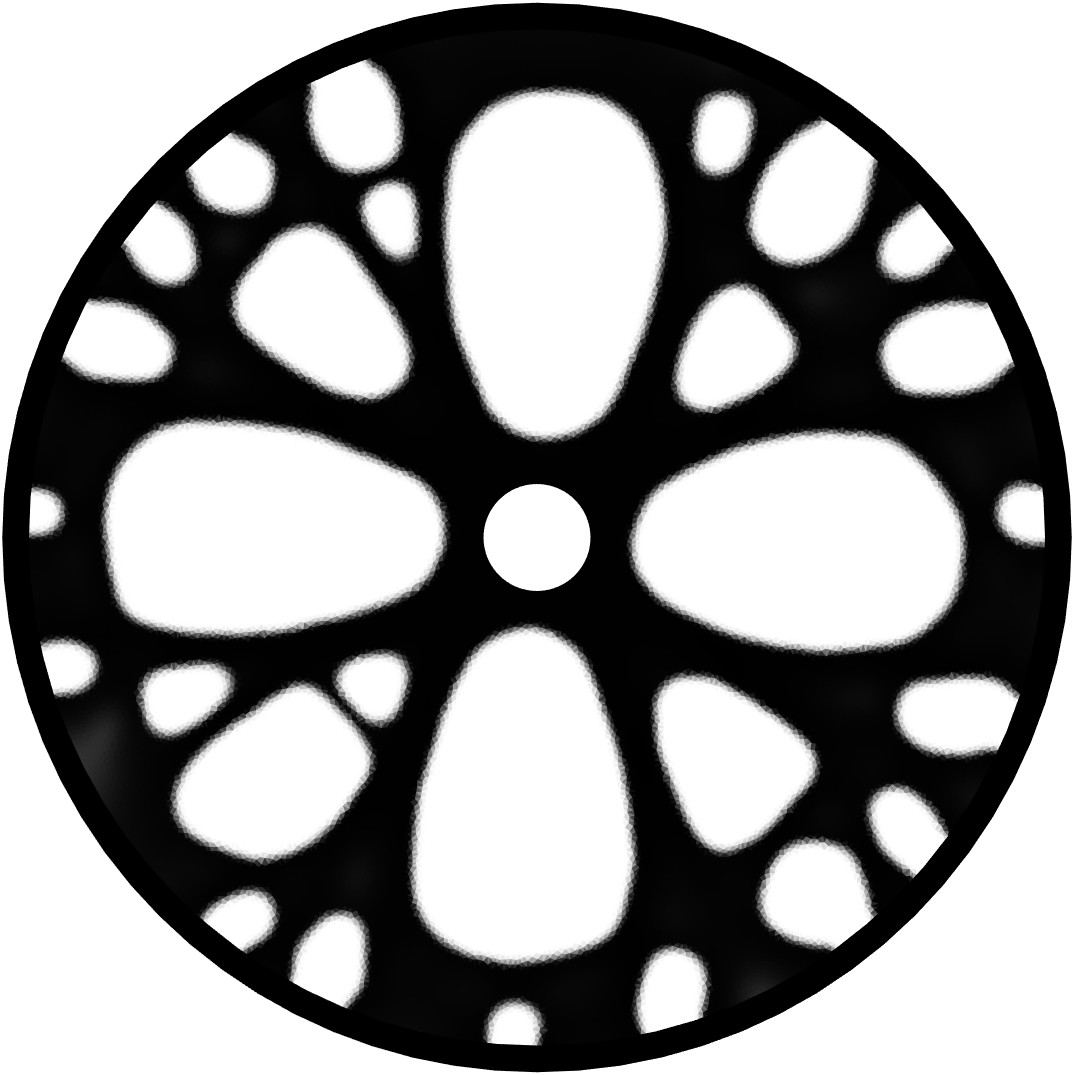} & \includegraphics[width=.135\textwidth,keepaspectratio,valign=c]{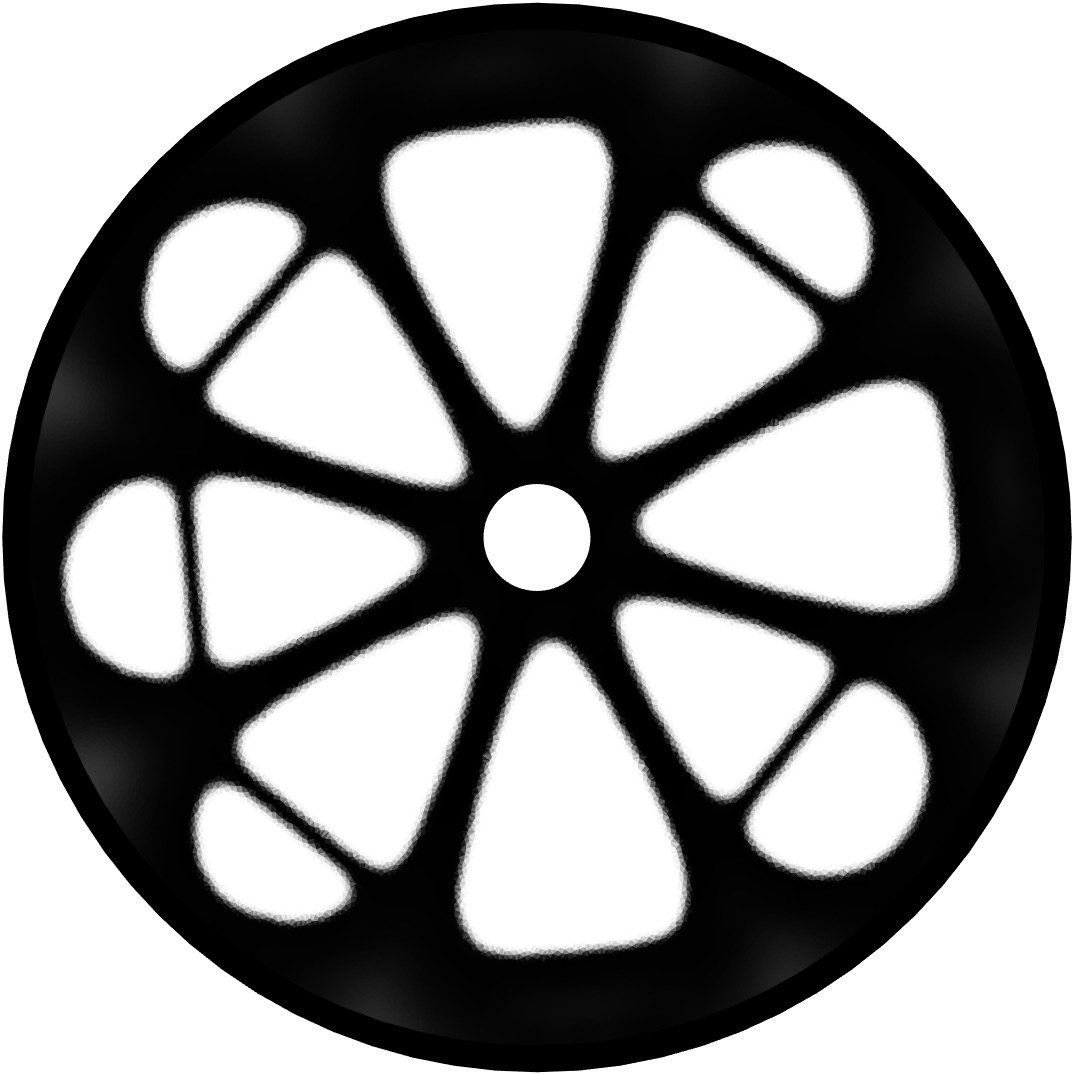} & \includegraphics[width=.135\textwidth,keepaspectratio,valign=c]{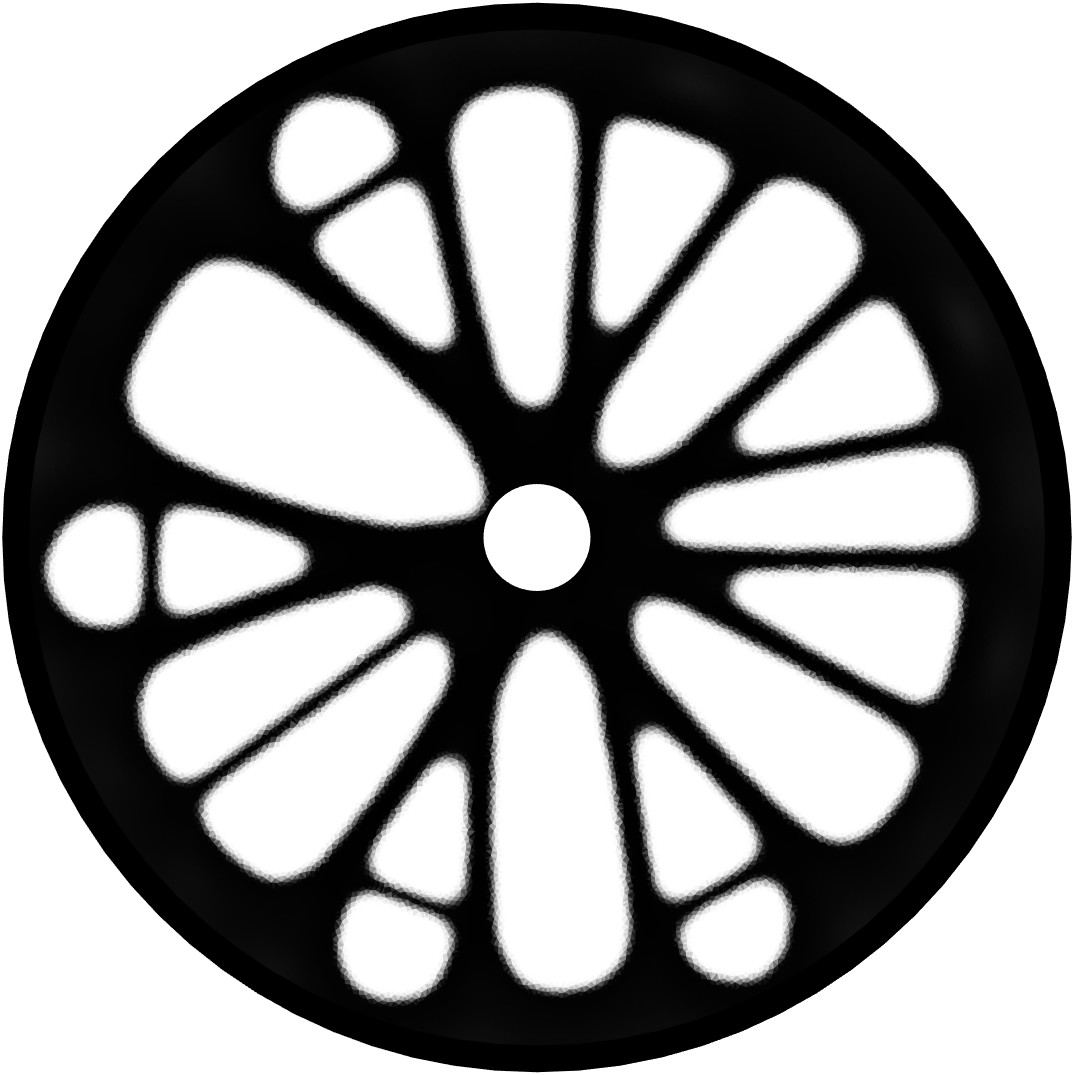} & 
         \includegraphics[width=.135\textwidth,keepaspectratio,valign=c]{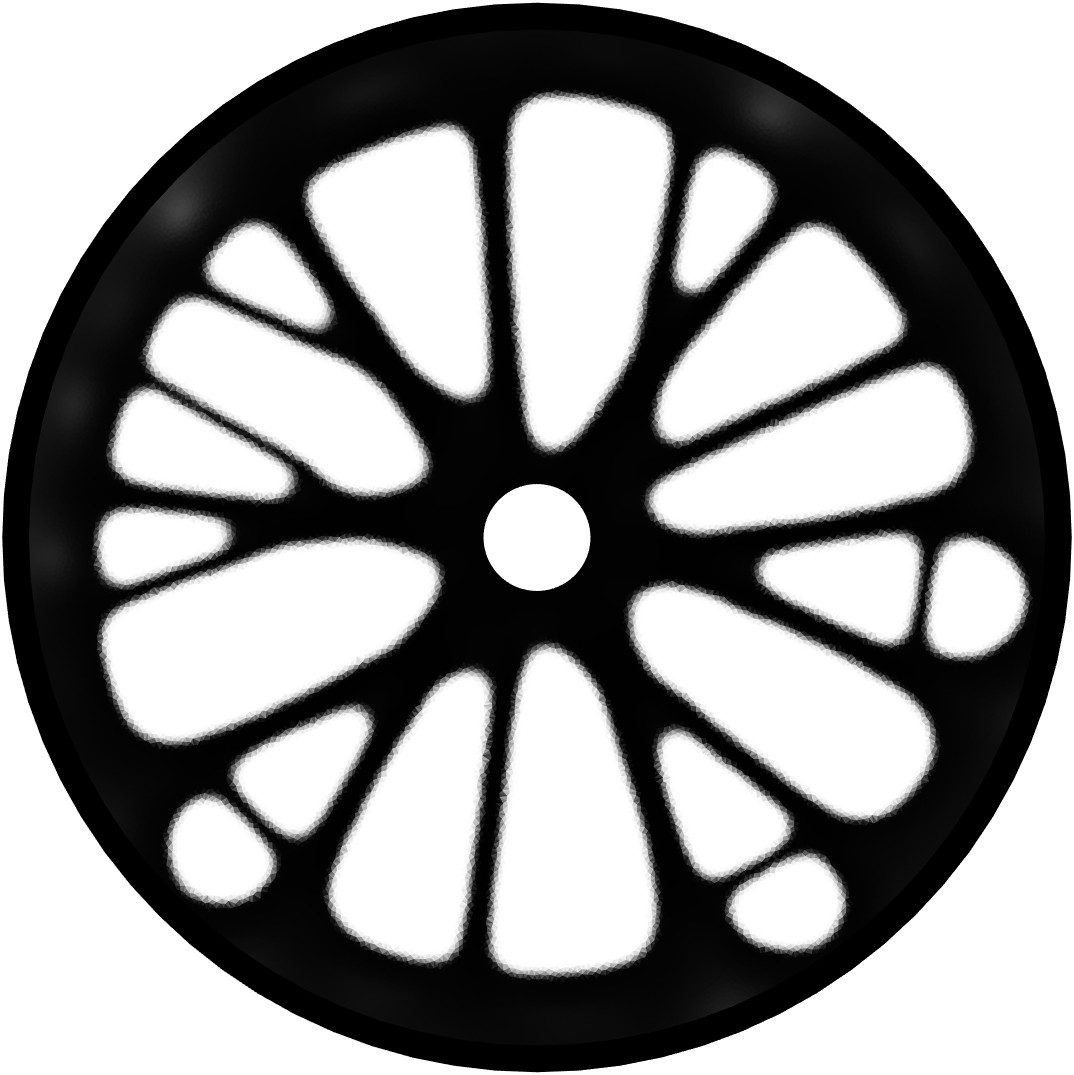} & \includegraphics[width=.135\textwidth,keepaspectratio,valign=c]{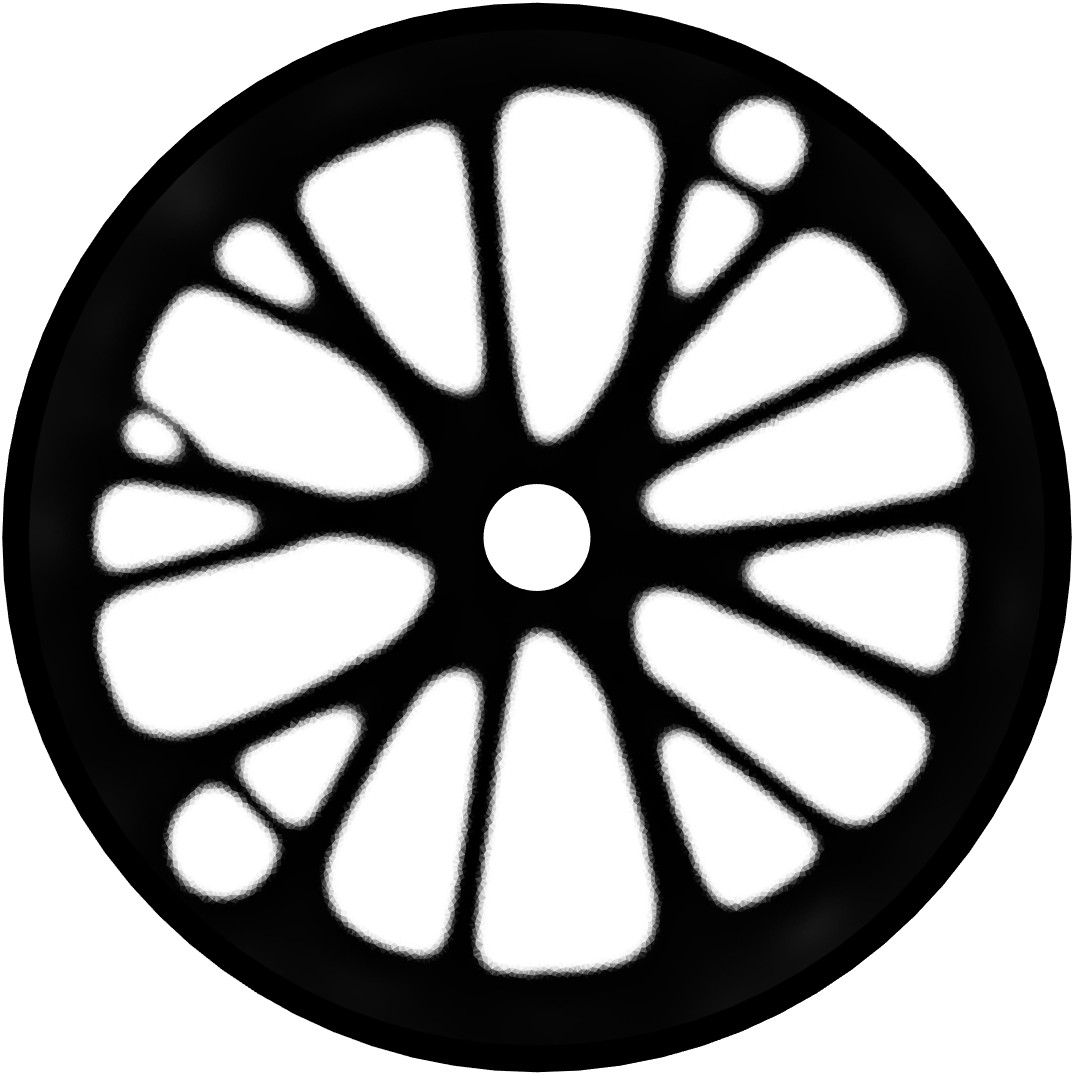} \\
    \end{tabular}
    \caption{Final designs for sMMA after 400 iterations with different batch sizes $\mathcal{B}$ (left to right) and move limits $\tau$ (top to bottom). Pseudo-densities ${\rho\in[0,1]}$ are depicted on a grayscale with $\rho=1$ and $\rho=0$ corresponding to black and white, respectively.}
    \label{tab:wheel_designs_mcmsa}
\end{table}
\begin{table}
    \centering
    \begin{tabular}{l|ccccc}
    & $\mathcal{B}=4$ & $\mathcal{B}=8$ & $\mathcal{B}=16$ & $\mathcal{B}=32$ & $\mathcal{B}=64$ \\[1em] \hline
        $\tau=1$ & \includegraphics[width=.135\textwidth,keepaspectratio,valign=c]{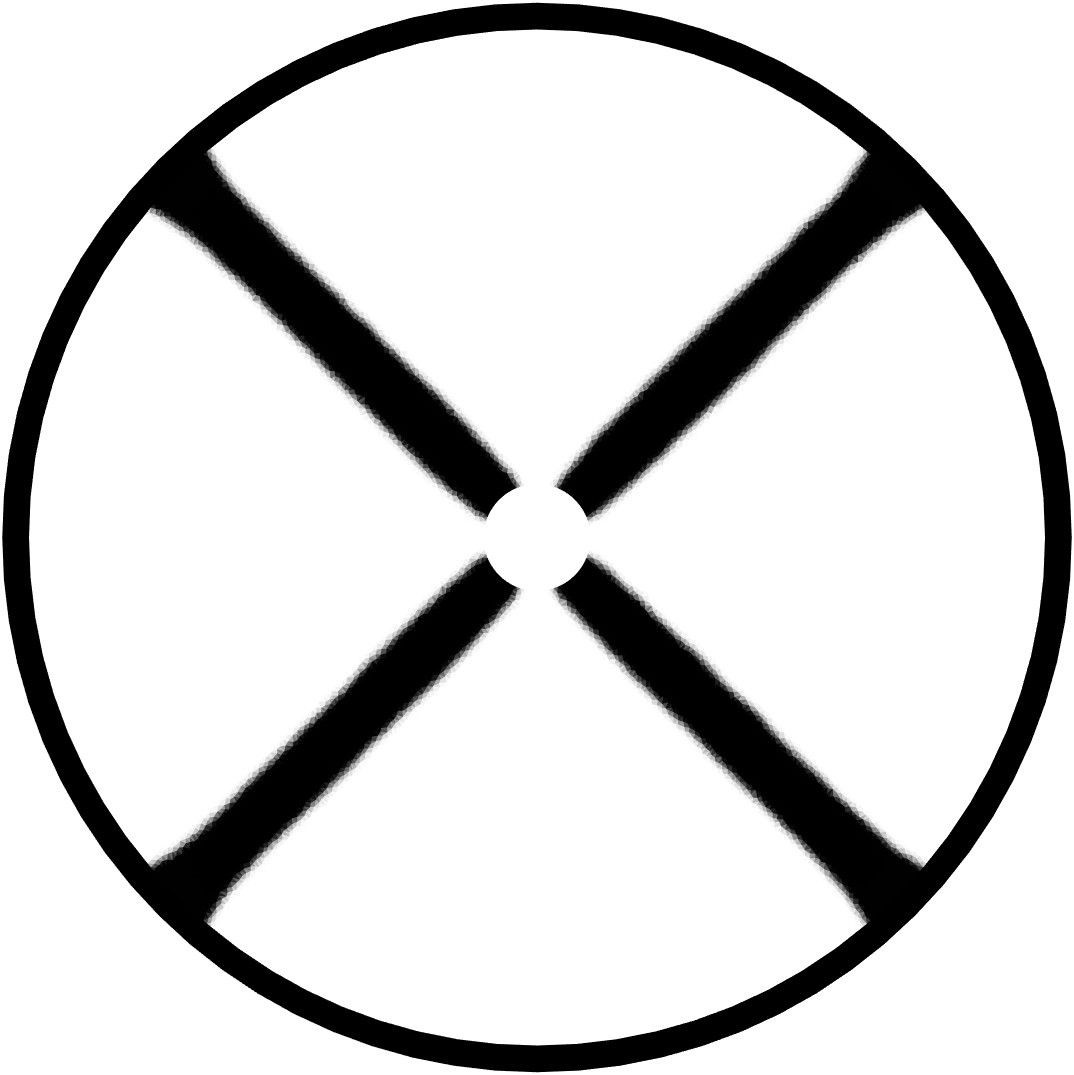} & \includegraphics[width=.135\textwidth,keepaspectratio,valign=c]{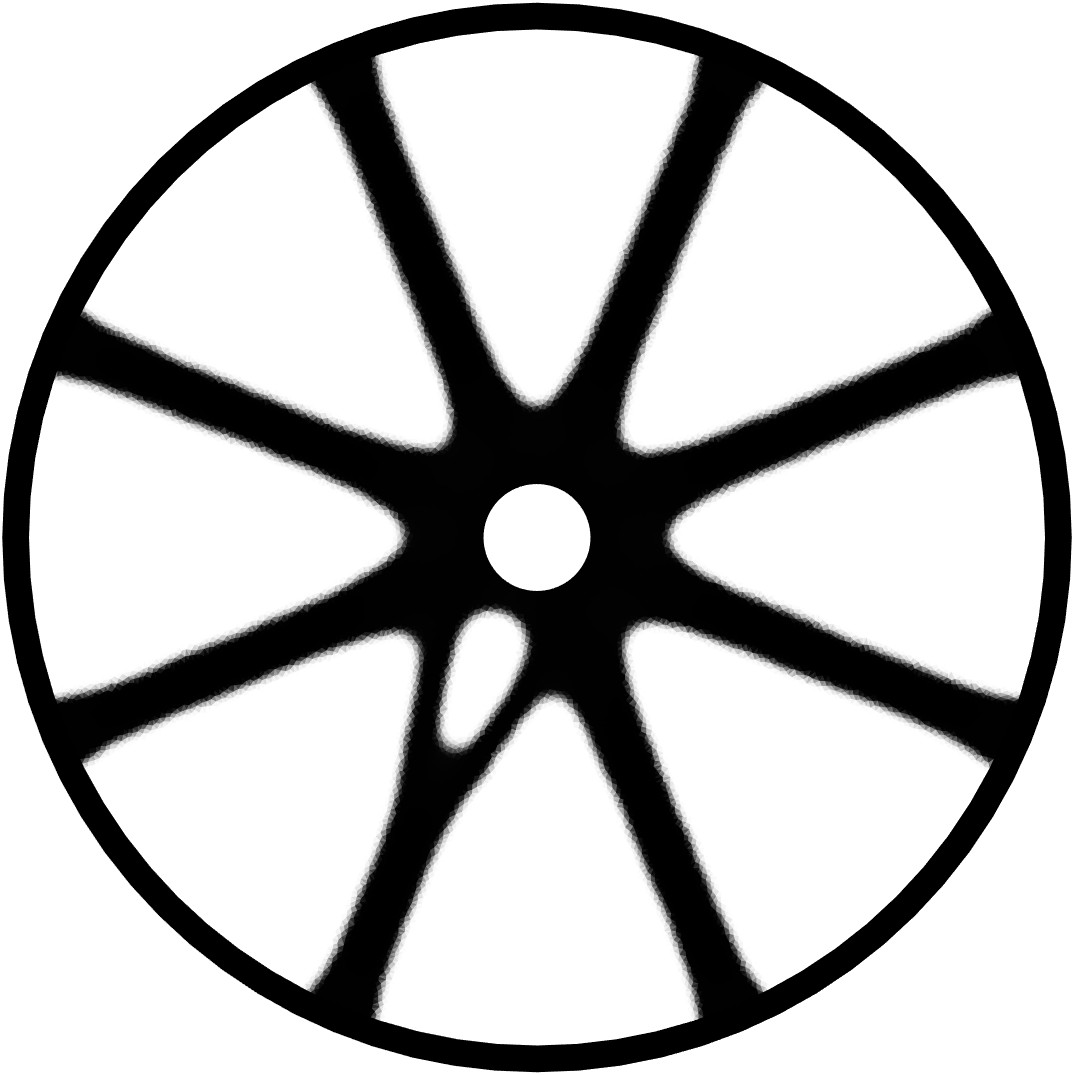} & \includegraphics[width=.135\textwidth,keepaspectratio,valign=c]{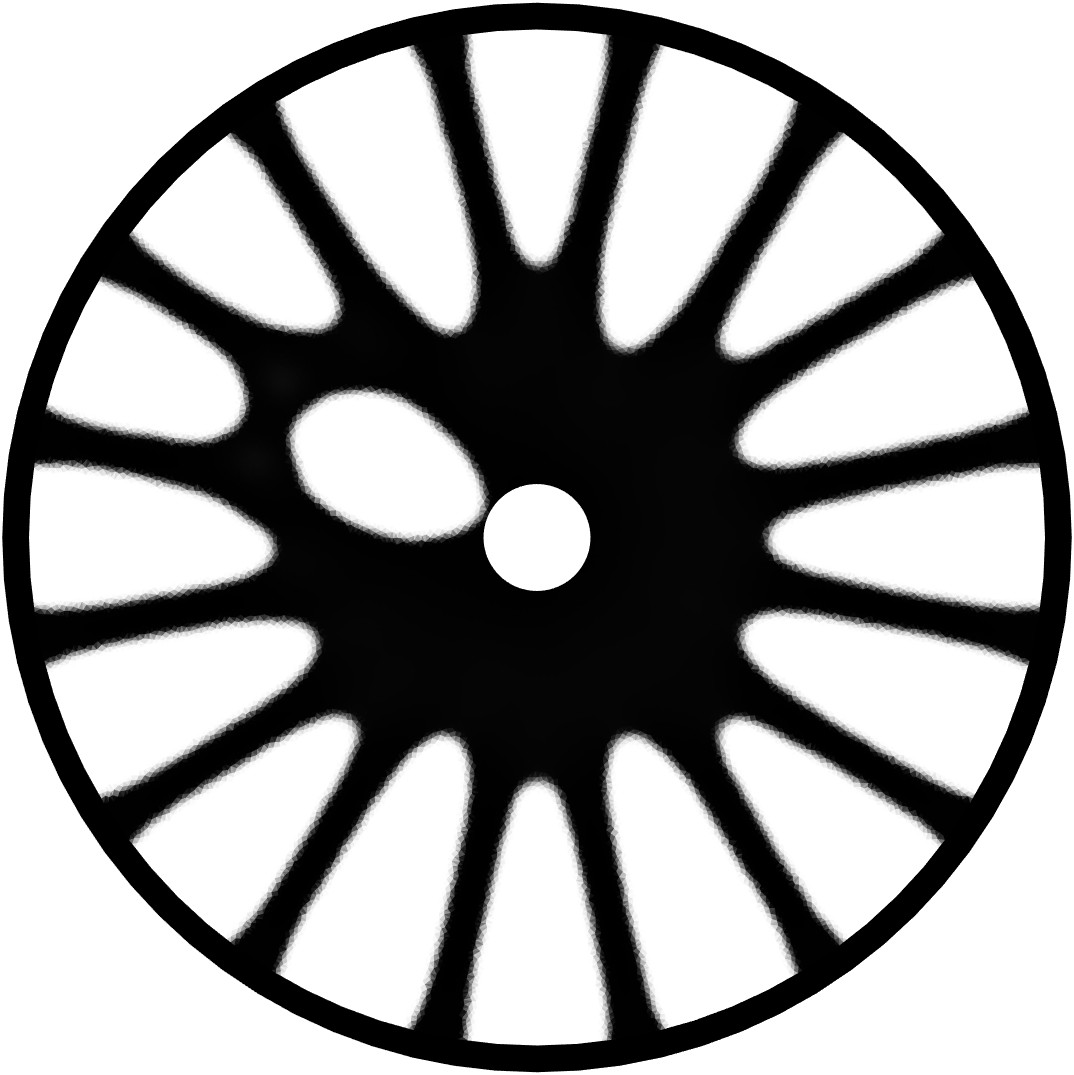} & 
         \includegraphics[width=.135\textwidth,keepaspectratio,valign=c]{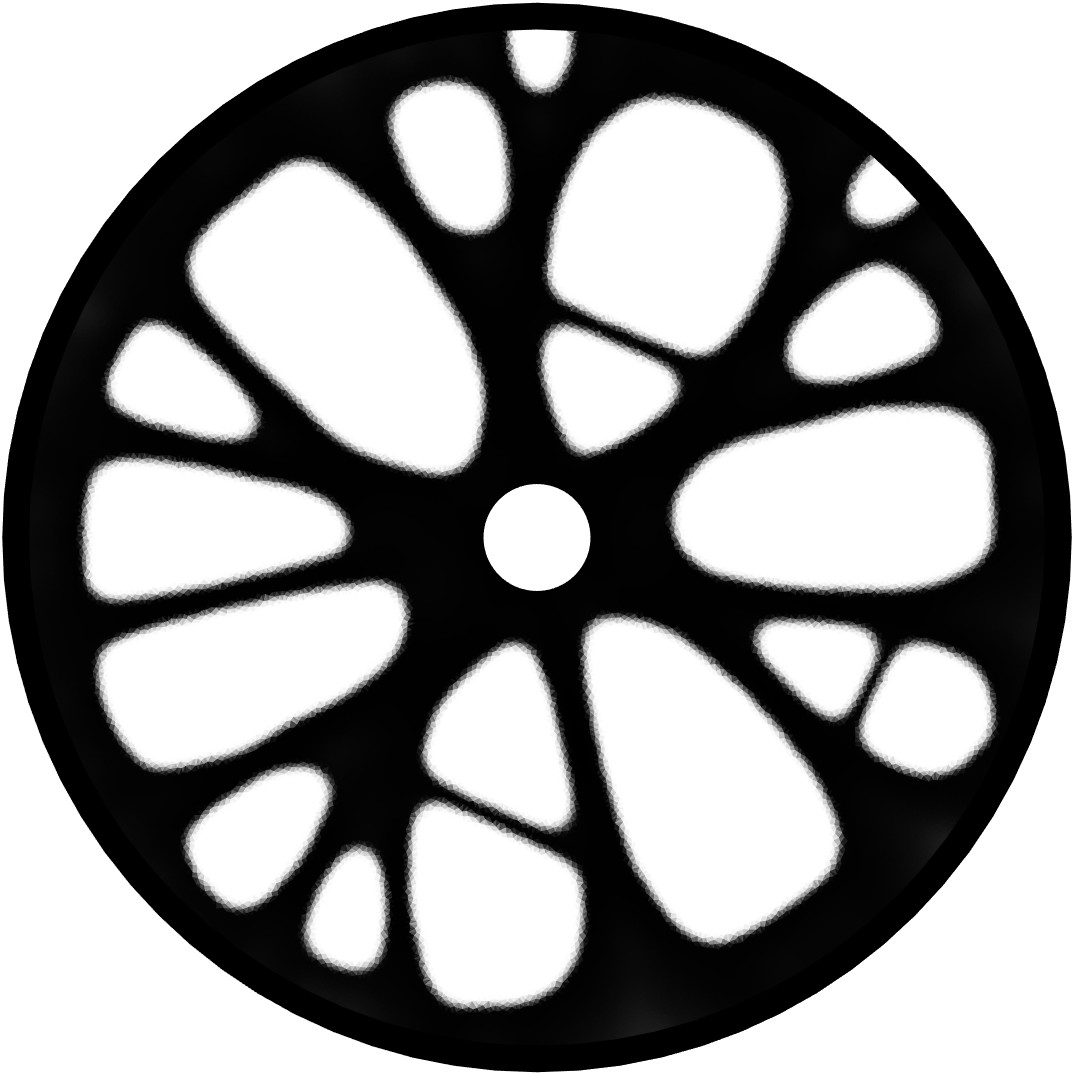} & \includegraphics[width=.135\textwidth,keepaspectratio,valign=c]{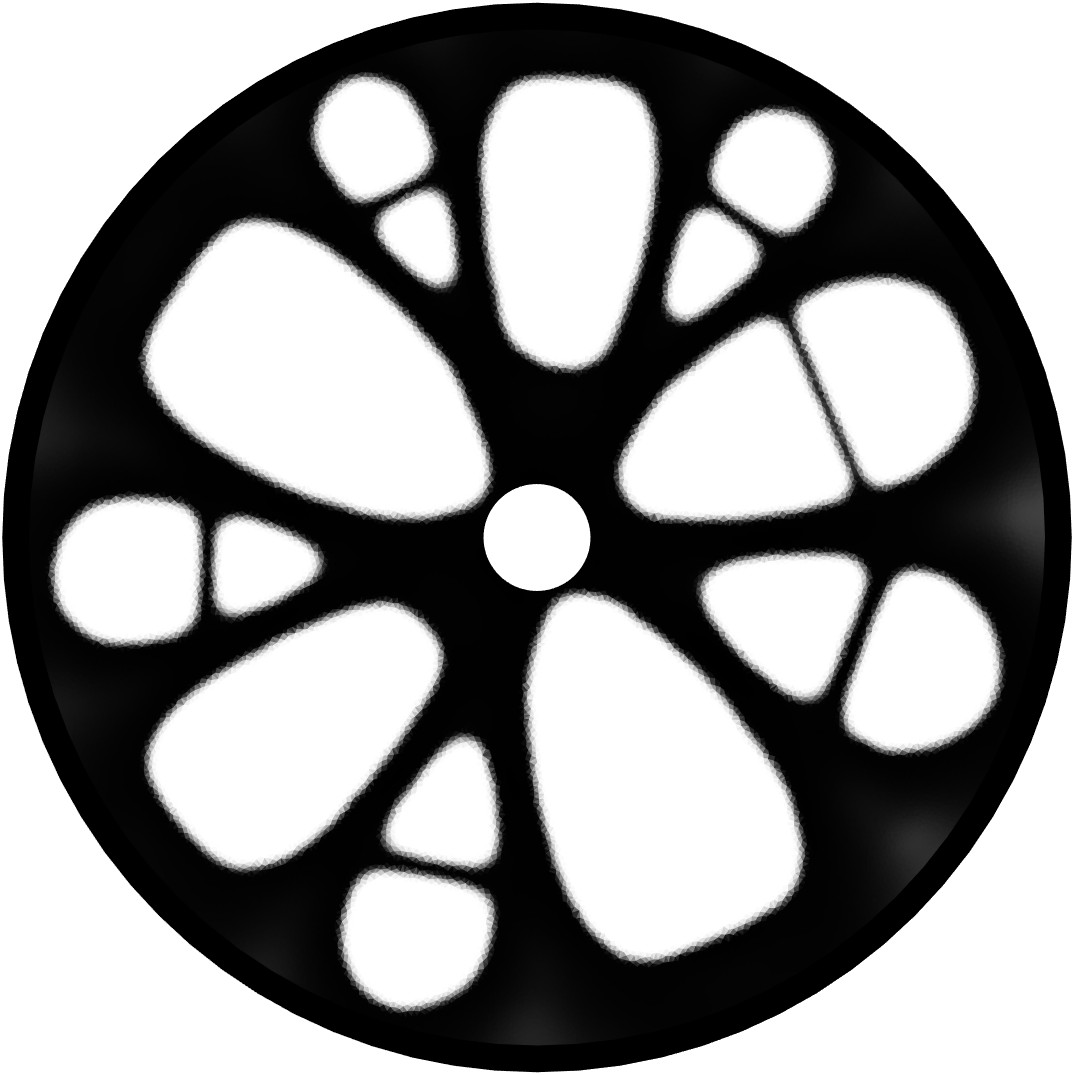}\\
        $\tau=\tfrac{3}{4}$ & \includegraphics[width=.135\textwidth,keepaspectratio,valign=c]{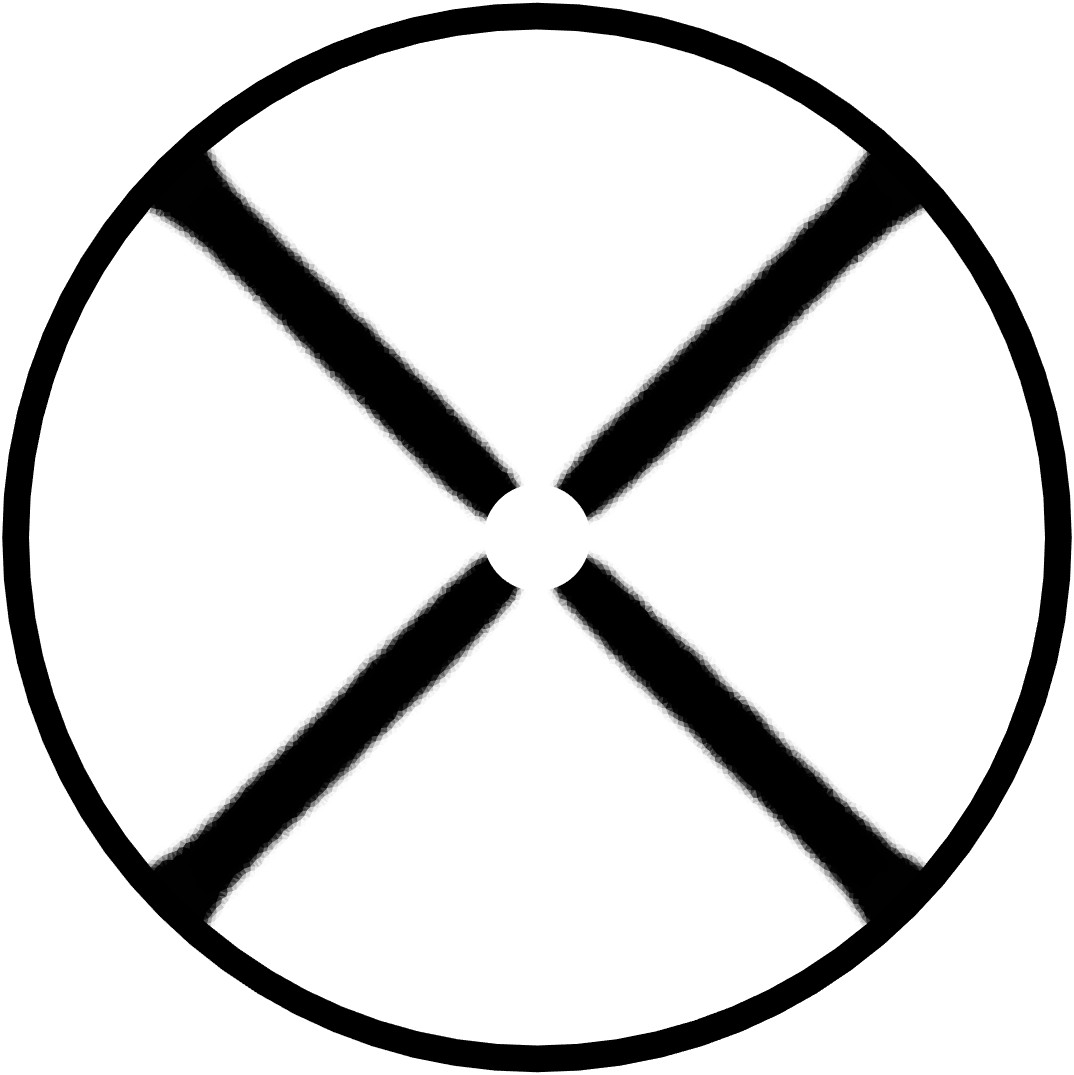} & \includegraphics[width=.135\textwidth,keepaspectratio,valign=c]{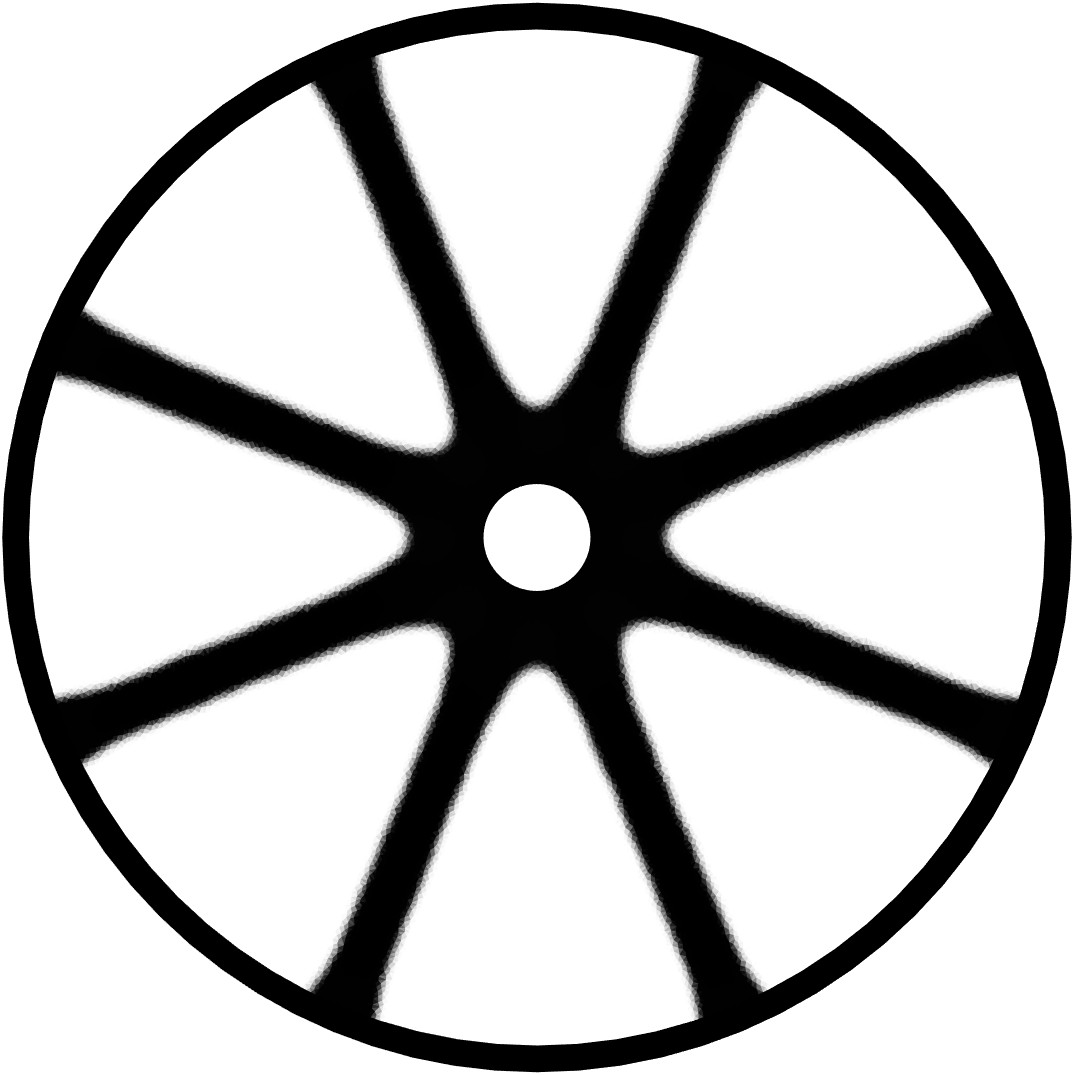} & \includegraphics[width=.135\textwidth,keepaspectratio,valign=c]{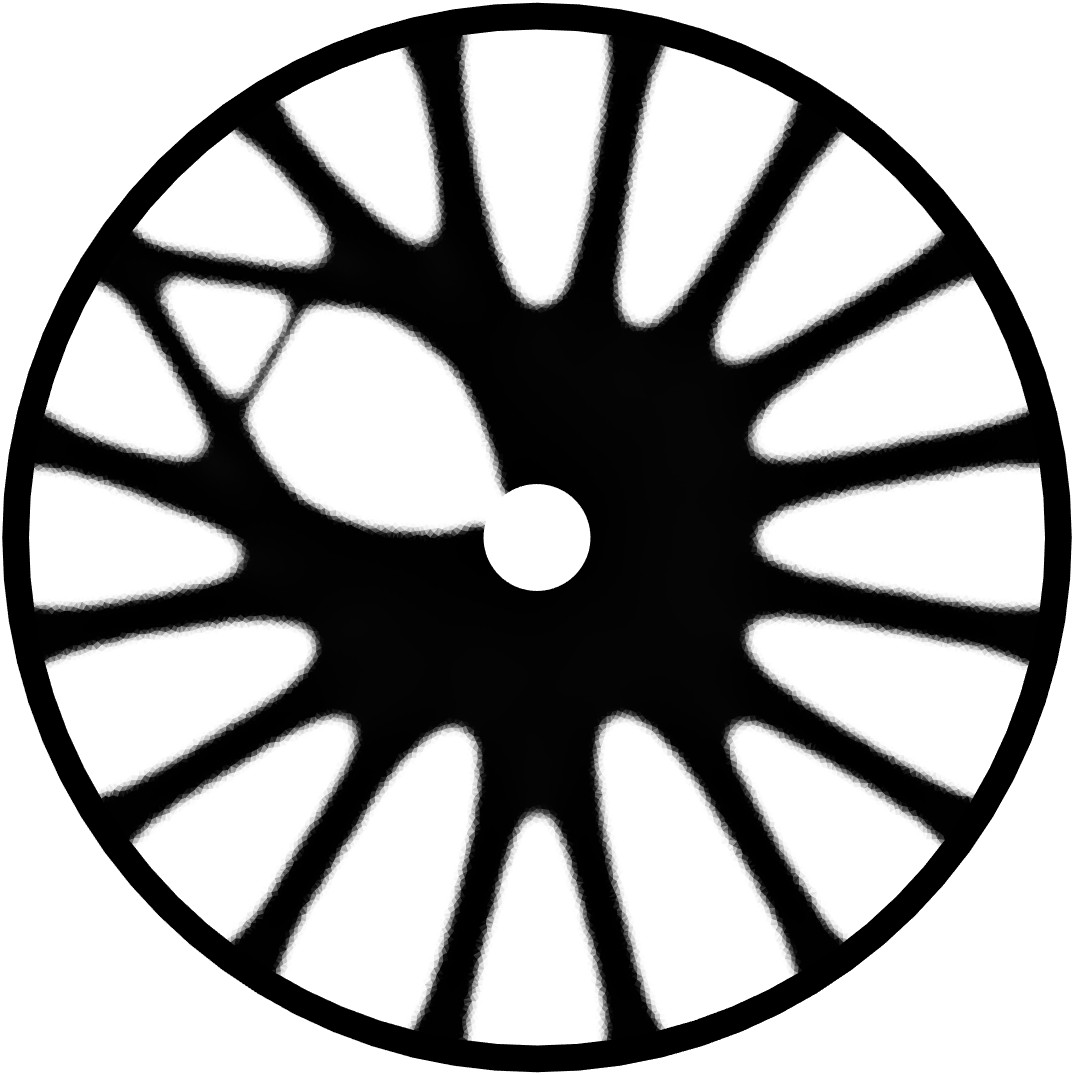} & 
         \includegraphics[width=.135\textwidth,keepaspectratio,valign=c]{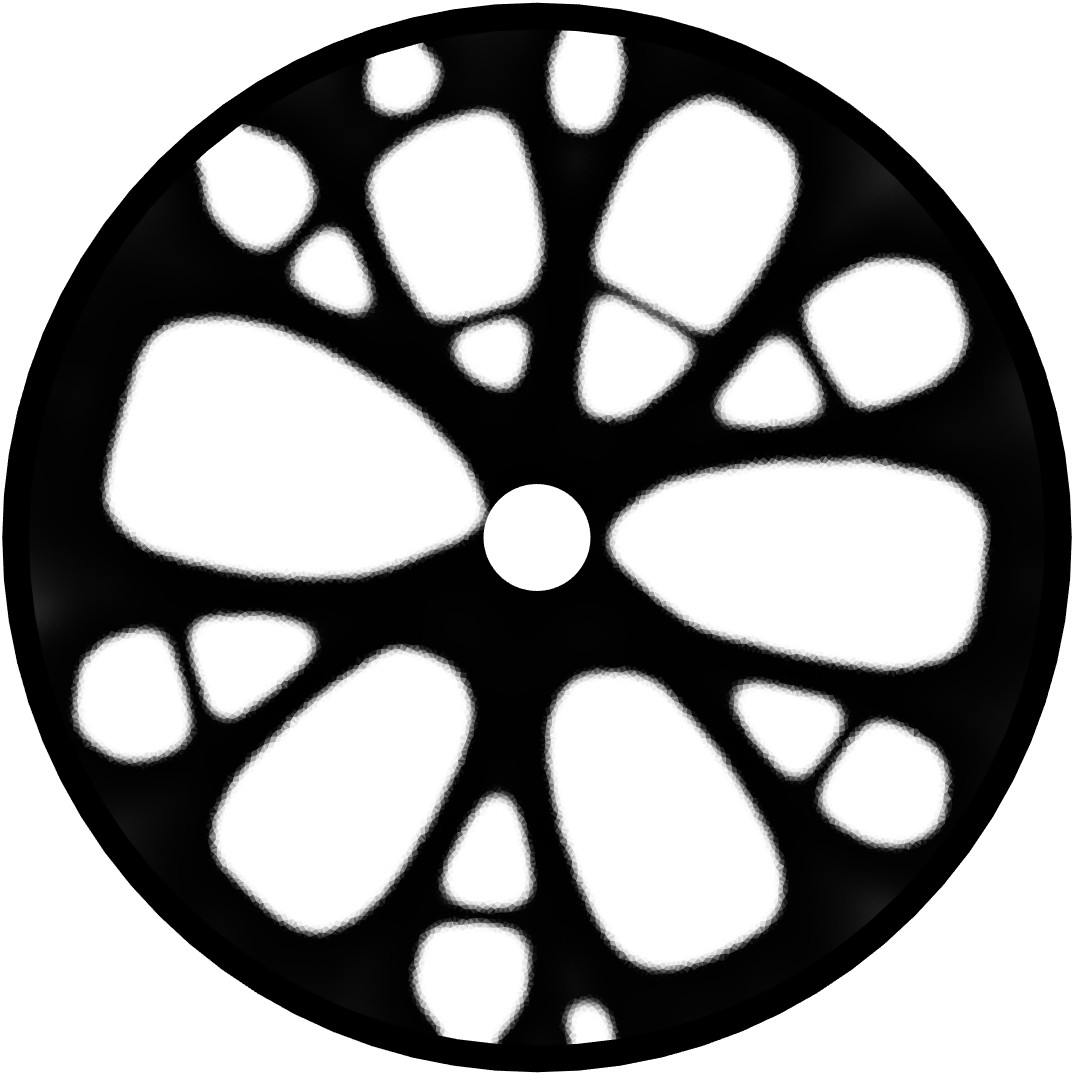} & \includegraphics[width=.135\textwidth,keepaspectratio,valign=c]{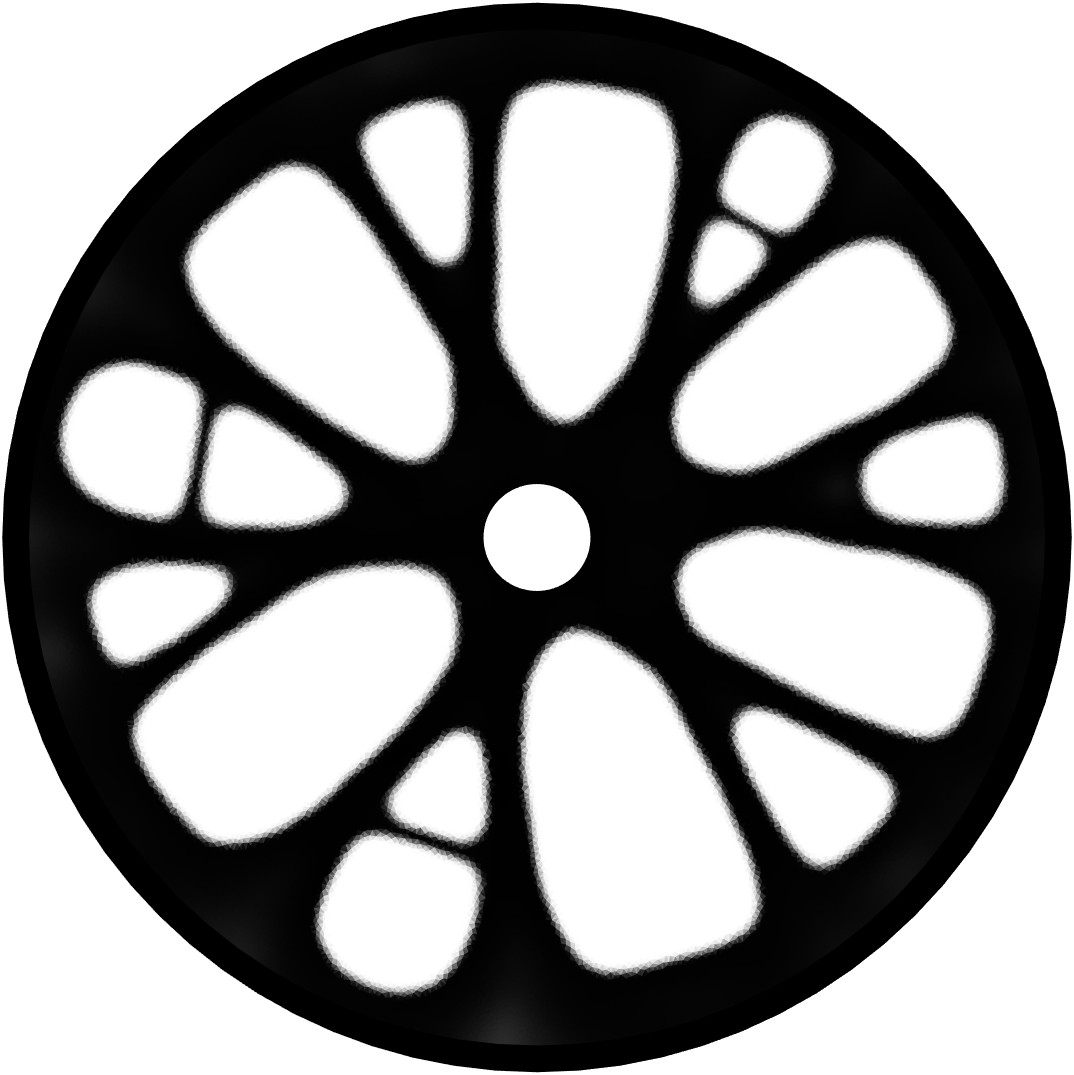} \\
         $\tau=\tfrac{1}{2}$ & \includegraphics[width=.135\textwidth,keepaspectratio,valign=c]{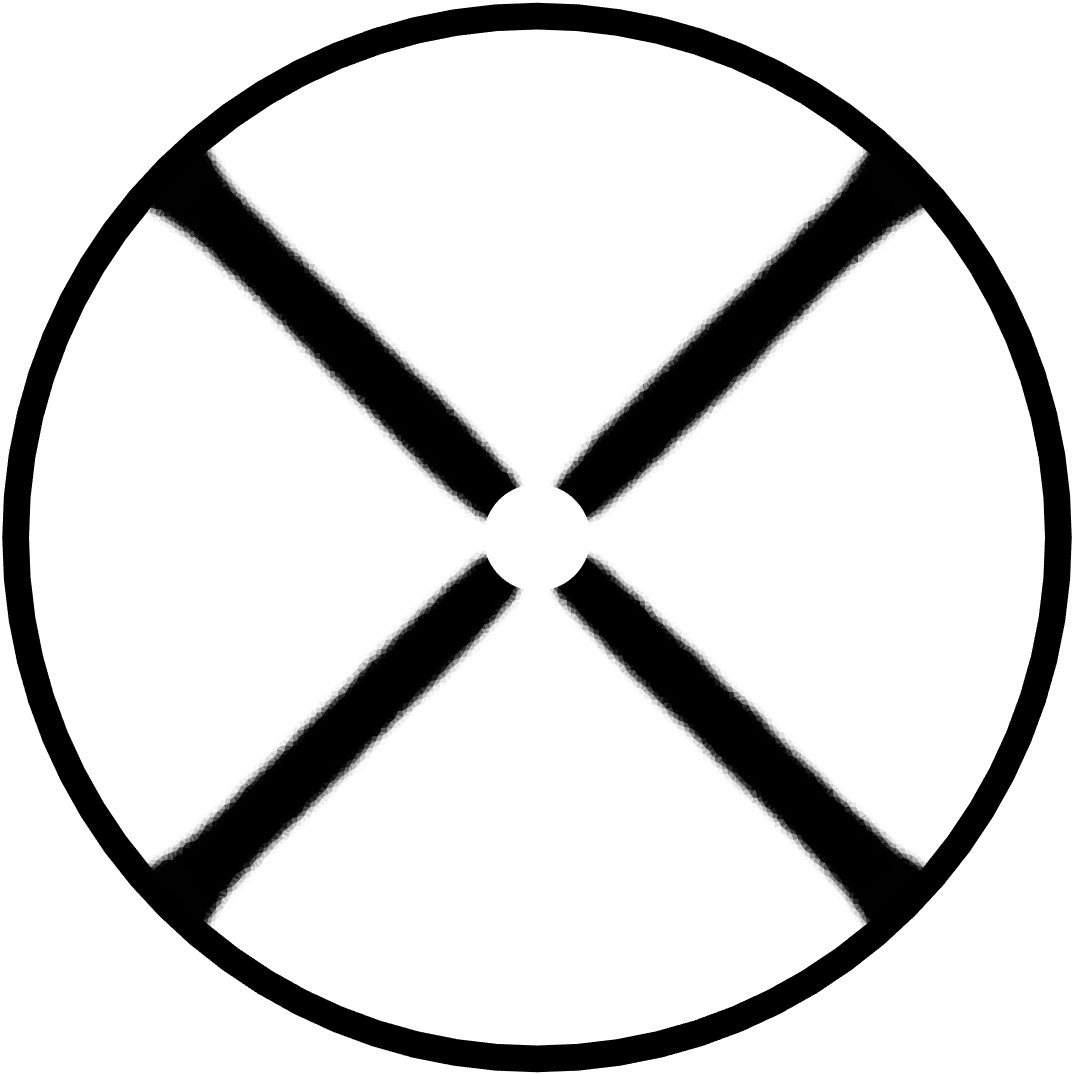} & \includegraphics[width=.135\textwidth,keepaspectratio,valign=c]{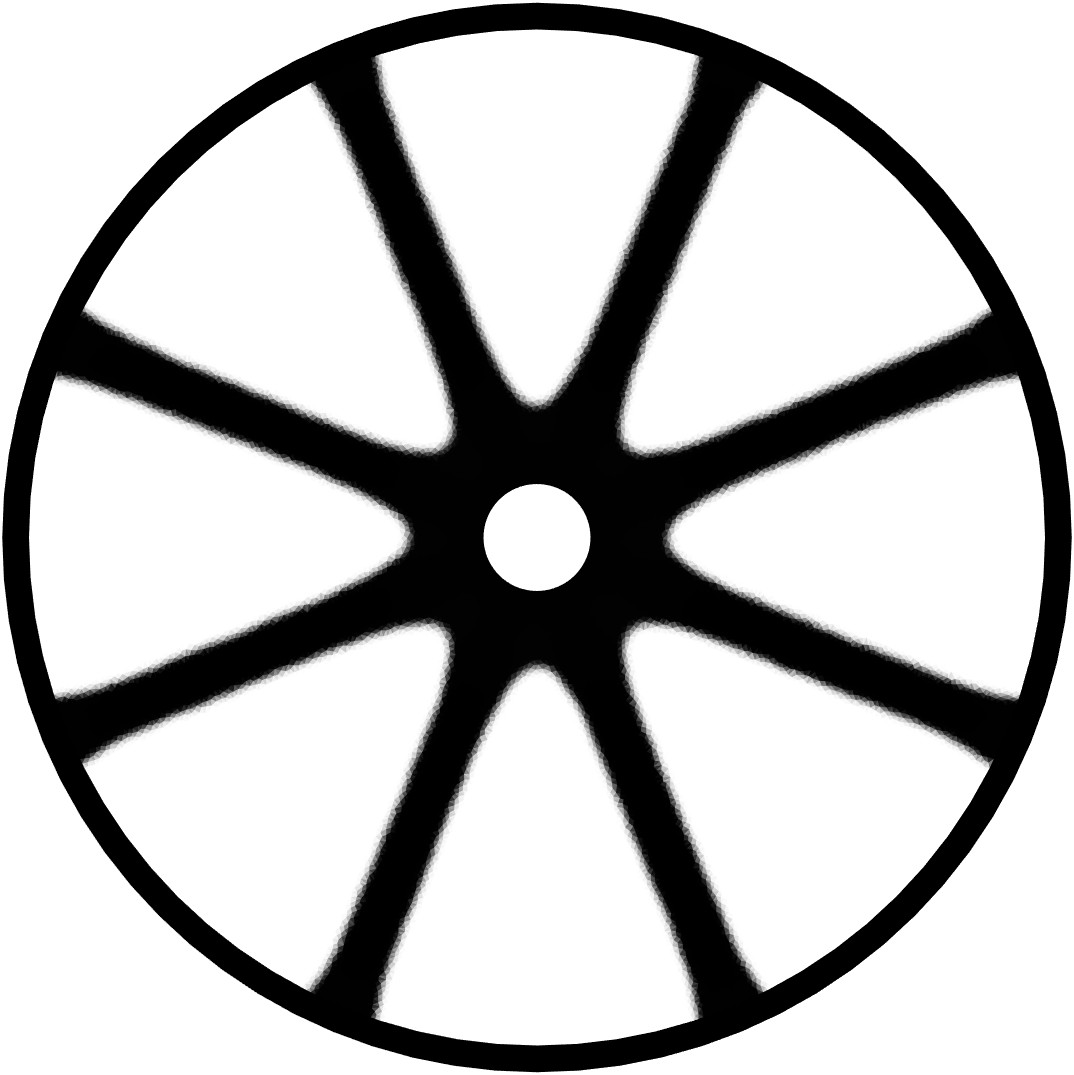} & \includegraphics[width=.135\textwidth,keepaspectratio,valign=c]{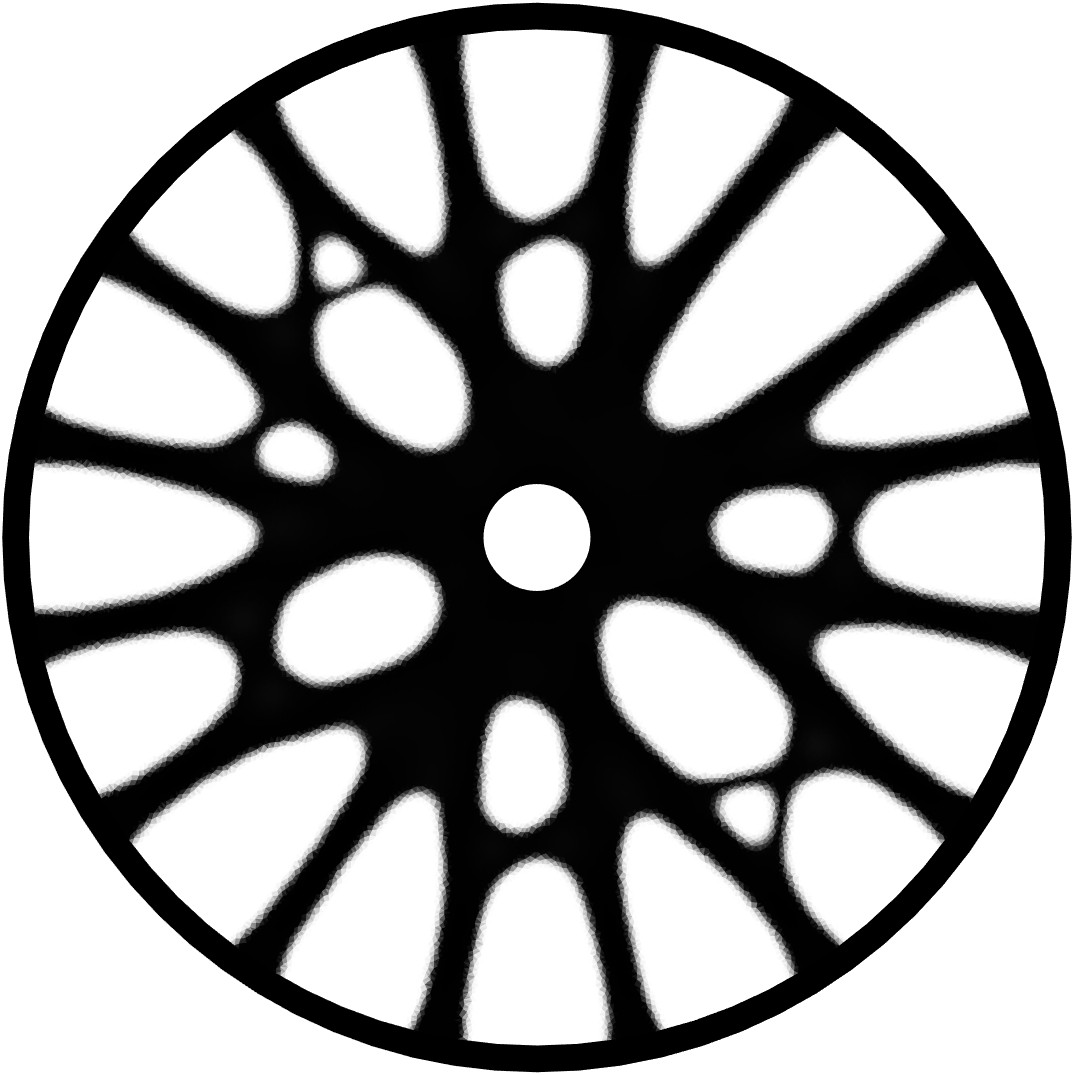} & 
         \includegraphics[width=.135\textwidth,keepaspectratio,valign=c]{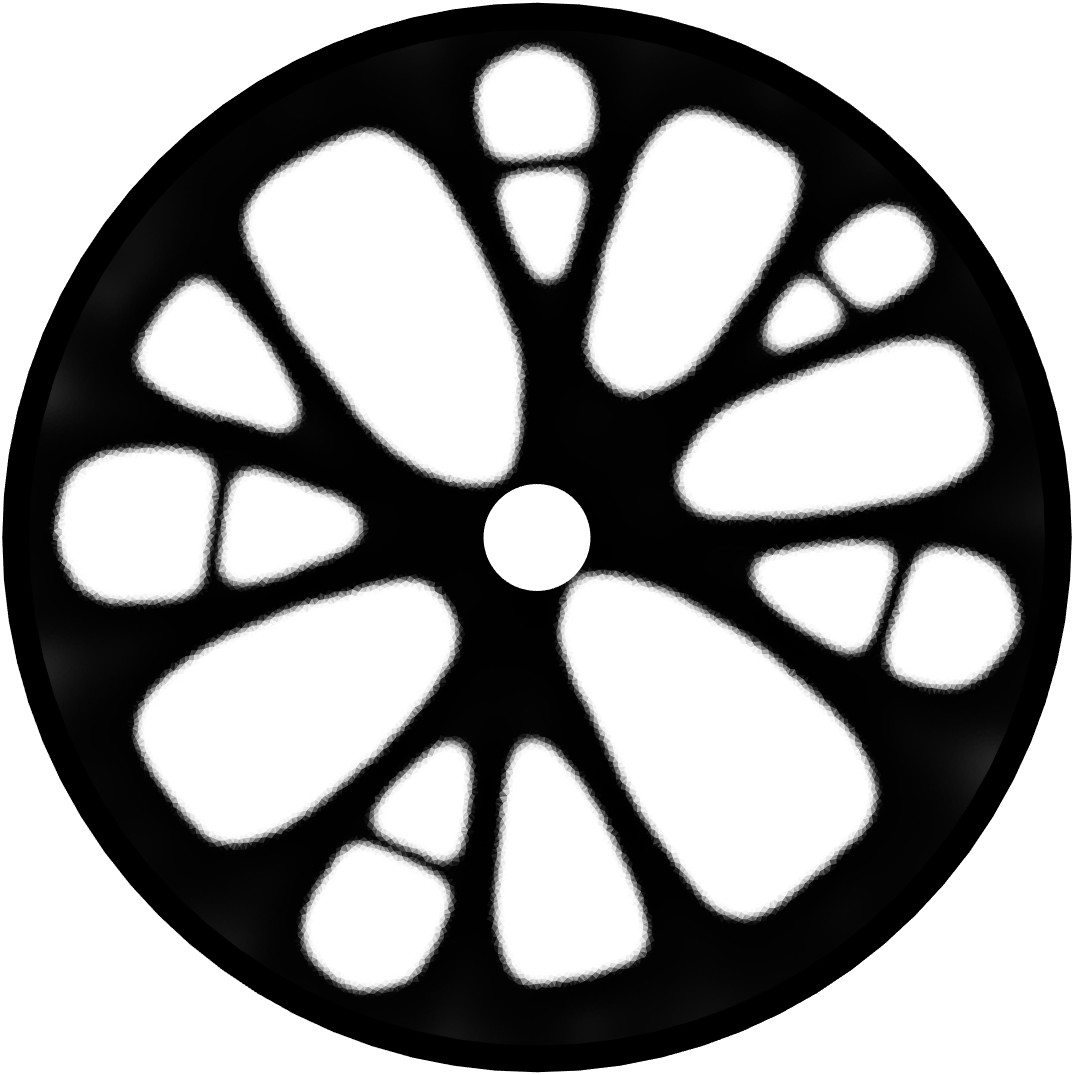} & \includegraphics[width=.135\textwidth,keepaspectratio,valign=c]{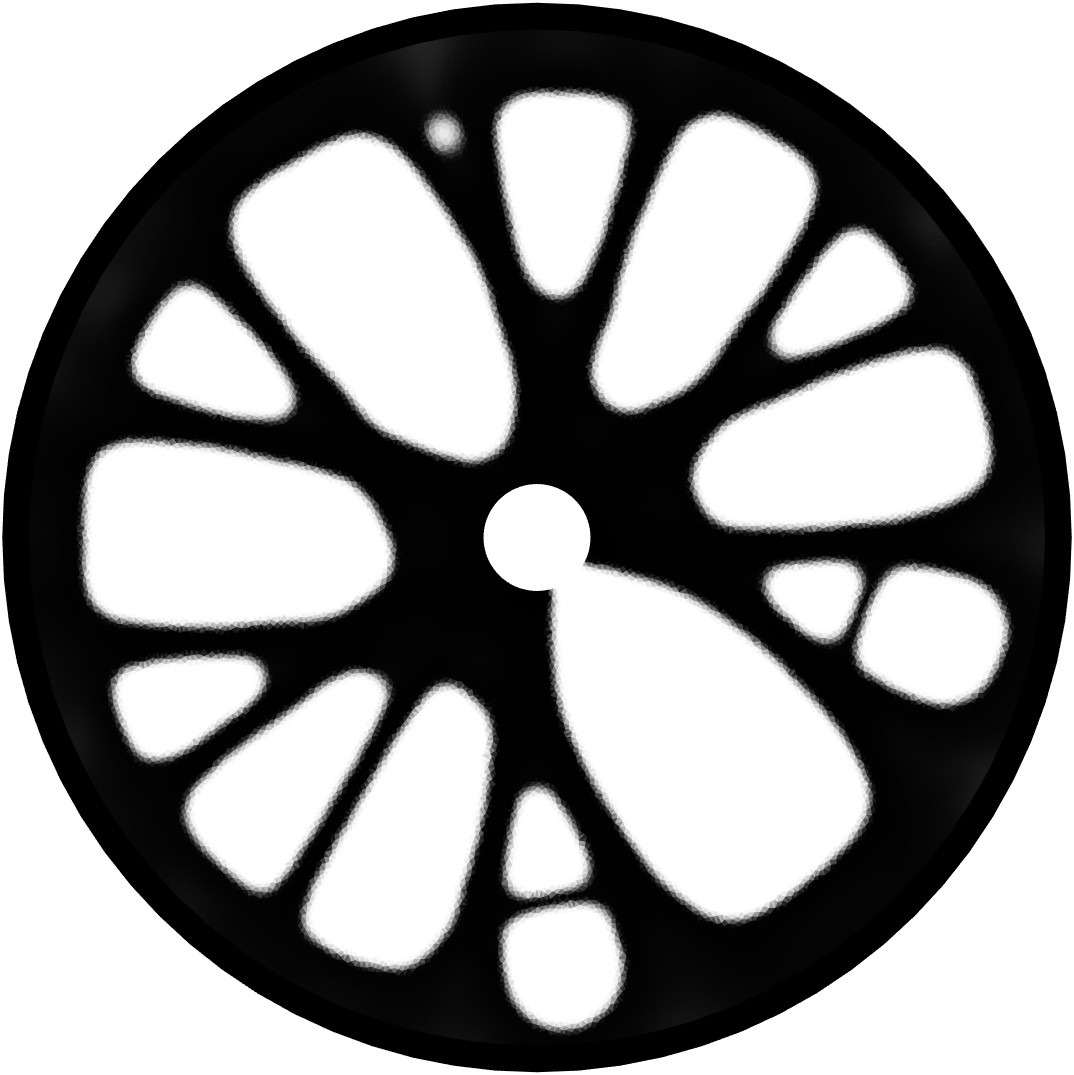}\\
         $\tau=\tfrac{1}{4}$ & \includegraphics[width=.135\textwidth,keepaspectratio,valign=c]{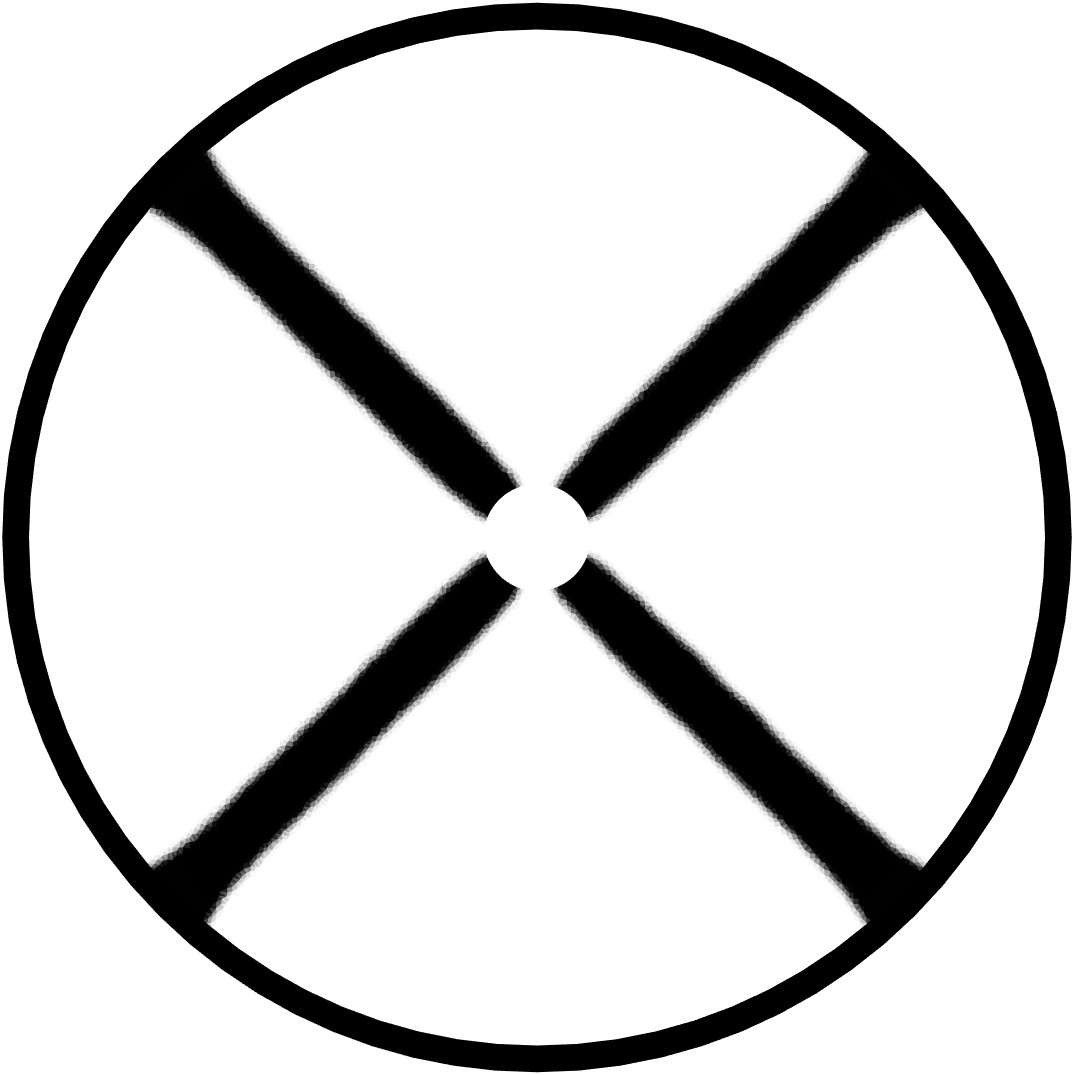} & \includegraphics[width=.135\textwidth,keepaspectratio,valign=c]{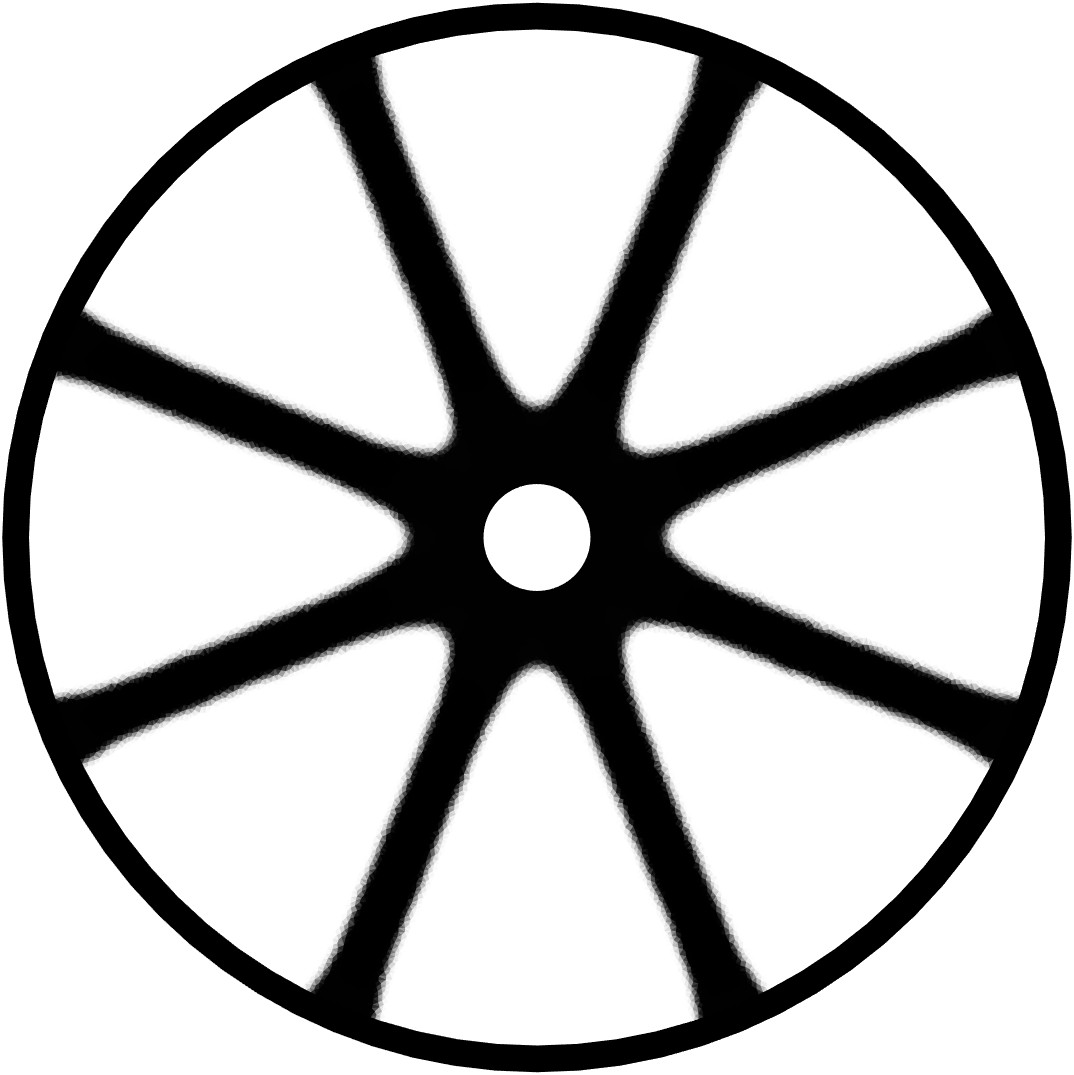} & \includegraphics[width=.135\textwidth,keepaspectratio,valign=c]{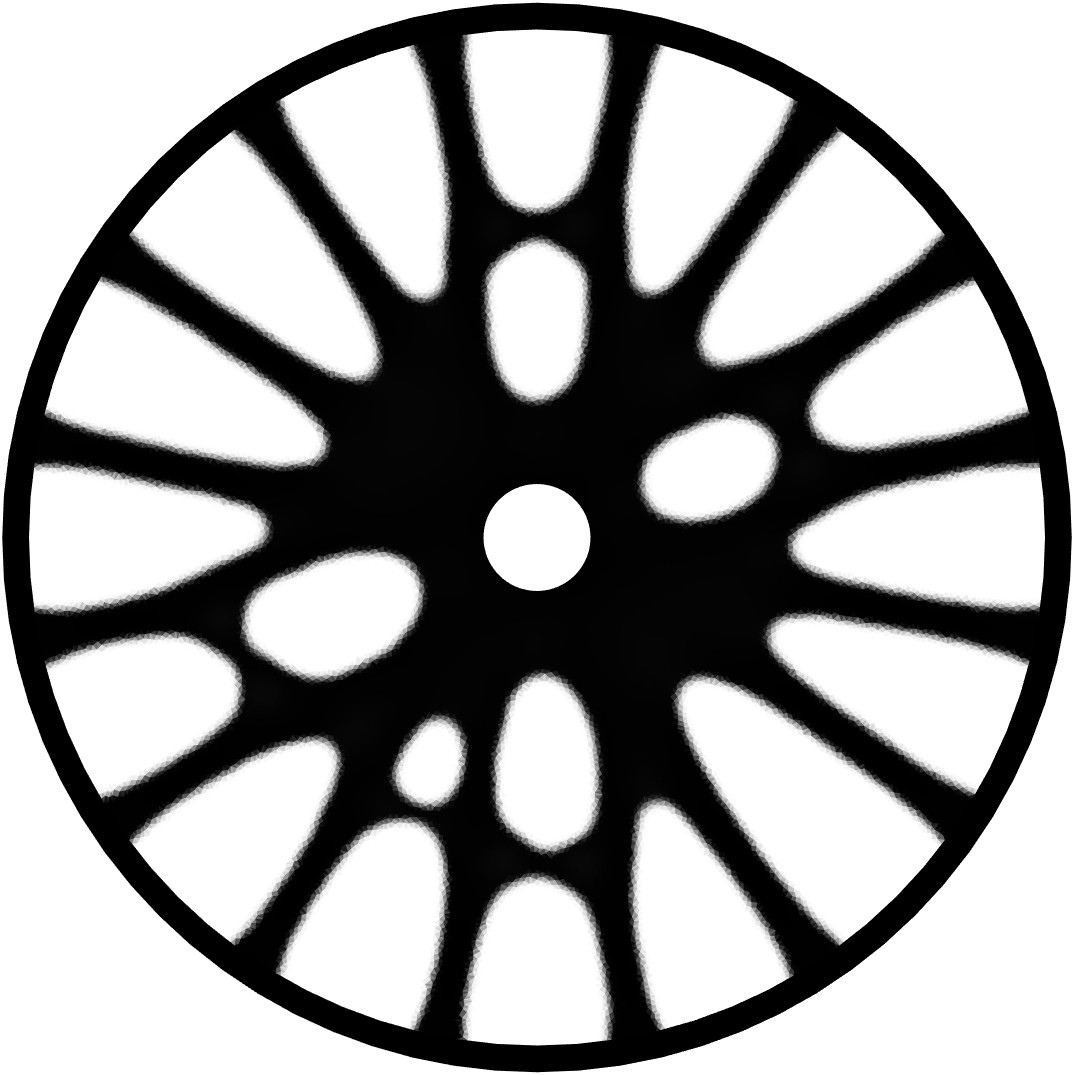} & 
         \includegraphics[width=.135\textwidth,keepaspectratio,valign=c]{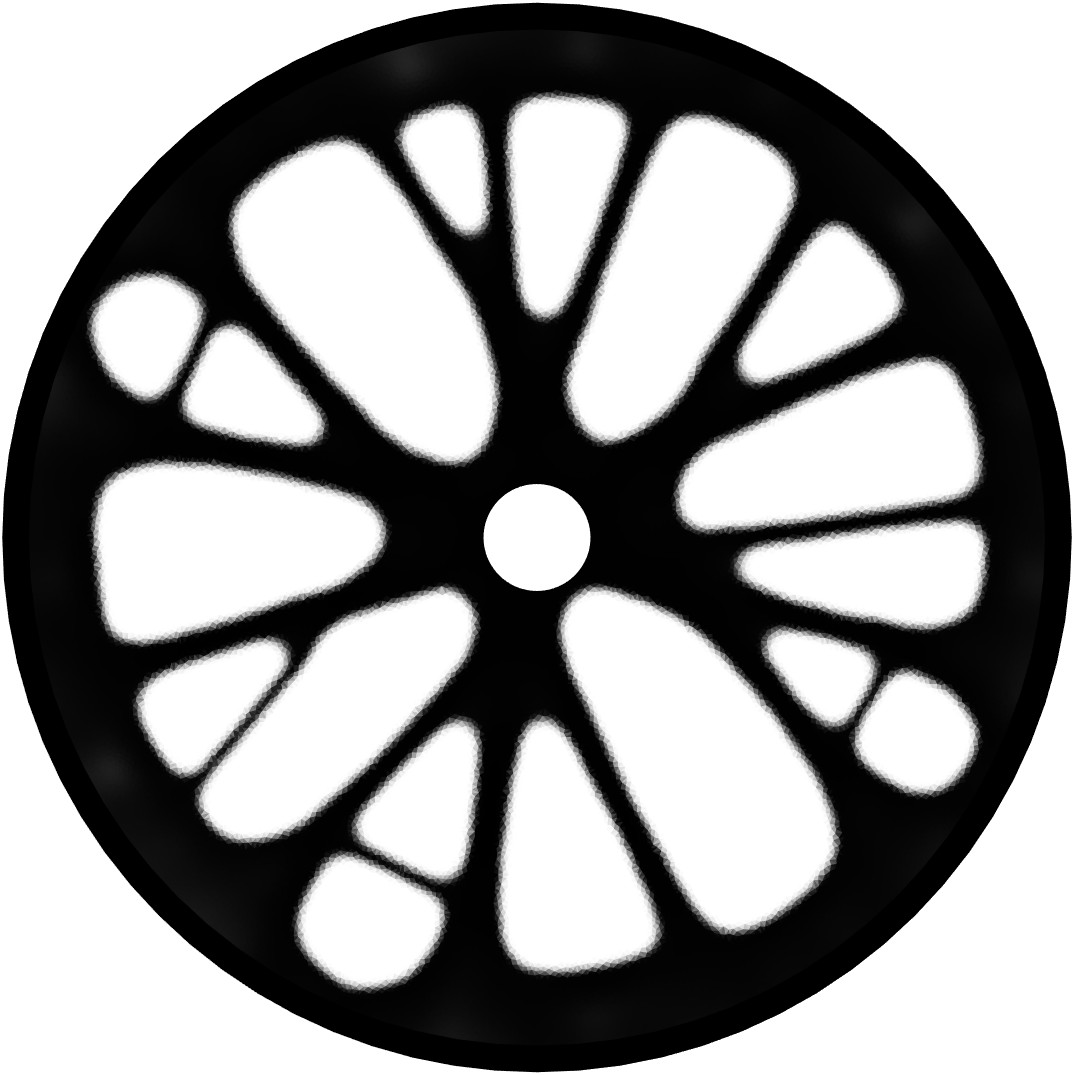} & \includegraphics[width=.135\textwidth,keepaspectratio,valign=c]{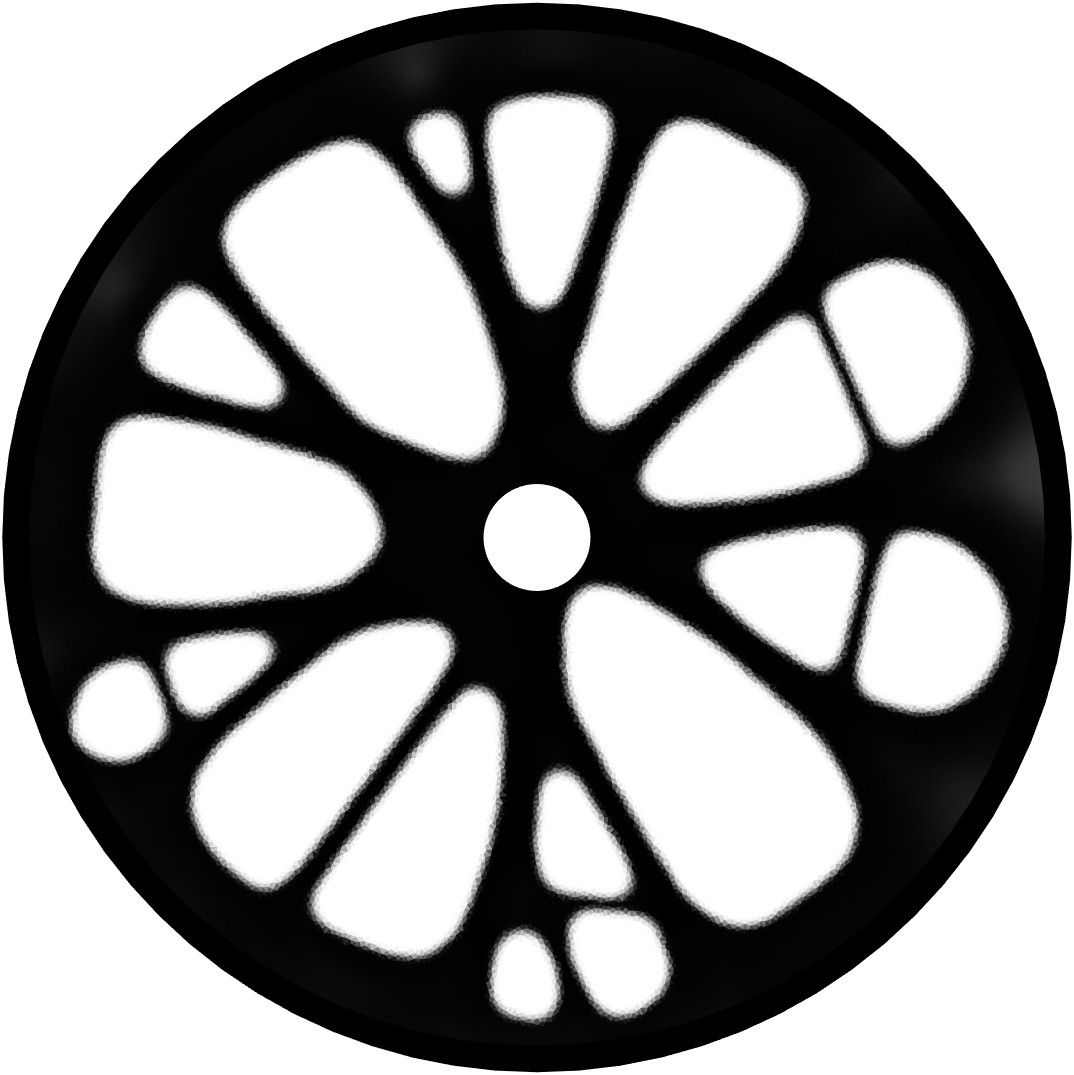} \\
    \end{tabular}
    \caption{Final designs for MMA after 400 iterations with different batch sizes $\mathcal{B}$ (left to right) and move limits $\tau$ (top to bottom). Pseudo-densities ${\rho\in[0,1]}$ are depicted on a grayscale with $\rho=1$ and $\rho=0$ corresponding to black and white, respectively. The beam-type structures for designs corresponding to ${\mathcal{B}\in\{4,8,16\}}$ is a result of the poor approximation to the full chance constraint integral, when using only $4,8,16$ integration points.}
    \label{tab:wheel_designs_mma}
\end{table}
%#############################################################################################
%################################### 2d Table ################################################
%#############################################################################################
%\FloatBarrier
\subsection{2d compliance with load and material uncertainty}
\begin{figure}
  \begin{minipage}[c]{0.5\textwidth}
    \input{2dtable_setup}
  \end{minipage}\hfill
  \begin{minipage}[c]{0.45\textwidth}
    \caption{Rectangular design domain $\mathscr{D}$ (light grey) with height $\ell$, width $2\ell$ and support at the bottom (dark grey). At the top, the structure is loaded by a distributed force (blue) with a range of possible angles (black dotted lines). Possible midpoints for the loaded area are contained in $\Omega$ (solid red line). Within $\mathscr{D}$, we model a ``hole" induced by the manufacturing process by significantly weakening the material parameters inside the blue circle. The midpoint of this circle can lie anywhere within $\Xi$ (green dashed rectangle). The illustration is true to scale with regard to the parameters selected for our numerical experiments.}
    \label{fig:2dtable_setup}
  \end{minipage}
\end{figure}
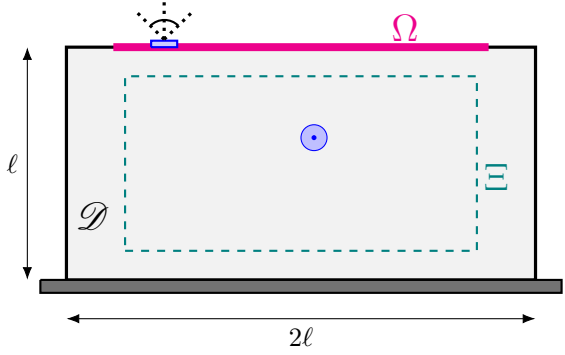

For our second example, we consider a rectangular design domain $\mathscr{D}$, firmly supported at the bottom and loaded at a random area from the top. Additionally, we model a possible manufacturing error by significantly weakening the material in a circular region within $\mathscr{D}$. Of course, the position of this region is assumed to be random. An illustration is given in~\Cref{fig:2dtable_setup}.

% \begin{figure}
%     \centering
%     \input{2dtable_setup}
%     \caption{Rectangular design domain $\mathscr{D}$ (light grey) with height $\ell$, width $2\ell$ and support at the bottom (dark grey). At the top, the structure is loaded by a distributed force (blue) with a range of possible angles (black dotted lines). Possible midpoints for the loaded area are contained in $\Omega$ (solid red line). Within $\mathscr{D}$, we model a ``hole" induced by the manufacturing process by significantly weakening the material parameters inside the blue circle. The midpoint of this circle can lie anywhere within $\Xi$ (green dashed rectangle). The illustration is true to scale with regard to the parameters selected for our numerical experiments.}
%     \label{fig:2dtable_setup}
% \end{figure}

To be precise, for $\ell>0$, we define ${\mathscr{D}:=[0,2\ell]\times[0,\ell]}$, ${\Omega:=\left[\tfrac{\ell}{5},\tfrac{4\ell}{5}\right]}$ as well as ${\Xi:= \left[ \tfrac{\ell}{4},\tfrac{7\ell}{4} \right]\times\left[ \tfrac{\ell}{8},\tfrac{7\ell}{8} \right]}$. Now, for ${\omega\in\Omega}$, the force is applied at a region of width $\tfrac{\ell}{9}$, centered around $\omega$. The intensity of the force is modeled by the smooth bump function
\begin{equation*}
    f_\omega(t):= \begin{cases}
        \exp\left( -\frac{\frac{1}{10}}{1-\left(\frac{18}{\ell}\right)^2(t-\omega)^2} \right), & \vert t-\omega\vert < \tfrac{\ell}{18}, \\
        0 & \text{else.}
    \end{cases}
\end{equation*}

Moreover, the force acts downward from an angle ${\alpha\in\mathcal{A}:=\left[\tfrac{\pi}{4},\tfrac{3\pi}{4}\right]}$ and the \emph{expected} compliance (w.r.t. $\alpha$) will be used in our calculations. Given ${\xi\in\Xi}$, the material is weakened in the circular domain ${\left\{x\in\mathscr{D}\,:\, \Vert x-\xi\Vert_2 < \tfrac{\ell}{18} \right\}\subset\mathscr{D}}$. This is done by reducing the stiffness parameter according to the two-dimensional bump function
\begin{equation*}
    g_\xi(x):=\begin{cases}
        \frac{99}{100}\cdot\exp\left( -\frac{\Vert x-\xi\Vert_2^2}{r^2(r^2-\Vert x-\xi\Vert_2^2)} \right), & \Vert x-\xi\Vert_2 < \tfrac{\ell}{18} \\
        0, & \text{else.}
    \end{cases}
\end{equation*}
By construction, the function is much steeper than $f_\omega$ and scaled such that the parameter is lowered by 99\% at the center.

For the regularization parameters, we fixed ${a_1=35}$, ${a_2=\tfrac{1}{20}}$ and ${a_3=5}$. The design domain was discretized into ${360\times180}$ square finite elements and we again chose 1 and $10^{-4}$ as stiffness parameters for material and void. The SIMP parameter was kept fixed at 5 during the optimization process. Both sMMA and MMA were initialized with ${\rho_1\equiv 0.65}$ and performed $5,000$ iterations. The move limits were chosen as ${\tau=1}$ and cut in half each $1,000$ iterations. Lastly, we picked ${p=0.05}$ and $\alpha$, $\xi$ as well as $\omega$ were assumed to follow uniform probability distributions.

The resulting optimization problem is of the following structure:
\begin{align*}
    \min_{\rho\in\R^d}\quad & \rvol(\rho),\\
    \text{s.t.}\quad & \frac{1}{\vert \Xi\vert\cdot\vert\Omega\vert}\int_\Xi\int_\Omega h_{a_1,a_2,a_3}\left(\frac{1}{\vert\mathcal{A}\vert}\int_\mathcal{A}\mathbf{F}(\alpha,\omega)^\top\mathbf{U}_{\alpha,\xi,\omega}\dd{\alpha}-\cm\right)\dd{\omega}\dd{\xi}\le p,\\
        &\mathbf{K}(\rho,\xi)\mathbf{U}_{\alpha,\xi,\omega}=\mathbf{F}(\alpha,\omega),\\
        &0\le \rho\le 1,   
\end{align*}
In this setup, uncertainty enters the state equation at two points: the right hand side $\mathbf{F}$ and the system matrix $\mathbf{K}$. As discussed for the previous example, uncertainty in $\mathbf{F}$ is much easier to handle than uncertainty in $\mathbf{K}$, since the moderate design dimension still allows for a simultaneous solution for multiple right hand sides (decomposing the system matrix requires ${\sim 125}$GB of memory). Thus, to reduce the numerical cost as much as possible, we first deal with the integrals over $\mathcal{A}$ and $\Omega$. 

First, note that ${\mathbf{F}^\top\mathbf{U}}$ is linear in $\mathbf{F}$. Additionally, for fixed ${\omega\in\Omega}$, the angle $\alpha$ has no impact on which finite elements are loaded, only on how large the load is. Thus, fixing ${\omega\in\Omega}$, for every ${\alpha\in\mathcal{A}}$, we can reformulate ${\mathbf{F}(\alpha,\omega)}$ as a linear combination of the corresponding forces acing in $x$ and $y$ direction:
\begin{equation*}
    \mathbf{F}(\alpha,\omega) = c_x(\alpha)\mathbf{F}_x(\omega) + c_y(\alpha)\mathbf{F}_y(\omega).
\end{equation*}
As a result, for fixed ${\omega\in\Omega}$, ${\xi\in\Xi}$, we have
\begin{align*}
    \int_\mathcal{A}\mathbf{F}(\alpha)^\top\mathbf{U}\dd{\alpha} &= \int_\mathcal{A}\mathbf{F}(\alpha)^\top\mathbf{K}^{-1}\mathbf{F}(\alpha)\dd{\alpha} \\
    &= \mathbf{F}_x^\top\mathbf{K}^{-1}\mathbf{F}_x\int_{\mathcal{A}}c_x(\alpha)^2\dd{\alpha} + \mathbf{F}_y^\top\mathbf{K}^{-1}\mathbf{F}_y\int_{\mathcal{A}}c_y(\alpha)^2\dd{\alpha} \\
    &\qquad + \left(\mathbf{F}_x^\top\mathbf{K}^{-1}\mathbf{F}_y + \mathbf{F}_y^\top\mathbf{K}^{-1}\mathbf{F}_x\right)\int_{\mathcal{A}}c_x(\alpha)c_y(\alpha)\dd{\alpha}.
\end{align*}
Since $c_x$ and $c_x$ are simple trigonometric functions, the appearing integrals can be calculated analytically. Consequentially, given ${\omega\in\Omega}$ and ${\xi\in\Xi}$, we can calculate the full integral over angles
\begin{equation*}
    \frac{1}{\vert\mathcal{A}\vert}\int_\mathcal{A}\mathbf{F}(\alpha,\omega)^\top\mathbf{U}_{\alpha,\xi,\omega}\dd{\alpha}
\end{equation*}
by solving the state equation for only two right hand sides ($\mathbf{F}_x$ and $\mathbf{F}_y$).

Now, since integrating over $\mathcal{A}$ poses no problem, the integral over $\Omega$ can be carried out using a simple trapezoidal rule with high enough resolution. For our numerical experiments, we chose 32 equidistant integration points. All in all, this means that, given ${\xi\in\Xi}$, integrating over $\mathcal{A}$ and $\Omega$ is efficiently realized by solving the state equation with fixed system matrix ${\mathbf{K}(\rho,\xi)}$ for 64 right hand sides simultaneously.

In contrast, $\xi$ directly impacts the system matrix ${\mathbf{K}(\rho,\xi)}$. As a consequence, each realization of ${\xi\in\Xi}$ requires us to solve a new state equation, making the integration over $\Xi$ much more expensive. For sMMA, we chose a batch size of 1 for the CSG-type approximation, meaning that only a single realization ${\xi\in\Xi}$ is evaluated in each iteration. In contrast, MMA uses a trapezoidal rule on a ${5\times5}$ space grid to numerically approximate the integral. Thus, when comparing our numerical results, it is important to keep in mind that each MMA iteration is 25 times as expensive as an sMMA iteration, since a system solve ($\widehat{=}$ solving the state equation for fixed ${\xi\in\Xi}$ and 64 right hand sides as discussed above) is by far the most time-consuming task in each optimization step.

For our analysis, we also evaluate the chance constraint value of intermediate designs using a ${50\times50}$ grid on $\Xi$ and use this as a comparison to the internal approximations of sMMA and MMA. Due to the tremendous numerical effort, this is only done every 100th iteration.

\begin{figure} 
    \centering
    \begin{minipage}[b][][c]{0.48\textwidth}
        \centering
        \includegraphics[width=\textwidth,keepaspectratio]{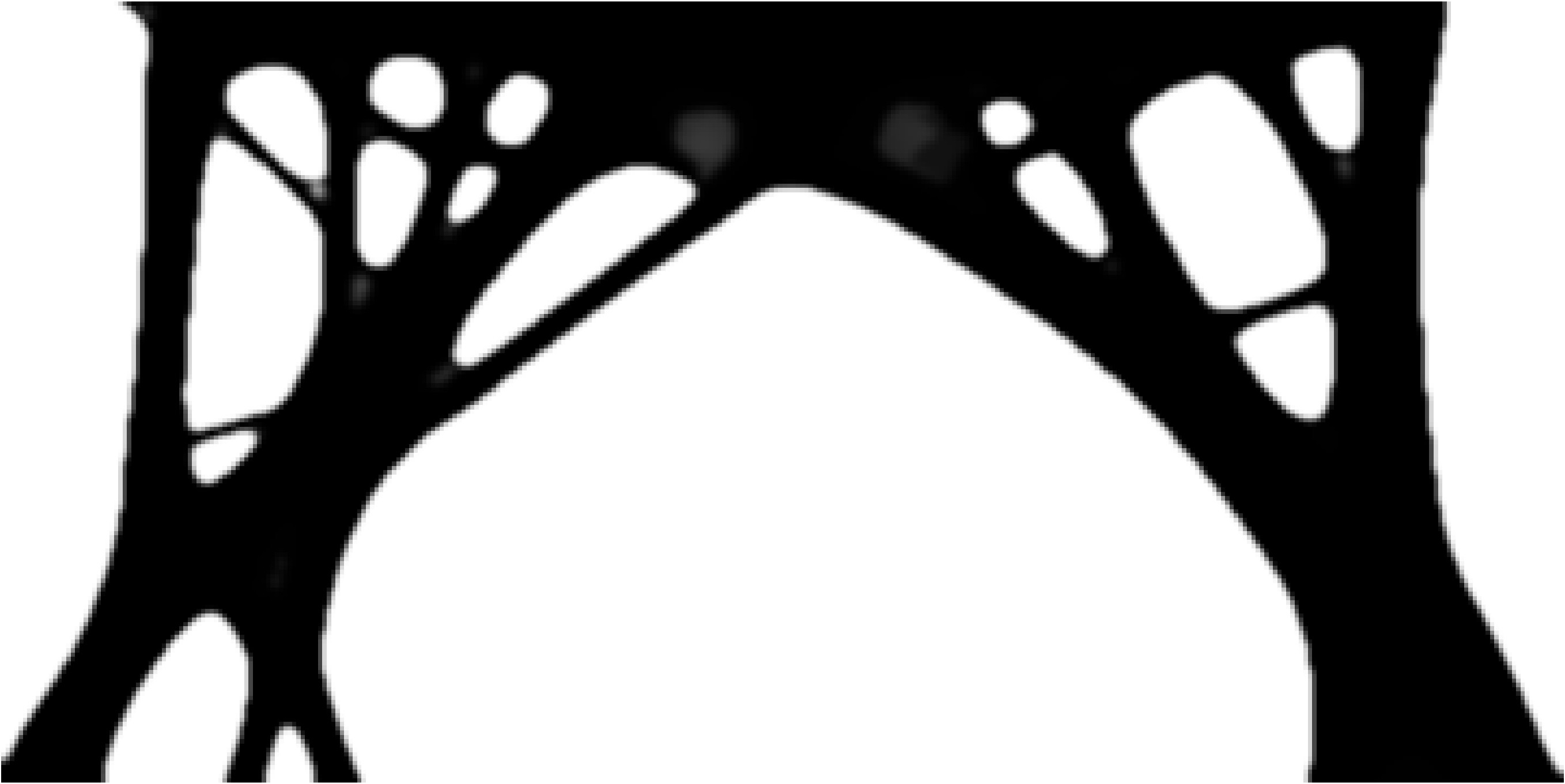}
    \end{minipage}\hfill
    \begin{minipage}[b][][c]{0.48\textwidth}
        \centering
        \includegraphics[width=\textwidth,keepaspectratio]{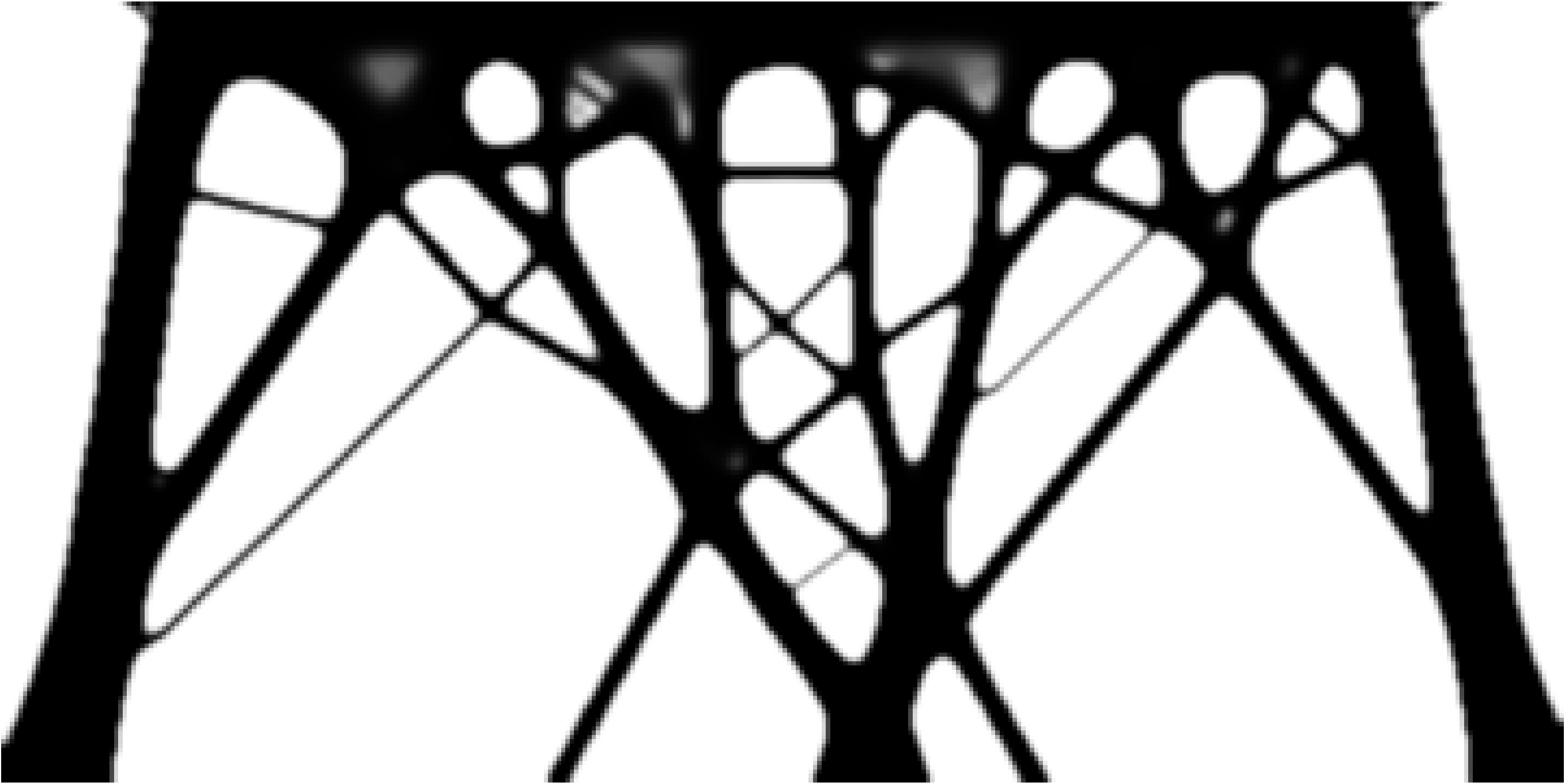}
    \end{minipage}
    \par
    \caption{Final designs obtained by sMMA (left) and MMA (right).}
    \label{fig:2d_designs}
\end{figure}

\begin{figure}
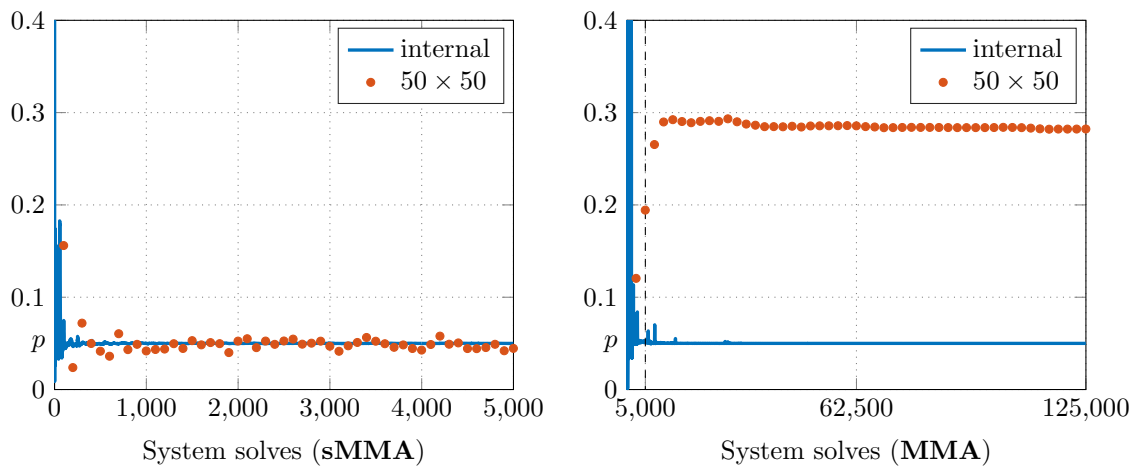
 
    \centering
    \begin{subfigure}{0.5\textwidth}
        \centering
        \input{MCMSA_CC}
    \end{subfigure}%
    \begin{subfigure}{0.5\textwidth}
        \centering
        \input{MMA_CC}%
    \end{subfigure}
   \caption{Chance constraint value evolution for sMMA (left) and MMA (right) during the optimization process. Solid lines correspond to the internal approximations using the CSG-type integration (sMMA) or a trapezoidal rule with 25 integration points (MMA). Every 100 iterations, intermediate designs are also analyzed using $2,500$ integration points (red dots) for comparison. Despite MMA requiring 25 times the numerical effort of sMMA, the produced design is infeasible, due to the large discretization error. For sMMA, we observe a close approximation to the true chance constraint value even for early iterations.}
    \label{fig:2d_cc}
\end{figure}

\begin{figure}
  \begin{minipage}[c]{0.5\textwidth}
    \includegraphics[width=\textwidth]{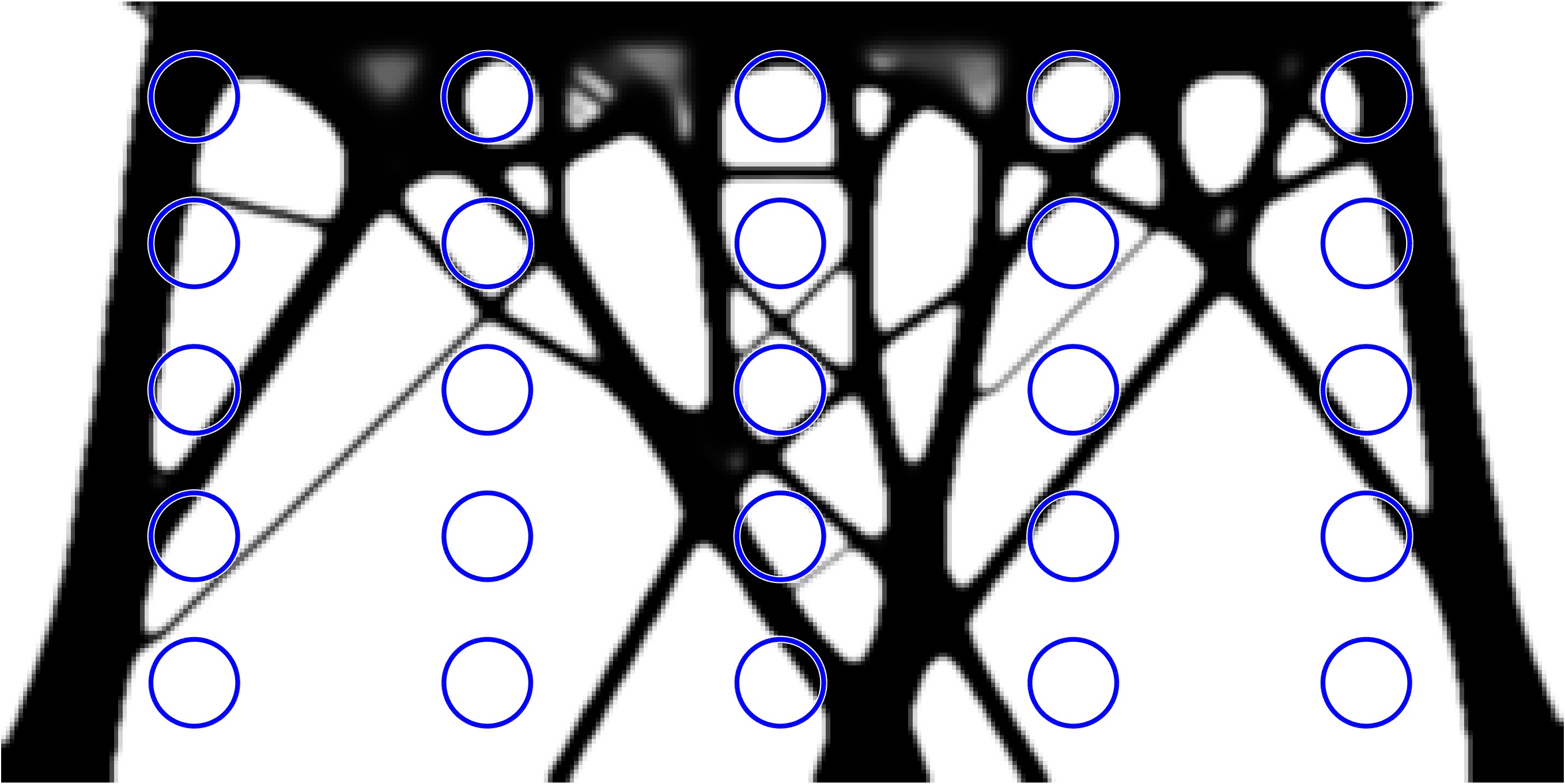}
  \end{minipage}\hfill
  \begin{minipage}[c]{0.45\textwidth}
    \caption{Final MMA design and discretization grid for the integration over $\Xi$. The 25 blue circles indicate the individual realizations for $\xi$, appearing in the trapezoidal rule used for the MMA optimization approach. As expected, the optimizer avoids to place material inside these circles, but constructs small structures in the rest of the design region.}
    \label{fig:2d_disc}
  \end{minipage}
\end{figure}

\begin{figure}
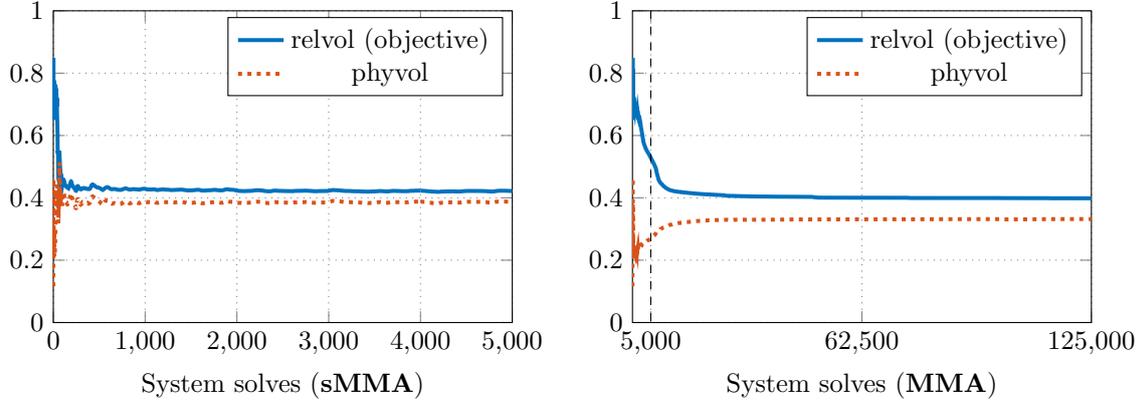
 
    \centering
    \begin{subfigure}{0.5\textwidth}
        \input{MCMSA_Obj}
    \end{subfigure}%
    \begin{subfigure}{0.5\textwidth}
        \input{MMA_Obj}%
    \end{subfigure}
    \caption{Objective function value ($\rvol$) during the optimization process for sMMA (left) and MMA (right), depicted by solid blue curves. The volumes associated to the physical interpretations via the SIMP method ($\pvol$) are indicated by the dotted red lines. Both optimization approaches yield comparable objective function values, despite only the sMMA design being feasible.}
    \label{fig:2d_objective}
\end{figure}

The final structures for sMMA and MMA are shown in~\Cref{fig:2d_designs}. Again, we observe that sMMA yields a feasible solution, while MMA, despite the computational costs being 25 times as high, fails to do so (\Cref{fig:2d_cc}). In~\Cref{fig:2d_disc}, it can be seen how the design found by MMA avoids the regions corresponding to integration points for the integration over $\Xi$, indicating that an even finer (even more expensive) quadrature scheme would be needed to obtain a feasible design. This, however, is an unreasonable approach, since~\Cref{fig:2d_objective} (bottom) already shows a superior efficiency of sMMA (with ${\mathcal{B}=1}$) with respect to system solves. Therefore, when deciding on a batch size for MMA, we face the following problem: If we want to find a feasible solution, the batch size must be chosen so large that the required amount of system solves is several orders of magnitude higher than for sMMA. If, on the other hand, we want the numerical cost of MMA to be comparable to sMMA, we no longer obtain a feasible solution.

\begin{figure} 
    \centering
    \begin{minipage}[b][][c]{0.48\textwidth}
        \centering
        \includegraphics[width=\textwidth,keepaspectratio]{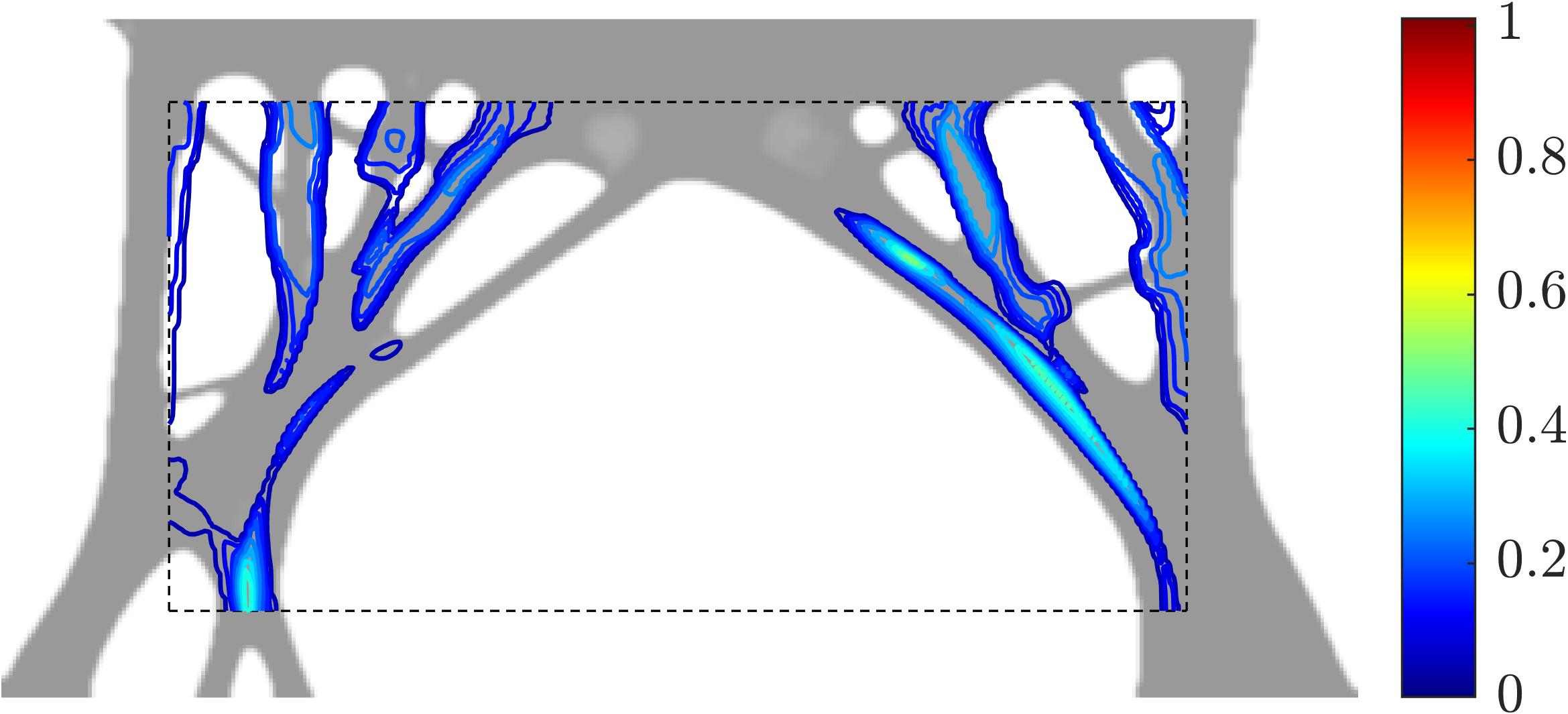}
    \end{minipage}\hfill
    \begin{minipage}[b][][c]{0.48\textwidth}
        \centering
        \includegraphics[width=\textwidth,keepaspectratio]{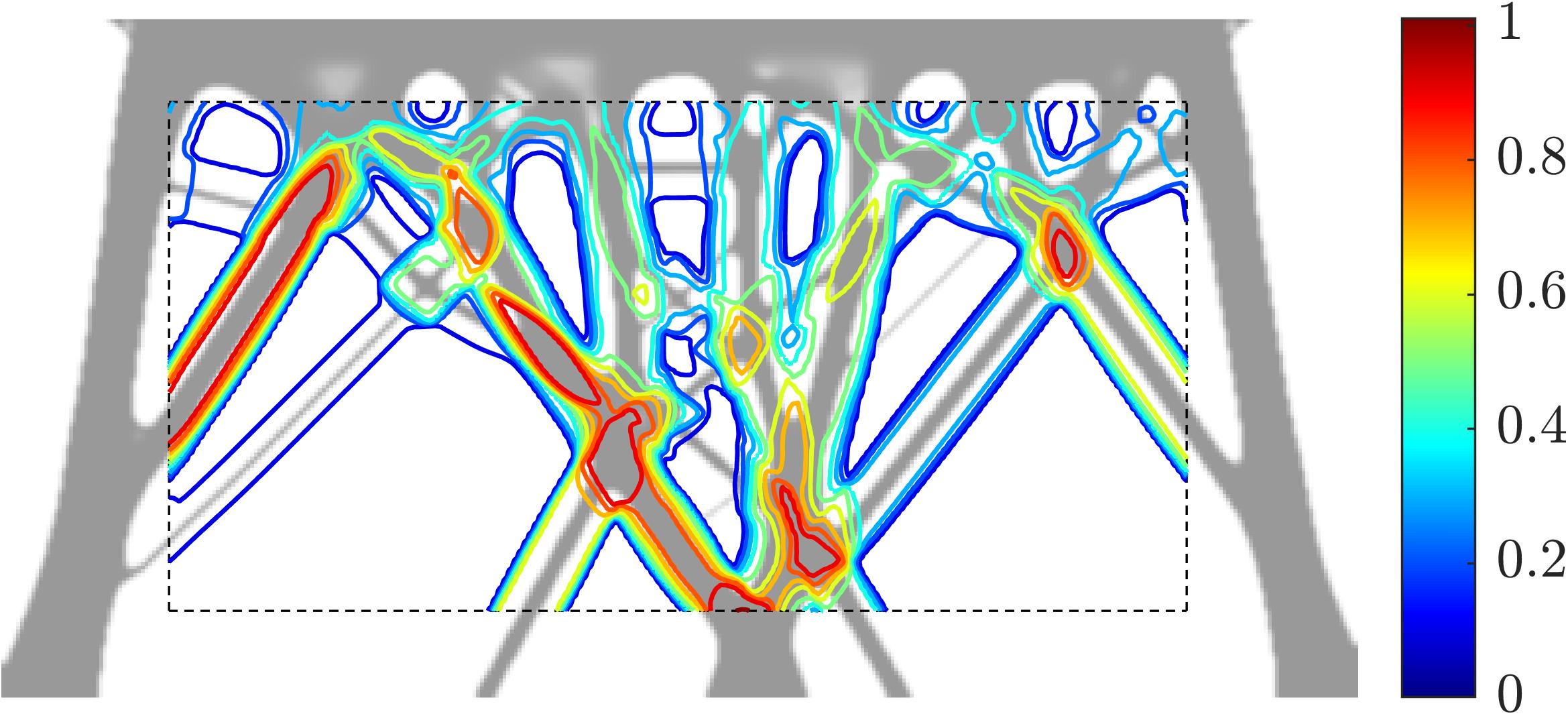}
    \end{minipage}
    \par
    \caption{Final designs for sMMA (left) and MM (right) with indicators for $\Xi$ (dashed rectangle). The colored contour lines depict values of the smoothed chance constraint $H_1$, as defined in~\eqref{eq:2d_h1}. For better visibility, regions with ${H_1(\xi)<0}$ are ignored. As we can see, the compliance bound is violated the most, if the material weakness is placed such that it essentially cuts a beam in the design. For sMMA, not many such fragile structures exist. In contrast, the MMA design consists of several very thin connections, most of which drastically reduce the structure's compliance when cut.}
    \label{fig:2d_analysis}
\end{figure}

In contrast to our previous example, it is no longer obvious at which realizations ${(\xi,\omega)\in\Xi\times\Omega}$ the compliance bound ${\mathbf{F}(\omega)^\top\mathbf{U}_{\xi,\omega}\le\cm}$ is violated. Thus, the $\xi$-dependent smoothed chance constraint values
\begin{equation}\label{eq:2d_h1}
    H_1(\xi):=\int_\Omega h_{a_1,a_2,a_3}\left(\frac{1}{\vert\mathcal{A}\vert}\int_\mathcal{A}\mathbf{F}(\alpha,\omega)^\top\mathbf{U}_{\alpha,\xi,\omega}\dd{\alpha}-\cm\right)\dd{\omega}
\end{equation}
for both final designs are given in~\Cref{fig:2d_analysis}. Intuitively, $H_1$ relates to how strong the compliance bound is violated, if the material is weakened at position ${\xi\in\Xi}$. By construction, ${H_1\le\tfrac{1}{2}}$ indicates that the compliance bound is satisfied even for the damaged structure, while ${H_1>\tfrac{1}{2}}$ means the compliance bound is violated for the particular choice of $\xi$. If we take another step back, we may also define the relative compliance
\begin{equation*}
    H_2(\xi,\omega):= \frac{\tfrac{1}{\vert\mathcal{A}\vert}\int_{\mathcal{A}}\mathbf{F}(\alpha,\omega)^\top\mathbf{U}_{\alpha,\xi,\omega}\dd{\alpha}}{\cm}.
\end{equation*}
% Now, for any ${0\le q_1<q_2\le 1}$, we can consider the quantile sets 
% \begin{equation*}
%     P_{q_1,q_2}(\omega):= [r_1(\omega),r_2(\omega)]\subset\R, 
% \end{equation*}
% with
% \begin{equation}\label{eq:2d_quantDef}
%     \P_\xi\big[ H_2(\xi,\omega)\le r_1(\omega) \big] = q_1\quad\text{and}\quad \P_\xi\big[ H_2(\xi,\omega)\ge r_2(\omega) \big] = 1-q_2.
% \end{equation}
% For example, given ${\omega\in\Omega}$, the set $P_{0.2,0.75}(\omega)$ gives the range of all values for ${H_2(\cdot,\omega)}$, if the smallest 20\% and largest 25\% of values are discarded. Different quantile sets for both final designs are shown in~\Cref{fig:2d_quants}. 

Notably, although only the sMMA design is feasible, the objective function value of both results are very close. This can be seen in~\Cref{fig:2d_objective}, where the values of $\rvol$ and $\pvol$ during the optimization are shown. The large difference between these values for MMA is a consequence of the many thin structures incorporated in the final MMA design, as interfaces between material and void are the dominant source of fictitious material, due to the applied density filter.

% \begin{figure}
%   \begin{minipage}[c]{0.5\textwidth}
%     \input{2d_quants}
%   \end{minipage}\hfill
%   \begin{minipage}[c]{0.45\textwidth}
%     \caption{Quantile regions of relative compliance for sMMA (blue) and MMA (red) as defined in~\eqref{eq:2d_quantDef}. Lightly shaded areas correspond to ${P_{0.05,0.95}(\omega)}$, whereas dark areas indicate ${P_{0.25,0.75}(\omega)}$. The solid lines depict the median of ${H_2(\cdot,\omega)}$, depending on $\omega$.}
%     \label{fig:2d_quants}
%   \end{minipage}
% \end{figure}

%#############################################################################################
%################################### 3d Table ################################################
%#############################################################################################
%\FloatBarrier
\subsection{3d compliance with load uncertainty}
\begin{figure}
  \begin{minipage}[c]{0.5\textwidth}
    \input{3dSetupSpat}
  \end{minipage}\hfill
  \begin{minipage}[c]{0.45\textwidth}
    \caption{Cubic design region $\mathscr{D}$ with three fixed supports at the bottom (grey squares). A downward force is applied at $\mathscr{F}$ (blue square), where the midpoint is drawn at random within ${\Omega\times\{\ell\}}$ (red region). The illustration is true to scale.}
    \label{fig:3d_setup}
  \end{minipage}
\end{figure}
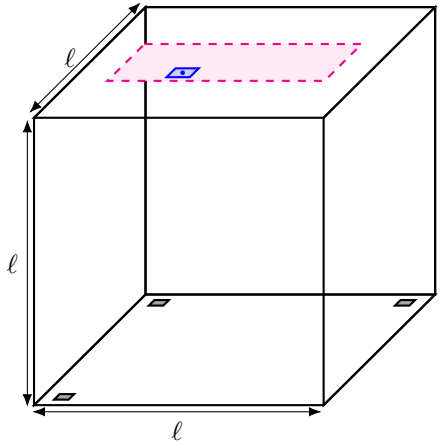

For our final example, we switch to a three-dimensional setup. Specifically, given ${\ell>0}$, we define a cubic design region ${\mathscr{D}:=[0,\ell]^3}$, which is again loaded from the top in a random area. This time, the structure is held in place at three squares on the bottom of the design region. Each square has a side length of $\tfrac{\ell}{20}$ and is spaced with a distance of $\tfrac{\ell}{20}$ to the boundary of $\mathscr{D}$. Defining ${\Omega:= \left[ \tfrac{\ell}{8},\tfrac{7\ell}{8} \right] \times \left[ \tfrac{\ell}{3},\tfrac{2\ell}{3} \right]}$, the force is applied in negative $z$-direction at 
\begin{equation*}
    \mathscr{F}:= \left[ \omega_1-\tfrac{\ell}{25},\omega_1+\tfrac{\ell}{25} \right] \times \left[ \omega_2-\tfrac{\ell}{25},\omega_2+\tfrac{\ell}{25} \right] \times \{\ell\},
\end{equation*}
where ${(\omega_1,\omega_2):=\omega\in\Omega}$. An illustration is given in~\Cref{fig:3d_setup}.

For our numerical experiments, the design region was discretized into $60^3$ cubic finite elements. The stiffness parameters for material and void were once again chosen as 1 and $10^{-4}$, respectively. The SIMP parameter was held fixed at a value of 5 during the optimization and the smoothing parameters were ${a_1=20}$, ${a_2=\tfrac{1}{50}}$ as well as ${a_3=5}$. Move limits were initialized as ${\tau=1}$ and halved each 250 iterations. Lastly, we chose the probability level ${p=0.05}$, initial design ${\rho_1\equiv 0.75}$ and assumed a uniform distribution for $\omega$.

As was the case for our first example, uncertainty enters the optimization problem
% \begin{align*}
%         \min_{\rho\in\R^d}\quad & \rvol(\rho),\\
%         \text{s.t.}\quad & \frac{1}{\vert\Omega\vert}\int_\Omega h_{a_1,a_2,a_3}\left(\mathbf{F}(\omega)^\top\mathbf{U}_{\omega}-\cm\right)\dd{\omega}\le p,\\
%             &\mathbf{K}(\rho)\mathbf{U}_{\omega}=\mathbf{F}(\omega),\\
%             &0\le \rho\le 1,
% \end{align*}
exclusively through the right hand side of the state equation. This time, however, the system matrix is far too large to apply a direct solver (this would require $\sim3.5$TB of memory), meaning that the evaluation for every individual $\omega$ is costly. Thus, for MMA, we discretized the integral using ${4\times4=16}$ integration points per iteration. Since the large design dimension might also slow down the sMMA integration weight computation, we used limited memory sMMA to solve the problem. To be precise, we limited the number of stored gradient information to $3,000$ samples. For a fair comparison, limited memory sMMA was run with a batch size of 16 as well. Both optimizers were allowed to perform $1,000$ iterations. For our analysis, the final designs were evaluated on a regular grid of ${100\times50=5,000}$ load cases.

\begin{figure} 
    \centering
    \begin{subfigure}{0.4\textwidth}
        \includegraphics[width=\textwidth,keepaspectratio]{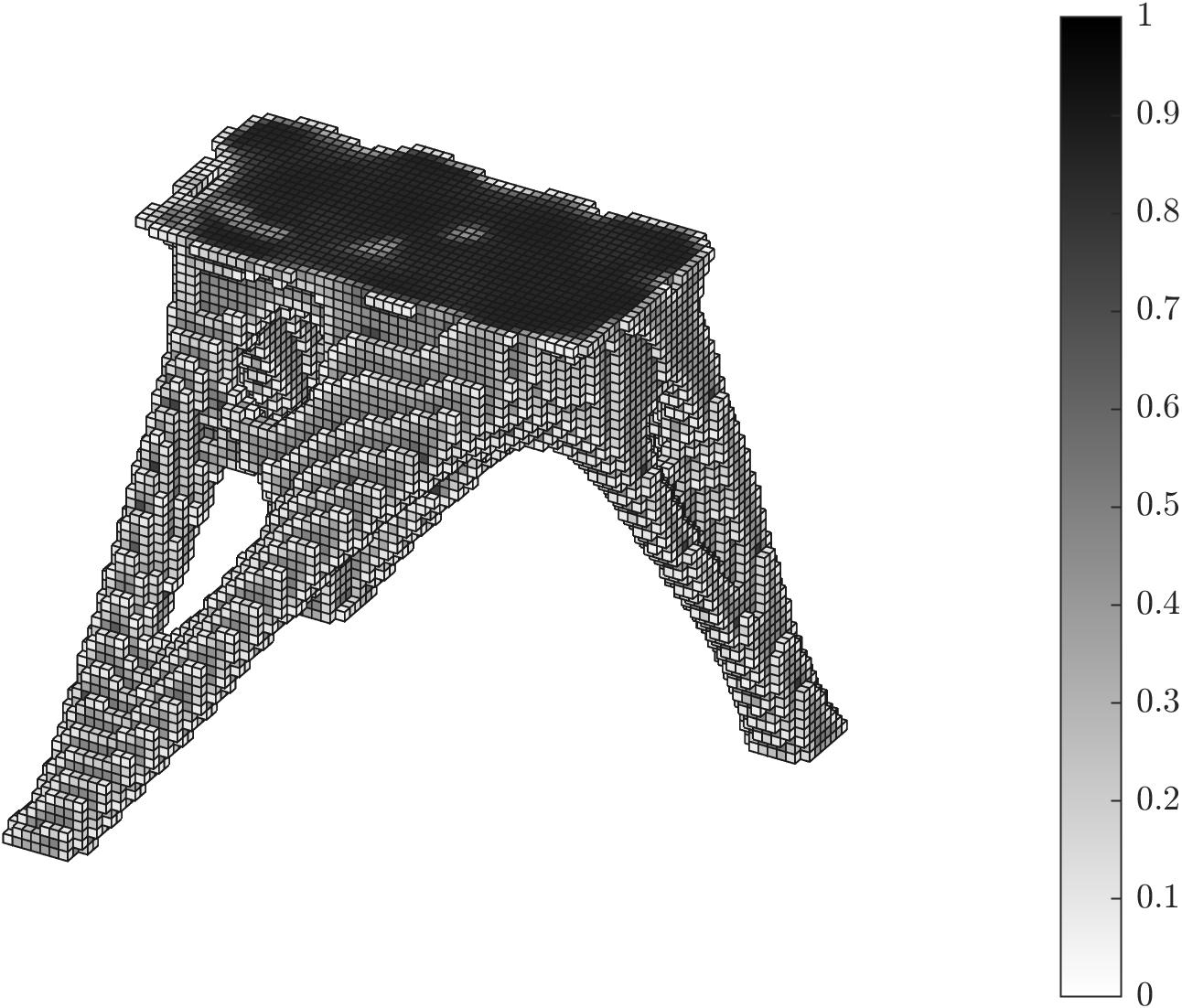}
    \end{subfigure}\hfill
    \begin{subfigure}{0.4\textwidth}
        \includegraphics[width=\textwidth,keepaspectratio]{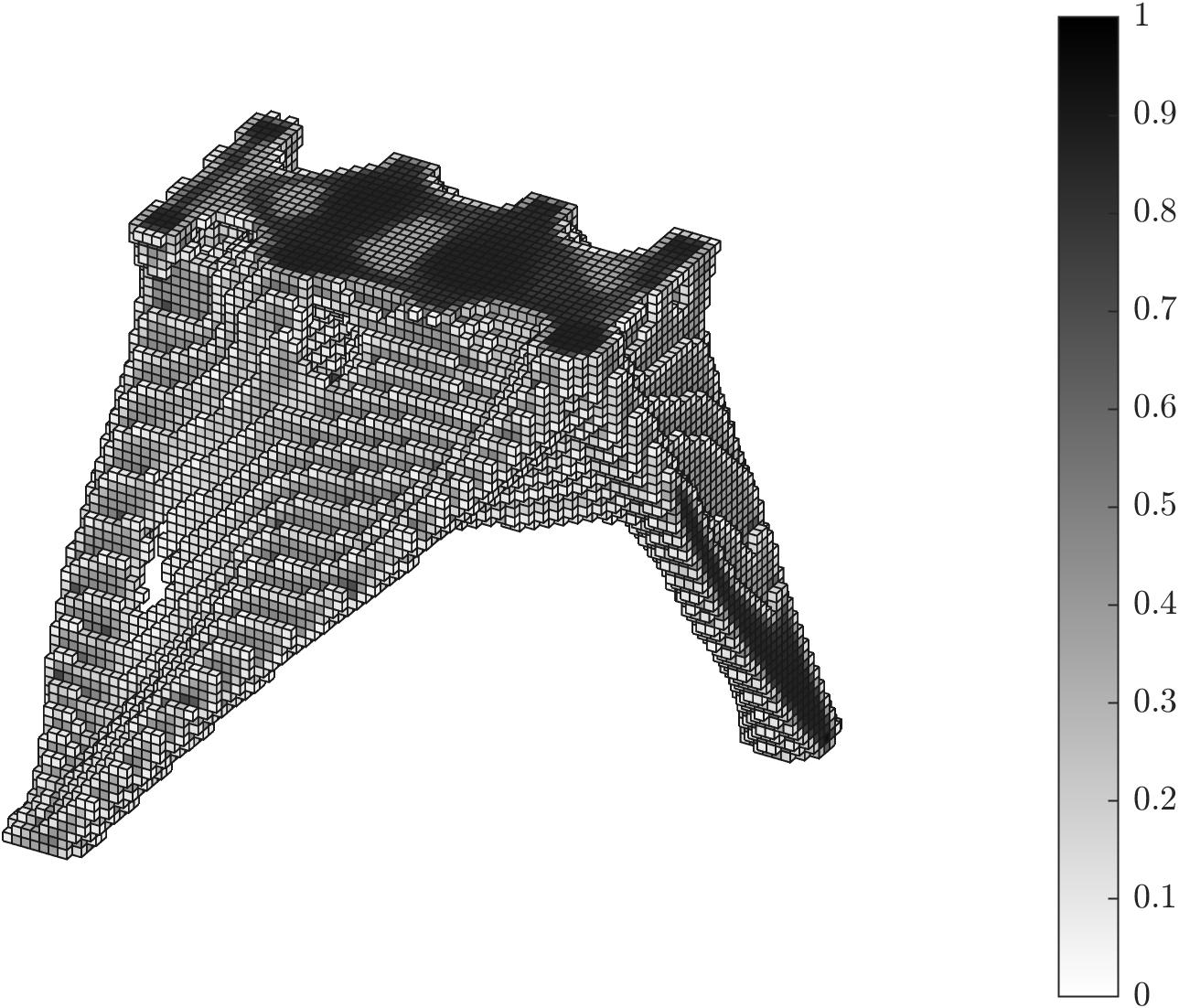}
    \end{subfigure}
    \caption{Final designs produced by limited memory sMMA (left) and MMA (right). Note that, in contrast to sMMA, MMA does not cover the whole region ${\Omega\times\{\ell\}}$ with material. As it turns out, the 16 integration points used for the discretization of $\Omega$ in MMA are spaced too far apart, resulting in the optimizer not placing material in between.}
    \label{fig:3d_designs}
\end{figure}

\begin{figure} 
    \centering
    \begin{subfigure}{0.48\textwidth}
        \includegraphics[width=\textwidth,keepaspectratio]{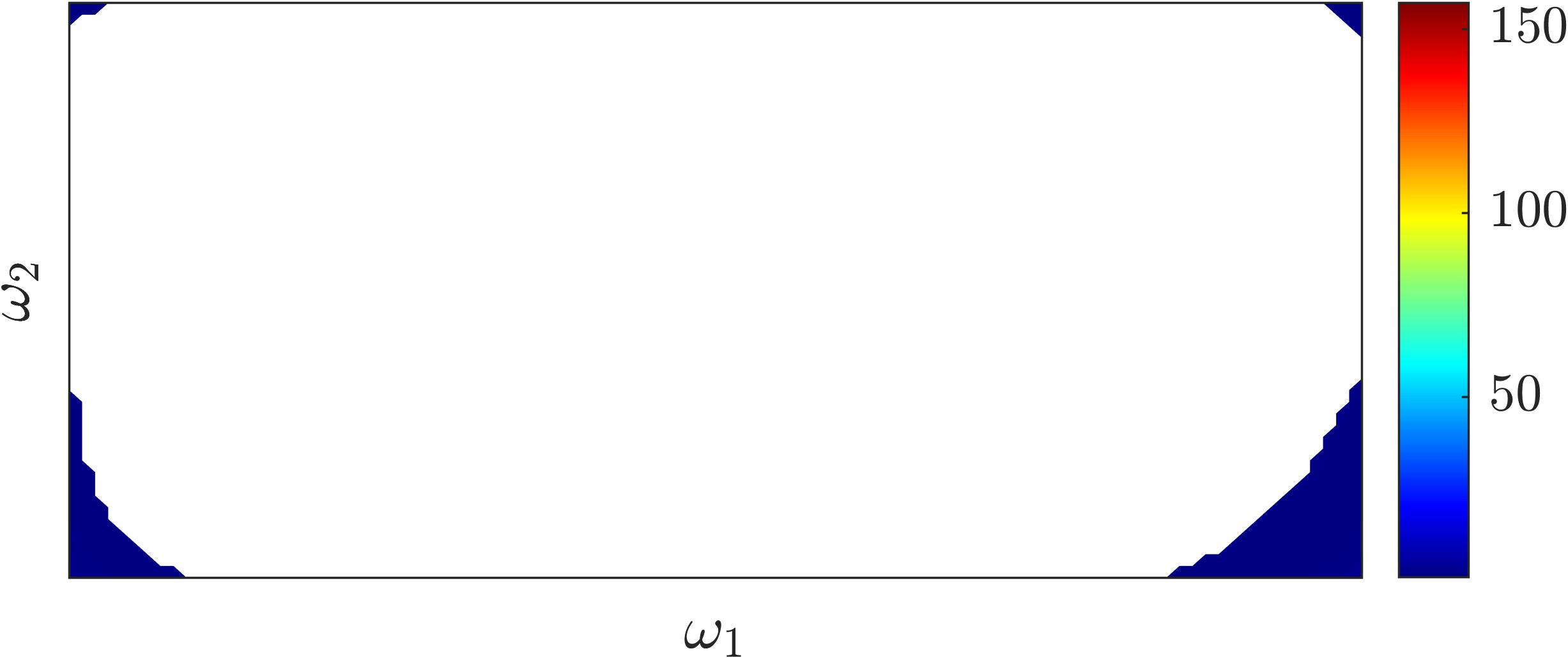}
    \end{subfigure}\hfill
    \begin{subfigure}{0.48\textwidth}
        \includegraphics[width=\textwidth,keepaspectratio]{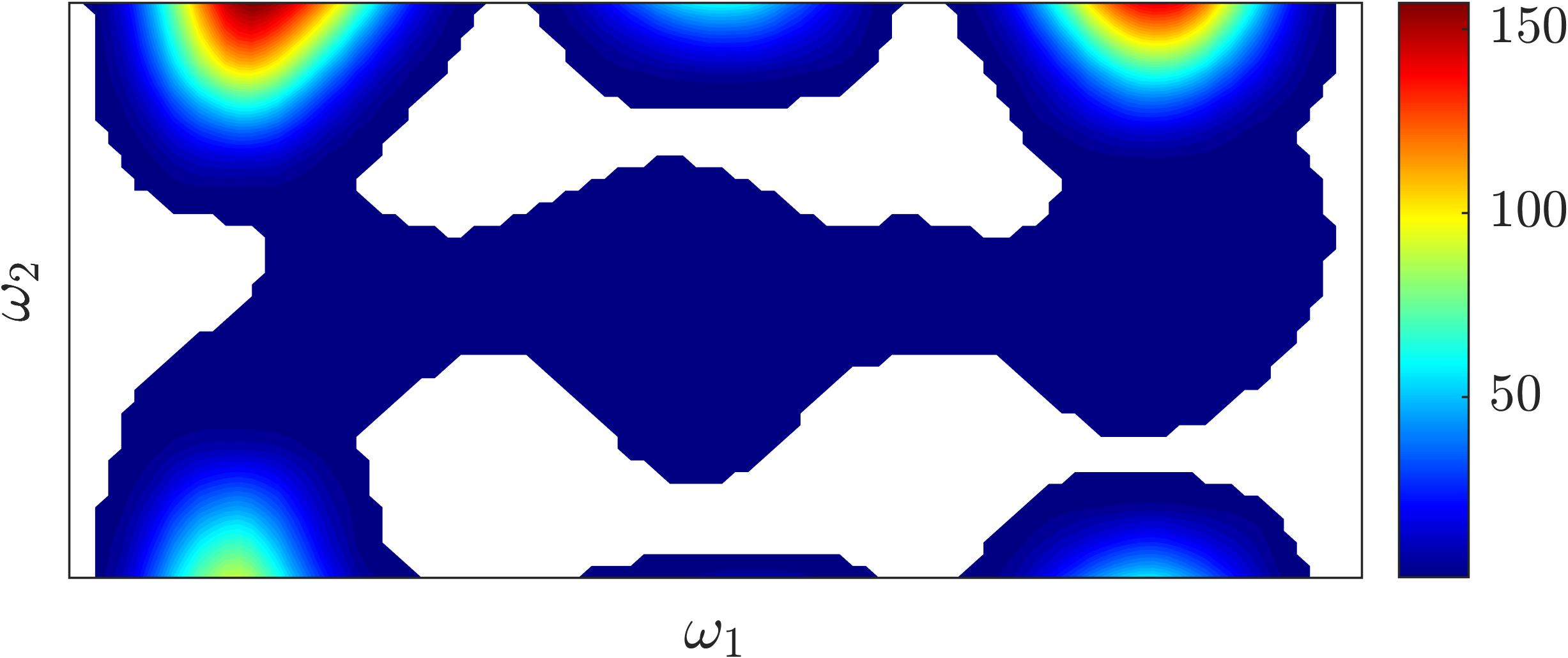}
    \end{subfigure}
    \caption{Relative compliance ${H_3(\omega)}$ for the final designs obtained by limited memory sMMA (left) and MMA (right). The values were calculated on the ${100\times50}$ verification grid. Negative values of $H_1$ were omitted. As we can see, compliance violations for the MMA design are especially high in between the discretization points used for the integral approximation. For limited memory sMMA, the relative compliance values are located at the extreme points (corners) of $\Omega$.}
    \label{fig:3d_viols}
\end{figure}

\begin{figure} 
    \centering
    \begin{minipage}[b][][c]{0.48\textwidth}
        \centering
        \input{MCMSA_rel_comp_histo}
    \end{minipage}\hfill
    \begin{minipage}[b][][c]{0.48\textwidth}
        \centering
        \input{MMA_rel_comp_histo}
    \end{minipage}
    \par
    \caption{Relative compliance values ${H_3(\omega)}$ for the final designs of limited memory sMMA (left) and MMA (right), evaluated on the ${100\times50}$ verification grid. Values smaller than 1 (left of dashed line) correspond to realizations of $\omega$ that satisfy the compliance bound. For MMA, values larger than 4 have been grouped together in the rightmost bar for the sake of better readability. For limited memory sMMA, no values larger than $3.5$ were observed.}
    \label{fig:3d_comp_histo}
\end{figure}
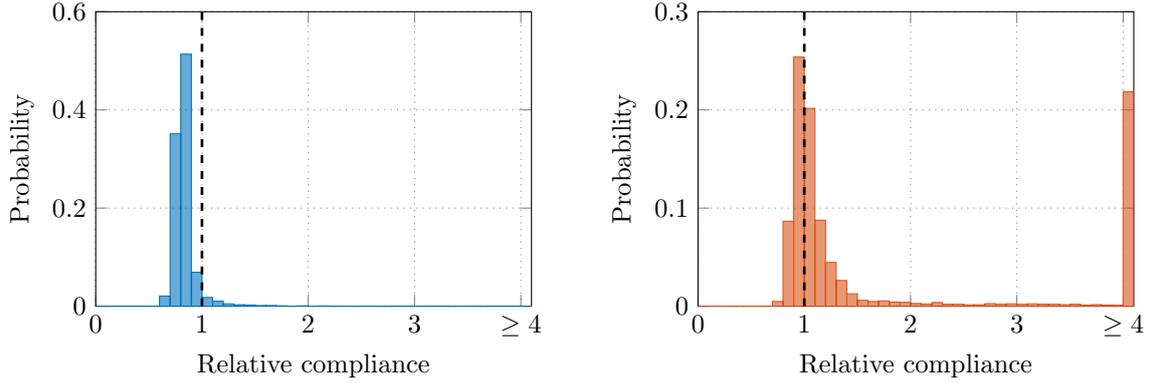

The final designs obtained by both methods are shown in~\Cref{fig:3d_designs}. Therein, we can already see artifacts of the discretization for the MMA design. As it turns out, MMA again fails to produce a feasible design, yielding a chance constraint value of ${p_{\text{veri}}=1.44}$ in the verification process (note that values larger than 1 are possible due to ${a_2>0}$). In contrast, for the final design obtained by limited memory sMMA, we have ${p_{\text{veri}}=0.046}$. The relative compliance
\begin{equation*}
    H_3(\omega):=\frac{\mathbf{F}(\omega)^\top\mathbf{U}_\omega}{\cm},
\end{equation*}
evaluated on the ${100\times50}$ verification grid, can be found in~\Cref{fig:3d_viols}. Additionally, histograms of $H_3$ for both final designs are shown in~\Cref{fig:3d_comp_histo}.

\begin{figure}
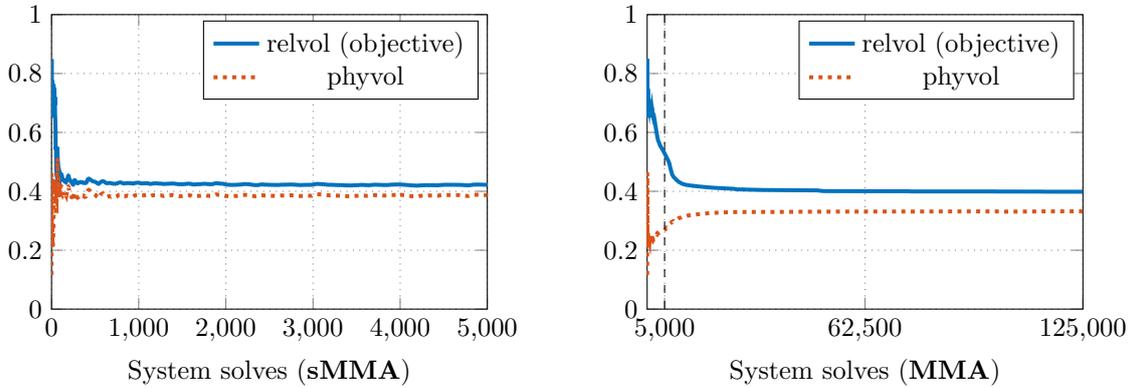
 
    \centering
    \begin{minipage}[b][][c]{0.48\textwidth}
        \centering
        \input{MCMSA_Obj}
    \end{minipage}\hfill
    \begin{minipage}[b][][c]{0.48\textwidth}
        \centering
        \input{MMA_Obj}
    \end{minipage}
    \par
    \caption{Objective function values (solid blue curves) during the optimization process for limited memory sMMA (left) and MMA (right). Volumes corresponding to the physical SIMP interpretation of designs ($\pvol$) are shown by the dotted red line.}
    \label{fig:3d_objective}
\end{figure}

Surprisingly, the feasible design produced by limited memory sMMA also admits a (slightly) lower objective function value than the infeasible design obtained by MMA. Moreover, as we can see in the objective function evolution (\Cref{fig:3d_objective}), the convergence of limited memory sMMA is also faster.

%#############################################################################################
%################################## Conclusion ###############################################
%#############################################################################################
%\FloatBarrier
\section{Conclusion}
To tackle topology optimization problems involving uncertainty, e.g., chance constraints or reliability-based topology optimization, we proposed the stochastic method of moving asymptotes (sMMA). By construction, the sample-based strategy requires much less simulations of the system than traditional approaches from literature. Nonetheless, the gradient approximation error of sMMA almost surely converges to zero during the optimization process, guaranteeing that the lower computational cost has no negative impact on the quality of solutions obtained. This is achieved by adaptively recombining old sample information to construct a nearest neighbor surrogate model for the gradient. Abusing the special structure of this model, the integration can be reduced to simple distance calculations, meaning that the associated numerical effort is negligible in contrast to rest of the iteration, especially the system simulation.

For several two- and three-dimensional applications from structural design optimization, the performance of sMMA was analyzed and compared to a standard deterministic MMA approach. In each example, we observed one the following outcomes:
\begin{itemize}
    \item If the numerical cost of sMMA and MMA were roughly equal, only sMMA produced a feasible design.
    \item If MMA and sMMA both produced a feasible design, sMMA required far less system solves to do so.
\end{itemize}

% \section*{Declaration of interests}
% The authors declare that they have no known competing financial interests or personal relationships that could have appeared to influence the work reported in this paper.
\section*{Acknowledgments}
The research was funded by the Deutsche Forschungsgemeinschaft (DFG, German Research Foundation) under Project-ID 416229255 (CRC 1411) and Project-ID 239904186 (TRR 154).
%#############################################################################################
%################################## Literature ###############################################
%#############################################################################################

\bibliographystyle{siamplain}
\bibliography{Biblio}
\end{document}

%% file: ex_shared.tex
\usepackage[a4paper, total={6in, 9in}]{geometry}
\usepackage{amsfonts}
\usepackage{graphicx}
\usepackage{algorithmic}

\usepackage{mathrsfs}
\usepackage{amsmath,amsfonts,amsthm}
\usepackage{multirow}
\usepackage{xcolor}
\usepackage{tikz}
\usepackage{pgfplots}
\usepackage{placeins}
\usepackage{algorithm}
\usepackage[export]{adjustbox}
\usepackage[title]{appendix}
\usepackage{subcaption}

\usepackage{hyperref}
\hypersetup{
    colorlinks=true,
    linkcolor=blue,
    citecolor=blue,
}
\usepackage{cleveref}
\usepackage{pgfplots}
\pgfplotsset{compat=newest}
\pgfplotsset{grid style={dotted,gray}}
\usetikzlibrary{plotmarks,pgfplots.fillbetween,shapes.geometric, arrows.meta}

\newcommand{\R}{\mathbb{R}}
\newcommand{\N}{\mathbb{N}}
\newcommand{\E}{\mathbb{E}}
\renewcommand{\P}{\mathbb{P}}
\newcommand{\U}{\mathcal{U}}
\newcommand{\X}{\mathcal{X}}

\newcommand{\dd}{\,\mathrm{d}}
\newcommand{\cm}{c_{\textnormal{max}}}

\DeclareMathOperator{\vol}{vol}
\DeclareMathOperator{\rvol}{relvol}
\DeclareMathOperator{\pvol}{phyvol}

\DeclareMathOperator{\supp}{supp}
\DeclareMathOperator*{\argmin}{arg\,min}
\DeclareMathOperator{\Dim}{dim}

\newtheorem{theorem}{Theorem}[section]
\newtheorem{proposition}[theorem]{Proposition}
\newtheorem{assumption}{Assumption}
\crefname{theorem}{Theorem}{Theorem}
\crefname{assumption}{Assumption}{Assumption}
\crefname{proposition}{Proposition}{Proposition}
\theoremstyle{definition}

\title{A stochastic method of moving asymptotes for topology optimization under uncertainty\thanks{Corresponding author: andrian.uihlein@fau.de}}%\thanks{Submitted to the editors DATE.
\date{}

\author{Lukas Pflug\thanks{FAU Competence Center Scientific Computing, Friedrich-Alexander-Universität Erlangen-Nürnberg.}
\and Michael Stingl\thanks{Department of Mathematics, Chair of Applied Mathematics, Friedrich-Alexander-Universität Erlangen-Nürnberg.}
\and Andrian Uihlein\footnotemark[3]}

\usepackage{amsopn}

%% file: cc_smooth.tex
\definecolor{mycolor1}{rgb}{0.00000,0.44700,0.74100}%
\definecolor{mycolor2}{rgb}{0.85000,0.32500,0.09800}%
\definecolor{mycolor3}{rgb}{0.92900,0.69400,0.12500}%

\begin{tikzpicture}
\begin{axis}[%
width=\textwidth,
height=.75\textwidth,
%scale only axis,
xmin=-2,
xmax=2,
ymin=-.2,
ymax=1.2,
xmajorgrids,
ymajorgrids,
legend style={legend cell align=left, align=left, draw=white!15!black},
legend pos = north west
]
\addplot[domain=-2:0,smooth,samples=101,mycolor1,line width=2.5pt,line join=round]{0};
\addlegendentry{$\chi_{(0,\infty)}$};
\addplot[domain=0:2,smooth,samples=101,mycolor1,line width=2.5pt,line join=round,forget plot]{1};
\addplot[domain=-2:2,smooth,samples=101,mycolor2,line width=1pt]{.5*(tanh(35*x)+1)};
\addlegendentry{$h_{a_1}$};
\addplot[domain=-2:2,smooth,samples=101,mycolor3,line width=1pt,dashed]{.5*(tanh(35*x)+1)+5/100*(x-x/(1+exp(5*x))};
\addlegendentry{$h_{a_1,a_2,a_3}$};
\end{axis}        
\end{tikzpicture}

%% file: weights_L2.tex
\definecolor{mycolor6}{rgb}{0.30100,0.74500,0.93300}%
\definecolor{mycolor3}{rgb}{0.46600,0.67400,0.18800}%
\definecolor{mycolor4}{rgb}{0.85000,0.32500,0.09800}%
\definecolor{mycolor2}{rgb}{0.00000,0.44700,0.74100}%
\definecolor{mycolor5}{rgb}{0.49400,0.18400,0.55600}%
\definecolor{mycolor1}{rgb}{0.92900,0.69400,0.12500}%
\begin{tikzpicture}

\begin{axis}[%
width=\textwidth,
height=\textwidth,
at={(0.758in,0.481in)},
%scale only axis,
xmin=0,
xmax=1,
xtick={0.00897241721128805,0.108636483506987,0.177325192604266,0.350531015462249,0.410696618759897,0.475129885991535,0.53829413184652,0.743905921003457,0.894972342700167,0.977299811263042},
xticklabels={{$x_3$},{$x_{10}$},{$x_9$},{$x_8$},{$x_4$},{$x_2$},{$x_1$},{$x_7$},{$x_6$},{$x_5$}},
xlabel style={font=\color{white!15!black}},
xlabel={$\mathcal X$},
ymin=0,
ymax=1,
ytick={0.0843594963324943,0.134465272753184,0.207614999142686,0.273810181102767,0.328223085934377,0.378356401350255,0.519949634436407,0.58730070907623,0.732493861020181,0.820529287580996},
yticklabels={{$u_2$},{$u_{10}$},{$u_1$},{$u_8$},{$u_5$},{$u_6$},{$u_4$},{$u_3$},{$u_9$},{$u_7$}},
ylabel style={font=\color{white!15!black}},
ylabel={$\mathcal U$},
axis background/.style={fill=white},
x post scale=-1, % reverse axis direction
y post scale=-1
]
\addplot [color=white!60!black, only marks, mark=*, mark size = 2.5, mark options={solid, fill=white!60!black, white!60!black}, forget plot]
  table[row sep=crcr]{%
0.53829413184652	0.207614999142686\\
0.475129885991535	0.0843594963324943\\
0.00897241721128805	0.58730070907623\\
0.410696618759897	0.519949634436407\\
0.977299811263042	0.328223085934377\\
0.894972342700167	0.378356401350255\\
0.743905921003457	0.820529287580996\\
0.350531015462249	0.273810181102767\\
0.177325192604266	0.732493861020181\\
0.108636483506987	0.134465272753184\\
};
\addplot [color=black, dashed, line width=1.0pt, forget plot]
  table[row sep=crcr]{%
0	0.134465272753184\\
1	0.134465272753184\\
};
\addplot [color=mycolor1, line width=2.0pt, forget plot]
  table[row sep=crcr]{%
0	0.134465272753184\\
0.2698	0.134465272753184\\
};
\addplot [color=mycolor2, line width=2.0pt, forget plot]
  table[row sep=crcr]{%
0.2696	0.134465272753184\\
0.345	0.134465272753184\\
};
\addplot [color=mycolor3, line width=2.0pt, forget plot]
  table[row sep=crcr]{%
0.3448	0.134465272753184\\
0.5292	0.134465272753184\\
};
\addplot [color=mycolor4, line width=2.0pt, forget plot]
  table[row sep=crcr]{%
0.529	0.134465272753184\\
0.7926	0.134465272753184\\
};
\addplot [color=mycolor5, line width=2.0pt, forget plot]
  table[row sep=crcr]{%
0.7924	0.134465272753184\\
0.803	0.134465272753184\\
};
\addplot [color=mycolor6, line width=2.0pt, forget plot]
  table[row sep=crcr]{%
0.8028	0.134465272753184\\
1	0.134465272753184\\
};

\addplot[area legend, draw=white!70!black, fill=white!10!black, fill opacity=0.2, forget plot]
table[row sep=crcr] {%
x	y\\
0.191313138760448	0.443152008253133\\
0.205760196793843	0.529323355692337\\
-4.4257077073758	5.89955236009779\\
-9.97380095830437	6.66133814775094e-16\\
-8.38886821693578	-1.4983562408777\\
0.130233581282631	0.376603751855166\\
}--cycle;

\addplot[area legend, draw=white!70!black, fill=white!10!black, fill opacity=0.2, forget plot]
table[row sep=crcr] {%
x	y\\
0.490039018510314	0.370132333705433\\
0.648679251655518	0.434941346019338\\
0.68316377793438	0.552884969202744\\
0.456283037232469	0.804394918432142\\
0.205760196793843	0.529323355692337\\
0.191313138760448	0.443152008253133\\
}--cycle;

\addplot[area legend, draw=white!70!black, fill=white!10!black, fill opacity=0.2, forget plot]
table[row sep=crcr] {%
x	y\\
0.68316377793438	0.552884969202744\\
1.15606421903481	0.714449324866466\\
7.19135423645013	3.57567636466009\\
-1.11022302462516e-16	10.827058656007\\
-0.931615644360846	9.73666998860743\\
0.456283037232469	0.804394918432142\\
}--cycle;

\addplot[area legend, draw=white!70!black, fill=white!10!black, fill opacity=0.2, forget plot]
table[row sep=crcr] {%
x	y\\
0.205760196793843	0.529323355692337\\
-4.4257077073758	5.89955236009779\\
-0.931615644360846	9.73666998860743\\
0.456283037232469	0.804394918432142\\
}--cycle;
\addplot [color=mycolor4, only marks, mark=*, mark size = 2.5, mark options={solid, fill=mycolor4, mycolor4}, forget plot]
  table[row sep=crcr]{%
0.53829413184652	0.207614999142686\\
};

\addplot[area legend, draw=white!70!black, fill=mycolor4, fill opacity=0.2, forget plot]
table[row sep=crcr] {%
x	y\\
0.425660854787906	0.187523202803647\\
0.837927013728227	-0.0237491576774736\\
0.797079135971164	0.12493449036376\\
0.648679251655518	0.434941346019338\\
0.490039018510314	0.370132333705433\\
}--cycle;
\addplot [color=mycolor3, only marks, mark=*, mark size = 2.5, mark options={solid, fill=mycolor3, mycolor3}, forget plot]
  table[row sep=crcr]{%
0.475129885991535	0.0843594963324943\\
};

\addplot[area legend, draw=white!70!black, fill=mycolor3, fill opacity=0.2, forget plot]
table[row sep=crcr] {%
x	y\\
0.425660854787906	0.187523202803647\\
0.837927013728227	-0.0237491576774736\\
3.86402170560819	-6.25515806856365\\
5.55111512312578e-17	-9.90480035640677\\
-0.943269227939864	-8.92497930939605\\
0.290400057561397	0.0985642075814856\\
}--cycle;
\addplot [color=mycolor6, only marks, mark=*, mark size = 2.5, mark options={solid, fill=mycolor6, mycolor6}, forget plot]
  table[row sep=crcr]{%
0.977299811263042	0.328223085934377\\
};

\addplot[area legend, draw=white!70!black, fill=mycolor6, fill opacity=0.2, forget plot]
table[row sep=crcr] {%
x	y\\
0.797079135971164	0.12493449036376\\
0.837927013728227	-0.0237491576774736\\
3.86402170560819	-6.25515806856365\\
11.024258288727	0\\
7.19135423645013	3.57567636466009\\
1.15606421903481	0.714449324866466\\
}--cycle;
\addplot [color=mycolor5, only marks, mark=*, mark size = 2.5, mark options={solid, fill=mycolor5, mycolor5}, forget plot]
  table[row sep=crcr]{%
0.894972342700167	0.378356401350255\\
};

\addplot[area legend, draw=white!70!black, fill=mycolor5, fill opacity=0.2, forget plot]
table[row sep=crcr] {%
x	y\\
0.648679251655518	0.434941346019338\\
0.797079135971164	0.12493449036376\\
1.15606421903481	0.714449324866466\\
0.68316377793438	0.552884969202744\\
}--cycle;
\addplot [color=mycolor2, only marks, mark=*, mark size = 2.5, mark options={solid, fill=mycolor2, mycolor2}, forget plot]
  table[row sep=crcr]{%
0.350531015462249	0.273810181102767\\
};

\addplot[area legend, draw=white!70!black, fill=mycolor2, fill opacity=0.2, forget plot]
table[row sep=crcr] {%
x	y\\
0.425660854787906	0.187523202803647\\
0.490039018510314	0.370132333705433\\
0.191313138760448	0.443152008253133\\
0.130233581282631	0.376603751855166\\
0.290400057561397	0.0985642075814856\\
}--cycle;
\addplot [color=mycolor1, only marks, mark=*, mark size = 2.5, mark options={solid, fill=mycolor1, mycolor1}, forget plot]
  table[row sep=crcr]{%
0.108636483506987	0.134465272753184\\
};

\addplot[area legend, draw=white!70!black, fill=mycolor1, fill opacity=0.2, forget plot]
table[row sep=crcr] {%
x	y\\
0.290400057561397	0.0985642075814856\\
-0.943269227939864	-8.92497930939605\\
-8.38886821693578	-1.4983562408777\\
0.130233581282631	0.376603751855166\\
}--cycle;
\addplot [color=black, dotted, forget plot]
  table[row sep=crcr]{%
0.53829413184652	1\\
0.53829413184652	0.207614999142686\\
};
\addplot [color=black, dotted, forget plot]
  table[row sep=crcr]{%
1	0.207614999142686\\
0.53829413184652	0.207614999142686\\
};
\addplot [color=black, dotted, forget plot]
  table[row sep=crcr]{%
0.475129885991535	1\\
0.475129885991535	0.0843594963324943\\
};
\addplot [color=black, dotted, forget plot]
  table[row sep=crcr]{%
1	0.0843594963324943\\
0.475129885991535	0.0843594963324943\\
};
\addplot [color=black, dotted, forget plot]
  table[row sep=crcr]{%
0.00897241721128805	1\\
0.00897241721128805	0.58730070907623\\
};
\addplot [color=black, dotted, forget plot]
  table[row sep=crcr]{%
1	0.58730070907623\\
0.00897241721128805	0.58730070907623\\
};
\addplot [color=black, dotted, forget plot]
  table[row sep=crcr]{%
0.410696618759897	1\\
0.410696618759897	0.519949634436407\\
};
\addplot [color=black, dotted, forget plot]
  table[row sep=crcr]{%
1	0.519949634436407\\
0.410696618759897	0.519949634436407\\
};
\addplot [color=black, dotted, forget plot]
  table[row sep=crcr]{%
0.977299811263042	1\\
0.977299811263042	0.328223085934377\\
};
\addplot [color=black, dotted, forget plot]
  table[row sep=crcr]{%
1	0.328223085934377\\
0.977299811263042	0.328223085934377\\
};
\addplot [color=black, dotted, forget plot]
  table[row sep=crcr]{%
0.894972342700167	1\\
0.894972342700167	0.378356401350255\\
};
\addplot [color=black, dotted, forget plot]
  table[row sep=crcr]{%
1	0.378356401350255\\
0.894972342700167	0.378356401350255\\
};
\addplot [color=black, dotted, forget plot]
  table[row sep=crcr]{%
0.743905921003457	1\\
0.743905921003457	0.820529287580996\\
};
\addplot [color=black, dotted, forget plot]
  table[row sep=crcr]{%
1	0.820529287580996\\
0.743905921003457	0.820529287580996\\
};
\addplot [color=black, dotted, forget plot]
  table[row sep=crcr]{%
0.350531015462249	1\\
0.350531015462249	0.273810181102767\\
};
\addplot [color=black, dotted, forget plot]
  table[row sep=crcr]{%
1	0.273810181102767\\
0.350531015462249	0.273810181102767\\
};
\addplot [color=black, dotted, forget plot]
  table[row sep=crcr]{%
0.177325192604266	1\\
0.177325192604266	0.732493861020181\\
};
\addplot [color=black, dotted, forget plot]
  table[row sep=crcr]{%
1	0.732493861020181\\
0.177325192604266	0.732493861020181\\
};
\addplot [color=black, dotted, forget plot]
  table[row sep=crcr]{%
0.108636483506987	1\\
0.108636483506987	0.134465272753184\\
};
\addplot [color=black, dotted, forget plot]
  table[row sep=crcr]{%
1	0.134465272753184\\
0.108636483506987	0.134465272753184\\
};
\end{axis}

% \begin{axis}[%
% width=5.833in,
% height=4.375in,
% at={(0in,0in)},
% scale only axis,
% xmin=0,
% xmax=1,
% ymin=0,
% ymax=1,
% axis line style={draw=none},
% ticks=none,
% axis x line*=bottom,
% axis y line*=left
% ]
% \end{axis}
\end{tikzpicture}%

%% file: weights_pseudoexact.tex
\definecolor{mycolor6}{rgb}{0.30100,0.74500,0.93300}%
\definecolor{mycolor3}{rgb}{0.46600,0.67400,0.18800}%
\definecolor{mycolor4}{rgb}{0.85000,0.32500,0.09800}%
\definecolor{mycolor2}{rgb}{0.00000,0.44700,0.74100}%
\definecolor{mycolor5}{rgb}{0.49400,0.18400,0.55600}%
\definecolor{mycolor1}{rgb}{0.92900,0.69400,0.12500}%
\begin{tikzpicture}

\begin{axis}[%
width=\textwidth,
height=\textwidth,
at={(0.758in,0.481in)},
%scale only axis,
xmin=0,
xmax=1,
xtick={0.00897241721128805,0.108636483506987,0.177325192604266,0.350531015462249,0.410696618759897,0.475129885991535,0.53829413184652,0.743905921003457,0.894972342700167,0.977299811263042},
xticklabels={{$x_3$},{$x_{10}$},{$x_9$},{$x_8$},{$x_4$},{$x_2$},{$x_1$},{$x_7$},{$x_6$},{$x_5$}},
xlabel style={font=\color{white!15!black}},
xlabel={$\mathcal X$},
ymin=0,
ymax=1,
ytick={0.0843594963324943,0.134465272753184,0.207614999142686,0.273810181102767,0.328223085934377,0.378356401350255,0.519949634436407,0.58730070907623,0.732493861020181,0.820529287580996},
yticklabels={{$u_2$},{$u_{10}$},{$u_1$},{$u_8$},{$u_5$},{$u_6$},{$u_4$},{$u_3$},{$u_9$},{$u_7$}},
ylabel style={font=\color{white!15!black}},
ylabel={$\mathcal U$},
axis background/.style={fill=white},
x post scale=-1, % reverse axis direction
y post scale=-1
]
\addplot [color=white!60!black, only marks, mark=*, mark size = 2.5, mark options={solid, fill=white!60!black, white!60!black}, forget plot]
  table[row sep=crcr]{%
0.53829413184652	0.207614999142686\\
0.475129885991535	0.0843594963324943\\
0.00897241721128805	0.58730070907623\\
0.410696618759897	0.519949634436407\\
0.977299811263042	0.328223085934377\\
0.894972342700167	0.378356401350255\\
0.743905921003457	0.820529287580996\\
0.350531015462249	0.273810181102767\\
0.177325192604266	0.732493861020181\\
0.108636483506987	0.134465272753184\\
};
\addplot [color=black, dashed, line width=1.0pt, forget plot]
  table[row sep=crcr]{%
0	0.134465272753184\\
1	0.134465272753184\\
};
\addplot [color=black, line width=2.0pt, forget plot]
  table[row sep=crcr]{%
0	0.134465272753184\\
1	0.134465272753184\\
};
% \addplot [color=mycolor2, line width=2.0pt, forget plot]
%   table[row sep=crcr]{%
% 0.2696	0.134465272753184\\
% 0.345	0.134465272753184\\
% };
% \addplot [color=mycolor3, line width=2.0pt, forget plot]
%   table[row sep=crcr]{%
% 0.3448	0.134465272753184\\
% 0.5292	0.134465272753184\\
% };
% \addplot [color=mycolor4, line width=2.0pt, forget plot]
%   table[row sep=crcr]{%
% 0.529	0.134465272753184\\
% 0.7926	0.134465272753184\\
% };
% \addplot [color=mycolor5, line width=2.0pt, forget plot]
%   table[row sep=crcr]{%
% 0.7924	0.134465272753184\\
% 0.803	0.134465272753184\\
% };
% \addplot [color=mycolor6, line width=2.0pt, forget plot]
%   table[row sep=crcr]{%
% 0.8028	0.134465272753184\\
% 1	0.134465272753184\\
% };

\addplot[area legend, draw=white!70!black, fill=white!10!black, fill opacity=0.2, forget plot]
table[row sep=crcr] {%
x	y\\
0.191313138760448	0.443152008253133\\
0.205760196793843	0.529323355692337\\
-4.4257077073758	5.89955236009779\\
-9.97380095830437	6.66133814775094e-16\\
-8.38886821693578	-1.4983562408777\\
0.130233581282631	0.376603751855166\\
}--cycle;

\addplot[area legend, draw=white!70!black, fill=white!10!black, fill opacity=0.2, forget plot]
table[row sep=crcr] {%
x	y\\
0.490039018510314	0.370132333705433\\
0.648679251655518	0.434941346019338\\
0.68316377793438	0.552884969202744\\
0.456283037232469	0.804394918432142\\
0.205760196793843	0.529323355692337\\
0.191313138760448	0.443152008253133\\
}--cycle;

\addplot[area legend, draw=white!70!black, fill=white!10!black, fill opacity=0.2, forget plot]
table[row sep=crcr] {%
x	y\\
0.68316377793438	0.552884969202744\\
1.15606421903481	0.714449324866466\\
7.19135423645013	3.57567636466009\\
-1.11022302462516e-16	10.827058656007\\
-0.931615644360846	9.73666998860743\\
0.456283037232469	0.804394918432142\\
}--cycle;

\addplot[area legend, draw=white!70!black, fill=white!10!black, fill opacity=0.2, forget plot]
table[row sep=crcr] {%
x	y\\
0.205760196793843	0.529323355692337\\
-4.4257077073758	5.89955236009779\\
-0.931615644360846	9.73666998860743\\
0.456283037232469	0.804394918432142\\
}--cycle;
\addplot [color=mycolor4, only marks, mark=*, mark size = 2.5, mark options={solid, fill=mycolor4, mycolor4}, forget plot]
  table[row sep=crcr]{%
0.53829413184652	0.207614999142686\\
};

\addplot[area legend, draw=white!70!black, fill=mycolor4, fill opacity=0.2, forget plot]
table[row sep=crcr] {%
x	y\\
0.425660854787906	0.187523202803647\\
0.837927013728227	-0.0237491576774736\\
0.797079135971164	0.12493449036376\\
0.648679251655518	0.434941346019338\\
0.490039018510314	0.370132333705433\\
}--cycle;
\addplot [color=mycolor3, only marks, mark=*, mark size = 2.5, mark options={solid, fill=mycolor3, mycolor3}, forget plot]
  table[row sep=crcr]{%
0.475129885991535	0.0843594963324943\\
};

\addplot[area legend, draw=white!70!black, fill=mycolor3, fill opacity=0.2, forget plot]
table[row sep=crcr] {%
x	y\\
0.425660854787906	0.187523202803647\\
0.837927013728227	-0.0237491576774736\\
3.86402170560819	-6.25515806856365\\
5.55111512312578e-17	-9.90480035640677\\
-0.943269227939864	-8.92497930939605\\
0.290400057561397	0.0985642075814856\\
}--cycle;
\addplot [color=mycolor6, only marks, mark=*, mark size = 2.5, mark options={solid, fill=mycolor6, mycolor6}, forget plot]
  table[row sep=crcr]{%
0.977299811263042	0.328223085934377\\
};

\addplot[area legend, draw=white!70!black, fill=mycolor6, fill opacity=0.2, forget plot]
table[row sep=crcr] {%
x	y\\
0.797079135971164	0.12493449036376\\
0.837927013728227	-0.0237491576774736\\
3.86402170560819	-6.25515806856365\\
11.024258288727	0\\
7.19135423645013	3.57567636466009\\
1.15606421903481	0.714449324866466\\
}--cycle;
\addplot [color=mycolor5, only marks, mark=*, mark size = 2.5, mark options={solid, fill=mycolor5, mycolor5}, forget plot]
  table[row sep=crcr]{%
0.894972342700167	0.378356401350255\\
};

\addplot[area legend, draw=white!70!black, fill=mycolor5, fill opacity=0.2, forget plot]
table[row sep=crcr] {%
x	y\\
0.648679251655518	0.434941346019338\\
0.797079135971164	0.12493449036376\\
1.15606421903481	0.714449324866466\\
0.68316377793438	0.552884969202744\\
}--cycle;
\addplot [color=mycolor2, only marks, mark=*, mark size = 2.5, mark options={solid, fill=mycolor2, mycolor2}, forget plot]
  table[row sep=crcr]{%
0.350531015462249	0.273810181102767\\
};

\addplot[area legend, draw=white!70!black, fill=mycolor2, fill opacity=0.2, forget plot]
table[row sep=crcr] {%
x	y\\
0.425660854787906	0.187523202803647\\
0.490039018510314	0.370132333705433\\
0.191313138760448	0.443152008253133\\
0.130233581282631	0.376603751855166\\
0.290400057561397	0.0985642075814856\\
}--cycle;
\addplot [color=mycolor1, only marks, mark=*, mark size = 2.5, mark options={solid, fill=mycolor1, mycolor1}, forget plot]
  table[row sep=crcr]{%
0.108636483506987	0.134465272753184\\
};

\addplot[area legend, draw=white!70!black, fill=mycolor1, fill opacity=0.2, forget plot]
table[row sep=crcr] {%
x	y\\
0.290400057561397	0.0985642075814856\\
-0.943269227939864	-8.92497930939605\\
-8.38886821693578	-1.4983562408777\\
0.130233581282631	0.376603751855166\\
}--cycle;
\addplot [color=black, dotted, forget plot]
  table[row sep=crcr]{%
0.53829413184652	1\\
0.53829413184652	0.207614999142686\\
};
\addplot [color=black, dotted, forget plot]
  table[row sep=crcr]{%
1	0.207614999142686\\
0.53829413184652	0.207614999142686\\
};
\addplot [color=black, dotted, forget plot]
  table[row sep=crcr]{%
0.475129885991535	1\\
0.475129885991535	0.0843594963324943\\
};
\addplot [color=black, dotted, forget plot]
  table[row sep=crcr]{%
1	0.0843594963324943\\
0.475129885991535	0.0843594963324943\\
};
\addplot [color=black, dotted, forget plot]
  table[row sep=crcr]{%
0.00897241721128805	1\\
0.00897241721128805	0.58730070907623\\
};
\addplot [color=black, dotted, forget plot]
  table[row sep=crcr]{%
1	0.58730070907623\\
0.00897241721128805	0.58730070907623\\
};
\addplot [color=black, dotted, forget plot]
  table[row sep=crcr]{%
0.410696618759897	1\\
0.410696618759897	0.519949634436407\\
};
\addplot [color=black, dotted, forget plot]
  table[row sep=crcr]{%
1	0.519949634436407\\
0.410696618759897	0.519949634436407\\
};
\addplot [color=black, dotted, forget plot]
  table[row sep=crcr]{%
0.977299811263042	1\\
0.977299811263042	0.328223085934377\\
};
\addplot [color=black, dotted, forget plot]
  table[row sep=crcr]{%
1	0.328223085934377\\
0.977299811263042	0.328223085934377\\
};
\addplot [color=black, dotted, forget plot]
  table[row sep=crcr]{%
0.894972342700167	1\\
0.894972342700167	0.378356401350255\\
};
\addplot [color=black, dotted, forget plot]
  table[row sep=crcr]{%
1	0.378356401350255\\
0.894972342700167	0.378356401350255\\
};
\addplot [color=black, dotted, forget plot]
  table[row sep=crcr]{%
0.743905921003457	1\\
0.743905921003457	0.820529287580996\\
};
\addplot [color=black, dotted, forget plot]
  table[row sep=crcr]{%
1	0.820529287580996\\
0.743905921003457	0.820529287580996\\
};
\addplot [color=black, dotted, forget plot]
  table[row sep=crcr]{%
0.350531015462249	1\\
0.350531015462249	0.273810181102767\\
};
\addplot [color=black, dotted, forget plot]
  table[row sep=crcr]{%
1	0.273810181102767\\
0.350531015462249	0.273810181102767\\
};
\addplot [color=black, dotted, forget plot]
  table[row sep=crcr]{%
0.177325192604266	1\\
0.177325192604266	0.732493861020181\\
};
\addplot [color=black, dotted, forget plot]
  table[row sep=crcr]{%
1	0.732493861020181\\
0.177325192604266	0.732493861020181\\
};
\addplot [color=black, dotted, forget plot]
  table[row sep=crcr]{%
0.108636483506987	1\\
0.108636483506987	0.134465272753184\\
};
\addplot [color=black, dotted, forget plot]
  table[row sep=crcr]{%
1	0.134465272753184\\
0.108636483506987	0.134465272753184\\
};
\addplot [black,only marks, mark=diamond*, mark size = 2.0, mark options={solid, fill=mycolor1}, forget plot]
  table[row sep=crcr]{%
0	0.134465272753184\\
0.05	0.134465272753184\\
0.1	0.134465272753184\\
0.15	0.134465272753184\\
0.2	0.134465272753184\\
0.25	0.134465272753184\\
};
\addplot [black,only marks, mark=diamond*, mark size = 2.0, mark options={solid, fill=mycolor2}, forget plot]
  table[row sep=crcr]{%
0.3	0.134465272753184\\
};
\addplot [black,only marks, mark=diamond*, mark size = 2.0, mark options={solid, fill=mycolor3}, forget plot]
  table[row sep=crcr]{%
0.35	0.134465272753184\\
0.4	0.134465272753184\\
0.45 0.134465272753184\\
0.5	0.134465272753184\\
};
\addplot [black,only marks, mark=diamond*, mark size = 2.0, mark options={solid, fill=mycolor4}, forget plot]
  table[row sep=crcr]{%
0.55	0.134465272753184\\
0.6	0.134465272753184\\
0.65 0.134465272753184\\
0.7	0.134465272753184\\
0.75	0.134465272753184\\
};
\addplot [black,only marks, mark=diamond*, mark size = 2.0, mark options={solid, fill=mycolor5}, forget plot]
  table[row sep=crcr]{%
0.8	0.134465272753184\\
};
\addplot [black,only marks, mark=diamond*, mark size = 2.0, mark options={solid, fill=mycolor6}, forget plot]
  table[row sep=crcr]{%
0.85	0.134465272753184\\
0.9	0.134465272753184\\
0.95	0.134465272753184\\
1	0.134465272753184\\
};
\end{axis}
\end{tikzpicture}%

%% file: wheel_setup_2.tex
\begin{tikzpicture}
\begin{axis}[ 
    width = \textwidth,
    %height=.777\textwidth,
    ticks=none,
    axis line style={draw=none},
    ymin=-3.5,
    ymax=3.5,
    xmin=-3.5, 
    xmax=5.5,
    axis equal image
]
\filldraw[line width=2pt,fill=black!55,draw=magenta] (0,0) circle(3);
\filldraw[fill=black!10,draw=black!10] (0,0) circle(2.7);
\filldraw[line width = 1.5pt,fill=white,draw=teal] (0,0) circle(0.3);
%\draw (-.5,.5) node{\color{blue} \LARGE{$\Gamma_D$}};
%\draw (-2.6,2.6) node{\color{NatDarkGreen} \LARGE{$\Gamma_N$}};
\draw (-1,-1) node{\Large{$\mathscr{D}$}};
\draw[dashed, line width=1.5pt,black] (3.8939,1.5732) -- (0,0);
\draw[dashed, line width=1.5pt,black] (5.385,0) -- (0,0);
\draw[line width = 1.5pt] (2,0) arc (0:22:2);
\draw (1.5,0.3) node{$\omega$};
\draw[->, >= stealth,line width = 1.5pt,blue] (5,2) --node[above,midway,blue]{$F(\omega)$} (3.9866,1.6107);

%% force fun
\addplot [color=blue, line width=1.0pt, forget plot]
  table[row sep=crcr]{%
3.1	0\\
3.09975523703183	0.0389547236383931\\
3.09902098677823	0.0779032958743462\\
3.09779736518583	0.116839566276797\\
3.09608456547865	0.155757386357286\\
3.09388285812764	0.194650610540871\\
3.0911925908079	0.233513097136591\\
3.08801418834384	0.272338709307304\\
3.08434815264203	0.311121316038766\\
3.08019506261203	0.349854793107793\\
3.07555557407488	0.388533024049343\\
3.07043041965962	0.427149901122378\\
3.06482040868753	0.465699326274347\\
3.05872642704439	0.504175212104139\\
3.05214943704054	0.542571482823356\\
3.04509047725894	0.580882075215746\\
3.03755066239117	0.619100939594663\\
3.02953118306204	0.657222040758508\\
3.0210333056771	0.69523935895199\\
3.0120583737486	0.733146891190525\\
3.00260785495379	0.770938664449484\\
2.99268439278576	0.808609024103596\\
2.98230847755655	0.846157439439044\\
2.97170417510535	0.883644943674955\\
2.9627200959225	0.921650162742602\\
2.96590358199843	0.963680490919017\\
3.0195302226844	1.02322950477213\\
3.19171673096582	1.12648607617322\\
3.47967294544247	1.27751294558159\\
3.74591235358056	1.42892679308906\\
3.88226947346715	1.53709922224602\\
3.87071153493265	1.58907312265622\\
3.71346465958557	1.57933273517162\\
3.42709743164987	1.50867296220447\\
3.10891194430024	1.41550240815008\\
2.89770494651337	1.36355644532125\\
2.80833347514539	1.36486510153434\\
2.77401098489634	1.39154433132024\\
2.75356838240455	1.4248762958527\\
2.73513083438211	1.45919714506559\\
2.71655236949149	1.49343730305324\\
2.69756971141604	1.52745479953484\\
2.67816259557505	1.56123192639839\\
2.65833263425305	1.5947625547543\\
2.63808289356444	1.62804135288057\\
2.61741656905626	1.6610630644349\\
2.59633692413064	1.69382247487623\\
2.57484728750702	1.72631441111338\\
2.55295105264813	1.75853374229295\\
2.53065167722327	1.79047538060903\\
2.50795268256234	1.82213428210667\\
2.48485765309972	1.85350544747831\\
2.46137023580829	1.88458392285328\\
2.43749413962352	1.91536480058004\\
2.41323313485777	1.94584322000117\\
2.38859105260495	1.97601436822094\\
2.36357178413549	2.00587348086528\\
2.33817928028192	2.03541584283417\\
2.31241755081496	2.06463678904618\\
2.28629066381034	2.09353170517518\\
2.25980274500638	2.12209602837893\\
2.23295797715251	2.1503252480197\\
2.20576059934878	2.17821490637643\\
2.17821490637643	2.20576059934878\\
2.1503252480197	2.23295797715251\\
2.12209602837894	2.25980274500638\\
2.09353170517518	2.28629066381034\\
2.06463678904618	2.31241755081496\\
2.03541584283417	2.33817928028192\\
2.00587348086528	2.36357178413549\\
1.97601436822094	2.38859105260495\\
1.94584322000117	2.41323313485777\\
1.91536480058004	2.43749413962352\\
1.88458392285328	2.46137023580829\\
1.85350544747831	2.48485765309972\\
1.82213428210667	2.50795268256234\\
1.79047538060903	2.53065167722327\\
1.75853374229295	2.55295105264813\\
1.72631441111338	2.57484728750702\\
1.69382247487623	2.59633692413064\\
1.66106306443489	2.61741656905625\\
1.62804135288002	2.63808289356354\\
1.59476255472267	2.65833263420032\\
1.56123192507086	2.67816259329779\\
1.5274547587997	2.69756963947553\\
1.49343638971532	2.71655070813598\\
1.45918218971253	2.73510280194836\\
1.42469756792661	2.75322299132199\\
1.3899879698791	2.77090841486892\\
1.35505887661789	2.78815627985575\\
1.31991580385173	2.80496386264466\\
1.28456430107918	2.82132850912349\\
1.24900995071235	2.83724763512484\\
1.21325836719533	2.8527187268342\\
1.17731519611758	2.86773934118682\\
1.1411861133225	2.88230710625358\\
1.10487682401108	2.8964197216155\\
1.068393061841	2.91007495872701\\
1.03174058802126	2.92327066126788\\
0.994925190402349	2.93600474548371\\
0.957952682562337	2.94827520051498\\
0.920828902888808	2.96008008871459\\
0.883559713656926	2.97141754595386\\
0.846151000103709	2.98228578191687\\
0.80860866949868	2.99268308038315\\
0.77093865021105	3.00260779949876\\
0.733146890773547	3.01205837203549\\
0.695239358943082	3.02103330563839\\
0.65722204075837	3.0295311830614\\
0.619100939594662	3.03755066239117\\
0.580882075215746	3.04509047725894\\
0.542571482823356	3.05214943704054\\
0.50417521210414	3.05872642704439\\
0.465699326274347	3.06482040868753\\
0.427149901122378	3.07043041965962\\
0.388533024049343	3.07555557407488\\
0.349854793107793	3.08019506261203\\
0.311121316038766	3.08434815264203\\
0.272338709307304	3.08801418834384\\
0.233513097136592	3.0911925908079\\
0.194650610540872	3.09388285812764\\
0.155757386357286	3.09608456547865\\
0.116839566276797	3.09779736518583\\
0.0779032958743464	3.09902098677823\\
0.0389547236383936	3.09975523703183\\
1.8982025386784e-16	3.1\\
-0.0389547236383932	3.09975523703183\\
-0.077903295874346	3.09902098677823\\
-0.116839566276797	3.09779736518583\\
-0.155757386357285	3.09608456547865\\
-0.194650610540871	3.09388285812764\\
-0.233513097136592	3.0911925908079\\
-0.272338709307304	3.08801418834384\\
-0.311121316038766	3.08434815264203\\
-0.349854793107793	3.08019506261203\\
-0.388533024049343	3.07555557407488\\
-0.427149901122378	3.07043041965962\\
-0.465699326274347	3.06482040868753\\
-0.504175212104139	3.05872642704439\\
-0.542571482823355	3.05214943704054\\
-0.580882075215746	3.04509047725894\\
-0.619100939594661	3.03755066239117\\
-0.657222040758369	3.0295311830614\\
-0.695239358943082	3.02103330563839\\
-0.733146890773546	3.01205837203549\\
-0.770938650211049	3.00260779949876\\
-0.80860866949868	2.99268308038315\\
-0.846151000103708	2.98228578191687\\
-0.883559713656926	2.97141754595387\\
-0.920828902888808	2.96008008871459\\
-0.957952682562337	2.94827520051498\\
-0.994925190402349	2.93600474548371\\
-1.03174058802126	2.92327066126788\\
-1.068393061841	2.91007495872701\\
-1.10487682401108	2.8964197216155\\
-1.1411861133225	2.88230710625358\\
-1.17731519611758	2.86773934118682\\
-1.21325836719533	2.8527187268342\\
-1.24900995071235	2.83724763512484\\
-1.28456430107918	2.82132850912349\\
-1.31991580385173	2.80496386264466\\
-1.35505887661789	2.78815627985575\\
-1.3899879698791	2.77090841486892\\
-1.42469756792661	2.75322299132199\\
-1.45918218971253	2.73510280194836\\
-1.49343638971532	2.71655070813598\\
-1.5274547587997	2.69756963947553\\
-1.56123192507086	2.67816259329779\\
-1.59476255472267	2.65833263420032\\
-1.62804135288002	2.63808289356354\\
-1.66106306443489	2.61741656905625\\
-1.69382247487623	2.59633692413064\\
-1.72631441111338	2.57484728750702\\
-1.75853374229294	2.55295105264813\\
-1.79047538060903	2.53065167722327\\
-1.82213428210667	2.50795268256234\\
-1.85350544747831	2.48485765309972\\
-1.88458392285328	2.46137023580829\\
-1.91536480058004	2.43749413962352\\
-1.94584322000117	2.41323313485777\\
-1.97601436822094	2.38859105260495\\
-2.00587348086528	2.36357178413549\\
-2.03541584283416	2.33817928028192\\
-2.06463678904618	2.31241755081497\\
-2.09353170517518	2.28629066381034\\
-2.12209602837893	2.25980274500638\\
-2.1503252480197	2.23295797715251\\
-2.17821490637643	2.20576059934878\\
-2.20576059934878	2.17821490637643\\
-2.23295797715251	2.1503252480197\\
-2.25980274500638	2.12209602837893\\
-2.28629066381034	2.09353170517518\\
-2.31241755081496	2.06463678904618\\
-2.33817928028192	2.03541584283417\\
-2.36357178413549	2.00587348086528\\
-2.38859105260495	1.97601436822094\\
-2.41323313485777	1.94584322000117\\
-2.43749413962352	1.91536480058004\\
-2.46137023580829	1.88458392285328\\
-2.48485765309972	1.85350544747831\\
-2.50795268256234	1.82213428210667\\
-2.53065167722327	1.79047538060903\\
-2.55295105264813	1.75853374229295\\
-2.57484728750702	1.72631441111338\\
-2.59633692413064	1.69382247487623\\
-2.61741656905625	1.66106306443489\\
-2.63808289356354	1.62804135288002\\
-2.65833263420032	1.59476255472267\\
-2.67816259329779	1.56123192507086\\
-2.69756963947553	1.5274547587997\\
-2.71655070813598	1.49343638971532\\
-2.73510280194836	1.45918218971253\\
-2.75322299132199	1.42469756792661\\
-2.77090841486892	1.3899879698791\\
-2.78815627985575	1.35505887661789\\
-2.80496386264466	1.31991580385173\\
-2.82132850912349	1.28456430107918\\
-2.83724763512484	1.24900995071235\\
-2.8527187268342	1.21325836719533\\
-2.86773934118682	1.17731519611758\\
-2.88230710625358	1.1411861133225\\
-2.8964197216155	1.10487682401108\\
-2.91007495872701	1.068393061841\\
-2.92327066126788	1.03174058802126\\
-2.93600474548371	0.99492519040235\\
-2.94827520051498	0.957952682562337\\
-2.96008008871459	0.920828902888808\\
-2.97141754595386	0.883559713656927\\
-2.98228578191687	0.846151000103709\\
-2.99268308038315	0.808608669498681\\
-3.00260779949876	0.77093865021105\\
-3.01205837203549	0.733146890773546\\
-3.02103330563839	0.695239358943081\\
-3.0295311830614	0.65722204075837\\
-3.03755066239117	0.619100939594662\\
-3.04509047725894	0.580882075215746\\
-3.05214943704054	0.542571482823357\\
-3.05872642704439	0.50417521210414\\
-3.06482040868753	0.465699326274347\\
-3.07043041965962	0.427149901122379\\
-3.07555557407488	0.388533024049344\\
-3.08019506261203	0.349854793107792\\
-3.08434815264203	0.311121316038766\\
-3.08801418834384	0.272338709307304\\
-3.0911925908079	0.233513097136591\\
-3.09388285812764	0.194650610540872\\
-3.09608456547865	0.155757386357286\\
-3.09779736518583	0.116839566276797\\
-3.09902098677823	0.0779032958743472\\
-3.09975523703183	0.0389547236383938\\
-3.1	3.79640507735679e-16\\
-3.09975523703183	-0.038954723638393\\
-3.09902098677823	-0.0779032958743465\\
-3.09779736518583	-0.116839566276798\\
-3.09608456547865	-0.155757386357285\\
-3.09388285812764	-0.194650610540871\\
-3.0911925908079	-0.233513097136592\\
-3.08801418834384	-0.272338709307303\\
-3.08434815264203	-0.311121316038766\\
-3.08019506261203	-0.349854793107793\\
-3.07555557407488	-0.388533024049342\\
-3.07043041965962	-0.427149901122377\\
-3.06482040868753	-0.465699326274348\\
-3.05872642704439	-0.504175212104139\\
-3.05214943704054	-0.542571482823356\\
-3.04509047725893	-0.580882075215747\\
-3.03755066239117	-0.619100939594661\\
-3.0295311830614	-0.657222040758369\\
-3.02103330563839	-0.695239358943082\\
-3.01205837203549	-0.733146890773546\\
-3.00260779949876	-0.770938650211049\\
-2.99268308038315	-0.80860866949868\\
-2.98228578191687	-0.846151000103708\\
-2.97141754595386	-0.883559713656926\\
-2.96008008871459	-0.920828902888809\\
-2.94827520051498	-0.957952682562337\\
-2.93600474548371	-0.994925190402349\\
-2.92327066126788	-1.03174058802126\\
-2.91007495872701	-1.068393061841\\
-2.8964197216155	-1.10487682401108\\
-2.88230710625358	-1.1411861133225\\
-2.86773934118682	-1.17731519611758\\
-2.8527187268342	-1.21325836719533\\
-2.83724763512484	-1.24900995071236\\
-2.82132850912349	-1.28456430107918\\
-2.80496386264466	-1.31991580385173\\
-2.78815627985575	-1.35505887661789\\
-2.77090841486892	-1.3899879698791\\
-2.75322299132199	-1.42469756792661\\
-2.73510280194836	-1.45918218971253\\
-2.71655070813598	-1.49343638971532\\
-2.69756963947553	-1.5274547587997\\
-2.67816259329779	-1.56123192507086\\
-2.65833263420032	-1.59476255472267\\
-2.63808289356354	-1.62804135288002\\
-2.61741656905625	-1.66106306443489\\
-2.59633692413064	-1.69382247487623\\
-2.57484728750702	-1.72631441111338\\
-2.55295105264812	-1.75853374229295\\
-2.53065167722327	-1.79047538060903\\
-2.50795268256234	-1.82213428210667\\
-2.48485765309972	-1.85350544747831\\
-2.46137023580829	-1.88458392285328\\
-2.43749413962352	-1.91536480058004\\
-2.41323313485777	-1.94584322000117\\
-2.38859105260495	-1.97601436822094\\
-2.36357178413549	-2.00587348086528\\
-2.33817928028192	-2.03541584283417\\
-2.31241755081497	-2.06463678904618\\
-2.28629066381034	-2.09353170517518\\
-2.25980274500638	-2.12209602837894\\
-2.23295797715251	-2.15032524801969\\
-2.20576059934879	-2.17821490637643\\
-2.17821490637643	-2.20576059934878\\
-2.1503252480197	-2.23295797715251\\
-2.12209602837893	-2.25980274500638\\
-2.09353170517518	-2.28629066381034\\
-2.06463678904618	-2.31241755081496\\
-2.03541584283417	-2.33817928028192\\
-2.00587348086528	-2.36357178413549\\
-1.97601436822094	-2.38859105260495\\
-1.94584322000117	-2.41323313485777\\
-1.91536480058004	-2.43749413962352\\
-1.88458392285328	-2.46137023580829\\
-1.85350544747831	-2.48485765309972\\
-1.82213428210667	-2.50795268256234\\
-1.79047538060903	-2.53065167722327\\
-1.75853374229295	-2.55295105264813\\
-1.72631441111338	-2.57484728750702\\
-1.69382247487623	-2.59633692413064\\
-1.66106306443489	-2.61741656905625\\
-1.62804135288002	-2.63808289356354\\
-1.59476255472267	-2.65833263420032\\
-1.56123192507086	-2.67816259329779\\
-1.5274547587997	-2.69756963947553\\
-1.49343638971532	-2.71655070813598\\
-1.45918218971253	-2.73510280194836\\
-1.42469756792661	-2.75322299132199\\
-1.3899879698791	-2.77090841486892\\
-1.35505887661789	-2.78815627985575\\
-1.31991580385173	-2.80496386264466\\
-1.28456430107918	-2.82132850912349\\
-1.24900995071235	-2.83724763512484\\
-1.21325836719533	-2.8527187268342\\
-1.17731519611758	-2.86773934118682\\
-1.1411861133225	-2.88230710625358\\
-1.10487682401108	-2.8964197216155\\
-1.068393061841	-2.91007495872701\\
-1.03174058802126	-2.92327066126788\\
-0.99492519040235	-2.93600474548371\\
-0.957952682562337	-2.94827520051498\\
-0.920828902888808	-2.96008008871459\\
-0.883559713656929	-2.97141754595386\\
-0.846151000103708	-2.98228578191687\\
-0.80860866949868	-2.99268308038315\\
-0.770938650211051	-3.00260779949876\\
-0.733146890773545	-3.01205837203549\\
-0.695239358943083	-3.02103330563839\\
-0.65722204075837	-3.0295311830614\\
-0.619100939594662	-3.03755066239117\\
-0.580882075215746	-3.04509047725894\\
-0.542571482823358	-3.05214943704054\\
-0.504175212104139	-3.05872642704439\\
-0.465699326274346	-3.06482040868753\\
-0.427149901122379	-3.07043041965962\\
-0.388533024049342	-3.07555557407488\\
-0.349854793107794	-3.08019506261203\\
-0.311121316038767	-3.08434815264203\\
-0.272338709307304	-3.08801418834384\\
-0.233513097136591	-3.0911925908079\\
-0.194650610540874	-3.09388285812764\\
-0.155757386357285	-3.09608456547865\\
-0.116839566276799	-3.09779736518583\\
-0.0779032958743474	-3.09902098677823\\
-0.0389547236383912	-3.09975523703183\\
-5.69460761603519e-16	-3.1\\
0.0389547236383928	-3.09975523703183\\
0.0779032958743463	-3.09902098677823\\
0.116839566276797	-3.09779736518583\\
0.155757386357284	-3.09608456547865\\
0.194650610540873	-3.09388285812764\\
0.23351309713659	-3.0911925908079\\
0.272338709307303	-3.08801418834384\\
0.311121316038765	-3.08434815264203\\
0.349854793107793	-3.08019506261203\\
0.388533024049343	-3.07555557407488\\
0.427149901122378	-3.07043041965962\\
0.465699326274348	-3.06482040868753\\
0.504175212104137	-3.05872642704439\\
0.542571482823357	-3.05214943704054\\
0.580882075215745	-3.04509047725894\\
0.619100939594661	-3.03755066239117\\
0.657222040758369	-3.0295311830614\\
0.695239358943082	-3.02103330563839\\
0.733146890773544	-3.01205837203549\\
0.77093865021105	-3.00260779949876\\
0.808608669498681	-2.99268308038315\\
0.846151000103707	-2.98228578191687\\
0.883559713656928	-2.97141754595386\\
0.920828902888807	-2.96008008871459\\
0.957952682562336	-2.94827520051498\\
0.994925190402349	-2.93600474548371\\
1.03174058802126	-2.92327066126788\\
1.068393061841	-2.91007495872701\\
1.10487682401108	-2.8964197216155\\
1.1411861133225	-2.88230710625358\\
1.17731519611758	-2.86773934118682\\
1.21325836719533	-2.8527187268342\\
1.24900995071235	-2.83724763512484\\
1.28456430107918	-2.82132850912349\\
1.31991580385173	-2.80496386264466\\
1.35505887661789	-2.78815627985575\\
1.3899879698791	-2.77090841486892\\
1.42469756792661	-2.75322299132199\\
1.45918218971253	-2.73510280194836\\
1.49343638971532	-2.71655070813598\\
1.52745475879971	-2.69756963947553\\
1.56123192507086	-2.67816259329779\\
1.59476255472267	-2.65833263420032\\
1.62804135288002	-2.63808289356354\\
1.66106306443489	-2.61741656905625\\
1.69382247487623	-2.59633692413064\\
1.72631441111338	-2.57484728750702\\
1.75853374229294	-2.55295105264813\\
1.79047538060903	-2.53065167722327\\
1.82213428210667	-2.50795268256234\\
1.85350544747831	-2.48485765309972\\
1.88458392285328	-2.46137023580829\\
1.91536480058004	-2.43749413962352\\
1.94584322000117	-2.41323313485777\\
1.97601436822094	-2.38859105260495\\
2.00587348086528	-2.36357178413549\\
2.03541584283416	-2.33817928028192\\
2.06463678904618	-2.31241755081497\\
2.09353170517518	-2.28629066381034\\
2.12209602837893	-2.25980274500638\\
2.15032524801969	-2.23295797715251\\
2.17821490637643	-2.20576059934878\\
2.20576059934879	-2.17821490637643\\
2.23295797715251	-2.1503252480197\\
2.25980274500638	-2.12209602837893\\
2.28629066381034	-2.09353170517518\\
2.31241755081496	-2.06463678904618\\
2.33817928028192	-2.03541584283417\\
2.36357178413549	-2.00587348086528\\
2.38859105260495	-1.97601436822094\\
2.41323313485777	-1.94584322000117\\
2.43749413962352	-1.91536480058004\\
2.46137023580829	-1.88458392285328\\
2.48485765309972	-1.85350544747831\\
2.50795268256234	-1.82213428210667\\
2.53065167722327	-1.79047538060903\\
2.55295105264813	-1.75853374229295\\
2.57484728750702	-1.72631441111338\\
2.59633692413064	-1.69382247487624\\
2.61741656905625	-1.66106306443489\\
2.63808289356354	-1.62804135288002\\
2.65833263420032	-1.59476255472267\\
2.67816259329779	-1.56123192507086\\
2.69756963947553	-1.5274547587997\\
2.71655070813598	-1.49343638971532\\
2.73510280194836	-1.45918218971253\\
2.75322299132199	-1.42469756792661\\
2.77090841486892	-1.3899879698791\\
2.78815627985575	-1.35505887661789\\
2.80496386264466	-1.31991580385173\\
2.82132850912349	-1.28456430107918\\
2.83724763512484	-1.24900995071235\\
2.8527187268342	-1.21325836719533\\
2.86773934118682	-1.17731519611759\\
2.88230710625358	-1.1411861133225\\
2.8964197216155	-1.10487682401108\\
2.91007495872701	-1.068393061841\\
2.92327066126788	-1.03174058802126\\
2.93600474548371	-0.994925190402351\\
2.94827520051498	-0.957952682562338\\
2.96008008871459	-0.920828902888809\\
2.97141754595386	-0.883559713656926\\
2.98228578191687	-0.846151000103711\\
2.99268308038315	-0.80860866949868\\
3.00260779949876	-0.770938650211052\\
3.01205837203549	-0.733146890773548\\
3.02103330563839	-0.69523935894308\\
3.0295311830614	-0.65722204075837\\
3.03755066239117	-0.619100939594662\\
3.04509047725894	-0.580882075215747\\
3.05214943704054	-0.542571482823356\\
3.05872642704439	-0.504175212104141\\
3.06482040868753	-0.465699326274346\\
3.07043041965962	-0.42714990112238\\
3.07555557407488	-0.388533024049344\\
3.08019506261203	-0.349854793107791\\
3.08434815264203	-0.311121316038767\\
3.08801418834384	-0.272338709307304\\
3.0911925908079	-0.233513097136591\\
3.09388285812764	-0.194650610540871\\
3.09608456547865	-0.155757386357288\\
3.09779736518583	-0.116839566276796\\
3.09902098677823	-0.0779032958743476\\
3.09975523703183	-0.0389547236383942\\
3.1	-7.59281015471359e-16\\
};
\end{axis}
\end{tikzpicture}

%% file: wheel_pareto.tex
\definecolor{mycolor1}{rgb}{0.00000,0.44700,0.74100}%
\definecolor{mycolor2}{rgb}{0.85000,0.32500,0.09800}%
% \definecolor{mycolor3}{rgb}{0.92900,0.69400,0.12500}%
% \definecolor{mycolor4}{rgb}{0.49400,0.18400,0.55600}%
% \definecolor{mycolor5}{rgb}{0.46600,0.67400,0.18800}%
% \definecolor{mycolor6}{rgb}{0.30100,0.74500,0.93300}%
% \definecolor{mycolor7}{rgb}{0.63500,0.07800,0.18400}%
%
\begin{tikzpicture}

\begin{axis}[%
width=\textwidth,
height=.75\textwidth,
at={(1.517in,0.962in)},
%scale only axis,
xmin=0,%0
xmax=1.2,%1.2
xlabel style={font=\color{white!15!black}},
xlabel={Chance constraint value},
ymin=0.1,%0.1
ymax=0.5,%0.5
ylabel style={font=\color{white!15!black}},
ylabel={$\pvol$},
xmajorgrids,
xminorgrids,
ymajorgrids,
yminorgrids,
legend style={legend cell align=left, align=left, draw=white!15!black},
axis background/.style={fill=white}
]

\addplot[dashed,forget plot] table[row sep=crcr]{%
x   y\\
0.025 0.1\\
0.025 0.5\\
};

\addplot[only marks, mark=*, mark options={}, mark size=1.2pt, color=mycolor1, fill=mycolor1] table[row sep=crcr]{%
x	y\\
0.0247840261148748	0.473020674728141\\
0.0255529933059216	0.45731504469348\\
0.0255806903160329	0.474355995938482\\
0.0257791162402105	0.471709486307729\\
0.0254273623830777	0.472211656813625\\
};
\addlegendentry{sMMA};

\addplot[only marks, mark=triangle*, mark options={}, mark size=1.5pt, color=mycolor2, fill=mycolor2] table[row sep=crcr]{%
x	y\\
1.10432823066256	0.136781278344474\\
0.838808802047539	0.262969091162073\\
0.490136209101736	0.423705376807686\\
0.0251878531556103	0.47824990594896\\
0.0250281230019199	0.477878891251179\\
};
\addlegendentry{MMA};

\addplot[only marks, mark=*, mark options={}, mark size=1.2pt, color=mycolor1, fill=mycolor1, forget plot] table[row sep=crcr]{%
x	y\\
0.0243599322343047	0.464628110485212\\
0.0258068407011863	0.464073302598008\\
0.0245814649510011	0.477515773444111\\
0.0253702917785468	0.464321474083927\\
0.0252741445973414	0.466773208489206\\
};
\addplot[only marks, mark=triangle*, mark options={}, mark size=1.5pt, color=mycolor2, fill=mycolor2, forget plot] table[row sep=crcr]{%
x	y\\
1.10434323745219	0.13678291064281\\
0.840699057279076	0.262338164899652\\
0.478525152779	0.418234790118707\\
0.0252463066588944	0.474127366173422\\
0.0250321997137167	0.476729867755772\\
};
\addplot[only marks, mark=*, mark options={}, mark size=1.2pt, color=mycolor1, fill=mycolor1, forget plot] table[row sep=crcr]{%
x	y\\
0.0250206256821411	0.465782079088995\\
0.0251563947116631	0.481055950346529\\
0.024753546527105	0.474656965949866\\
0.025072877567134	0.468792756947\\
0.0254152397019789	0.467762111571891\\
};
\addplot[only marks, mark=triangle*, mark options={}, mark size=1.5pt, color=mycolor2, fill=mycolor2, forget plot] table[row sep=crcr]{%
x	y\\
1.10427438360777	0.136778958240685\\
0.839804655898168	0.262354323283541\\
0.452259844315991	0.418265830998063\\
0.0250081703988064	0.475463467807045\\
0.0250082366042228	0.478775300297812\\
};
\addplot[only marks, mark=*, mark options={}, mark size=1.2pt, color=mycolor1, fill=mycolor1, forget plot] table[row sep=crcr]{%
x	y\\
0.0249361483567161	0.478031668688297\\
0.0248202333459136	0.497055101561615\\
0.0248915070600666	0.469659288146445\\
0.0254888046491379	0.47084898908096\\
0.0244935536677588	0.481381488794859\\
};
\addplot[only marks, mark=triangle*, mark options={}, mark size=1.5pt, color=mycolor2, fill=mycolor2, forget plot] table[row sep=crcr]{%
x	y\\
1.10419223179461	0.136773530630534\\
0.840156256820081	0.262389191853929\\
0.471928478464833	0.42145607473311\\
0.0249934481031141	0.472745339668167\\
0.0250312399467903	0.474901453032746\\
};
\end{axis}
\end{tikzpicture}%

%% file: wheel_pareto_zoom.tex
\definecolor{mycolor1}{rgb}{0.00000,0.44700,0.74100}%
\definecolor{mycolor2}{rgb}{0.85000,0.32500,0.09800}%
% \definecolor{mycolor3}{rgb}{0.92900,0.69400,0.12500}%
% \definecolor{mycolor4}{rgb}{0.49400,0.18400,0.55600}%
% \definecolor{mycolor5}{rgb}{0.46600,0.67400,0.18800}%
% \definecolor{mycolor6}{rgb}{0.30100,0.74500,0.93300}%
% \definecolor{mycolor7}{rgb}{0.63500,0.07800,0.18400}%
%
\begin{tikzpicture}

\begin{axis}[%
width=\textwidth,
height=.75\textwidth,
at={(1.517in,0.962in)},
%scale only axis,
xmin=0,%0
xmax=0.1,%1.2
xlabel style={font=\color{white!15!black}},
xlabel={Chance constraint value},
ymin=0.45,%0.1
ymax=0.5,%0.5
ylabel style={font=\color{white!15!black}},
ylabel={$\pvol$},
xmajorgrids,
xminorgrids,
ymajorgrids,
yminorgrids,
xtick={0,0.02,0.04,0.06,0.08,0.1},
xticklabels={0,0.02,0.04,0.06,0.08,0.1},
legend style={legend cell align=left, align=left, draw=white!15!black},
axis background/.style={fill=white}
]

\addplot[dashed,forget plot] table[row sep=crcr]{%
x   y\\
0.025 0.1\\
0.025 0.5\\
};

\addplot[only marks, mark=*, mark options={}, mark size=1.2pt, color=mycolor1, fill=mycolor1] table[row sep=crcr]{%
x	y\\
0.0247840261148748	0.473020674728141\\
0.0255529933059216	0.45731504469348\\
0.0255806903160329	0.474355995938482\\
0.0257791162402105	0.471709486307729\\
0.0254273623830777	0.472211656813625\\
};
\addlegendentry{sMMA};

\addplot[only marks, mark=triangle*, mark options={}, mark size=1.5pt, color=mycolor2, fill=mycolor2] table[row sep=crcr]{%
x	y\\
1.10432823066256	0.136781278344474\\
0.838808802047539	0.262969091162073\\
0.490136209101736	0.423705376807686\\
0.0251878531556103	0.47824990594896\\
0.0250281230019199	0.477878891251179\\
};
\addlegendentry{MMA};

\addplot[only marks, mark=*, mark options={}, mark size=1.2pt, color=mycolor1, fill=mycolor1, forget plot] table[row sep=crcr]{%
x	y\\
0.0243599322343047	0.464628110485212\\
0.0258068407011863	0.464073302598008\\
0.0245814649510011	0.477515773444111\\
0.0253702917785468	0.464321474083927\\
0.0252741445973414	0.466773208489206\\
};
\addplot[only marks, mark=triangle*, mark options={}, mark size=1.5pt, color=mycolor2, fill=mycolor2, forget plot] table[row sep=crcr]{%
x	y\\
1.10434323745219	0.13678291064281\\
0.840699057279076	0.262338164899652\\
0.478525152779	0.418234790118707\\
0.0252463066588944	0.474127366173422\\
0.0250321997137167	0.476729867755772\\
};
\addplot[only marks, mark=*, mark options={}, mark size=1.2pt, color=mycolor1, fill=mycolor1, forget plot] table[row sep=crcr]{%
x	y\\
0.0250206256821411	0.465782079088995\\
0.0251563947116631	0.481055950346529\\
0.024753546527105	0.474656965949866\\
0.025072877567134	0.468792756947\\
0.0254152397019789	0.467762111571891\\
};
\addplot[only marks, mark=triangle*, mark options={}, mark size=1.5pt, color=mycolor2, fill=mycolor2, forget plot] table[row sep=crcr]{%
x	y\\
1.10427438360777	0.136778958240685\\
0.839804655898168	0.262354323283541\\
0.452259844315991	0.418265830998063\\
0.0250081703988064	0.475463467807045\\
0.0250082366042228	0.478775300297812\\
};
\addplot[only marks, mark=*, mark options={}, mark size=1.2pt, color=mycolor1, fill=mycolor1, forget plot] table[row sep=crcr]{%
x	y\\
0.0249361483567161	0.478031668688297\\
0.0248202333459136	0.497055101561615\\
0.0248915070600666	0.469659288146445\\
0.0254888046491379	0.47084898908096\\
0.0244935536677588	0.481381488794859\\
};
\addplot[only marks, mark=triangle*, mark options={}, mark size=1.5pt, color=mycolor2, fill=mycolor2, forget plot] table[row sep=crcr]{%
x	y\\
1.10419223179461	0.136773530630534\\
0.840156256820081	0.262389191853929\\
0.471928478464833	0.42145607473311\\
0.0249934481031141	0.472745339668167\\
0.0250312399467903	0.474901453032746\\
};
\end{axis}
\end{tikzpicture}%

%% file: mcmsa_16_cc.tex
\definecolor{mycolor1}{rgb}{0.00000,0.44700,0.74100}%
\definecolor{mycolor2}{rgb}{0.85000,0.32500,0.09800}%
\definecolor{mycolor3}{rgb}{0.49400,0.18400,0.55600}%
\begin{tikzpicture}

\begin{axis}[%
width=\textwidth,
height=.75\textwidth,
at={(0.758in,0.481in)},
%scale only axis,
xmin=0,
xmax=400,
ymin=0,
ymax=1,
axis background/.style={fill=white},
title style={font=\bfseries},
title={sMMA, $\mathcal{B}=16$},
every axis title/.style={at={(0,1)}, anchor=north west, draw=black, fill=white},
xminorgrids,
yminorgrids,
xmajorgrids,
ymajorgrids,
xlabel={\vphantom{Iteration}},
ytick={0.025,0.2,0.4,0.6,0.8,1},
yticklabels={$p$,0.2,0.4,0.6,0.8,1},
legend style={legend cell align=left, align=left, draw=white!15!black}
]
\addplot [color=mycolor1, line width=1.5pt,line join=round]
  table[row sep=crcr]{%
1	1.00332906196697\\
2	0.867400565076876\\
3	0.456897448131919\\
4	0.0274465257467138\\
5	0.00164380206016538\\
6	0.119717376619768\\
7	0.0367810233688647\\
8	0.0137983107491279\\
9	0.0203655084946325\\
10	0.010528021846066\\
11	0.0533048597990082\\
12	0.020052886450357\\
13	0.023369981256629\\
14	0.0226736163772046\\
15	0.0241806445525247\\
16	0.0260451560476189\\
17	0.0279067253863797\\
18	0.0278447328767913\\
19	0.0270805898939073\\
20	0.0271204464012097\\
21	0.0260608123494101\\
22	0.0242979841671138\\
23	0.0234491547462394\\
24	0.0224357507189816\\
25	0.0241780768029679\\
26	0.024389960720692\\
27	0.0242785547794594\\
28	0.0248033934709488\\
29	0.0252856801855033\\
30	0.0256545655056504\\
31	0.0255925398184676\\
32	0.0250961380819776\\
33	0.0244291051094407\\
34	0.0233439234780474\\
35	0.0208210082409589\\
36	0.0176167978714286\\
37	0.0151845724713245\\
38	0.0137941620920443\\
39	0.0153174753638683\\
40	0.0187049124998194\\
41	0.0243197555166458\\
42	0.0274105466995123\\
43	0.0277157728412415\\
44	0.027744023453068\\
45	0.0248632805688524\\
46	0.0227317527829796\\
47	0.0234175431087903\\
48	0.0238316584975218\\
49	0.0235127814551839\\
50	0.0238629587470109\\
51	0.024248476433464\\
52	0.0252744603959931\\
53	0.0249373053808842\\
54	0.0251388749866212\\
55	0.0252260705786018\\
56	0.0252086144363254\\
57	0.0249205665198857\\
58	0.0243539468662982\\
59	0.0240552733122405\\
60	0.0241769440726828\\
61	0.0244157011446518\\
62	0.0244762769346296\\
63	0.0246773835527369\\
64	0.0248306139702383\\
65	0.0250445843942278\\
66	0.0250050406413273\\
67	0.0250820108909692\\
68	0.0250099253291523\\
69	0.0249608141757827\\
70	0.0248499340396589\\
71	0.0247486911088107\\
72	0.0246690232193388\\
73	0.0246466256650915\\
74	0.024594391689508\\
75	0.0245983120365988\\
76	0.0247055490144859\\
77	0.0248358360831651\\
78	0.0248426285228879\\
79	0.0248920895207793\\
80	0.0249939987192069\\
81	0.0250585537553264\\
82	0.0249876335849816\\
83	0.0249115700570205\\
84	0.0247993332417612\\
85	0.0247435777180359\\
86	0.0246519266310042\\
87	0.0245861685637258\\
88	0.0246516474784186\\
89	0.0247395817448208\\
90	0.0247780313823073\\
91	0.0248489821208802\\
92	0.0249113616218679\\
93	0.0250436971948379\\
94	0.0250888346796541\\
95	0.0250637171925279\\
96	0.0250955084921438\\
97	0.0250904958501722\\
98	0.0250121862118932\\
99	0.0249297751695269\\
100	0.0248416837062121\\
101	0.0247986619385992\\
102	0.0247117251003454\\
103	0.0246374604256708\\
104	0.0246513741682115\\
105	0.0247471355427747\\
106	0.0247933202322811\\
107	0.0248740972938855\\
108	0.024932276315439\\
109	0.0250175535413596\\
110	0.0250487404343151\\
111	0.0250441377920791\\
112	0.0250951101273501\\
113	0.0251191541491672\\
114	0.0250334298689122\\
115	0.0249599819803426\\
116	0.0248647301092347\\
117	0.0248174358225933\\
118	0.0247416601440371\\
119	0.0247473096569534\\
120	0.024784080153641\\
121	0.0248485543248742\\
122	0.024814205079955\\
123	0.0248477316050759\\
124	0.0249027314790106\\
125	0.0249995164277623\\
126	0.0250551608292549\\
127	0.0250493010136797\\
128	0.0250890649689267\\
129	0.0251154981078451\\
130	0.0250587021134347\\
131	0.0250036048146331\\
132	0.0249220797257793\\
133	0.0249333001756123\\
134	0.0248551667797495\\
135	0.0247732675473489\\
136	0.0247626957394904\\
137	0.024816499525487\\
138	0.0248290932595875\\
139	0.0248860986918849\\
140	0.0249274213112937\\
141	0.0250201454109723\\
142	0.0250504056333226\\
143	0.0250310115076874\\
144	0.02508261970671\\
145	0.0250946289467106\\
146	0.0250597237299627\\
147	0.0250196709877763\\
148	0.0250008033145716\\
149	0.0249734997176342\\
150	0.0249322082639901\\
151	0.0248785379867519\\
152	0.0248765021982212\\
153	0.024902967743661\\
154	0.0248771237973795\\
155	0.0248855755630389\\
156	0.0248759209062156\\
157	0.0249049225479232\\
158	0.0248921547798091\\
159	0.0248846900994771\\
160	0.0249361234974984\\
161	0.0249780255941671\\
162	0.0249795095604343\\
163	0.0250163387586605\\
164	0.0250306699527285\\
165	0.025089128151488\\
166	0.0250795366872316\\
167	0.0250730092151195\\
168	0.0250686693601547\\
169	0.025054192778038\\
170	0.0250051521606692\\
171	0.0249505587234862\\
172	0.024926135567287\\
173	0.0248851258274146\\
174	0.0248441919832262\\
175	0.0247980794634921\\
176	0.0248188927444286\\
177	0.024843988260067\\
178	0.0248455482113145\\
179	0.0248668684091376\\
180	0.0249208259279639\\
181	0.0249819081894844\\
182	0.0250171002596042\\
183	0.0250299661613628\\
184	0.0251001886088286\\
185	0.0251296989338178\\
186	0.0251002812697563\\
187	0.0250951750914631\\
188	0.0250697272764381\\
189	0.0250617884580193\\
190	0.0249995908999215\\
191	0.024943657907779\\
192	0.0249489500794256\\
193	0.0249178964314686\\
194	0.0248952645713768\\
195	0.0248987556917872\\
196	0.0249077595026637\\
197	0.0248630295090137\\
198	0.0248404914656179\\
199	0.0248198687962403\\
200	0.0261359582463627\\
201	0.0263168342077945\\
202	0.0267854755799052\\
203	0.0271803346708307\\
204	0.027103370140116\\
205	0.0268923129190313\\
206	0.0262392253016462\\
207	0.0252925293425264\\
208	0.024292968341514\\
209	0.0230855924257187\\
210	0.023099420169079\\
211	0.02310451323972\\
212	0.0245264649172139\\
213	0.0250773563137764\\
214	0.025151808319586\\
215	0.025396890174923\\
216	0.0254761940426033\\
217	0.0251185494996559\\
218	0.0249108699569536\\
219	0.02477687633812\\
220	0.0247051446095232\\
221	0.0249624992261981\\
222	0.024904248337762\\
223	0.0249589415625565\\
224	0.0250134349202057\\
225	0.0253585033946666\\
226	0.0254003802261857\\
227	0.0254001374528597\\
228	0.0252856494936078\\
229	0.025119382972424\\
230	0.0249540498180583\\
231	0.0248092312028504\\
232	0.0246733321281778\\
233	0.024627599070881\\
234	0.0246048053849339\\
235	0.0246190291931286\\
236	0.0246800968132847\\
237	0.0247528007320262\\
238	0.0248188297189008\\
239	0.0248765077188444\\
240	0.0249805564662186\\
241	0.0250850403559466\\
242	0.0251336603514576\\
243	0.0251879301097872\\
244	0.0252504397103696\\
245	0.0252995060575691\\
246	0.0253144090336319\\
247	0.0252491989358129\\
248	0.0251498575285519\\
249	0.025110132177358\\
250	0.0250764427863509\\
251	0.0250187777820851\\
252	0.0249324548476615\\
253	0.0248553616060804\\
254	0.0247825292651541\\
255	0.0246945541281058\\
256	0.024668445707191\\
257	0.0246489387272897\\
258	0.0247052978354739\\
259	0.0248168782315327\\
260	0.0249090743163089\\
261	0.0250401248374724\\
262	0.0251340640630421\\
263	0.0252222747405692\\
264	0.0253232380501659\\
265	0.0252478409018536\\
266	0.0251718179829886\\
267	0.0250816785129461\\
268	0.0249989512267859\\
269	0.0249119725463919\\
270	0.0248509959932337\\
271	0.0248180870488934\\
272	0.0248290838431951\\
273	0.0248373654156204\\
274	0.0248388117347435\\
275	0.0248580208614705\\
276	0.0248753389206246\\
277	0.024903744976031\\
278	0.0249159193763768\\
279	0.0249385288179053\\
280	0.0249718187598078\\
281	0.024996208220123\\
282	0.0250111752550893\\
283	0.0250374043045681\\
284	0.0250808753851731\\
285	0.0251207174519103\\
286	0.025122855400589\\
287	0.0251415363437803\\
288	0.0251730090423838\\
289	0.0251616252163165\\
290	0.0251293296545646\\
291	0.0250732194498832\\
292	0.0250142545081447\\
293	0.0249324491355196\\
294	0.0248601163502456\\
295	0.0248012686427045\\
296	0.0247396688289343\\
297	0.0247289569775684\\
298	0.0247226674446811\\
299	0.0247431116950528\\
300	0.0248591469845545\\
301	0.0249333502621531\\
302	0.0249758232796457\\
303	0.0250269865415377\\
304	0.0251242577543989\\
305	0.0251576469832471\\
306	0.0251522325474181\\
307	0.025167043289291\\
308	0.0251562060308003\\
309	0.0251445823796241\\
310	0.0251130574731901\\
311	0.0250863461181801\\
312	0.0250433311073052\\
313	0.0250342089303923\\
314	0.0250141272494487\\
315	0.0249780620430087\\
316	0.0249380094489313\\
317	0.0248960573668299\\
318	0.0248618704139187\\
319	0.0248203278309107\\
320	0.0247999801023612\\
321	0.0248025848464226\\
322	0.0248126912594545\\
323	0.0248518481136972\\
324	0.0248837952222801\\
325	0.0249519942839831\\
326	0.0249928569332435\\
327	0.0250445636456845\\
328	0.0251001335648907\\
329	0.0251241016094371\\
330	0.0251211857422295\\
331	0.0251164966028269\\
332	0.0251126527080273\\
333	0.0251022004916828\\
334	0.0250843799816087\\
335	0.0250611496865552\\
336	0.0250430730075937\\
337	0.0250090087957634\\
338	0.0249761593153864\\
339	0.0249474302355938\\
340	0.0249140057718254\\
341	0.0248913726146256\\
342	0.0248716158166777\\
343	0.0248647345162845\\
344	0.0248649696813768\\
345	0.024863767013362\\
346	0.0248664938162278\\
347	0.024878219604386\\
348	0.0249106299365663\\
349	0.0249512432871476\\
350	0.0249847329007482\\
351	0.0250242440276376\\
352	0.0250751929995616\\
353	0.0251274223054084\\
354	0.0251460408881006\\
355	0.0251709532876716\\
356	0.0251685619839807\\
357	0.0251639050306992\\
358	0.0251352477485805\\
359	0.0250579975029085\\
360	0.0249997323551179\\
361	0.0249451329316488\\
362	0.0248970690888613\\
363	0.0248612676322185\\
364	0.0248242149389488\\
365	0.0247920742405196\\
366	0.0247724051442857\\
367	0.0247688082869027\\
368	0.0248424660580856\\
369	0.0249140745795478\\
370	0.0249421362317895\\
371	0.0249870770895306\\
372	0.0250246217554467\\
373	0.0250729220479118\\
374	0.0250935216256072\\
375	0.0250914331664256\\
376	0.0250919709715343\\
377	0.0251052199286053\\
378	0.0251263331564036\\
379	0.0251440386994303\\
380	0.0251591687023975\\
381	0.0251357734718603\\
382	0.0251077433216968\\
383	0.0250553444168894\\
384	0.0250247971952255\\
385	0.0249681284580506\\
386	0.0249595266463999\\
387	0.0249419373531285\\
388	0.024920272539579\\
389	0.0248694460742992\\
390	0.0248431021152485\\
391	0.0248315837087407\\
392	0.0248938069618949\\
393	0.0248836081900529\\
394	0.0248773924233975\\
395	0.0248794519001776\\
396	0.0248904813619335\\
397	0.0249074531833906\\
398	0.0249219140894006\\
399	0.0249346720267004\\
400	0.0249713356348633\\
};
\addlegendentry{internal}

\addplot [color=mycolor2, line width=1pt,line join = round]
  table[row sep=crcr]{%
1	1.00331519624177\\
2	0.817277239280557\\
3	0.0175739347990537\\
4	0.00191584040642272\\
5	0.00124893130820033\\
6	0.0708163371202274\\
7	0.0093743677773168\\
8	0.00539448610246812\\
9	0.0101341423346233\\
10	0.0124134164456548\\
11	0.0392843662872431\\
12	0.0195125906587405\\
13	0.0244905544798148\\
14	0.0261172127361559\\
15	0.028977638848721\\
16	0.029983232719482\\
17	0.0288485923745001\\
18	0.02614007534322\\
19	0.0238415828438225\\
20	0.0222523852293939\\
21	0.0206327650018054\\
22	0.0198200598739183\\
23	0.0204160445466995\\
24	0.0217900050712054\\
25	0.024421720594235\\
26	0.0253007507532161\\
27	0.025950812220769\\
28	0.0267100479376679\\
29	0.0268368261682328\\
30	0.0263395664970972\\
31	0.0253058206960862\\
32	0.0241303236955129\\
33	0.0229259313833535\\
34	0.0209554328060933\\
35	0.0171719378010867\\
36	0.0146565123513448\\
37	0.0119227665050947\\
38	0.0123689395518029\\
39	0.0185847630564384\\
40	0.028422660922743\\
41	0.0356510630674017\\
42	0.0343917198542466\\
43	0.0291067686151274\\
44	0.0240990412911144\\
45	0.0210702127439931\\
46	0.0207391018467359\\
47	0.022363355785673\\
48	0.0234053684250568\\
49	0.0243278798494925\\
50	0.0253510559882913\\
51	0.0259921780033776\\
52	0.0263745929545902\\
53	0.0257702275912129\\
54	0.0254490714227378\\
55	0.0249515378434057\\
56	0.024386707821839\\
57	0.0238863701560423\\
58	0.023659336297278\\
59	0.0239493340103837\\
60	0.0245574234252014\\
61	0.0250445095570516\\
62	0.025323532734849\\
63	0.0255728760318653\\
64	0.0256353772349363\\
65	0.0255764197742215\\
66	0.0252995389447065\\
67	0.0250725431254606\\
68	0.0247745773403051\\
69	0.0245663512344158\\
70	0.0244134695121892\\
71	0.0243714811994683\\
72	0.0244260109275675\\
73	0.024553626845828\\
74	0.0247159574876338\\
75	0.0249298892577418\\
76	0.0251281401606863\\
77	0.0252156103433981\\
78	0.0251747579913421\\
79	0.0251400117499301\\
80	0.0250642474471538\\
81	0.0248884891890003\\
82	0.0246554799922717\\
83	0.0245017345477521\\
84	0.0244168454770398\\
85	0.0244440840305711\\
86	0.0245364088298419\\
87	0.0247314475662143\\
88	0.0249955411070363\\
89	0.0251844309401427\\
90	0.0252964100940318\\
91	0.0253848960824397\\
92	0.0254087380738708\\
93	0.0253775352165378\\
94	0.0252080382589195\\
95	0.0250113223915681\\
96	0.0248455641202944\\
97	0.0246602112701183\\
98	0.0244721963248314\\
99	0.024372401488599\\
100	0.0243482320449974\\
101	0.0244099085450842\\
102	0.024501160619092\\
103	0.0246807948944066\\
104	0.0249359324624219\\
105	0.0251758653734846\\
106	0.0253252728300336\\
107	0.0254307864885233\\
108	0.0254601738965687\\
109	0.0254413521961566\\
110	0.0253335061542132\\
111	0.0252004708345473\\
112	0.0250743603117161\\
113	0.0248968340234645\\
114	0.0246903142587658\\
115	0.024571244380478\\
116	0.024522889033997\\
117	0.0245694155588752\\
118	0.0246616698404873\\
119	0.0248308599833817\\
120	0.0249974438740835\\
121	0.025118414233898\\
122	0.02518006558275\\
123	0.0252824489065737\\
124	0.025359164372657\\
125	0.0253870256564381\\
126	0.0253176952179094\\
127	0.025200232067708\\
128	0.0250928278349705\\
129	0.0249548766276821\\
130	0.02478395205551\\
131	0.0246741109820979\\
132	0.0246160390924755\\
133	0.0246367720930028\\
134	0.0246424036771059\\
135	0.0247255059700657\\
136	0.0248891427330082\\
137	0.0250610004185703\\
138	0.0251819632983014\\
139	0.0252941705214576\\
140	0.0253504727673473\\
141	0.0253708023212892\\
142	0.0252996084976007\\
143	0.0252023688521531\\
144	0.0251244966557485\\
145	0.0249925766516324\\
146	0.0248527041615231\\
147	0.0247500668023986\\
148	0.0246861205420633\\
149	0.0246413361993066\\
150	0.0246214696239675\\
151	0.0246442246009677\\
152	0.0247203978608165\\
153	0.0247969466113776\\
154	0.0248483393035512\\
155	0.0249295477751711\\
156	0.0250021477979886\\
157	0.0250846796890676\\
158	0.0251361721241365\\
159	0.0252003023522896\\
160	0.0252719482755589\\
161	0.0252907842580929\\
162	0.0252681783467097\\
163	0.0252479309486102\\
164	0.0251923097300118\\
165	0.0251240870159756\\
166	0.0249966320935931\\
167	0.0248805706155674\\
168	0.024771696092184\\
169	0.0246655485251763\\
170	0.024574770755415\\
171	0.0245344093139423\\
172	0.024548491962773\\
173	0.0245875212058292\\
174	0.0246651961934797\\
175	0.0247819048099962\\
176	0.0249437999111253\\
177	0.0250831376943235\\
178	0.0251974687666397\\
179	0.0253132358906177\\
180	0.0254071590210196\\
181	0.0254464489797401\\
182	0.0254227498966781\\
183	0.0253652756554492\\
184	0.0252950424195792\\
185	0.0251542445556947\\
186	0.0249873963295952\\
187	0.0248523924412165\\
188	0.0247238159543982\\
189	0.0246212304442035\\
190	0.0245258423010339\\
191	0.0244934225576607\\
192	0.024517054268499\\
193	0.0245371831347859\\
194	0.0245865369572037\\
195	0.0246608624947897\\
196	0.0247341264101514\\
197	0.0248028608743292\\
198	0.0249111300540069\\
199	0.0250416132049074\\
200	0.0361809251732859\\
201	0.03414685190299\\
202	0.0320748755577277\\
203	0.0294566889124462\\
204	0.0265130550679192\\
205	0.023896461295185\\
206	0.0217165488287877\\
207	0.0203428179695585\\
208	0.0199882412438962\\
209	0.0206010066422027\\
210	0.0223339863269817\\
211	0.0241095846087596\\
212	0.0259253208612034\\
213	0.0263251199601443\\
214	0.0261944951028044\\
215	0.0260089134254423\\
216	0.0255776657060302\\
217	0.0250628363909282\\
218	0.0249116462657398\\
219	0.024973120567169\\
220	0.0251696085304145\\
221	0.025435772776436\\
222	0.0254408097345794\\
223	0.0255075407934558\\
224	0.0255206366993947\\
225	0.025481971994475\\
226	0.0250972319205702\\
227	0.0246700340908191\\
228	0.0242398214255812\\
229	0.0239242094897115\\
230	0.0237751450733433\\
231	0.0237928587273251\\
232	0.0239537372830815\\
233	0.0242447218332026\\
234	0.0245781716342944\\
235	0.0249333944052107\\
236	0.0252737025561707\\
237	0.0255527675280495\\
238	0.0257610701324356\\
239	0.0259072430802671\\
240	0.0259993635402377\\
241	0.0259891344051465\\
242	0.0258744634889801\\
243	0.0257135451956922\\
244	0.0254982911485343\\
245	0.0252206921767348\\
246	0.0248944180416609\\
247	0.0245548275540257\\
248	0.0242847109553673\\
249	0.0241160056258928\\
250	0.02398821881169\\
251	0.0238957566457502\\
252	0.0238620668338042\\
253	0.0239132304516917\\
254	0.0240389334240212\\
255	0.024234446112017\\
256	0.0245140072443754\\
257	0.024818754887604\\
258	0.0251412077583453\\
259	0.0254095560304156\\
260	0.0255698213564401\\
261	0.0256425957995897\\
262	0.0255861903051036\\
263	0.0254382552766275\\
264	0.0252016045278857\\
265	0.0248619905398711\\
266	0.0245965504189784\\
267	0.0244068196355161\\
268	0.024307327131346\\
269	0.0242909415635311\\
270	0.0243601543673725\\
271	0.024489131717282\\
272	0.0246497361615535\\
273	0.0247972788662952\\
274	0.0249373173811499\\
275	0.0250767201395173\\
276	0.0251974940896595\\
277	0.0253023901415635\\
278	0.0253804813506234\\
279	0.0254481070042312\\
280	0.025494468200421\\
281	0.0255083368019746\\
282	0.0254970248670104\\
283	0.0254711098280684\\
284	0.0254186903173655\\
285	0.0253224174191158\\
286	0.0251852086242783\\
287	0.0250462530279251\\
288	0.0248889608925041\\
289	0.0247003148280183\\
290	0.024522956258558\\
291	0.0243791723529841\\
292	0.0242919695267246\\
293	0.0242642851867038\\
294	0.0243174277222833\\
295	0.0244422138035524\\
296	0.0246242262713974\\
297	0.0248655425686441\\
298	0.0251155891585098\\
299	0.0253715218901229\\
300	0.0256065101928544\\
301	0.0257273965407613\\
302	0.0257744505863573\\
303	0.0257800356169215\\
304	0.0257347255251848\\
305	0.0255913692366087\\
306	0.0254152036568312\\
307	0.0252447172723743\\
308	0.0250601907366271\\
309	0.0248869084803638\\
310	0.0247258473203504\\
311	0.0245978118734763\\
312	0.0244978470097753\\
313	0.0244417932665122\\
314	0.0243945232648147\\
315	0.0243676237941496\\
316	0.0243762535235398\\
317	0.0244237736711677\\
318	0.0245111972846405\\
319	0.0246312882203033\\
320	0.0247913930662151\\
321	0.0249715271287967\\
322	0.0251483580876269\\
323	0.0253173048283762\\
324	0.0254484121376228\\
325	0.0255492198456862\\
326	0.025582317803206\\
327	0.0255757996762383\\
328	0.0255183524783173\\
329	0.0254055790125061\\
330	0.0252674130311826\\
331	0.0251321695437606\\
332	0.0250013050745283\\
333	0.0248736409047367\\
334	0.0247560514190076\\
335	0.0246563191081163\\
336	0.02457988717461\\
337	0.0245211936305468\\
338	0.0244973320865471\\
339	0.0245072446263741\\
340	0.0245464528663698\\
341	0.024619121885337\\
342	0.0247141464537987\\
343	0.024828869121066\\
344	0.0249507100043395\\
345	0.0250727245824618\\
346	0.0251949962283886\\
347	0.0253145168327588\\
348	0.0254220562513573\\
349	0.0254970054726964\\
350	0.0255314572851515\\
351	0.0255326209663454\\
352	0.0254941318145029\\
353	0.025404722767387\\
354	0.0252631740091667\\
355	0.025104250777677\\
356	0.0249199631503873\\
357	0.0247380043002409\\
358	0.0245607881843103\\
359	0.0244133454584809\\
360	0.0243445463292446\\
361	0.0243338680070678\\
362	0.0243770356313078\\
363	0.0244676618003355\\
364	0.0245928967158603\\
365	0.024753827861689\\
366	0.0249454619257077\\
367	0.0251556924400444\\
368	0.0253687632983227\\
369	0.0255092034141651\\
370	0.0255803375996634\\
371	0.0256250943297699\\
372	0.0256262528814588\\
373	0.0255908888177711\\
374	0.0255078343652331\\
375	0.0254049907793243\\
376	0.0253042273060513\\
377	0.0252028586596574\\
378	0.0250873901869002\\
379	0.0249511342352637\\
380	0.0247971704684054\\
381	0.024628002873373\\
382	0.0244821006535146\\
383	0.0243646430122268\\
384	0.0242992998120156\\
385	0.024265583669021\\
386	0.0242873988882529\\
387	0.0243190980592372\\
388	0.0243681364844379\\
389	0.0244385229513239\\
390	0.0245577651217799\\
391	0.0247024312885103\\
392	0.0248578168645656\\
393	0.0249527277583832\\
394	0.0250571130453371\\
395	0.0251674532751548\\
396	0.0252750741235033\\
397	0.0253713846125678\\
398	0.0254503124358527\\
399	0.0255147785843419\\
400	0.025566451874783\\
401	0.0255806903160329\\
};
\addlegendentry{exact}

\addplot [color=black, dashed, forget plot]
  table[row sep=crcr]{%
1	0.025\\
400	0.025\\
};
\end{axis}
\end{tikzpicture}%

%% file: mma_16_cc.tex
\definecolor{mycolor1}{rgb}{0.00000,0.44700,0.74100}%
\definecolor{mycolor2}{rgb}{0.85000,0.32500,0.09800}%
\definecolor{mycolor3}{rgb}{0.49400,0.18400,0.55600}%
\begin{tikzpicture}

\begin{axis}[%
width=\textwidth,
height=.75\textwidth,
at={(0.758in,0.481in)},
%scale only axis,
xmin=0,
xmax=400,
ymin=0,
ymax=1,
axis background/.style={fill=white},
title style={font=\bfseries},
title={MMA, $\mathcal{B}=16$},
every axis title/.style={at={(0,1)}, anchor=north west, draw=black, fill=white},
xmajorgrids,
ymajorgrids,
xminorgrids,
yminorgrids,
xlabel={\vphantom{Iteration}},
ytick={0.025,0.2,0.4,0.6,0.8,1},
yticklabels={$p$,0.2,0.4,0.6,0.8,1},
legend style={legend cell align=left, align=left, draw=white!15!black}
]
\addplot [color=mycolor1, line width=1.5pt,line join=round]
  table[row sep=crcr]{%
1	1.0033299237069\\
2	0.809205179229764\\
3	0.0665748341454685\\
4	0.0371780467220003\\
5	0.0404712198193883\\
6	0.337111967610145\\
7	0.323403622740422\\
8	0.279066026994625\\
9	0.151046196473375\\
10	0.0799838689364968\\
11	0.0405422503955068\\
12	0.0267208497202702\\
13	0.0248584399753426\\
14	0.0248097804785039\\
15	0.0247641426380739\\
16	0.0246880727135889\\
17	0.0245795339991631\\
18	0.0244262435171393\\
19	0.0242203917621281\\
20	0.0239824330469772\\
21	0.0236955581724109\\
22	0.0233820062826247\\
23	0.0229967327411453\\
24	0.0223494070571863\\
25	0.0214954378148756\\
26	0.0206498387545893\\
27	0.0200547989128094\\
28	0.0199861149404129\\
29	0.020158821196476\\
30	0.020137555651898\\
31	0.0203585364010856\\
32	0.0208910681578921\\
33	0.0214550991889072\\
34	0.0222494585679393\\
35	0.0230284443012997\\
36	0.0239086301270689\\
37	0.0242914117663171\\
38	0.0244276986806826\\
39	0.0245538372692184\\
40	0.0246726423757367\\
41	0.0247624771000718\\
42	0.0248134920599015\\
43	0.0248468703505193\\
44	0.0248706764480855\\
45	0.0248908968551159\\
46	0.0249000561602967\\
47	0.0249113788237063\\
48	0.0249200719850403\\
49	0.0249288363380647\\
50	0.0249361022086803\\
51	0.0249410933115545\\
52	0.0249421988945107\\
53	0.0249495432298355\\
54	0.024959160547484\\
55	0.024965735703551\\
56	0.0249690700909125\\
57	0.0249725490034583\\
58	0.0249751848038373\\
59	0.0249785957554657\\
60	0.0249806472304361\\
61	0.024982347083688\\
62	0.0249841484506521\\
63	0.0249841395085583\\
64	0.0249842340083456\\
65	0.024985613172997\\
66	0.0249870636000633\\
67	0.0249883328509288\\
68	0.0249895603904819\\
69	0.0249904529742106\\
70	0.0249908758880481\\
71	0.0249913210566994\\
72	0.0249918898888656\\
73	0.0249927073760086\\
74	0.0249931753059924\\
75	0.0249933083669656\\
76	0.024993045491478\\
77	0.0249935217086599\\
78	0.0249941925370577\\
79	0.0249946764852046\\
80	0.0249950564001786\\
81	0.0249956268300401\\
82	0.0249961778605106\\
83	0.0249964289816296\\
84	0.0249965303564184\\
85	0.0249968029645782\\
86	0.0249969518094871\\
87	0.0249967452090894\\
88	0.0249968767099395\\
89	0.0249972008201045\\
90	0.0249972961143554\\
91	0.0249973512398484\\
92	0.0249974410728686\\
93	0.0249974203062248\\
94	0.0249976107992372\\
95	0.024997761664019\\
96	0.0249977682329786\\
97	0.0249978419648266\\
98	0.0249978571885313\\
99	0.0249978271579198\\
100	0.0249978706210419\\
101	0.0249978558137701\\
102	0.0249979002268465\\
103	0.0249978594975718\\
104	0.0249978974832864\\
105	0.0249978314766137\\
106	0.0249977760555984\\
107	0.0249977472951687\\
108	0.0249979131773309\\
109	0.0249979852824317\\
110	0.024997939059158\\
111	0.0249978305973196\\
112	0.0249976326259886\\
113	0.0249976248406135\\
114	0.0249977331142039\\
115	0.0249979543324184\\
116	0.0249979823199654\\
117	0.0249980931916095\\
118	0.0249981119490022\\
119	0.0249982326911061\\
120	0.0249981731765441\\
121	0.0249982271272657\\
122	0.0249982276564718\\
123	0.0249982572818911\\
124	0.0249983324855029\\
125	0.024998322821564\\
126	0.0249982477726624\\
127	0.0249981803105583\\
128	0.024997984308444\\
129	0.0249980591637857\\
130	0.0249982269903576\\
131	0.0249982116302565\\
132	0.0249982895375389\\
133	0.0249983739342824\\
134	0.0249984424199522\\
135	0.0249984097747104\\
136	0.0249983550125984\\
137	0.0249983288638954\\
138	0.0249983362936011\\
139	0.024998280819554\\
140	0.0249982623406704\\
141	0.0249982953881826\\
142	0.0249982707747138\\
143	0.0249982287208255\\
144	0.0249981648851721\\
145	0.0249981009212351\\
146	0.0249979364574598\\
147	0.0249979251122487\\
148	0.0249978528750461\\
149	0.0249979538815513\\
150	0.0249982192745704\\
151	0.0249985074957038\\
152	0.0249985413110232\\
153	0.0249986082550416\\
154	0.0249987113717428\\
155	0.0249988199412299\\
156	0.0249988279991682\\
157	0.0249987857042953\\
158	0.0249987556563473\\
159	0.0249987641611107\\
160	0.0249987268593637\\
161	0.0249986845779545\\
162	0.0249986560107952\\
163	0.0249986079550189\\
164	0.0249986120379135\\
165	0.0249986089388716\\
166	0.0249985711431547\\
167	0.0249984942841681\\
168	0.0249984505295721\\
169	0.0249984714620695\\
170	0.0249984895851593\\
171	0.0249985262705971\\
172	0.0249986086199349\\
173	0.0249987075548835\\
174	0.0249987006721615\\
175	0.024998685664374\\
176	0.024998677680432\\
177	0.024998703856945\\
178	0.0249987692083074\\
179	0.0249988458884706\\
180	0.0249988846058991\\
181	0.0249989510879356\\
182	0.0249989684571766\\
183	0.0249989152958328\\
184	0.0249988585586001\\
185	0.0249988117999146\\
186	0.024998756897333\\
187	0.0249987334436245\\
188	0.0249987025217438\\
189	0.0249986698684281\\
190	0.0249986409968748\\
191	0.0249986517004235\\
192	0.0249985711054883\\
193	0.0249986611156568\\
194	0.0249987816960042\\
195	0.0249988341475946\\
196	0.0249989105266842\\
197	0.0249989587242301\\
198	0.0249989598624069\\
199	0.0249989352428804\\
200	0.0322807327290019\\
201	0.0234495348338493\\
202	0.0248868895450728\\
203	0.0249602702966741\\
204	0.0249786960687626\\
205	0.0249875119995285\\
206	0.0249917414267033\\
207	0.0249933884172256\\
208	0.0249945489042169\\
209	0.0249957038957999\\
210	0.0249966293411197\\
211	0.0249971454809374\\
212	0.0249975112932402\\
213	0.024997766599514\\
214	0.024997911051244\\
215	0.0249980462596652\\
216	0.0249981362864632\\
217	0.0249981459188765\\
218	0.0249981896316083\\
219	0.024998404805545\\
220	0.0249985974578426\\
221	0.0249986641518168\\
222	0.024998715657645\\
223	0.0249987504084216\\
224	0.024998842649333\\
225	0.0249988845868248\\
226	0.0249989386789039\\
227	0.0249989680156047\\
228	0.0249989695385129\\
229	0.0249989469397776\\
230	0.0249989637554708\\
231	0.0249989924710172\\
232	0.0249989947445041\\
233	0.0249990897081929\\
234	0.0249991289664064\\
235	0.0249991283248171\\
236	0.0249991054205186\\
237	0.0249991050056106\\
238	0.0249990328325943\\
239	0.0249990267860946\\
240	0.0249989697462683\\
241	0.0249989202635936\\
242	0.0249988741890486\\
243	0.0249988745012548\\
244	0.0249988452219671\\
245	0.0249987817657711\\
246	0.0249987265864046\\
247	0.0249986687153657\\
248	0.0249985727421861\\
249	0.0249984769725747\\
250	0.0249983919964782\\
251	0.0249983440666496\\
252	0.0249983971932917\\
253	0.0249983723026776\\
254	0.0249984630487083\\
255	0.0249985904736856\\
256	0.0249986360427035\\
257	0.0249987757625155\\
258	0.0249988728904844\\
259	0.02499889597058\\
260	0.0249988283854492\\
261	0.0249988110453582\\
262	0.024998832623072\\
263	0.0249989921981744\\
264	0.0249990512479797\\
265	0.0249991044023089\\
266	0.024999111060144\\
267	0.0249992007932821\\
268	0.024999304783623\\
269	0.0249993330117868\\
270	0.0249993943737527\\
271	0.0249994067403895\\
272	0.0249993917528981\\
273	0.0249993749539752\\
274	0.0249993913878457\\
275	0.0249994033213952\\
276	0.0249994147786656\\
277	0.0249994153616731\\
278	0.0249993808619562\\
279	0.0249993497821329\\
280	0.0249993039553676\\
281	0.0249992395162302\\
282	0.0249991946547816\\
283	0.0249992060709174\\
284	0.0249992374658075\\
285	0.0249992504121502\\
286	0.0249992446120978\\
287	0.0249992328874103\\
288	0.0249992331607271\\
289	0.0249992030502222\\
290	0.0249991363059737\\
291	0.0249990709307074\\
292	0.0249990055345749\\
293	0.0249989915397818\\
294	0.0249989868438148\\
295	0.0249989501260402\\
296	0.0249988975576344\\
297	0.0249988682184342\\
298	0.0249988358582276\\
299	0.0249987529497902\\
300	0.0249987121225724\\
301	0.0249987273760004\\
302	0.0249987359274607\\
303	0.024998733681438\\
304	0.0249987546186261\\
305	0.0249988223643324\\
306	0.024998858974625\\
307	0.0249988958615926\\
308	0.024998847168686\\
309	0.0249987961022532\\
310	0.0249987225995188\\
311	0.0249985881279306\\
312	0.0249984929066968\\
313	0.0249984071032275\\
314	0.0249983379211635\\
315	0.0249983746097882\\
316	0.0249983453634952\\
317	0.0249984507537756\\
318	0.0249984056910579\\
319	0.0249983150529039\\
320	0.0249981724821648\\
321	0.0249980323378417\\
322	0.0249979674958199\\
323	0.0249979550381001\\
324	0.0249978641104723\\
325	0.0249978463178174\\
326	0.0249977835518956\\
327	0.0249979485886859\\
328	0.0249980935308948\\
329	0.0249982014499909\\
330	0.0249983535824686\\
331	0.0249984047729857\\
332	0.0249985681431312\\
333	0.0249987850358026\\
334	0.0249989620127359\\
335	0.0249990274313203\\
336	0.0249991207827863\\
337	0.0249991820415494\\
338	0.0249992489764419\\
339	0.0249993107007852\\
340	0.0249993533890572\\
341	0.0249993700303028\\
342	0.024999383301407\\
343	0.0249993808054541\\
344	0.0249993682609524\\
345	0.0249993674871876\\
346	0.0249993709455873\\
347	0.024999357069742\\
348	0.0249993495731539\\
349	0.0249993333099113\\
350	0.0249993138102872\\
351	0.0249992964842357\\
352	0.0249992725749763\\
353	0.0249992282853959\\
354	0.0249992393209869\\
355	0.0249992316057263\\
356	0.0249991838045253\\
357	0.0249991242959977\\
358	0.0249990283663264\\
359	0.0249989405574306\\
360	0.0249989116896598\\
361	0.024998961173774\\
362	0.0249989631353065\\
363	0.0249989227957655\\
364	0.024998889551987\\
365	0.0249988253822327\\
366	0.0249987770649064\\
367	0.0249987794432709\\
368	0.0249988203312238\\
369	0.0249988067644371\\
370	0.0249988066308351\\
371	0.0249987677349662\\
372	0.0249986688410253\\
373	0.024998664917281\\
374	0.0249986514037147\\
375	0.024998694741665\\
376	0.0249987464819886\\
377	0.0249987875233087\\
378	0.024998858377764\\
379	0.0249988650949621\\
380	0.0249988504143055\\
381	0.0249988324522763\\
382	0.0249988149776301\\
383	0.0249988322253215\\
384	0.0249988865470408\\
385	0.0249989156758766\\
386	0.0249989144928096\\
387	0.0249989034070666\\
388	0.0249988617650708\\
389	0.0249987853261431\\
390	0.0249987420140147\\
391	0.0249987554857193\\
392	0.0249987557093255\\
393	0.0249987572671932\\
394	0.0249988046337091\\
395	0.0249988590633539\\
396	0.0249988878303324\\
397	0.0249988736939698\\
398	0.0249988085950697\\
399	0.0249987126810864\\
400	0.0249986182386424\\
};
\addlegendentry{internal}

\addplot [color=mycolor2, line width=1pt,line join=round]
  table[row sep=crcr]{%
1	1.00331519624177\\
2	0.817370799108866\\
3	0.0764244368174531\\
4	0.0464899835293668\\
5	0.0539698981841624\\
6	0.366063649102888\\
7	0.387604829764022\\
8	0.356212550904089\\
9	0.238714645061228\\
10	0.139494488300636\\
11	0.077329298401878\\
12	0.0529751465002873\\
13	0.0514859270878881\\
14	0.0536745511464473\\
15	0.0565424095450421\\
16	0.0603200933033912\\
17	0.0654700525742283\\
18	0.0726138031714211\\
19	0.0825544132119864\\
20	0.096262616223276\\
21	0.114927289500552\\
22	0.140579729884613\\
23	0.173738426375319\\
24	0.210310522353227\\
25	0.24525446931072\\
26	0.277916306858356\\
27	0.310435251142246\\
28	0.345510494151035\\
29	0.380778354351074\\
30	0.411340088007886\\
31	0.4390425494392\\
32	0.462766461391578\\
33	0.480651028520871\\
34	0.494975298211559\\
35	0.506595667797223\\
36	0.516743729648207\\
37	0.523556602257478\\
38	0.527553828768996\\
39	0.530457803772513\\
40	0.532291689180238\\
41	0.5332116518923\\
42	0.533186383935337\\
43	0.532501369569433\\
44	0.531511498620934\\
45	0.530424694902067\\
46	0.529218860578108\\
47	0.527974837313493\\
48	0.526733229012271\\
49	0.525533939403907\\
50	0.524332884105127\\
51	0.52313389523935\\
52	0.52192873839633\\
53	0.520834527920333\\
54	0.519822412458239\\
55	0.518907861019538\\
56	0.518032753016818\\
57	0.517201437930459\\
58	0.516368097012245\\
59	0.515541652307606\\
60	0.514715501357748\\
61	0.513915424190366\\
62	0.51312818477952\\
63	0.512334464652988\\
64	0.51156075447855\\
65	0.510820506056653\\
66	0.510101404605392\\
67	0.509408934225777\\
68	0.508762105467271\\
69	0.508158807527688\\
70	0.50760298082331\\
71	0.507076020296835\\
72	0.50656697555498\\
73	0.506086147885151\\
74	0.50561798989834\\
75	0.505162578587994\\
76	0.504719653334308\\
77	0.504321196194286\\
78	0.503971750159063\\
79	0.503652946123143\\
80	0.503352323836632\\
81	0.503064860251211\\
82	0.502794048676054\\
83	0.502539078356226\\
84	0.502295732908105\\
85	0.502062894964049\\
86	0.501840525684617\\
87	0.50161765654614\\
88	0.501398069033748\\
89	0.501186278095539\\
90	0.500977895684917\\
91	0.500770202279256\\
92	0.50056519102213\\
93	0.500361389591099\\
94	0.500168854332021\\
95	0.499986117138421\\
96	0.499809675479575\\
97	0.499639833318022\\
98	0.499475841578724\\
99	0.499314984882266\\
100	0.499162708840358\\
101	0.499017042923924\\
102	0.498878009122952\\
103	0.498740100368135\\
104	0.498603782783459\\
105	0.498466366705715\\
106	0.498327275020962\\
107	0.498183394867278\\
108	0.498040058818539\\
109	0.497894806959399\\
110	0.497747991297215\\
111	0.497597851161258\\
112	0.497439140124414\\
113	0.497277931359742\\
114	0.497119856385398\\
115	0.496965913417047\\
116	0.496813029773263\\
117	0.496665159219517\\
118	0.496519636205487\\
119	0.496380344389115\\
120	0.496242476819067\\
121	0.496108909986551\\
122	0.495979423259889\\
123	0.49585230587647\\
124	0.495728474336696\\
125	0.495604703812822\\
126	0.495479518929113\\
127	0.495351134857202\\
128	0.495216686791193\\
129	0.495082188226519\\
130	0.494950310069806\\
131	0.494817558920627\\
132	0.494687392901837\\
133	0.494561275535067\\
134	0.494439579989548\\
135	0.49431965453163\\
136	0.494199259876083\\
137	0.494081182027718\\
138	0.493967248116357\\
139	0.493854275170192\\
140	0.493742026544602\\
141	0.493635108156932\\
142	0.493532580282527\\
143	0.493432660958852\\
144	0.493337636338233\\
145	0.493246723267198\\
146	0.4931567786742\\
147	0.493069336002872\\
148	0.49298021344741\\
149	0.492898146634272\\
150	0.492826934324258\\
151	0.492761260987626\\
152	0.4926935305754\\
153	0.492628197398441\\
154	0.492569666052668\\
155	0.492518876340807\\
156	0.492471264610757\\
157	0.492422967928594\\
158	0.492373403042797\\
159	0.492325158061838\\
160	0.492278555865337\\
161	0.49223454807007\\
162	0.492194194560029\\
163	0.492152592860447\\
164	0.492108749956887\\
165	0.492066806369313\\
166	0.492024020287193\\
167	0.491977825315834\\
168	0.491930124954966\\
169	0.491885109555925\\
170	0.491843792892997\\
171	0.491804832361076\\
172	0.491766865542653\\
173	0.49173079658841\\
174	0.491695105338134\\
175	0.491660706184246\\
176	0.491628142243566\\
177	0.491596530255456\\
178	0.491566746295247\\
179	0.491538807491565\\
180	0.491510442276172\\
181	0.491482710847189\\
182	0.491455863910527\\
183	0.491426813622046\\
184	0.491394630010751\\
185	0.491359233085014\\
186	0.491319156441382\\
187	0.491275951303126\\
188	0.49122962422419\\
189	0.491179177921651\\
190	0.491125473087922\\
191	0.491073149936439\\
192	0.491019891548979\\
193	0.490969122703656\\
194	0.490922909857268\\
195	0.490876416964921\\
196	0.490829384822616\\
197	0.490781741158958\\
198	0.490730409124908\\
199	0.490674176845007\\
200	0.546796174045338\\
201	0.482178222001184\\
202	0.493859884404957\\
203	0.494890167601727\\
204	0.495112750132716\\
205	0.49510587685155\\
206	0.494996130607518\\
207	0.494845245253436\\
208	0.494706387489011\\
209	0.494608276523651\\
210	0.494548442551017\\
211	0.494507117493036\\
212	0.494473706799692\\
213	0.494447758636456\\
214	0.494421166439048\\
215	0.494392507791162\\
216	0.494356157250024\\
217	0.49431598164641\\
218	0.494278856468427\\
219	0.494253765610498\\
220	0.494246226879648\\
221	0.494247776966027\\
222	0.494258059277325\\
223	0.49427371189844\\
224	0.494287774210335\\
225	0.494295574771084\\
226	0.494298143542698\\
227	0.494298298473545\\
228	0.494295546425425\\
229	0.494288500146642\\
230	0.494275942949113\\
231	0.494260228271585\\
232	0.494242906694747\\
233	0.49422510322267\\
234	0.494206518603026\\
235	0.494186465989467\\
236	0.494163479958214\\
237	0.494136488908049\\
238	0.494105382938271\\
239	0.494072247775153\\
240	0.494039548828497\\
241	0.494005107580848\\
242	0.493967057339058\\
243	0.493928063512634\\
244	0.493892064157045\\
245	0.493859667456316\\
246	0.493829139670427\\
247	0.49379770651339\\
248	0.493765006863867\\
249	0.493732227540235\\
250	0.493699155505955\\
251	0.493663349703469\\
252	0.493624437312726\\
253	0.493581329084374\\
254	0.493540402520617\\
255	0.49350384301824\\
256	0.493467936297175\\
257	0.493431422746\\
258	0.493394418425449\\
259	0.493357152004025\\
260	0.4933168078448\\
261	0.493276285375227\\
262	0.493235491322874\\
263	0.493198740183726\\
264	0.493163769539813\\
265	0.493129124397378\\
266	0.493094449876784\\
267	0.493062473733085\\
268	0.493035346006498\\
269	0.493010966641029\\
270	0.492992078911006\\
271	0.492976837673199\\
272	0.492963343354191\\
273	0.492949489335561\\
274	0.492934762337459\\
275	0.492918824200627\\
276	0.492901732152146\\
277	0.492884854130897\\
278	0.492868547732444\\
279	0.492851801576618\\
280	0.492834100648566\\
281	0.492814741166969\\
282	0.492794030826182\\
283	0.492774428924558\\
284	0.49275620885984\\
285	0.492736960884291\\
286	0.492715138553313\\
287	0.492693825648327\\
288	0.492676771292898\\
289	0.492662767526476\\
290	0.492649747817715\\
291	0.492639103110519\\
292	0.492632806820245\\
293	0.492630942204279\\
294	0.492629400842387\\
295	0.49262524947543\\
296	0.49261707993731\\
297	0.492603856881293\\
298	0.492585824454916\\
299	0.492565379921063\\
300	0.492543919424339\\
301	0.492524571524015\\
302	0.492508981696506\\
303	0.492499012308829\\
304	0.492495099281769\\
305	0.492497358438387\\
306	0.4925016235416\\
307	0.492502085205725\\
308	0.492495658176886\\
309	0.492481988864084\\
310	0.492464758691395\\
311	0.492444866484513\\
312	0.49242485388904\\
313	0.492405262284453\\
314	0.49238383778019\\
315	0.492359808556601\\
316	0.492333888449385\\
317	0.492313475454012\\
318	0.492294850797347\\
319	0.492273139868944\\
320	0.492247602939199\\
321	0.492217097650514\\
322	0.492179909976868\\
323	0.492142713371287\\
324	0.492105220878325\\
325	0.492068203365225\\
326	0.492033886323144\\
327	0.492010873371103\\
328	0.492000370894198\\
329	0.491998346054331\\
330	0.492005329534847\\
331	0.492010537527616\\
332	0.492007011743257\\
333	0.491994782650224\\
334	0.491978342691041\\
335	0.491958745161961\\
336	0.491939015480823\\
337	0.491920348203086\\
338	0.491902327563802\\
339	0.491882239158113\\
340	0.491858069184586\\
341	0.491831290043357\\
342	0.491803531867586\\
343	0.491775131487909\\
344	0.491747037430946\\
345	0.491719568700851\\
346	0.491692705595271\\
347	0.491667044794666\\
348	0.491642375643369\\
349	0.491616772841183\\
350	0.49158820967186\\
351	0.491556438762003\\
352	0.491521711614717\\
353	0.491484827155487\\
354	0.491448043589249\\
355	0.491410557426988\\
356	0.491370157634369\\
357	0.491326303617227\\
358	0.491277911087033\\
359	0.491225719773231\\
360	0.491172229575347\\
361	0.491119189937637\\
362	0.491065453426323\\
363	0.491012125921631\\
364	0.49096102436319\\
365	0.490911521667245\\
366	0.490862943964864\\
367	0.490815482566459\\
368	0.490770665661236\\
369	0.490727447639008\\
370	0.490685584446485\\
371	0.490646176071934\\
372	0.490607905473077\\
373	0.490571668500241\\
374	0.490538192696119\\
375	0.490506932120435\\
376	0.490474858589221\\
377	0.490443602290547\\
378	0.490416615095355\\
379	0.490392228058729\\
380	0.490369005982919\\
381	0.490347580097516\\
382	0.490326021078915\\
383	0.490306099536659\\
384	0.490287197037812\\
385	0.490269327854389\\
386	0.490254543020963\\
387	0.490243592510589\\
388	0.490236429074075\\
389	0.490232722051334\\
390	0.490230802564742\\
391	0.490226554288736\\
392	0.49021777208096\\
393	0.4902058045226\\
394	0.490192155899377\\
395	0.490178252059623\\
396	0.490167125779257\\
397	0.490160904747784\\
398	0.490157526541866\\
399	0.490153896415267\\
400	0.490147476534243\\
401	0.490136209101736\\
};
\addlegendentry{exact}

\addplot [color=black, dashed, forget plot]
  table[row sep=crcr]{%
1	0.025\\
400	0.025\\
};

\end{axis}
\end{tikzpicture}%

%% file: mcmsa_32_cc.tex
\definecolor{mycolor1}{rgb}{0.00000,0.44700,0.74100}%
\definecolor{mycolor2}{rgb}{0.85000,0.32500,0.09800}%
\definecolor{mycolor3}{rgb}{0.49400,0.18400,0.55600}%
\begin{tikzpicture}

\begin{axis}[%
width=\textwidth,
height=.75\textwidth,
at={(0.758in,0.481in)},
%scale only axis,
xmin=0,
xmax=400,
ymin=0,
ymax=1,
axis background/.style={fill=white},
title style={font=\bfseries},
title={sMMA, $\mathcal{B}=32$},
every axis title/.style={at={(0,1)}, anchor=north west, draw=black, fill=white},
xmajorgrids,
ymajorgrids,
xminorgrids,
yminorgrids,
xlabel={Iteration},
ytick={0.025,0.2,0.4,0.6,0.8,1},
yticklabels={$p$,0.2,0.4,0.6,0.8,1},
legend style={legend cell align=left, align=left, draw=white!15!black}
]
\addplot [color=mycolor1, line width=1.5pt,line join=round]
  table[row sep=crcr]{%
1	1.00332949283693\\
2	0.858993887433342\\
3	0.4632412322807\\
4	0.0336904826145071\\
5	0.00125903237163015\\
6	0.104084911479211\\
7	0.00192807194689496\\
8	0.180296963815484\\
9	0.00159934175979515\\
10	0.0156773935044056\\
11	0.0389303339660486\\
12	0.00233111104106581\\
13	0.0197743197611827\\
14	0.0153828447420685\\
15	0.0398992002468299\\
16	0.00986053542578846\\
17	0.0181294805304316\\
18	0.0197966890272074\\
19	0.0242356248150721\\
20	0.0242914957118459\\
21	0.025985009043102\\
22	0.0234335212581055\\
23	0.0264454060800398\\
24	0.0239134475654103\\
25	0.0249541258530592\\
26	0.0246536263456764\\
27	0.0247269638673015\\
28	0.0248437665976006\\
29	0.0251871173700864\\
30	0.0253163596726571\\
31	0.0258582346081626\\
32	0.0258349418543717\\
33	0.0262513079260525\\
34	0.0262413930966152\\
35	0.0262703604135758\\
36	0.0258124378034438\\
37	0.0252808600263847\\
38	0.0249016538454834\\
39	0.0245602639986977\\
40	0.024365813301117\\
41	0.0249271034176973\\
42	0.0247591117737398\\
43	0.0246376536503301\\
44	0.0245276522950777\\
45	0.0244328689573829\\
46	0.0244155689199041\\
47	0.0244100617343379\\
48	0.0244386650864634\\
49	0.0245631514277844\\
50	0.0247034283927828\\
51	0.0248154781491486\\
52	0.024951595463894\\
53	0.0250575725601193\\
54	0.025146884664888\\
55	0.0251425735474189\\
56	0.0251493295348248\\
57	0.0251686517506657\\
58	0.0251354856895302\\
59	0.0250905468383168\\
60	0.0250515080345295\\
61	0.0250667695316071\\
62	0.024997702602595\\
63	0.0249269094634452\\
64	0.0248933058107315\\
65	0.0248999123316211\\
66	0.0249262774930445\\
67	0.0250227947894541\\
68	0.0251118714626345\\
69	0.025182547823658\\
70	0.0252769336708304\\
71	0.0252325208620417\\
72	0.0251512201980841\\
73	0.0249393872161457\\
74	0.0242801286493739\\
75	0.0226490053675449\\
76	0.0202128090259311\\
77	0.0186501972964093\\
78	0.0172193865711984\\
79	0.0166707968505273\\
80	0.017865139688065\\
81	0.0191196021025441\\
82	0.0241248165171448\\
83	0.0255585798303643\\
84	0.0276769017073622\\
85	0.0258559298540855\\
86	0.0236133812413179\\
87	0.0224231287708415\\
88	0.0223173928667843\\
89	0.0224516522118565\\
90	0.0233375841466255\\
91	0.0237512753554555\\
92	0.023676846923384\\
93	0.023814114640927\\
94	0.0240357144559522\\
95	0.0241446091185468\\
96	0.024427186737725\\
97	0.0246554388669977\\
98	0.0245148372674248\\
99	0.024290895817975\\
100	0.0240669676919893\\
101	0.0240209043487397\\
102	0.0240378432572353\\
103	0.0240135821981035\\
104	0.0240289923010943\\
105	0.0241346318946789\\
106	0.0242667028014595\\
107	0.0243164676326302\\
108	0.0244555188133049\\
109	0.0244434214747704\\
110	0.0244229951873325\\
111	0.0244629429323495\\
112	0.0245047493904973\\
113	0.0245417094000955\\
114	0.0245590065426974\\
115	0.0245630571351127\\
116	0.0245367893915943\\
117	0.0244714311831371\\
118	0.0243835584807803\\
119	0.0242657579303277\\
120	0.0241925319027932\\
121	0.024157258217484\\
122	0.02421465385481\\
123	0.0243305165402073\\
124	0.0245016957852658\\
125	0.0246000002701577\\
126	0.0246949336727845\\
127	0.0248693911400064\\
128	0.0249500137797898\\
129	0.0250215522535313\\
130	0.0248941682041828\\
131	0.0247525794599018\\
132	0.0245792465340546\\
133	0.0245868846868398\\
134	0.0246268478333353\\
135	0.0246682263915319\\
136	0.0246794236235476\\
137	0.0247214273414895\\
138	0.0247941141081042\\
139	0.024873176514549\\
140	0.0249344094719723\\
141	0.0249357255339964\\
142	0.0249099565569503\\
143	0.0248600473382285\\
144	0.0248158573466098\\
145	0.0247756227329935\\
146	0.0247494066438332\\
147	0.0247167900987393\\
148	0.0247181974739653\\
149	0.0247525828002195\\
150	0.0248069985761815\\
151	0.0248797657244746\\
152	0.0249527660821752\\
153	0.024973311507445\\
154	0.0249768657411301\\
155	0.0249603024311674\\
156	0.0249344587383823\\
157	0.0249157409358654\\
158	0.0248834659220374\\
159	0.024862435696502\\
160	0.0248392592810787\\
161	0.0248466646208801\\
162	0.0248584295468219\\
163	0.0248766798538908\\
164	0.024900428021299\\
165	0.0249117820721596\\
166	0.0249280022826425\\
167	0.0249462572189715\\
168	0.0249651464922969\\
169	0.0249610513224378\\
170	0.0249577223445527\\
171	0.024944709792421\\
172	0.0249323344807423\\
173	0.0249280290596745\\
174	0.0249236927063724\\
175	0.0249251824640143\\
176	0.0249287016259275\\
177	0.0249291070502337\\
178	0.0249301145628738\\
179	0.0249306912391034\\
180	0.0249322175669122\\
181	0.0249304076208132\\
182	0.024930527123569\\
183	0.0249329994184579\\
184	0.0249379213310264\\
185	0.0249423035208326\\
186	0.0249476104926304\\
187	0.0249518097989562\\
188	0.0249567394419499\\
189	0.024963927373598\\
190	0.0249689022607563\\
191	0.0249742013158496\\
192	0.0249749406181194\\
193	0.0249759368550129\\
194	0.0249709844772137\\
195	0.0249619657720738\\
196	0.0249524592774162\\
197	0.0249373725601543\\
198	0.0249240041526532\\
199	0.0249129940029546\\
200	0.026377560542691\\
201	0.0265276716309876\\
202	0.0273480477283557\\
203	0.0277090406623251\\
204	0.0267695899494299\\
205	0.0260536578817536\\
206	0.0253503996716597\\
207	0.0244613265980158\\
208	0.0241623277957287\\
209	0.0233451156385175\\
210	0.0234215708002955\\
211	0.0240628207095425\\
212	0.0242932471807327\\
213	0.0250879842103254\\
214	0.0252954656644951\\
215	0.0256881766465287\\
216	0.0257043573877114\\
217	0.0252169084928435\\
218	0.0249558191390271\\
219	0.0247798606832509\\
220	0.0249472099418135\\
221	0.0249948075253021\\
222	0.024854786745126\\
223	0.0248389726964705\\
224	0.0248137491642655\\
225	0.0249991165204298\\
226	0.0250577212255652\\
227	0.0250600022345764\\
228	0.0250213227943391\\
229	0.0249364061671241\\
230	0.0248662520210888\\
231	0.0248192494363087\\
232	0.0248161799889788\\
233	0.0248481845845888\\
234	0.0249003191565865\\
235	0.0249588930106707\\
236	0.0250118660016223\\
237	0.0250661350350276\\
238	0.0251153684157744\\
239	0.0251441976884587\\
240	0.0251495448118106\\
241	0.0251311150854549\\
242	0.0251030694010026\\
243	0.025067973689027\\
244	0.0250307510815099\\
245	0.0249933292349168\\
246	0.0249566208170149\\
247	0.0249241872675249\\
248	0.0248978983987904\\
249	0.0248723516002371\\
250	0.0248478545028952\\
251	0.0248299278451476\\
252	0.0248123963913491\\
253	0.0248131988799596\\
254	0.0248339750318656\\
255	0.0248614542673754\\
256	0.0248797140532126\\
257	0.0249108405464054\\
258	0.0249353258531514\\
259	0.0249795414441216\\
260	0.0250208146284768\\
261	0.0251149936655542\\
262	0.0251854220493137\\
263	0.0252245638764826\\
264	0.0252332591970332\\
265	0.0251974972451983\\
266	0.0251624575037715\\
267	0.0251219949090668\\
268	0.0250834188887409\\
269	0.0250269473633052\\
270	0.0249715032452529\\
271	0.0249267007699521\\
272	0.0248911656673862\\
273	0.0248524009705983\\
274	0.0248284214166304\\
275	0.0248164763667349\\
276	0.0248162146062983\\
277	0.0248246291779361\\
278	0.0248449347983968\\
279	0.024876997794611\\
280	0.0249145515416064\\
281	0.0249254229912096\\
282	0.0249491287786296\\
283	0.0249856976311902\\
284	0.0250329135466353\\
285	0.0250650775644396\\
286	0.0250937602979082\\
287	0.0251215849162987\\
288	0.0251395798510504\\
289	0.0251599719370086\\
290	0.02516201321713\\
291	0.0251455456937377\\
292	0.02511264402\\
293	0.0250598899344385\\
294	0.0249937194744861\\
295	0.0249248937865788\\
296	0.0248628909593497\\
297	0.0248139506265836\\
298	0.0247826119869707\\
299	0.0247667855543974\\
300	0.024842659020069\\
301	0.0248674345025195\\
302	0.0248861390806065\\
303	0.024915216227564\\
304	0.0249700247892795\\
305	0.0250088706936876\\
306	0.0250303413731139\\
307	0.0250622162910833\\
308	0.0250878091224405\\
309	0.0250980578482947\\
310	0.0250853996677391\\
311	0.0250805489567038\\
312	0.0250675862915215\\
313	0.0250913083780743\\
314	0.0251072709352612\\
315	0.0251086596525976\\
316	0.0250984890976846\\
317	0.0250689438557443\\
318	0.0250373502264983\\
319	0.0249881031667105\\
320	0.0249455479862427\\
321	0.0249216446880044\\
322	0.024908685665288\\
323	0.0249012052357622\\
324	0.0249003450873675\\
325	0.0248985248550135\\
326	0.0249017521870247\\
327	0.0249141395896163\\
328	0.024933016337521\\
329	0.0249372776472329\\
330	0.0249412455014711\\
331	0.0249464222789491\\
332	0.0249532507013033\\
333	0.0249618099031248\\
334	0.0249711403892035\\
335	0.0249813486802105\\
336	0.0249926371966763\\
337	0.0250083933006682\\
338	0.0250212422872574\\
339	0.0250289395447056\\
340	0.0250360495339411\\
341	0.0250405628449496\\
342	0.0250411893486466\\
343	0.0250404476952086\\
344	0.0250383723627474\\
345	0.0250398471534683\\
346	0.025038291289503\\
347	0.0250355567859802\\
348	0.025033266591342\\
349	0.0250335760437918\\
350	0.0250330904766659\\
351	0.0250338966366163\\
352	0.0250320705865917\\
353	0.0250223791487515\\
354	0.0250101453015234\\
355	0.024995868746038\\
356	0.0249814696266631\\
357	0.0249663279557593\\
358	0.024952180568827\\
359	0.0249396067791131\\
360	0.0249296663158347\\
361	0.0249180898821844\\
362	0.0249106563489356\\
363	0.0249069099843823\\
364	0.024907364449035\\
365	0.0249123679886203\\
366	0.0249202088470097\\
367	0.0249330736427317\\
368	0.0249479426224701\\
369	0.0249651083937933\\
370	0.0249816892648028\\
371	0.0249974230246595\\
372	0.0250113729982838\\
373	0.0250244585579767\\
374	0.0250347872612652\\
375	0.0250402027866768\\
376	0.0250452968331465\\
377	0.0250625406049874\\
378	0.0250770439196249\\
379	0.0250869394566439\\
380	0.0250902420082394\\
381	0.0250940844826649\\
382	0.025093193853528\\
383	0.0250870889328775\\
384	0.0250790947507349\\
385	0.0250708327463922\\
386	0.0250590391735738\\
387	0.0250389367133824\\
388	0.0250199797654413\\
389	0.0249694128208186\\
390	0.0249638934635998\\
391	0.0249403539263651\\
392	0.0249216284693383\\
393	0.0248815237549514\\
394	0.0248468057771664\\
395	0.0248195172708245\\
396	0.0248039450652603\\
397	0.0247963755779273\\
398	0.0248035103106741\\
399	0.0248204167602069\\
400	0.0248518720414913\\
};
\addlegendentry{internal}

\addplot [color=mycolor2, line width=1pt,line join=round]
  table[row sep=crcr]{%
1	1.00331519624177\\
2	0.815894145238307\\
3	0.0199065616588609\\
4	0.00230271975981887\\
5	0.000279384620234748\\
6	0.0401164087778137\\
7	0.00221271377037118\\
8	0.0371509941415653\\
9	0.000100333999432385\\
10	0.00492784737774984\\
11	0.00820839921754774\\
12	0.0024377149649826\\
13	0.0107626577110193\\
14	0.0131392363788318\\
15	0.0205566371975766\\
16	0.00731407612398377\\
17	0.0152109462211563\\
18	0.0200520129461631\\
19	0.0249739671271875\\
20	0.0257523634979252\\
21	0.0265016773980169\\
22	0.0255065728700989\\
23	0.0272175960910993\\
24	0.0257220609809751\\
25	0.0268732586095858\\
26	0.0269274338637986\\
27	0.0273100883346172\\
28	0.0276200658119455\\
29	0.0278041085571138\\
30	0.0276049532051191\\
31	0.027262474442808\\
32	0.0263479472030969\\
33	0.0254924228888421\\
34	0.0242810281014069\\
35	0.0231480763817731\\
36	0.0220565256730249\\
37	0.0213810842432433\\
38	0.0211469974874576\\
39	0.0212257537922383\\
40	0.0215968945234777\\
41	0.022152459271195\\
42	0.0222157125012861\\
43	0.0224294440192451\\
44	0.0227579512300996\\
45	0.0231981186542934\\
46	0.0237442088011866\\
47	0.0243261115169637\\
48	0.0249328159490175\\
49	0.0255270001275726\\
50	0.0259981045477754\\
51	0.0263194474351839\\
52	0.0265168433405202\\
53	0.0265688321628421\\
54	0.0265137736592444\\
55	0.0263658324740145\\
56	0.0262150180856944\\
57	0.0260503660298395\\
58	0.0258691322081459\\
59	0.0257299834703584\\
60	0.0256440688079778\\
61	0.0255998231405053\\
62	0.0255373913728708\\
63	0.0255356144516941\\
64	0.0256034065267614\\
65	0.0257069273556207\\
66	0.0257992413315905\\
67	0.0258579741640583\\
68	0.0258051897147859\\
69	0.0256465998807982\\
70	0.0254113376269176\\
71	0.0250148475776117\\
72	0.0245275418197473\\
73	0.0237575260191869\\
74	0.022359662717324\\
75	0.0194972131437748\\
76	0.0171778519516366\\
77	0.016317176231702\\
78	0.0147258795903719\\
79	0.0160699665770317\\
80	0.0209223730466656\\
81	0.0266395449282198\\
82	0.0321817254811864\\
83	0.0314232172575393\\
84	0.0272737528890154\\
85	0.0229215965656655\\
86	0.0208709087521969\\
87	0.0205965877504944\\
88	0.0214420218273714\\
89	0.0227728153841885\\
90	0.0242266568208583\\
91	0.0245959352141174\\
92	0.0246226125888071\\
93	0.0247133752408225\\
94	0.0247673648213651\\
95	0.0247595066276113\\
96	0.0246976060456202\\
97	0.024475062296968\\
98	0.0241270313928631\\
99	0.0238232845900694\\
100	0.023712879463039\\
101	0.0238231572515837\\
102	0.0239706176402301\\
103	0.0240974047459988\\
104	0.0242385593082335\\
105	0.0244188758643754\\
106	0.0245175491653298\\
107	0.0245559106649395\\
108	0.0245503353683323\\
109	0.0244607551041207\\
110	0.0244629688549538\\
111	0.024547357126785\\
112	0.0245820497948582\\
113	0.0246106865360377\\
114	0.024596163733524\\
115	0.024564279713517\\
116	0.0244159333863392\\
117	0.0242018575679614\\
118	0.0240233716801731\\
119	0.023959571790115\\
120	0.024009588613786\\
121	0.024182083073353\\
122	0.024511362699857\\
123	0.0248549589839632\\
124	0.0251054574735255\\
125	0.0252090928717771\\
126	0.0252564147179973\\
127	0.0252412767468016\\
128	0.0250535165352005\\
129	0.0248182796959041\\
130	0.0245416949278607\\
131	0.0244059972957139\\
132	0.0243868743617761\\
133	0.0245371709740179\\
134	0.0247004453193578\\
135	0.0248325027034049\\
136	0.0249251177142095\\
137	0.0250046592202923\\
138	0.0250668398356296\\
139	0.0250585407114543\\
140	0.0249683044660525\\
141	0.0248224544469749\\
142	0.0247016999976028\\
143	0.0246291982796846\\
144	0.0246058806269341\\
145	0.0246196183727268\\
146	0.0246769662133379\\
147	0.024759515856637\\
148	0.0248796298102867\\
149	0.0249980340100901\\
150	0.0250858067251842\\
151	0.0251304318541296\\
152	0.0251122954096422\\
153	0.0250248993613113\\
154	0.0249297145016092\\
155	0.0248399306877928\\
156	0.0247763383251503\\
157	0.0247470015319081\\
158	0.0247316060692999\\
159	0.0247498099607626\\
160	0.0247903221356684\\
161	0.0248641062071422\\
162	0.0249267226036687\\
163	0.0249848788350497\\
164	0.0250285056156507\\
165	0.0250526144946969\\
166	0.0250551481415448\\
167	0.0250387813471677\\
168	0.0250061150921812\\
169	0.0249540224757686\\
170	0.0249212446128596\\
171	0.024891718206305\\
172	0.0248796528774224\\
173	0.0248859296352238\\
174	0.0248932102211952\\
175	0.0249045484726624\\
176	0.0249134081952263\\
177	0.0249169094363733\\
178	0.0249217801880442\\
179	0.0249308899570785\\
180	0.0249412346055796\\
181	0.0249510699512418\\
182	0.0249631788683154\\
183	0.0249771500526436\\
184	0.0249883824111847\\
185	0.0249946545891709\\
186	0.0250005152317608\\
187	0.025004166548124\\
188	0.0250057224483809\\
189	0.0250024829360024\\
190	0.0249886085421812\\
191	0.0249680958802633\\
192	0.0249413315420043\\
193	0.0249163643525008\\
194	0.0248888031939115\\
195	0.0248698898091335\\
196	0.0248622638989584\\
197	0.0248659807708632\\
198	0.0248797273606968\\
199	0.024903832293523\\
200	0.0358750608118796\\
201	0.033415801100153\\
202	0.0310883764500864\\
203	0.0278749478964072\\
204	0.0245517285492922\\
205	0.0225192592468248\\
206	0.0213526957181208\\
207	0.0209454564719446\\
208	0.0214070582981034\\
209	0.0221415928616338\\
210	0.023645977775633\\
211	0.0251370660306332\\
212	0.0259962701766395\\
213	0.026638267503204\\
214	0.0264830473885251\\
215	0.0261413527970966\\
216	0.0254169500720075\\
217	0.0246748971919239\\
218	0.0244247178754691\\
219	0.0244386216698097\\
220	0.0246281498112033\\
221	0.0246507086632348\\
222	0.0246248206880614\\
223	0.0247372302255975\\
224	0.0248639108022112\\
225	0.0250174295504031\\
226	0.0249888028054023\\
227	0.0249068379028387\\
228	0.0248238184758127\\
229	0.0247793488133743\\
230	0.0248184420300186\\
231	0.0249257667079414\\
232	0.0250786985321891\\
233	0.0252323805401262\\
234	0.0253541730046436\\
235	0.0254255697913096\\
236	0.0254407214212063\\
237	0.0254049576647335\\
238	0.0253153540587658\\
239	0.0251766650434363\\
240	0.0250092652880547\\
241	0.0248375998081622\\
242	0.02468577181241\\
243	0.0245624491686361\\
244	0.0244742064515903\\
245	0.0244230540468049\\
246	0.0244094848756858\\
247	0.0244334471524223\\
248	0.024490775596117\\
249	0.0245748136728506\\
250	0.0246836647826474\\
251	0.0248168516940655\\
252	0.0249668626036228\\
253	0.0251327032887242\\
254	0.0252962713120372\\
255	0.0254385977216508\\
256	0.0255531270082904\\
257	0.0256519061853057\\
258	0.0257187519775529\\
259	0.0257646508472741\\
260	0.0257674459656051\\
261	0.0257303644694299\\
262	0.0255975776891865\\
263	0.0253946847990195\\
264	0.0251530343586744\\
265	0.0249027813305161\\
266	0.0246891107692016\\
267	0.0245107718520798\\
268	0.0243734311659462\\
269	0.0242751170594089\\
270	0.0242319663165429\\
271	0.0242430330803971\\
272	0.024298058492457\\
273	0.0243878457446293\\
274	0.0245164978356256\\
275	0.0246682011364743\\
276	0.0248322495218235\\
277	0.0249979720053391\\
278	0.0251558125022761\\
279	0.0252940971120684\\
280	0.0254012008154685\\
281	0.0254720789094844\\
282	0.0255324288634625\\
283	0.0255692878111154\\
284	0.025569938249593\\
285	0.0255239905540741\\
286	0.0254457016710566\\
287	0.0253386615794255\\
288	0.0252034448786567\\
289	0.0250505129807515\\
290	0.024877544832549\\
291	0.0247030885636189\\
292	0.0245450449791696\\
293	0.0244200274735046\\
294	0.0243478409636093\\
295	0.0243420464174433\\
296	0.0244046735098923\\
297	0.0245283676431253\\
298	0.0246995326249051\\
299	0.0249006529527872\\
300	0.0251163138416082\\
301	0.0252586525731362\\
302	0.0253766797878631\\
303	0.0254766335073606\\
304	0.0255476995568413\\
305	0.0255642917524959\\
306	0.0255435327260706\\
307	0.0255018952151337\\
308	0.0254294899005131\\
309	0.0253323707379182\\
310	0.0252254052383137\\
311	0.0251309396039749\\
312	0.0250418257332879\\
313	0.0249658549968945\\
314	0.0248659098929355\\
315	0.0247498631603701\\
316	0.0246324119344505\\
317	0.0245253843213715\\
318	0.0244478557700937\\
319	0.024401739906072\\
320	0.0244039933773838\\
321	0.0244490463272353\\
322	0.0245176656879565\\
323	0.0245994928640149\\
324	0.0246888800715557\\
325	0.0247791693780335\\
326	0.0248711908412507\\
327	0.0249600672067371\\
328	0.0250367741423534\\
329	0.0250949092041412\\
330	0.0251486636066347\\
331	0.0251983852128975\\
332	0.0252426050022871\\
333	0.0252797296635788\\
334	0.0253080988366111\\
335	0.0253270198528596\\
336	0.0253357689402895\\
337	0.0253333112705653\\
338	0.0253160109239419\\
339	0.0252868039385529\\
340	0.0252505457803647\\
341	0.0252074404318557\\
342	0.025160007770472\\
343	0.0251122739148107\\
344	0.0250655745633924\\
345	0.0250211292068486\\
346	0.0249749716384948\\
347	0.024930139399425\\
348	0.0248879772150246\\
349	0.0248480788454237\\
350	0.0248077384316637\\
351	0.0247676300910589\\
352	0.0247264201685522\\
353	0.0246871624984027\\
354	0.0246578216611606\\
355	0.0246409731195569\\
356	0.0246383649380672\\
357	0.0246501026686689\\
358	0.0246767930117012\\
359	0.0247172891677\\
360	0.0247698947390733\\
361	0.024832094568652\\
362	0.0249054009571807\\
363	0.024986006494101\\
364	0.0250702939785312\\
365	0.0251542723372055\\
366	0.0252331932544884\\
367	0.0253040119362747\\
368	0.0253621303243495\\
369	0.0254055798409677\\
370	0.0254322373407312\\
371	0.0254426046783871\\
372	0.0254376366081143\\
373	0.0254191486533539\\
374	0.0253879339047836\\
375	0.025346788422801\\
376	0.025300673697522\\
377	0.0252497640166058\\
378	0.0251813056090524\\
379	0.0250980954110843\\
380	0.0250048316936505\\
381	0.0249082525637266\\
382	0.0248076441401091\\
383	0.0247076229600409\\
384	0.024613543219941\\
385	0.024527897650852\\
386	0.0244503729652805\\
387	0.024385037916645\\
388	0.0243399828945806\\
389	0.024314067011748\\
390	0.0243379103635122\\
391	0.0243677170555327\\
392	0.0244206035332849\\
393	0.0244920219957025\\
394	0.024602259219082\\
395	0.0247461871955045\\
396	0.0249162471666133\\
397	0.0251008683815346\\
398	0.0252922254889481\\
399	0.0254760519098752\\
400	0.0256431841926796\\
401	0.0257791162402105\\
};
\addlegendentry{exact}

\addplot [color=black, dashed, forget plot]
  table[row sep=crcr]{%
1	0.025\\
400	0.025\\
};

\end{axis}
\end{tikzpicture}%

%% file: mma_32_cc.tex
\definecolor{mycolor1}{rgb}{0.00000,0.44700,0.74100}%
\definecolor{mycolor2}{rgb}{0.85000,0.32500,0.09800}%
\definecolor{mycolor3}{rgb}{0.49400,0.18400,0.55600}%
\begin{tikzpicture}

\begin{axis}[%
width=\textwidth,
height=.75\textwidth,
at={(0.758in,0.481in)},
%scale only axis,
xmin=0,
xmax=400,
ymin=0,
ymax=1,
axis background/.style={fill=white},
title style={font=\bfseries},
title={MMA, $\mathcal{B}=32$},
every axis title/.style={at={(0,1)}, anchor=north west, draw=black, fill=white},
xmajorgrids,
ymajorgrids,
xminorgrids,
yminorgrids,
xlabel={Iteration},
ytick={0.025,0.2,0.4,0.6,0.8,1},
yticklabels={$p$,0.2,0.4,0.6,0.8,1},
legend style={legend cell align=left, align=left, draw=white!15!black}
]
\addplot [color=mycolor1, line width=1.5pt,line join = round]
  table[row sep=crcr]{%
1	1.00331710093093\\
2	0.813842026956059\\
3	0.0634856019003432\\
4	0.0352190031658391\\
5	0.0276155745722008\\
6	0.0376178631793569\\
7	0.362362486401247\\
8	0.328725300141496\\
9	0.245509766912179\\
10	0.152448686941307\\
11	0.072867930902628\\
12	0.0664709770260413\\
13	0.0652643215279396\\
14	0.0614310784378929\\
15	0.0235913669661241\\
16	0.0154630299619751\\
17	0.0239304079494352\\
18	0.0303021808819032\\
19	0.0262312146345507\\
20	0.024966211509688\\
21	0.024707382375658\\
22	0.0247841178897171\\
23	0.0247843419111666\\
24	0.0247496698730561\\
25	0.0246098925349461\\
26	0.0244236322356453\\
27	0.0245373558298069\\
28	0.0245852785000632\\
29	0.0245896752335312\\
30	0.024627115238196\\
31	0.0246381675869587\\
32	0.0245980526927745\\
33	0.0244499838261092\\
34	0.0241292723531832\\
35	0.0259344500526044\\
36	0.0269326514875352\\
37	0.0253473804465196\\
38	0.0244559350268363\\
39	0.0243862750358728\\
40	0.0241783471251992\\
41	0.0241145959798206\\
42	0.0241506742820867\\
43	0.0241194488569804\\
44	0.0241669151961835\\
45	0.0242374412688189\\
46	0.0241034393662183\\
47	0.0240346116629952\\
48	0.0241830455128307\\
49	0.0242580529676279\\
50	0.0242183260201542\\
51	0.0241742658356633\\
52	0.0242053920825722\\
53	0.0242649641181559\\
54	0.0241988541320681\\
55	0.0241067871446039\\
56	0.0240851755099831\\
57	0.0240861330520623\\
58	0.0241395421680988\\
59	0.0243580454781732\\
60	0.0246235068019951\\
61	0.0247040878832235\\
62	0.0246245194238517\\
63	0.0243823567532665\\
64	0.0242086374891387\\
65	0.0241409547683709\\
66	0.0241949506452751\\
67	0.0242522974168462\\
68	0.0244142644045754\\
69	0.0244803525162323\\
70	0.0244764110060519\\
71	0.02437746933311\\
72	0.024415654443001\\
73	0.0244893033768559\\
74	0.0245613188343552\\
75	0.0245621986039906\\
76	0.0245523010064049\\
77	0.0245567514476481\\
78	0.0245686785411053\\
79	0.0245743752866183\\
80	0.0245635512192074\\
81	0.024526462780194\\
82	0.0245061723051736\\
83	0.0245367730461843\\
84	0.0245890828586877\\
85	0.0246331793464341\\
86	0.0246390036203307\\
87	0.0246591536963851\\
88	0.0246960742539047\\
89	0.0247146207221019\\
90	0.0247271867577169\\
91	0.024755233666362\\
92	0.0247778936553962\\
93	0.0247964751166641\\
94	0.0248133294358179\\
95	0.0248298613459554\\
96	0.0248398279241458\\
97	0.0248512478628465\\
98	0.0248478743519096\\
99	0.0248534914380492\\
100	0.0248598902668787\\
101	0.0248599222419949\\
102	0.0248693788540864\\
103	0.0248669717646723\\
104	0.0248670761785464\\
105	0.0248678643245334\\
106	0.0248722583030568\\
107	0.024881067060822\\
108	0.0248920256904766\\
109	0.0248953929300526\\
110	0.0249012773974293\\
111	0.0249084491088597\\
112	0.0249154987570621\\
113	0.0249172210744648\\
114	0.0249182459428152\\
115	0.0249204661459823\\
116	0.0249212828341863\\
117	0.0249274661659949\\
118	0.0249309863907844\\
119	0.024934259416452\\
120	0.0249406659179076\\
121	0.0249432608329702\\
122	0.0249436070373654\\
123	0.0249490893801332\\
124	0.0249506381732167\\
125	0.0249521346362127\\
126	0.0249520769187835\\
127	0.0249535416559895\\
128	0.0249563872541622\\
129	0.0249583136096589\\
130	0.0249597496291709\\
131	0.0249631308527988\\
132	0.0249645202245958\\
133	0.0249679388004601\\
134	0.02496729102388\\
135	0.0249674320684235\\
136	0.0249673077322925\\
137	0.0249690970533957\\
138	0.0249731416272694\\
139	0.0249744666973843\\
140	0.0249747054337268\\
141	0.0249753468853426\\
142	0.0249764405229718\\
143	0.0249760993720146\\
144	0.0249767978835435\\
145	0.024977919637539\\
146	0.024979529485009\\
147	0.024981923134969\\
148	0.0249840060599159\\
149	0.0249835354221844\\
150	0.0249829781865408\\
151	0.0249824888893724\\
152	0.0249840436680738\\
153	0.0249855243506934\\
154	0.0249860420628727\\
155	0.0249870788320015\\
156	0.0249882809510686\\
157	0.0249887620314669\\
158	0.024988383772785\\
159	0.0249888061939178\\
160	0.0249892490574483\\
161	0.024990295011265\\
162	0.0249909960428655\\
163	0.0249911250500538\\
164	0.0249915676681334\\
165	0.024992185174268\\
166	0.0249922799555302\\
167	0.0249918681399194\\
168	0.0249916178125713\\
169	0.0249917564705226\\
170	0.0249916791539445\\
171	0.0249920363498129\\
172	0.0249921418118639\\
173	0.024991956371479\\
174	0.0249919343013539\\
175	0.0249923031902048\\
176	0.0249923238141781\\
177	0.0249926275388387\\
178	0.0249929605214899\\
179	0.0249933064603281\\
180	0.024993505935347\\
181	0.0249937803558193\\
182	0.0249939451511932\\
183	0.0249943112870825\\
184	0.024994514801299\\
185	0.024994686006028\\
186	0.024994630334098\\
187	0.0249945219211595\\
188	0.0249945797176579\\
189	0.0249947443567603\\
190	0.0249949392770733\\
191	0.0249949811563571\\
192	0.0249949228622673\\
193	0.024994926702848\\
194	0.0249951116364803\\
195	0.0249949292124002\\
196	0.0249949399038341\\
197	0.0249949625620313\\
198	0.0249947446934537\\
199	0.0249944881988739\\
200	0.0339130437751975\\
201	0.0245357706223141\\
202	0.0249181960400813\\
203	0.0249618011641458\\
204	0.0249745288879251\\
205	0.0249823013974068\\
206	0.0249869753026123\\
207	0.0249897599958017\\
208	0.0249915011230205\\
209	0.0249928453542151\\
210	0.0249937503650502\\
211	0.0249941934230619\\
212	0.0249941558300476\\
213	0.0249940076967787\\
214	0.024993639793566\\
215	0.0249936580708333\\
216	0.0249941204075019\\
217	0.0249944756657674\\
218	0.0249946115210915\\
219	0.0249947807726599\\
220	0.0249948478445088\\
221	0.0249950776544323\\
222	0.0249954808059944\\
223	0.0249959066867075\\
224	0.0249961315620331\\
225	0.0249960381860274\\
226	0.024995735117855\\
227	0.0249956205061473\\
228	0.0249959376548883\\
229	0.0249961871937063\\
230	0.0249965837018278\\
231	0.0249967834345029\\
232	0.0249969054592247\\
233	0.0249970958433001\\
234	0.0249969283533157\\
235	0.0249967935569464\\
236	0.0249968788679849\\
237	0.0249967672831299\\
238	0.0249968914135469\\
239	0.0249971086194826\\
240	0.024997397732436\\
241	0.0249973142555892\\
242	0.0249974797346276\\
243	0.0249974754196393\\
244	0.0249973530226861\\
245	0.02499751849065\\
246	0.0249975724557092\\
247	0.0249975437419419\\
248	0.0249976659058093\\
249	0.0249978226863109\\
250	0.0249976785640001\\
251	0.0249974961697324\\
252	0.0249973698533608\\
253	0.0249971033550932\\
254	0.0249967896644481\\
255	0.0249969851676062\\
256	0.0249971381239301\\
257	0.0249972417035465\\
258	0.0249972650496471\\
259	0.0249973469498553\\
260	0.0249973350368847\\
261	0.0249973609703729\\
262	0.0249975781197661\\
263	0.02499777677019\\
264	0.0249978965187291\\
265	0.0249979453175808\\
266	0.0249980272915217\\
267	0.024997952079892\\
268	0.0249977757748038\\
269	0.0249977233519291\\
270	0.0249978336719093\\
271	0.0249976759816337\\
272	0.0249975408102445\\
273	0.0249974731705111\\
274	0.0249974646945715\\
275	0.0249975356947144\\
276	0.024997707997095\\
277	0.0249976401567933\\
278	0.0249975094371338\\
279	0.0249974574453295\\
280	0.0249972721352637\\
281	0.0249972652951861\\
282	0.0249975634314785\\
283	0.024997715059064\\
284	0.0249976268443248\\
285	0.0249977011237389\\
286	0.0249977720543243\\
287	0.0249979313243637\\
288	0.0249979961123109\\
289	0.0249978571447386\\
290	0.0249977241411852\\
291	0.0249975913810249\\
292	0.0249974781487156\\
293	0.0249973593843395\\
294	0.0249971273185461\\
295	0.024996855144513\\
296	0.0249967788603958\\
297	0.024996590101791\\
298	0.0249962832137866\\
299	0.024996113844642\\
300	0.0249960677120278\\
301	0.0249963905038104\\
302	0.0249968402241211\\
303	0.0249969747613081\\
304	0.0249971940171607\\
305	0.0249972286549183\\
306	0.0249970459492622\\
307	0.0249969606948162\\
308	0.024997073960309\\
309	0.0249970945403203\\
310	0.0249971130254171\\
311	0.0249970859972183\\
312	0.0249969622326262\\
313	0.0249967263159828\\
314	0.0249966580304481\\
315	0.0249967892967085\\
316	0.0249965502419283\\
317	0.0249965511705063\\
318	0.0249966936333183\\
319	0.0249965481436145\\
320	0.0249964210671015\\
321	0.0249965508383899\\
322	0.0249964216861568\\
323	0.0249962908263494\\
324	0.0249965584281603\\
325	0.0249967930597585\\
326	0.0249969154818829\\
327	0.0249970544661966\\
328	0.0249970467535206\\
329	0.0249970203236756\\
330	0.0249973054234697\\
331	0.0249974383496499\\
332	0.0249974299827323\\
333	0.0249975668661003\\
334	0.0249978685423967\\
335	0.0249979213866789\\
336	0.0249978731850925\\
337	0.0249978335636414\\
338	0.0249977498211481\\
339	0.0249975556065141\\
340	0.0249975840308414\\
341	0.0249977878201468\\
342	0.0249977233402073\\
343	0.0249975974082892\\
344	0.0249974894449489\\
345	0.024997277264858\\
346	0.0249970727727276\\
347	0.0249970772765166\\
348	0.0249970279840387\\
349	0.0249970580668527\\
350	0.0249969933624937\\
351	0.0249971591787899\\
352	0.0249973264314387\\
353	0.024997675213045\\
354	0.0249978574079924\\
355	0.0249977567795188\\
356	0.0249977104438543\\
357	0.0249977848106025\\
358	0.0249979061291042\\
359	0.0249978986355051\\
360	0.0249979615426676\\
361	0.0249980490781595\\
362	0.0249980990873249\\
363	0.0249981607411641\\
364	0.0249981486189678\\
365	0.0249981578937809\\
366	0.0249982336918373\\
367	0.0249982642031915\\
368	0.0249983305206276\\
369	0.0249983967944042\\
370	0.0249984103197577\\
371	0.0249984958915185\\
372	0.0249985157326368\\
373	0.0249984574552831\\
374	0.0249983994776727\\
375	0.0249983438182869\\
376	0.0249983260438611\\
377	0.0249983453622993\\
378	0.0249983481312469\\
379	0.0249983596671581\\
380	0.0249983817940891\\
381	0.0249982933122239\\
382	0.0249981680611004\\
383	0.0249980555968226\\
384	0.0249979161403141\\
385	0.0249978039098572\\
386	0.0249977645224314\\
387	0.0249978627002802\\
388	0.0249977588046875\\
389	0.0249976286610468\\
390	0.0249975635186636\\
391	0.0249974379180947\\
392	0.0249972985054142\\
393	0.0249973231856387\\
394	0.0249974188385613\\
395	0.0249974280788983\\
396	0.0249974006003642\\
397	0.0249972587941684\\
398	0.0249973191788337\\
399	0.0249974837005006\\
400	0.0249975865744331\\
};
\addlegendentry{internal}

\addplot [color=mycolor2, line width=1pt,line join = round]
  table[row sep=crcr]{%
1	1.00331519624177\\
2	0.813744301851071\\
3	0.0634587099166002\\
4	0.0352235383939465\\
5	0.0275894574814653\\
6	0.0377013767072254\\
7	0.362029132751556\\
8	0.326226434319796\\
9	0.244484878407236\\
10	0.141365429635167\\
11	0.0845665103889118\\
12	0.0639394120850963\\
13	0.0573111341044894\\
14	0.0486095778548943\\
15	0.0273921452027234\\
16	0.0197692186675315\\
17	0.0259015037095096\\
18	0.0315110120707312\\
19	0.0270616642373874\\
20	0.025622746421994\\
21	0.0251930055715869\\
22	0.0251505574221936\\
23	0.0250432527392066\\
24	0.0249183304787433\\
25	0.0247164517475819\\
26	0.0245020226304072\\
27	0.0246040909579914\\
28	0.0246595474991766\\
29	0.0246679731546358\\
30	0.0247210323313541\\
31	0.0247491111250656\\
32	0.0247059606741015\\
33	0.0245595104520335\\
34	0.0242686105033605\\
35	0.0260763043573553\\
36	0.0271055095098585\\
37	0.0254440932120355\\
38	0.0245433124730791\\
39	0.0244695562481048\\
40	0.0242574829087\\
41	0.0241887907304345\\
42	0.0242210142499726\\
43	0.0241843073849858\\
44	0.0242279835635561\\
45	0.024298183581753\\
46	0.0241659590962183\\
47	0.0240956691102409\\
48	0.0242413233196046\\
49	0.0243120403132809\\
50	0.0242696569934975\\
51	0.0242248336552149\\
52	0.0242576000298165\\
53	0.0243196008732683\\
54	0.024256707781581\\
55	0.0241659980079601\\
56	0.0241437603860597\\
57	0.0241462410504875\\
58	0.0241999171933129\\
59	0.0244224425735199\\
60	0.0246836716083779\\
61	0.0247696989076156\\
62	0.0246873233175582\\
63	0.0244475111294464\\
64	0.0242745860012104\\
65	0.0242065683069885\\
66	0.0242630035358856\\
67	0.0243219423671966\\
68	0.0244852653123948\\
69	0.0245512890875403\\
70	0.0245468923509915\\
71	0.0244481153550374\\
72	0.0244875484988259\\
73	0.0245618083294634\\
74	0.0246331914763891\\
75	0.0246328996297163\\
76	0.0246209116306112\\
77	0.0246236099957965\\
78	0.02463403120141\\
79	0.0246394372400627\\
80	0.0246289881667694\\
81	0.024593056442578\\
82	0.0245748443555156\\
83	0.0246082826902486\\
84	0.0246618777379856\\
85	0.0247076866899013\\
86	0.0247139533392861\\
87	0.024736350925888\\
88	0.0247744645622983\\
89	0.0247951605940543\\
90	0.0248087487242135\\
91	0.0248396458406844\\
92	0.0248644289832564\\
93	0.0248834946218337\\
94	0.0249008309166395\\
95	0.0249184340896754\\
96	0.0249288604068672\\
97	0.0249400173259674\\
98	0.0249362706207024\\
99	0.0249411285137801\\
100	0.0249465937616097\\
101	0.0249450814519152\\
102	0.0249547386077688\\
103	0.0249526065166269\\
104	0.0249546770120624\\
105	0.0249573126333238\\
106	0.0249636759442268\\
107	0.024974093489946\\
108	0.0249867895801768\\
109	0.0249917373593846\\
110	0.0249986718521814\\
111	0.0250065786962955\\
112	0.025013885581748\\
113	0.0250157093689178\\
114	0.0250165827812334\\
115	0.0250185176669392\\
116	0.0250193136878136\\
117	0.0250252601285508\\
118	0.0250284349264308\\
119	0.0250314281146784\\
120	0.025037809875433\\
121	0.0250404363020391\\
122	0.0250408883245341\\
123	0.0250465177426692\\
124	0.0250482851615294\\
125	0.0250500030260448\\
126	0.0250500033586118\\
127	0.0250517302774162\\
128	0.0250552186208287\\
129	0.0250579369265666\\
130	0.0250602986722628\\
131	0.0250646856197426\\
132	0.0250669795492493\\
133	0.0250711436050173\\
134	0.0250714107050709\\
135	0.0250724518167057\\
136	0.025073193446957\\
137	0.0250760441976529\\
138	0.0250816016208158\\
139	0.0250843707758945\\
140	0.0250861098021749\\
141	0.0250882519078213\\
142	0.0250908040861461\\
143	0.0250917720983619\\
144	0.0250936972312383\\
145	0.0250960830805642\\
146	0.0250990849543528\\
147	0.0251032298862239\\
148	0.0251069838399617\\
149	0.0251081565738205\\
150	0.0251091524370122\\
151	0.025110035567032\\
152	0.0251127456754644\\
153	0.0251151061971385\\
154	0.0251164130582791\\
155	0.0251184914811274\\
156	0.0251208772593022\\
157	0.0251223736085815\\
158	0.0251229688245715\\
159	0.0251245203492112\\
160	0.0251264466678512\\
161	0.0251289745930996\\
162	0.0251308398411618\\
163	0.0251321117679689\\
164	0.0251335610329628\\
165	0.025135016570617\\
166	0.025135799718989\\
167	0.0251360144253811\\
168	0.0251364326237394\\
169	0.0251373778587099\\
170	0.0251381615568468\\
171	0.0251393104775645\\
172	0.0251400171855682\\
173	0.0251403142259166\\
174	0.0251406402269492\\
175	0.0251412623742703\\
176	0.0251414864939327\\
177	0.0251419628666948\\
178	0.0251422974845181\\
179	0.0251425739633742\\
180	0.0251429101796475\\
181	0.0251435102872271\\
182	0.0251439310617578\\
183	0.0251445874086597\\
184	0.0251450388546202\\
185	0.0251456120435296\\
186	0.0251459897924487\\
187	0.0251464416813064\\
188	0.0251472166630376\\
189	0.0251479985402842\\
190	0.0251486409287543\\
191	0.0251490518895777\\
192	0.0251494601920643\\
193	0.0251499892942418\\
194	0.0251506410750134\\
195	0.0251510081186638\\
196	0.0251516039713739\\
197	0.0251522216517808\\
198	0.0251526744366586\\
199	0.025153131836393\\
200	0.0341725670246206\\
201	0.0246919464616597\\
202	0.0250820119146997\\
203	0.025128206441805\\
204	0.0251423829557937\\
205	0.0251510408978122\\
206	0.0251556074151495\\
207	0.0251577909342281\\
208	0.0251589812686987\\
209	0.0251599739118457\\
210	0.0251606722059486\\
211	0.0251611665981815\\
212	0.0251613755638775\\
213	0.0251615825597469\\
214	0.0251615070799549\\
215	0.0251618216928665\\
216	0.0251625189347309\\
217	0.0251630373715445\\
218	0.0251633098203568\\
219	0.0251635535471584\\
220	0.0251636304927783\\
221	0.0251638317923068\\
222	0.0251642657704624\\
223	0.0251648664081234\\
224	0.0251655290101272\\
225	0.0251659758994455\\
226	0.0251662367012317\\
227	0.0251667654673711\\
228	0.025167807327291\\
229	0.0251688206905555\\
230	0.0251699879904375\\
231	0.0251709165889906\\
232	0.0251716162264232\\
233	0.0251722091023675\\
234	0.0251724114673532\\
235	0.0251726541379904\\
236	0.0251732040756169\\
237	0.0251736462574087\\
238	0.0251743162378956\\
239	0.0251748094078476\\
240	0.0251752186480701\\
241	0.0251750002467976\\
242	0.0251748232895207\\
243	0.0251742953028685\\
244	0.0251736505974647\\
245	0.0251738537102104\\
246	0.0251743215933054\\
247	0.0251747505976807\\
248	0.0251752830615695\\
249	0.0251756560660347\\
250	0.0251756147577859\\
251	0.0251757785014536\\
252	0.0251762179242286\\
253	0.0251766011711487\\
254	0.0251769376863106\\
255	0.0251777474524643\\
256	0.0251785200189667\\
257	0.0251792031700666\\
258	0.0251797002128223\\
259	0.0251801917686028\\
260	0.0251805846481069\\
261	0.0251810073459666\\
262	0.0251815506867738\\
263	0.0251820340657383\\
264	0.0251823925067069\\
265	0.0251826516924365\\
266	0.0251829140228286\\
267	0.0251830246701086\\
268	0.0251830354343917\\
269	0.0251831798034668\\
270	0.025183512859467\\
271	0.0251835844027298\\
272	0.0251836587520187\\
273	0.0251838043944825\\
274	0.0251840724301781\\
275	0.0251843180450648\\
276	0.0251845010132455\\
277	0.0251843685504244\\
278	0.0251841569098739\\
279	0.0251840383662032\\
280	0.0251837886784271\\
281	0.0251836417222214\\
282	0.025183613656441\\
283	0.0251831844793234\\
284	0.025182437895223\\
285	0.0251821853984604\\
286	0.0251823504809267\\
287	0.0251829099484578\\
288	0.0251833784416364\\
289	0.0251835710731947\\
290	0.0251837621401755\\
291	0.0251840153788709\\
292	0.0251843288756097\\
293	0.0251845774741291\\
294	0.0251846467020376\\
295	0.0251845969753882\\
296	0.0251848686790123\\
297	0.0251852031886456\\
298	0.0251853677646347\\
299	0.0251854622930473\\
300	0.0251854119486266\\
301	0.0251855113273616\\
302	0.0251856031621309\\
303	0.0251854408161537\\
304	0.0251854508965645\\
305	0.0251850872976133\\
306	0.0251843519960966\\
307	0.0251836242784728\\
308	0.0251831108306542\\
309	0.0251827896948372\\
310	0.0251827834974879\\
311	0.0251828758682545\\
312	0.025182869454074\\
313	0.0251827114702867\\
314	0.0251826788930027\\
315	0.0251829105292119\\
316	0.0251827716414372\\
317	0.025182795715602\\
318	0.0251829662150719\\
319	0.0251828380253488\\
320	0.025182726858766\\
321	0.0251828595977275\\
322	0.0251826798650948\\
323	0.0251825266672645\\
324	0.025182837450731\\
325	0.0251831546891976\\
326	0.0251834082269463\\
327	0.0251837635130027\\
328	0.0251841367505521\\
329	0.025184536892478\\
330	0.0251852000172791\\
331	0.0251856796484264\\
332	0.0251859956405286\\
333	0.0251864557000267\\
334	0.0251871439138307\\
335	0.0251875906287889\\
336	0.0251879232184995\\
337	0.0251881971106002\\
338	0.0251883454317913\\
339	0.0251883310811172\\
340	0.0251885200033336\\
341	0.0251888790897773\\
342	0.0251888997849273\\
343	0.0251887788123505\\
344	0.0251886332070321\\
345	0.0251884102955053\\
346	0.0251882554653768\\
347	0.0251882914117637\\
348	0.0251882055820434\\
349	0.025188163351583\\
350	0.0251879980100821\\
351	0.025187991790855\\
352	0.0251879091598044\\
353	0.0251880740190572\\
354	0.0251881727377597\\
355	0.0251879941742798\\
356	0.0251879064593824\\
357	0.0251880545994239\\
358	0.0251882865988735\\
359	0.0251882154907513\\
360	0.0251880958780345\\
361	0.0251879735459597\\
362	0.0251877927204716\\
363	0.025187591939306\\
364	0.0251873148240553\\
365	0.0251871492033622\\
366	0.0251871678485844\\
367	0.0251871273931428\\
368	0.0251870222083223\\
369	0.0251868561002723\\
370	0.0251867010597815\\
371	0.0251867215399535\\
372	0.0251867227725657\\
373	0.0251866594704739\\
374	0.0251866425895141\\
375	0.0251866845014954\\
376	0.0251868171003531\\
377	0.025186985474308\\
378	0.0251871140026483\\
379	0.0251872441310418\\
380	0.0251873780448487\\
381	0.0251874235820605\\
382	0.0251874608900383\\
383	0.0251875447277189\\
384	0.0251876243329447\\
385	0.025187739533943\\
386	0.0251879100832226\\
387	0.0251881783070206\\
388	0.0251882233811891\\
389	0.0251882108600238\\
390	0.0251882291072521\\
391	0.0251881464062324\\
392	0.0251880198579118\\
393	0.0251880215826824\\
394	0.0251880526567207\\
395	0.0251879733908633\\
396	0.0251878649749206\\
397	0.0251876127086051\\
398	0.0251875743914683\\
399	0.0251876699729469\\
400	0.0251877332531777\\
401	0.0251878531556103\\
};
\addlegendentry{exact}

\addplot [color=black, dashed, forget plot]
  table[row sep=crcr]{%
1	0.025\\
400	0.025\\
};

\end{axis}
\end{tikzpicture}%

%% file: mcmsa_16_1_smoothingeffect.tex
\definecolor{mycolor1}{rgb}{0.00000,0.44700,0.74100}%
\definecolor{mycolor2}{rgb}{0.85000,0.32500,0.09800}%
\definecolor{mycolor3}{rgb}{0.92900,0.69400,0.12500}%
\begin{tikzpicture}

\begin{axis}[%
width=\textwidth,
height=.85\textwidth,
at={(1.517in,0.962in)},
%scale only axis,
% title style={font=\bfseries},
% title={sMMA},
xmin=0,
xmax=400,
ymin=0,
ymax=0.06,
axis background/.style={fill=white},
xmajorgrids,
ymajorgrids,
legend style={legend cell align=left, align=left, draw=white!15!black},
xlabel={Iteration (\textbf{sMMA})},
yticklabel style={/pgf/number format/fixed},
scaled y ticks=false,
]
\addplot [color=mycolor1,line width=3pt,line join=round]
  table[row sep=crcr]{%
1	1\\
2	1\\
3	0\\
401	0\\
};
\addlegendentry{nonsmooth}

\addplot [color=mycolor2,line width=1.5pt,line join=round]
  table[row sep=crcr]{%
1	0.999156650065644\\
2	0.816502462687947\\
3	0.0193457379889739\\
4	0.0042730144777034\\
5	0.00366189521430814\\
6	0.0720138886293393\\
7	0.0113602565113961\\
8	0.00753945003053858\\
9	0.0120958297215644\\
10	0.0143083985378926\\
11	0.0407407427113716\\
12	0.0212459214243871\\
13	0.0261372926483857\\
14	0.0277386930442532\\
15	0.0305579593849642\\
16	0.0315496679670867\\
17	0.0304302173267338\\
18	0.0277606887583592\\
19	0.025498230203041\\
20	0.0239356684524695\\
21	0.0223444895734375\\
22	0.0215468685951184\\
23	0.0221320948075971\\
24	0.0234812423813105\\
25	0.026069036092393\\
26	0.0269342737143519\\
27	0.0275744722052684\\
28	0.0283224549352635\\
29	0.0284474224622537\\
30	0.0279575960145816\\
31	0.0269391700827661\\
32	0.0257823644367098\\
33	0.0245981019464081\\
34	0.0226616815961991\\
35	0.0189522260279323\\
36	0.0164943551663191\\
37	0.0138303253588005\\
38	0.0142656755837457\\
39	0.0203367147103799\\
40	0.0300147179372664\\
41	0.037156987297024\\
42	0.035916423090969\\
43	0.0306904677622891\\
44	0.0257518707378363\\
45	0.0227750468527405\\
46	0.0224502291475481\\
47	0.0240460107687295\\
48	0.0250705220939607\\
49	0.0259787769367642\\
50	0.0269843178535723\\
51	0.0276150540096917\\
52	0.0279917603486933\\
53	0.0273965787421904\\
54	0.0270801598267988\\
55	0.0265903048688391\\
56	0.0260344121139971\\
57	0.0255422843032704\\
58	0.0253191137309252\\
59	0.0256043278061908\\
60	0.0262025792551367\\
61	0.026681978430403\\
62	0.0269566961123906\\
63	0.0272022388378694\\
64	0.0272638248097194\\
65	0.0272057436680169\\
66	0.0269331180428646\\
67	0.0267096263459044\\
68	0.0264163361240532\\
69	0.0262114119110913\\
70	0.0260609523346084\\
71	0.0260196491043115\\
72	0.026073325245707\\
73	0.0261989095851836\\
74	0.0263586664537625\\
75	0.0265692336702163\\
76	0.0267644139214958\\
77	0.026850549451499\\
78	0.026810355941813\\
79	0.0267761798651099\\
80	0.0267016186771102\\
81	0.0265286133749205\\
82	0.0262992771497361\\
83	0.0261479711698423\\
84	0.0260644294419731\\
85	0.026091236191004\\
86	0.0261821094243055\\
87	0.026374087232864\\
88	0.0266340756342046\\
89	0.0268200634236734\\
90	0.0269303351166032\\
91	0.0270174799089387\\
92	0.0270409738862742\\
93	0.0270102675218175\\
94	0.0268434087804386\\
95	0.0266497688812712\\
96	0.0264866227013264\\
97	0.0263041915392165\\
98	0.0261191633787074\\
99	0.0260209674157289\\
100	0.0259971937832807\\
101	0.0260579009121641\\
102	0.0261477151837121\\
103	0.0263245298569105\\
104	0.0265756955301814\\
105	0.0268119256764978\\
106	0.0269590478030672\\
107	0.0270629680464651\\
108	0.0270919303721392\\
109	0.027073419884659\\
110	0.0269672563277007\\
111	0.0268363000482808\\
112	0.0267121720637743\\
113	0.02653743267141\\
114	0.0263341748067931\\
115	0.0262170101796888\\
116	0.0261694467961733\\
117	0.0262152689248369\\
118	0.026306097806317\\
119	0.0264726639061331\\
120	0.0266366806899079\\
121	0.0267558021783143\\
122	0.0268165241036514\\
123	0.0269173565904823\\
124	0.0269929187703497\\
125	0.0270203770291401\\
126	0.0269521388663847\\
127	0.0268365130493413\\
128	0.0267307979441032\\
129	0.0265950113971133\\
130	0.026426784738039\\
131	0.02631869737485\\
132	0.0262615692572608\\
133	0.0262820075090924\\
134	0.0262875786193778\\
135	0.0263694037395768\\
136	0.0265305070487906\\
137	0.0266997110216883\\
138	0.0268188169029876\\
139	0.026929312519887\\
140	0.0269847664913795\\
141	0.027004803697837\\
142	0.0269347249139657\\
143	0.026839005879753\\
144	0.0267623636400741\\
145	0.0266325160017134\\
146	0.0264948493834871\\
147	0.0263938457685535\\
148	0.0263309320485852\\
149	0.0262868804041034\\
150	0.0262673577882228\\
151	0.0262897876694353\\
152	0.026364798862889\\
153	0.0264401815025501\\
154	0.0264907995471882\\
155	0.0265707692508816\\
156	0.0266422627184336\\
157	0.0267235368486401\\
158	0.0267742524321868\\
159	0.0268374163372931\\
160	0.026907983538355\\
161	0.0269265544523689\\
162	0.0269043245739876\\
163	0.0268844161084771\\
164	0.0268296750755208\\
165	0.0267625276147267\\
166	0.0266370694889644\\
167	0.0265228389134977\\
168	0.0264156914811221\\
169	0.0263112377148365\\
170	0.0262219188214814\\
171	0.0261822229697307\\
172	0.0261961120707634\\
173	0.0262345521490506\\
174	0.0263110352326738\\
175	0.0264259454129331\\
176	0.0265853496540206\\
177	0.0267225576625183\\
178	0.0268351546382204\\
179	0.0269491700572251\\
180	0.027041678555998\\
181	0.0270803868691789\\
182	0.0270570681986478\\
183	0.0270004899652616\\
184	0.0269313530873228\\
185	0.0267927409781681\\
186	0.0266284932169536\\
187	0.0264956085565316\\
188	0.0263690576721995\\
189	0.0262680981075861\\
190	0.0261742271884882\\
191	0.0261423355440777\\
192	0.0261656109208409\\
193	0.0261854356216016\\
194	0.0262340304740403\\
195	0.0263072062043789\\
196	0.0263793378969659\\
197	0.026447011190602\\
198	0.0265536120671895\\
199	0.0266820899378948\\
200	0.0376742909474555\\
201	0.0356638541748447\\
202	0.0336172854070761\\
203	0.0310333041447172\\
204	0.0281313343890188\\
205	0.0255551698795821\\
206	0.0234118531520362\\
207	0.0220627685513152\\
208	0.0217147934044115\\
209	0.0223162982358444\\
210	0.0240188230654078\\
211	0.0257651701356366\\
212	0.0275527855202499\\
213	0.0279465926460357\\
214	0.0278179054060445\\
215	0.0276351059098313\\
216	0.0272103839901652\\
217	0.0267034809761238\\
218	0.0265546674581136\\
219	0.0266152122773168\\
220	0.0268086809807193\\
221	0.0270707808686322\\
222	0.0270757584816373\\
223	0.0271414908776219\\
224	0.0271544058750532\\
225	0.0271163433806778\\
226	0.0267375219755357\\
227	0.0263169886265757\\
228	0.0258935889796164\\
229	0.0255830475547719\\
230	0.0254364104299816\\
231	0.0254538596208159\\
232	0.0256121638836585\\
233	0.0258985135453369\\
234	0.0262267049710245\\
235	0.0265763903828193\\
236	0.0269114569627754\\
237	0.027186269993742\\
238	0.0273914297489951\\
239	0.0275354158363767\\
240	0.0276261704541815\\
241	0.0276161116955542\\
242	0.0275031797157548\\
243	0.0273447042488061\\
244	0.027132731242399\\
245	0.0268593933438551\\
246	0.0265381772695916\\
247	0.0262039093795802\\
248	0.0259380708722812\\
249	0.025772061835509\\
250	0.0256463317906307\\
251	0.0255553649340265\\
252	0.0255222269582339\\
253	0.0255725780766771\\
254	0.0256962769386282\\
255	0.0258886842655425\\
256	0.0261638371052814\\
257	0.0264638256067206\\
258	0.0267812983963655\\
259	0.0270455432864479\\
260	0.0272033782772042\\
261	0.0272750574031349\\
262	0.0272195148280835\\
263	0.0270738363035478\\
264	0.026840814845025\\
265	0.0265064594260158\\
266	0.0262451779208558\\
267	0.026058448907551\\
268	0.0259605491440276\\
269	0.0259444427474271\\
270	0.0260125809408402\\
271	0.0261395432380751\\
272	0.026297645160546\\
273	0.0264429031984441\\
274	0.0265807831504555\\
275	0.0267180460317185\\
276	0.0268369740465776\\
277	0.0269402730225872\\
278	0.0270171809859531\\
279	0.0270837847038255\\
280	0.0271294476690632\\
281	0.0271431120502094\\
282	0.027131980345996\\
283	0.0271064685129864\\
284	0.0270548575319477\\
285	0.0269600656091554\\
286	0.026824972512674\\
287	0.0266881704025253\\
288	0.0265333250217284\\
289	0.0263476270656078\\
290	0.0261730585191963\\
291	0.0260315491126739\\
292	0.0259457347136397\\
293	0.0259184985611426\\
294	0.025970810377626\\
295	0.0260936376811933\\
296	0.0262728053916176\\
297	0.0265103773955805\\
298	0.0267565767932067\\
299	0.0270086045274636\\
300	0.0272400374150042\\
301	0.0273591098386232\\
302	0.0274054662714091\\
303	0.0274109782198748\\
304	0.0273663598198036\\
305	0.0272251783323755\\
306	0.027051697717988\\
307	0.0268838253134703\\
308	0.0267021428828546\\
309	0.0265315469427245\\
310	0.0263729959931204\\
311	0.0262469665649696\\
312	0.0261485756210737\\
313	0.026093411098137\\
314	0.0260468949498442\\
315	0.026020429157297\\
316	0.0260289325348933\\
317	0.026075714317581\\
318	0.0261617739427807\\
319	0.0262799939636462\\
320	0.0264376124596131\\
321	0.0266149628864404\\
322	0.0267890775805029\\
323	0.0269554432385164\\
324	0.0270845584764453\\
325	0.0271838413634693\\
326	0.0272164427313916\\
327	0.0272100275319911\\
328	0.027153454559603\\
329	0.0270423976730092\\
330	0.0269063437723501\\
331	0.0267731777568749\\
332	0.0266443339175353\\
333	0.0265186506037984\\
334	0.0264028951287893\\
335	0.0263047252213451\\
336	0.0262294961207192\\
337	0.0261717304803606\\
338	0.0261482505680713\\
339	0.0261580129750521\\
340	0.0261966103889003\\
341	0.0262681452708835\\
342	0.0263616900253638\\
343	0.0264746320248675\\
344	0.0265945890456122\\
345	0.0267147249442671\\
346	0.0268351225303164\\
347	0.0269528188408722\\
348	0.0270587241366166\\
349	0.0271325396470345\\
350	0.027166474727329\\
351	0.027167626822685\\
352	0.0271297274417871\\
353	0.0270416830415424\\
354	0.0269022989957371\\
355	0.0267458157824877\\
356	0.0265643726282015\\
357	0.026385238933348\\
358	0.0262107907965942\\
359	0.0260656633921483\\
360	0.0259979498398091\\
361	0.0259874430805123\\
362	0.0260299351437396\\
363	0.0261191418507283\\
364	0.0262424227169743\\
365	0.0264008548597191\\
366	0.0265895312259592\\
367	0.0267965395193441\\
368	0.0270063694446962\\
369	0.0271446897247224\\
370	0.0272147582143962\\
371	0.027258847054451\\
372	0.0272599941888472\\
373	0.0272251687254291\\
374	0.0271433745252762\\
375	0.0270420943127877\\
376	0.0269428679652823\\
377	0.0268430512022523\\
378	0.0267293567353065\\
379	0.0265952005293082\\
380	0.0264436198037829\\
381	0.0262770836091713\\
382	0.0261334617812338\\
383	0.0260178478882364\\
384	0.0259535351384092\\
385	0.0259203520671565\\
386	0.0259418243115245\\
387	0.0259730242103582\\
388	0.0260212912374462\\
389	0.0260905728896152\\
390	0.0262079505010547\\
391	0.0263503659441473\\
392	0.0265033473265598\\
393	0.0265967975996799\\
394	0.0266995813847531\\
395	0.0268082343455391\\
396	0.026914216006762\\
397	0.0270090655450234\\
398	0.0270868012343996\\
399	0.0271502970800873\\
400	0.0272011944470444\\
401	0.0272152222343664\\
};
\addlegendentry{smooth}

\addplot [color=mycolor3,line width=1.5pt,line join=round]
  table[row sep=crcr]{%
1	1.00331519624177\\
2	0.817277239280557\\
3	0.0175739347990537\\
4	0.00191584040642272\\
5	0.00124893130820033\\
6	0.0708163371202274\\
7	0.0093743677773168\\
8	0.00539448610246812\\
9	0.0101341423346233\\
10	0.0124134164456548\\
11	0.0392843662872431\\
12	0.0195125906587405\\
13	0.0244905544798148\\
14	0.0261172127361559\\
15	0.028977638848721\\
16	0.029983232719482\\
17	0.0288485923745001\\
18	0.02614007534322\\
19	0.0238415828438225\\
20	0.0222523852293939\\
21	0.0206327650018054\\
22	0.0198200598739183\\
23	0.0204160445466995\\
24	0.0217900050712054\\
25	0.024421720594235\\
26	0.0253007507532161\\
27	0.025950812220769\\
28	0.0267100479376679\\
29	0.0268368261682328\\
30	0.0263395664970972\\
31	0.0253058206960862\\
32	0.0241303236955129\\
33	0.0229259313833535\\
34	0.0209554328060933\\
35	0.0171719378010867\\
36	0.0146565123513448\\
37	0.0119227665050947\\
38	0.0123689395518029\\
39	0.0185847630564384\\
40	0.028422660922743\\
41	0.0356510630674017\\
42	0.0343917198542466\\
43	0.0291067686151274\\
44	0.0240990412911144\\
45	0.0210702127439931\\
46	0.0207391018467359\\
47	0.022363355785673\\
48	0.0234053684250568\\
49	0.0243278798494925\\
50	0.0253510559882913\\
51	0.0259921780033776\\
52	0.0263745929545902\\
53	0.0257702275912129\\
54	0.0254490714227378\\
55	0.0249515378434057\\
56	0.024386707821839\\
57	0.0238863701560423\\
58	0.023659336297278\\
59	0.0239493340103837\\
60	0.0245574234252014\\
61	0.0250445095570516\\
62	0.025323532734849\\
63	0.0255728760318653\\
64	0.0256353772349363\\
65	0.0255764197742215\\
66	0.0252995389447065\\
67	0.0250725431254606\\
68	0.0247745773403051\\
69	0.0245663512344158\\
70	0.0244134695121892\\
71	0.0243714811994683\\
72	0.0244260109275675\\
73	0.024553626845828\\
74	0.0247159574876338\\
75	0.0249298892577418\\
76	0.0251281401606863\\
77	0.0252156103433981\\
78	0.0251747579913421\\
79	0.0251400117499301\\
80	0.0250642474471538\\
81	0.0248884891890003\\
82	0.0246554799922717\\
83	0.0245017345477521\\
84	0.0244168454770398\\
85	0.0244440840305711\\
86	0.0245364088298419\\
87	0.0247314475662143\\
88	0.0249955411070363\\
89	0.0251844309401427\\
90	0.0252964100940318\\
91	0.0253848960824397\\
92	0.0254087380738708\\
93	0.0253775352165378\\
94	0.0252080382589195\\
95	0.0250113223915681\\
96	0.0248455641202944\\
97	0.0246602112701183\\
98	0.0244721963248314\\
99	0.024372401488599\\
100	0.0243482320449974\\
101	0.0244099085450842\\
102	0.024501160619092\\
103	0.0246807948944066\\
104	0.0249359324624219\\
105	0.0251758653734846\\
106	0.0253252728300336\\
107	0.0254307864885233\\
108	0.0254601738965687\\
109	0.0254413521961566\\
110	0.0253335061542132\\
111	0.0252004708345473\\
112	0.0250743603117161\\
113	0.0248968340234645\\
114	0.0246903142587658\\
115	0.024571244380478\\
116	0.024522889033997\\
117	0.0245694155588752\\
118	0.0246616698404873\\
119	0.0248308599833817\\
120	0.0249974438740835\\
121	0.025118414233898\\
122	0.02518006558275\\
123	0.0252824489065737\\
124	0.025359164372657\\
125	0.0253870256564381\\
126	0.0253176952179094\\
127	0.025200232067708\\
128	0.0250928278349705\\
129	0.0249548766276821\\
130	0.02478395205551\\
131	0.0246741109820979\\
132	0.0246160390924755\\
133	0.0246367720930028\\
134	0.0246424036771059\\
135	0.0247255059700657\\
136	0.0248891427330082\\
137	0.0250610004185703\\
138	0.0251819632983014\\
139	0.0252941705214576\\
140	0.0253504727673473\\
141	0.0253708023212892\\
142	0.0252996084976007\\
143	0.0252023688521531\\
144	0.0251244966557485\\
145	0.0249925766516324\\
146	0.0248527041615231\\
147	0.0247500668023986\\
148	0.0246861205420633\\
149	0.0246413361993066\\
150	0.0246214696239675\\
151	0.0246442246009677\\
152	0.0247203978608165\\
153	0.0247969466113776\\
154	0.0248483393035512\\
155	0.0249295477751711\\
156	0.0250021477979886\\
157	0.0250846796890676\\
158	0.0251361721241365\\
159	0.0252003023522896\\
160	0.0252719482755589\\
161	0.0252907842580929\\
162	0.0252681783467097\\
163	0.0252479309486102\\
164	0.0251923097300118\\
165	0.0251240870159756\\
166	0.0249966320935931\\
167	0.0248805706155674\\
168	0.024771696092184\\
169	0.0246655485251763\\
170	0.024574770755415\\
171	0.0245344093139423\\
172	0.024548491962773\\
173	0.0245875212058292\\
174	0.0246651961934797\\
175	0.0247819048099962\\
176	0.0249437999111253\\
177	0.0250831376943235\\
178	0.0251974687666397\\
179	0.0253132358906177\\
180	0.0254071590210196\\
181	0.0254464489797401\\
182	0.0254227498966781\\
183	0.0253652756554492\\
184	0.0252950424195792\\
185	0.0251542445556947\\
186	0.0249873963295952\\
187	0.0248523924412165\\
188	0.0247238159543982\\
189	0.0246212304442035\\
190	0.0245258423010339\\
191	0.0244934225576607\\
192	0.024517054268499\\
193	0.0245371831347859\\
194	0.0245865369572037\\
195	0.0246608624947897\\
196	0.0247341264101514\\
197	0.0248028608743292\\
198	0.0249111300540069\\
199	0.0250416132049074\\
200	0.0361809251732859\\
201	0.03414685190299\\
202	0.0320748755577277\\
203	0.0294566889124462\\
204	0.0265130550679192\\
205	0.023896461295185\\
206	0.0217165488287877\\
207	0.0203428179695585\\
208	0.0199882412438962\\
209	0.0206010066422027\\
210	0.0223339863269817\\
211	0.0241095846087596\\
212	0.0259253208612034\\
213	0.0263251199601443\\
214	0.0261944951028044\\
215	0.0260089134254423\\
216	0.0255776657060302\\
217	0.0250628363909282\\
218	0.0249116462657398\\
219	0.024973120567169\\
220	0.0251696085304145\\
221	0.025435772776436\\
222	0.0254408097345794\\
223	0.0255075407934558\\
224	0.0255206366993947\\
225	0.025481971994475\\
226	0.0250972319205702\\
227	0.0246700340908191\\
228	0.0242398214255812\\
229	0.0239242094897115\\
230	0.0237751450733433\\
231	0.0237928587273251\\
232	0.0239537372830815\\
233	0.0242447218332026\\
234	0.0245781716342944\\
235	0.0249333944052107\\
236	0.0252737025561707\\
237	0.0255527675280495\\
238	0.0257610701324356\\
239	0.0259072430802671\\
240	0.0259993635402377\\
241	0.0259891344051465\\
242	0.0258744634889801\\
243	0.0257135451956922\\
244	0.0254982911485343\\
245	0.0252206921767348\\
246	0.0248944180416609\\
247	0.0245548275540257\\
248	0.0242847109553673\\
249	0.0241160056258928\\
250	0.02398821881169\\
251	0.0238957566457502\\
252	0.0238620668338042\\
253	0.0239132304516917\\
254	0.0240389334240212\\
255	0.024234446112017\\
256	0.0245140072443754\\
257	0.024818754887604\\
258	0.0251412077583453\\
259	0.0254095560304156\\
260	0.0255698213564401\\
261	0.0256425957995897\\
262	0.0255861903051036\\
263	0.0254382552766275\\
264	0.0252016045278857\\
265	0.0248619905398711\\
266	0.0245965504189784\\
267	0.0244068196355161\\
268	0.024307327131346\\
269	0.0242909415635311\\
270	0.0243601543673725\\
271	0.024489131717282\\
272	0.0246497361615535\\
273	0.0247972788662952\\
274	0.0249373173811499\\
275	0.0250767201395173\\
276	0.0251974940896595\\
277	0.0253023901415635\\
278	0.0253804813506234\\
279	0.0254481070042312\\
280	0.025494468200421\\
281	0.0255083368019746\\
282	0.0254970248670104\\
283	0.0254711098280684\\
284	0.0254186903173655\\
285	0.0253224174191158\\
286	0.0251852086242783\\
287	0.0250462530279251\\
288	0.0248889608925041\\
289	0.0247003148280183\\
290	0.024522956258558\\
291	0.0243791723529841\\
292	0.0242919695267246\\
293	0.0242642851867038\\
294	0.0243174277222833\\
295	0.0244422138035524\\
296	0.0246242262713974\\
297	0.0248655425686441\\
298	0.0251155891585098\\
299	0.0253715218901229\\
300	0.0256065101928544\\
301	0.0257273965407613\\
302	0.0257744505863573\\
303	0.0257800356169215\\
304	0.0257347255251848\\
305	0.0255913692366087\\
306	0.0254152036568312\\
307	0.0252447172723743\\
308	0.0250601907366271\\
309	0.0248869084803638\\
310	0.0247258473203504\\
311	0.0245978118734763\\
312	0.0244978470097753\\
313	0.0244417932665122\\
314	0.0243945232648147\\
315	0.0243676237941496\\
316	0.0243762535235398\\
317	0.0244237736711677\\
318	0.0245111972846405\\
319	0.0246312882203033\\
320	0.0247913930662151\\
321	0.0249715271287967\\
322	0.0251483580876269\\
323	0.0253173048283762\\
324	0.0254484121376228\\
325	0.0255492198456862\\
326	0.025582317803206\\
327	0.0255757996762383\\
328	0.0255183524783173\\
329	0.0254055790125061\\
330	0.0252674130311826\\
331	0.0251321695437606\\
332	0.0250013050745283\\
333	0.0248736409047367\\
334	0.0247560514190076\\
335	0.0246563191081163\\
336	0.02457988717461\\
337	0.0245211936305468\\
338	0.0244973320865471\\
339	0.0245072446263741\\
340	0.0245464528663698\\
341	0.024619121885337\\
342	0.0247141464537987\\
343	0.024828869121066\\
344	0.0249507100043395\\
345	0.0250727245824618\\
346	0.0251949962283886\\
347	0.0253145168327588\\
348	0.0254220562513573\\
349	0.0254970054726964\\
350	0.0255314572851515\\
351	0.0255326209663454\\
352	0.0254941318145029\\
353	0.025404722767387\\
354	0.0252631740091667\\
355	0.025104250777677\\
356	0.0249199631503873\\
357	0.0247380043002409\\
358	0.0245607881843103\\
359	0.0244133454584809\\
360	0.0243445463292446\\
361	0.0243338680070678\\
362	0.0243770356313078\\
363	0.0244676618003355\\
364	0.0245928967158603\\
365	0.024753827861689\\
366	0.0249454619257077\\
367	0.0251556924400444\\
368	0.0253687632983227\\
369	0.0255092034141651\\
370	0.0255803375996634\\
371	0.0256250943297699\\
372	0.0256262528814588\\
373	0.0255908888177711\\
374	0.0255078343652331\\
375	0.0254049907793243\\
376	0.0253042273060513\\
377	0.0252028586596574\\
378	0.0250873901869002\\
379	0.0249511342352637\\
380	0.0247971704684054\\
381	0.024628002873373\\
382	0.0244821006535146\\
383	0.0243646430122268\\
384	0.0242992998120156\\
385	0.024265583669021\\
386	0.0242873988882529\\
387	0.0243190980592372\\
388	0.0243681364844379\\
389	0.0244385229513239\\
390	0.0245577651217799\\
391	0.0247024312885103\\
392	0.0248578168645656\\
393	0.0249527277583832\\
394	0.0250571130453371\\
395	0.0251674532751548\\
396	0.0252750741235033\\
397	0.0253713846125678\\
398	0.0254503124358527\\
399	0.0255147785843419\\
400	0.025566451874783\\
401	0.0255806903160329\\
};
\addlegendentry{smooth \& steep}

% \addplot[dashed,line width=1pt,forget plot]
%   table[row sep=crcr]{%
% 0 0.025\\
% 401 0.025\\
% };

\end{axis}
\end{tikzpicture}%

%% file: mma_16_1_smoothingeffect.tex
\definecolor{mycolor1}{rgb}{0.00000,0.44700,0.74100}%
\definecolor{mycolor2}{rgb}{0.85000,0.32500,0.09800}%
\definecolor{mycolor3}{rgb}{0.92900,0.69400,0.12500}%
\begin{tikzpicture}

\begin{axis}[%
width=\textwidth,
height=.85\textwidth,
at={(1.517in,0.962in)},
%scale only axis,
% title style={font=\bfseries},
% title={MMA},
xmin=0,
xmax=400,
ymin=0.46,
ymax=0.54,
axis background/.style={fill=white},
xmajorgrids,
ymajorgrids,
legend style={legend cell align=left, align=left, draw=white!15!black},
xlabel={Iteration (\textbf{MMA})}
]

\addplot [color=mycolor1,line width=3pt,line join=round]
  table[row sep=crcr]{%
1	1\\
2	1\\
3	0\\
4	0\\
5	0\\
6	0.389814814814815\\
7	0.428703703703704\\
8	0.336111111111111\\
9	0.05\\
10	0\\
11	0\\
12	0\\
13	0\\
14	0\\
15	0\\
16	0\\
17	0\\
18	0\\
19	0\\
20	0\\
21	0.0037037037037037\\
22	0.037037037037037\\
23	0.0740740740740741\\
24	0.155555555555556\\
25	0.20462962962963\\
26	0.253703703703704\\
27	0.29537037037037\\
28	0.328703703703704\\
29	0.376851851851852\\
30	0.412962962962963\\
31	0.441666666666667\\
32	0.461111111111111\\
33	0.475925925925926\\
34	0.493518518518519\\
35	0.500925925925926\\
36	0.510185185185185\\
37	0.517592592592593\\
38	0.522222222222222\\
39	0.525\\
40	0.526851851851852\\
41	0.527777777777778\\
42	0.527777777777778\\
43	0.526851851851852\\
44	0.527777777777778\\
45	0.526851851851852\\
46	0.526851851851852\\
47	0.524074074074074\\
48	0.521296296296296\\
49	0.519444444444444\\
50	0.518518518518518\\
51	0.517592592592593\\
52	0.517592592592593\\
53	0.515740740740741\\
54	0.515740740740741\\
55	0.512962962962963\\
56	0.512962962962963\\
57	0.511111111111111\\
58	0.511111111111111\\
59	0.509259259259259\\
60	0.509259259259259\\
61	0.508333333333333\\
62	0.507407407407407\\
63	0.507407407407407\\
64	0.506481481481482\\
65	0.50462962962963\\
66	0.50462962962963\\
67	0.503703703703704\\
68	0.502777777777778\\
69	0.502777777777778\\
70	0.502777777777778\\
71	0.502777777777778\\
72	0.502777777777778\\
73	0.500925925925926\\
74	0.5\\
75	0.499074074074074\\
76	0.499074074074074\\
77	0.498148148148148\\
78	0.498148148148148\\
79	0.498148148148148\\
80	0.497222222222222\\
81	0.497222222222222\\
82	0.497222222222222\\
83	0.496296296296296\\
84	0.496296296296296\\
85	0.496296296296296\\
86	0.496296296296296\\
87	0.496296296296296\\
88	0.496296296296296\\
89	0.496296296296296\\
90	0.496296296296296\\
91	0.496296296296296\\
92	0.49537037037037\\
93	0.49537037037037\\
94	0.49537037037037\\
95	0.49537037037037\\
96	0.49537037037037\\
97	0.49537037037037\\
98	0.49537037037037\\
99	0.49537037037037\\
100	0.49537037037037\\
101	0.49537037037037\\
102	0.49537037037037\\
103	0.49537037037037\\
104	0.49537037037037\\
105	0.49537037037037\\
106	0.49537037037037\\
107	0.49537037037037\\
108	0.49537037037037\\
109	0.49537037037037\\
110	0.49537037037037\\
111	0.49537037037037\\
112	0.49537037037037\\
113	0.49537037037037\\
114	0.494444444444444\\
115	0.494444444444444\\
116	0.494444444444444\\
117	0.494444444444444\\
118	0.494444444444444\\
119	0.494444444444444\\
120	0.494444444444444\\
121	0.494444444444444\\
122	0.494444444444444\\
123	0.494444444444444\\
124	0.494444444444444\\
125	0.494444444444444\\
126	0.494444444444444\\
127	0.494444444444444\\
128	0.494444444444444\\
129	0.494444444444444\\
130	0.494444444444444\\
131	0.493518518518519\\
132	0.492592592592593\\
133	0.491666666666667\\
134	0.491666666666667\\
135	0.491666666666667\\
136	0.491666666666667\\
137	0.491666666666667\\
138	0.491666666666667\\
139	0.491666666666667\\
140	0.491666666666667\\
141	0.490740740740741\\
142	0.490740740740741\\
143	0.489814814814815\\
144	0.489814814814815\\
145	0.489814814814815\\
146	0.489814814814815\\
147	0.489814814814815\\
148	0.489814814814815\\
149	0.489814814814815\\
150	0.489814814814815\\
151	0.488888888888889\\
152	0.488888888888889\\
153	0.488888888888889\\
154	0.488888888888889\\
155	0.488888888888889\\
156	0.488888888888889\\
157	0.487962962962963\\
158	0.487962962962963\\
159	0.487962962962963\\
160	0.487962962962963\\
161	0.487962962962963\\
162	0.487962962962963\\
163	0.487962962962963\\
164	0.487962962962963\\
165	0.487962962962963\\
166	0.487037037037037\\
167	0.486111111111111\\
168	0.486111111111111\\
169	0.486111111111111\\
170	0.486111111111111\\
171	0.486111111111111\\
172	0.486111111111111\\
173	0.486111111111111\\
174	0.486111111111111\\
175	0.486111111111111\\
176	0.486111111111111\\
177	0.486111111111111\\
178	0.486111111111111\\
179	0.486111111111111\\
180	0.486111111111111\\
181	0.486111111111111\\
182	0.486111111111111\\
183	0.486111111111111\\
184	0.486111111111111\\
185	0.486111111111111\\
186	0.486111111111111\\
187	0.486111111111111\\
188	0.486111111111111\\
189	0.486111111111111\\
190	0.486111111111111\\
191	0.486111111111111\\
192	0.486111111111111\\
193	0.486111111111111\\
194	0.486111111111111\\
195	0.486111111111111\\
196	0.486111111111111\\
197	0.486111111111111\\
198	0.486111111111111\\
199	0.486111111111111\\
200	0.536111111111111\\
201	0.47962962962963\\
202	0.490740740740741\\
203	0.491666666666667\\
204	0.491666666666667\\
205	0.492592592592593\\
206	0.492592592592593\\
207	0.492592592592593\\
208	0.492592592592593\\
209	0.492592592592593\\
210	0.491666666666667\\
211	0.490740740740741\\
212	0.490740740740741\\
213	0.489814814814815\\
214	0.489814814814815\\
215	0.488888888888889\\
216	0.489814814814815\\
217	0.489814814814815\\
218	0.489814814814815\\
219	0.488888888888889\\
220	0.488888888888889\\
221	0.488888888888889\\
222	0.488888888888889\\
223	0.489814814814815\\
224	0.489814814814815\\
225	0.489814814814815\\
226	0.489814814814815\\
227	0.489814814814815\\
228	0.489814814814815\\
229	0.489814814814815\\
230	0.489814814814815\\
231	0.489814814814815\\
232	0.489814814814815\\
233	0.489814814814815\\
234	0.490740740740741\\
235	0.490740740740741\\
236	0.490740740740741\\
237	0.491666666666667\\
238	0.491666666666667\\
239	0.491666666666667\\
240	0.491666666666667\\
241	0.490740740740741\\
242	0.490740740740741\\
243	0.490740740740741\\
244	0.490740740740741\\
245	0.490740740740741\\
246	0.490740740740741\\
247	0.490740740740741\\
248	0.490740740740741\\
249	0.490740740740741\\
250	0.490740740740741\\
251	0.490740740740741\\
252	0.490740740740741\\
253	0.490740740740741\\
254	0.490740740740741\\
255	0.490740740740741\\
256	0.490740740740741\\
257	0.490740740740741\\
258	0.490740740740741\\
259	0.490740740740741\\
260	0.490740740740741\\
261	0.490740740740741\\
262	0.490740740740741\\
263	0.490740740740741\\
264	0.490740740740741\\
265	0.489814814814815\\
266	0.489814814814815\\
267	0.489814814814815\\
268	0.489814814814815\\
269	0.489814814814815\\
270	0.489814814814815\\
271	0.489814814814815\\
272	0.489814814814815\\
273	0.489814814814815\\
274	0.489814814814815\\
275	0.489814814814815\\
276	0.489814814814815\\
277	0.489814814814815\\
278	0.489814814814815\\
279	0.488888888888889\\
280	0.488888888888889\\
281	0.488888888888889\\
282	0.488888888888889\\
283	0.488888888888889\\
284	0.488888888888889\\
285	0.488888888888889\\
286	0.488888888888889\\
287	0.488888888888889\\
288	0.488888888888889\\
289	0.488888888888889\\
290	0.488888888888889\\
291	0.488888888888889\\
292	0.488888888888889\\
293	0.488888888888889\\
294	0.488888888888889\\
295	0.488888888888889\\
296	0.488888888888889\\
297	0.488888888888889\\
298	0.488888888888889\\
299	0.487962962962963\\
300	0.487962962962963\\
301	0.487962962962963\\
302	0.487962962962963\\
303	0.487962962962963\\
304	0.487962962962963\\
305	0.487962962962963\\
306	0.487962962962963\\
307	0.487962962962963\\
308	0.487962962962963\\
309	0.487962962962963\\
310	0.487962962962963\\
311	0.487037037037037\\
312	0.487037037037037\\
313	0.487037037037037\\
314	0.487037037037037\\
315	0.487037037037037\\
316	0.487037037037037\\
317	0.486111111111111\\
318	0.486111111111111\\
319	0.486111111111111\\
320	0.486111111111111\\
321	0.486111111111111\\
322	0.486111111111111\\
323	0.486111111111111\\
324	0.486111111111111\\
325	0.486111111111111\\
326	0.486111111111111\\
327	0.487037037037037\\
328	0.487037037037037\\
329	0.487037037037037\\
330	0.487037037037037\\
331	0.487037037037037\\
332	0.486111111111111\\
333	0.486111111111111\\
334	0.486111111111111\\
335	0.486111111111111\\
336	0.486111111111111\\
337	0.486111111111111\\
338	0.486111111111111\\
339	0.487037037037037\\
340	0.487037037037037\\
341	0.487037037037037\\
342	0.487037037037037\\
343	0.487037037037037\\
344	0.487037037037037\\
345	0.487037037037037\\
346	0.487037037037037\\
347	0.487037037037037\\
348	0.487037037037037\\
349	0.487037037037037\\
350	0.487962962962963\\
351	0.488888888888889\\
352	0.488888888888889\\
353	0.488888888888889\\
354	0.488888888888889\\
355	0.487962962962963\\
356	0.488888888888889\\
357	0.488888888888889\\
358	0.488888888888889\\
359	0.488888888888889\\
360	0.488888888888889\\
361	0.488888888888889\\
362	0.488888888888889\\
363	0.488888888888889\\
364	0.488888888888889\\
365	0.488888888888889\\
366	0.488888888888889\\
367	0.488888888888889\\
368	0.488888888888889\\
369	0.488888888888889\\
370	0.488888888888889\\
371	0.488888888888889\\
372	0.487962962962963\\
373	0.487037037037037\\
374	0.487037037037037\\
375	0.487962962962963\\
376	0.487962962962963\\
377	0.487962962962963\\
378	0.487962962962963\\
379	0.487962962962963\\
380	0.487962962962963\\
381	0.487037037037037\\
382	0.487037037037037\\
383	0.487037037037037\\
384	0.487037037037037\\
385	0.487037037037037\\
386	0.487037037037037\\
387	0.487037037037037\\
388	0.487037037037037\\
389	0.487037037037037\\
390	0.487037037037037\\
391	0.487962962962963\\
392	0.487962962962963\\
393	0.487962962962963\\
394	0.488888888888889\\
395	0.488888888888889\\
396	0.488888888888889\\
397	0.488888888888889\\
398	0.488888888888889\\
399	0.488888888888889\\
400	0.487962962962963\\
401	0.487962962962963\\
};
\addlegendentry{nonsmooth}

\addplot [color=mycolor2,line width=1.5pt,line join=round]
  table[row sep=crcr]{%
1	0.999156650065644\\
2	0.816595756057913\\
3	0.0775882304872165\\
4	0.0478862381993567\\
5	0.0553756194984354\\
6	0.366603258881778\\
7	0.387954888412419\\
8	0.356584889949596\\
9	0.239343518027135\\
10	0.140408072664929\\
11	0.078528360412846\\
12	0.0543415415216163\\
13	0.0528675685147501\\
14	0.0550429593078466\\
15	0.057894474742636\\
16	0.0616518964808358\\
17	0.0667763206530777\\
18	0.0738879747459651\\
19	0.0837886076790918\\
20	0.0974472104689641\\
21	0.116051008604862\\
22	0.141628266883015\\
23	0.174698236788813\\
24	0.211177051236828\\
25	0.246031182947741\\
26	0.278606566589521\\
27	0.31103513333901\\
28	0.346004128890065\\
29	0.381155616906308\\
30	0.411608686662539\\
31	0.439203720415992\\
32	0.462823066120177\\
33	0.480619413867379\\
34	0.494865840911432\\
35	0.50641966109416\\
36	0.516506119375991\\
37	0.52327602929204\\
38	0.52724631799475\\
39	0.530131119945597\\
40	0.53195275267781\\
41	0.532866988322803\\
42	0.53284255901119\\
43	0.532163058752744\\
44	0.5311807553055\\
45	0.530101996725694\\
46	0.528904803820913\\
47	0.527669546156809\\
48	0.526436576467566\\
49	0.525245568572938\\
50	0.524052676020065\\
51	0.522861724319827\\
52	0.521664580934756\\
53	0.520577563824453\\
54	0.519571986587759\\
55	0.518663333037728\\
56	0.51779380685771\\
57	0.516967747435568\\
58	0.516139664589222\\
59	0.515318432956215\\
60	0.514497498617641\\
61	0.513702434662589\\
62	0.512920120657616\\
63	0.512131316396731\\
64	0.511362289588002\\
65	0.510626430620174\\
66	0.509911524572945\\
67	0.509223044080988\\
68	0.508579889364341\\
69	0.50797999886415\\
70	0.507427322685018\\
71	0.506903358679042\\
72	0.506397219246429\\
73	0.505919154294284\\
74	0.50545369534171\\
75	0.505000921172581\\
76	0.504560568767321\\
77	0.504164422910319\\
78	0.503817023182927\\
79	0.50350009106399\\
80	0.503201231003024\\
81	0.502915452989934\\
82	0.502646228671331\\
83	0.502392760853718\\
84	0.502150861968036\\
85	0.501919413057029\\
86	0.501698365752452\\
87	0.501476817781761\\
88	0.501258532743594\\
89	0.501047996036285\\
90	0.50084084380738\\
91	0.500634372895596\\
92	0.500430564787215\\
93	0.500227953865523\\
94	0.500036534604521\\
95	0.499854848826097\\
96	0.499679413545083\\
97	0.499510532158805\\
98	0.499347460519044\\
99	0.499187499202284\\
100	0.499036061136902\\
101	0.498891189810231\\
102	0.498752910809335\\
103	0.498615753239168\\
104	0.498480180670408\\
105	0.498343520122737\\
106	0.498205198880854\\
107	0.498062120628175\\
108	0.49791958415331\\
109	0.497775143843248\\
110	0.4976291481341\\
111	0.497479846934396\\
112	0.497322024211315\\
113	0.497161717585366\\
114	0.497004527353167\\
115	0.496851446819858\\
116	0.496699421382731\\
117	0.496552379963706\\
118	0.496407671614525\\
119	0.496269157765489\\
120	0.496132057065206\\
121	0.495999231088572\\
122	0.495870460189634\\
123	0.49574404317214\\
124	0.495620890800782\\
125	0.495497795632217\\
126	0.495373291927739\\
127	0.495245607344114\\
128	0.495111892848734\\
129	0.494978127148955\\
130	0.494846965231149\\
131	0.494714932982597\\
132	0.494585467556207\\
133	0.494460023672475\\
134	0.494338974752257\\
135	0.494219685422881\\
136	0.494099931473997\\
137	0.493982484099057\\
138	0.493869158731368\\
139	0.493756790029192\\
140	0.493645142000468\\
141	0.493538791861522\\
142	0.493436803313569\\
143	0.493337407875575\\
144	0.493242879526314\\
145	0.493152439196206\\
146	0.493062961658619\\
147	0.49297597486469\\
148	0.492887319335308\\
149	0.492805678293583\\
150	0.492734832598745\\
151	0.492669500344967\\
152	0.492602127142073\\
153	0.492537139437819\\
154	0.492478913615473\\
155	0.492428383638885\\
156	0.492381011858981\\
157	0.492332958151215\\
158	0.492283642969717\\
159	0.49223563945174\\
160	0.492189269494375\\
161	0.492145480242288\\
162	0.49210532570927\\
163	0.492063932935662\\
164	0.492020314740211\\
165	0.491978585618632\\
166	0.491936019341438\\
167	0.491890064599201\\
168	0.491842614418166\\
169	0.491797835458452\\
170	0.491756735607141\\
171	0.491717978840104\\
172	0.491680210312883\\
173	0.491644328994794\\
174	0.491608822900172\\
175	0.491574599983102\\
176	0.491542200045969\\
177	0.491510745235526\\
178	0.491481106965567\\
179	0.491453302620466\\
180	0.491425074018564\\
181	0.491397474693264\\
182	0.491370753940297\\
183	0.491341842340195\\
184	0.491309815315451\\
185	0.491274592858472\\
186	0.491234716125862\\
187	0.491191726419452\\
188	0.491145629312775\\
189	0.491095434195383\\
190	0.491041998135515\\
191	0.490989933225131\\
192	0.490936935683599\\
193	0.490886412475443\\
194	0.490840419251047\\
195	0.490794149384955\\
196	0.490747345236782\\
197	0.490699932162129\\
198	0.490648850659303\\
199	0.490592895208685\\
200	0.546381125929401\\
201	0.482140041286137\\
202	0.493760267430956\\
203	0.494785080024608\\
204	0.49500669927587\\
205	0.495000203016274\\
206	0.494891364386318\\
207	0.494741577839522\\
208	0.494603695915668\\
209	0.494506279269031\\
210	0.494446891366009\\
211	0.494405894319241\\
212	0.494372760787206\\
213	0.494347040594765\\
214	0.494320677272223\\
215	0.49429225111379\\
216	0.49425616491534\\
217	0.494216261248453\\
218	0.494179378322302\\
219	0.494154451401297\\
220	0.494146967546013\\
221	0.494148513869065\\
222	0.494158737701979\\
223	0.494174297286211\\
224	0.494188272797309\\
225	0.494196019604783\\
226	0.494198560953747\\
227	0.494198699315716\\
228	0.49419594573587\\
229	0.494188920707791\\
230	0.494176413949328\\
231	0.4941607663188\\
232	0.494143520048365\\
233	0.494125794146487\\
234	0.494107293097536\\
235	0.494087335352151\\
236	0.494064463232899\\
237	0.494037611113095\\
238	0.494006669556159\\
239	0.493973711500989\\
240	0.493941188306233\\
241	0.493906935018372\\
242	0.493869095903423\\
243	0.49383032038116\\
244	0.493794523240708\\
245	0.493762309758717\\
246	0.49373195663264\\
247	0.493700705055521\\
248	0.493668195291965\\
249	0.493635606952747\\
250	0.493602727652897\\
251	0.493567128784058\\
252	0.493528439269209\\
253	0.493485576091785\\
254	0.493444882396384\\
255	0.493408531984139\\
256	0.493372831831923\\
257	0.493336529108914\\
258	0.493299738580749\\
259	0.493262687217834\\
260	0.493222574932402\\
261	0.493182284611086\\
262	0.493141723866832\\
263	0.493105183463153\\
264	0.493070413924918\\
265	0.493035968516735\\
266	0.493001493895915\\
267	0.492969702075482\\
268	0.492942731337099\\
269	0.492918493473211\\
270	0.492899716158886\\
271	0.492884564900497\\
272	0.492871150962489\\
273	0.492857379284173\\
274	0.492842739368849\\
275	0.49282689487935\\
276	0.492809903115623\\
277	0.49279312469082\\
278	0.492776915053681\\
279	0.492760268413494\\
280	0.492742672617994\\
281	0.492723428435565\\
282	0.492702842226601\\
283	0.492683359006998\\
284	0.492665250501756\\
285	0.492646120573894\\
286	0.492624431395507\\
287	0.492603247906707\\
288	0.492586296768809\\
289	0.492572377136425\\
290	0.492559434192681\\
291	0.492548850515048\\
292	0.492542588133232\\
293	0.49254073190271\\
294	0.492539198313321\\
295	0.492535068807179\\
296	0.492526941430877\\
297	0.49251378667823\\
298	0.492495847609587\\
299	0.492475509260829\\
300	0.492454161334659\\
301	0.492434914320427\\
302	0.492419403432757\\
303	0.492409479123075\\
304	0.492405573778585\\
305	0.492407802647956\\
306	0.492412025158533\\
307	0.492412465854865\\
308	0.492406056932562\\
309	0.492392447454424\\
310	0.492375297995439\\
311	0.492355501468946\\
312	0.492335584131229\\
313	0.492316084293956\\
314	0.492294761208873\\
315	0.492270845429103\\
316	0.492245045982982\\
317	0.492224720155292\\
318	0.492206171805353\\
319	0.492184554930051\\
320	0.492159135624851\\
321	0.492128776221537\\
322	0.492091770693424\\
323	0.49205475639487\\
324	0.492017445573854\\
325	0.491980603601966\\
326	0.491946443627623\\
327	0.491923521882947\\
328	0.491913035740666\\
329	0.491910974553433\\
330	0.49191786715266\\
331	0.491922995522084\\
332	0.491919450670939\\
333	0.491907262494967\\
334	0.49189089107494\\
335	0.491871384790691\\
336	0.491851751783651\\
337	0.491833176904002\\
338	0.491815245063358\\
339	0.491795256823457\\
340	0.491771210388658\\
341	0.491744569807261\\
342	0.491716955213805\\
343	0.491688701336418\\
344	0.491660751475015\\
345	0.49163342338638\\
346	0.49160669858181\\
347	0.491581171159882\\
348	0.491556631761609\\
349	0.491531164729712\\
350	0.491502754152566\\
351	0.491471153741974\\
352	0.49143661445613\\
353	0.491399931316897\\
354	0.491363349453792\\
355	0.491326069319504\\
356	0.491285892690614\\
357	0.491242282001322\\
358	0.491194158622437\\
359	0.491142257686153\\
360	0.491089064556246\\
361	0.491036317648197\\
362	0.490982874978478\\
363	0.490929838286091\\
364	0.490879015019764\\
365	0.490829782212865\\
366	0.490781469222641\\
367	0.49073426404373\\
368	0.490689686400131\\
369	0.490646697126868\\
370	0.490605052745915\\
371	0.490565848566563\\
372	0.490527775314467\\
373	0.490491722941048\\
374	0.49045841417714\\
375	0.490427305430915\\
376	0.490395385439146\\
377	0.490364276445224\\
378	0.49033741117604\\
379	0.490313131191459\\
380	0.490290008686454\\
381	0.490268670704036\\
382	0.490247198237902\\
383	0.490227353216849\\
384	0.490208523153371\\
385	0.490190723332028\\
386	0.490175992558722\\
387	0.49016507658611\\
388	0.490157927496902\\
389	0.490154214932407\\
390	0.490152278822155\\
391	0.49014802499403\\
392	0.490139261543377\\
393	0.490127330999845\\
394	0.49011372912221\\
395	0.490099872886767\\
396	0.490088777810317\\
397	0.490082560698468\\
398	0.490079170254663\\
399	0.490075528333774\\
400	0.490069110980036\\
401	0.490057873003851\\
};
\addlegendentry{smooth}

\addplot [color=mycolor3,line width=1.5pt,line join=round]
  table[row sep=crcr]{%
1	1.00331519624177\\
2	0.817370799108866\\
3	0.0764244368174531\\
4	0.0464899835293668\\
5	0.0539698981841624\\
6	0.366063649102888\\
7	0.387604829764022\\
8	0.356212550904089\\
9	0.238714645061228\\
10	0.139494488300636\\
11	0.077329298401878\\
12	0.0529751465002873\\
13	0.0514859270878881\\
14	0.0536745511464473\\
15	0.0565424095450421\\
16	0.0603200933033912\\
17	0.0654700525742283\\
18	0.0726138031714211\\
19	0.0825544132119864\\
20	0.096262616223276\\
21	0.114927289500552\\
22	0.140579729884613\\
23	0.173738426375319\\
24	0.210310522353227\\
25	0.24525446931072\\
26	0.277916306858356\\
27	0.310435251142246\\
28	0.345510494151035\\
29	0.380778354351074\\
30	0.411340088007886\\
31	0.4390425494392\\
32	0.462766461391578\\
33	0.480651028520871\\
34	0.494975298211559\\
35	0.506595667797223\\
36	0.516743729648207\\
37	0.523556602257478\\
38	0.527553828768996\\
39	0.530457803772513\\
40	0.532291689180238\\
41	0.5332116518923\\
42	0.533186383935337\\
43	0.532501369569433\\
44	0.531511498620934\\
45	0.530424694902067\\
46	0.529218860578108\\
47	0.527974837313493\\
48	0.526733229012271\\
49	0.525533939403907\\
50	0.524332884105127\\
51	0.52313389523935\\
52	0.52192873839633\\
53	0.520834527920333\\
54	0.519822412458239\\
55	0.518907861019538\\
56	0.518032753016818\\
57	0.517201437930459\\
58	0.516368097012245\\
59	0.515541652307606\\
60	0.514715501357748\\
61	0.513915424190366\\
62	0.51312818477952\\
63	0.512334464652988\\
64	0.51156075447855\\
65	0.510820506056653\\
66	0.510101404605392\\
67	0.509408934225777\\
68	0.508762105467271\\
69	0.508158807527688\\
70	0.50760298082331\\
71	0.507076020296835\\
72	0.50656697555498\\
73	0.506086147885151\\
74	0.50561798989834\\
75	0.505162578587994\\
76	0.504719653334308\\
77	0.504321196194286\\
78	0.503971750159063\\
79	0.503652946123143\\
80	0.503352323836632\\
81	0.503064860251211\\
82	0.502794048676054\\
83	0.502539078356226\\
84	0.502295732908105\\
85	0.502062894964049\\
86	0.501840525684617\\
87	0.50161765654614\\
88	0.501398069033748\\
89	0.501186278095539\\
90	0.500977895684917\\
91	0.500770202279256\\
92	0.50056519102213\\
93	0.500361389591099\\
94	0.500168854332021\\
95	0.499986117138421\\
96	0.499809675479575\\
97	0.499639833318022\\
98	0.499475841578724\\
99	0.499314984882266\\
100	0.499162708840358\\
101	0.499017042923924\\
102	0.498878009122952\\
103	0.498740100368135\\
104	0.498603782783459\\
105	0.498466366705715\\
106	0.498327275020962\\
107	0.498183394867278\\
108	0.498040058818539\\
109	0.497894806959399\\
110	0.497747991297215\\
111	0.497597851161258\\
112	0.497439140124414\\
113	0.497277931359742\\
114	0.497119856385398\\
115	0.496965913417047\\
116	0.496813029773263\\
117	0.496665159219517\\
118	0.496519636205487\\
119	0.496380344389115\\
120	0.496242476819067\\
121	0.496108909986551\\
122	0.495979423259889\\
123	0.49585230587647\\
124	0.495728474336696\\
125	0.495604703812822\\
126	0.495479518929113\\
127	0.495351134857202\\
128	0.495216686791193\\
129	0.495082188226519\\
130	0.494950310069806\\
131	0.494817558920627\\
132	0.494687392901837\\
133	0.494561275535067\\
134	0.494439579989548\\
135	0.49431965453163\\
136	0.494199259876083\\
137	0.494081182027718\\
138	0.493967248116357\\
139	0.493854275170192\\
140	0.493742026544602\\
141	0.493635108156932\\
142	0.493532580282527\\
143	0.493432660958852\\
144	0.493337636338233\\
145	0.493246723267198\\
146	0.4931567786742\\
147	0.493069336002872\\
148	0.49298021344741\\
149	0.492898146634272\\
150	0.492826934324258\\
151	0.492761260987626\\
152	0.4926935305754\\
153	0.492628197398441\\
154	0.492569666052668\\
155	0.492518876340807\\
156	0.492471264610757\\
157	0.492422967928594\\
158	0.492373403042797\\
159	0.492325158061838\\
160	0.492278555865337\\
161	0.49223454807007\\
162	0.492194194560029\\
163	0.492152592860447\\
164	0.492108749956887\\
165	0.492066806369313\\
166	0.492024020287193\\
167	0.491977825315834\\
168	0.491930124954966\\
169	0.491885109555925\\
170	0.491843792892997\\
171	0.491804832361076\\
172	0.491766865542653\\
173	0.49173079658841\\
174	0.491695105338134\\
175	0.491660706184246\\
176	0.491628142243566\\
177	0.491596530255456\\
178	0.491566746295247\\
179	0.491538807491565\\
180	0.491510442276172\\
181	0.491482710847189\\
182	0.491455863910527\\
183	0.491426813622046\\
184	0.491394630010751\\
185	0.491359233085014\\
186	0.491319156441382\\
187	0.491275951303126\\
188	0.49122962422419\\
189	0.491179177921651\\
190	0.491125473087922\\
191	0.491073149936439\\
192	0.491019891548979\\
193	0.490969122703656\\
194	0.490922909857268\\
195	0.490876416964921\\
196	0.490829384822616\\
197	0.490781741158958\\
198	0.490730409124908\\
199	0.490674176845007\\
200	0.546796174045338\\
201	0.482178222001184\\
202	0.493859884404957\\
203	0.494890167601727\\
204	0.495112750132716\\
205	0.49510587685155\\
206	0.494996130607518\\
207	0.494845245253436\\
208	0.494706387489011\\
209	0.494608276523651\\
210	0.494548442551017\\
211	0.494507117493036\\
212	0.494473706799692\\
213	0.494447758636456\\
214	0.494421166439048\\
215	0.494392507791162\\
216	0.494356157250024\\
217	0.49431598164641\\
218	0.494278856468427\\
219	0.494253765610498\\
220	0.494246226879648\\
221	0.494247776966027\\
222	0.494258059277325\\
223	0.49427371189844\\
224	0.494287774210335\\
225	0.494295574771084\\
226	0.494298143542698\\
227	0.494298298473545\\
228	0.494295546425425\\
229	0.494288500146642\\
230	0.494275942949113\\
231	0.494260228271585\\
232	0.494242906694747\\
233	0.49422510322267\\
234	0.494206518603026\\
235	0.494186465989467\\
236	0.494163479958214\\
237	0.494136488908049\\
238	0.494105382938271\\
239	0.494072247775153\\
240	0.494039548828497\\
241	0.494005107580848\\
242	0.493967057339058\\
243	0.493928063512634\\
244	0.493892064157045\\
245	0.493859667456316\\
246	0.493829139670427\\
247	0.49379770651339\\
248	0.493765006863867\\
249	0.493732227540235\\
250	0.493699155505955\\
251	0.493663349703469\\
252	0.493624437312726\\
253	0.493581329084374\\
254	0.493540402520617\\
255	0.49350384301824\\
256	0.493467936297175\\
257	0.493431422746\\
258	0.493394418425449\\
259	0.493357152004025\\
260	0.4933168078448\\
261	0.493276285375227\\
262	0.493235491322874\\
263	0.493198740183726\\
264	0.493163769539813\\
265	0.493129124397378\\
266	0.493094449876784\\
267	0.493062473733085\\
268	0.493035346006498\\
269	0.493010966641029\\
270	0.492992078911006\\
271	0.492976837673199\\
272	0.492963343354191\\
273	0.492949489335561\\
274	0.492934762337459\\
275	0.492918824200627\\
276	0.492901732152146\\
277	0.492884854130897\\
278	0.492868547732444\\
279	0.492851801576618\\
280	0.492834100648566\\
281	0.492814741166969\\
282	0.492794030826182\\
283	0.492774428924558\\
284	0.49275620885984\\
285	0.492736960884291\\
286	0.492715138553313\\
287	0.492693825648327\\
288	0.492676771292898\\
289	0.492662767526476\\
290	0.492649747817715\\
291	0.492639103110519\\
292	0.492632806820245\\
293	0.492630942204279\\
294	0.492629400842387\\
295	0.49262524947543\\
296	0.49261707993731\\
297	0.492603856881293\\
298	0.492585824454916\\
299	0.492565379921063\\
300	0.492543919424339\\
301	0.492524571524015\\
302	0.492508981696506\\
303	0.492499012308829\\
304	0.492495099281769\\
305	0.492497358438387\\
306	0.4925016235416\\
307	0.492502085205725\\
308	0.492495658176886\\
309	0.492481988864084\\
310	0.492464758691395\\
311	0.492444866484513\\
312	0.49242485388904\\
313	0.492405262284453\\
314	0.49238383778019\\
315	0.492359808556601\\
316	0.492333888449385\\
317	0.492313475454012\\
318	0.492294850797347\\
319	0.492273139868944\\
320	0.492247602939199\\
321	0.492217097650514\\
322	0.492179909976868\\
323	0.492142713371287\\
324	0.492105220878325\\
325	0.492068203365225\\
326	0.492033886323144\\
327	0.492010873371103\\
328	0.492000370894198\\
329	0.491998346054331\\
330	0.492005329534847\\
331	0.492010537527616\\
332	0.492007011743257\\
333	0.491994782650224\\
334	0.491978342691041\\
335	0.491958745161961\\
336	0.491939015480823\\
337	0.491920348203086\\
338	0.491902327563802\\
339	0.491882239158113\\
340	0.491858069184586\\
341	0.491831290043357\\
342	0.491803531867586\\
343	0.491775131487909\\
344	0.491747037430946\\
345	0.491719568700851\\
346	0.491692705595271\\
347	0.491667044794666\\
348	0.491642375643369\\
349	0.491616772841183\\
350	0.49158820967186\\
351	0.491556438762003\\
352	0.491521711614717\\
353	0.491484827155487\\
354	0.491448043589249\\
355	0.491410557426988\\
356	0.491370157634369\\
357	0.491326303617227\\
358	0.491277911087033\\
359	0.491225719773231\\
360	0.491172229575347\\
361	0.491119189937637\\
362	0.491065453426323\\
363	0.491012125921631\\
364	0.49096102436319\\
365	0.490911521667245\\
366	0.490862943964864\\
367	0.490815482566459\\
368	0.490770665661236\\
369	0.490727447639008\\
370	0.490685584446485\\
371	0.490646176071934\\
372	0.490607905473077\\
373	0.490571668500241\\
374	0.490538192696119\\
375	0.490506932120435\\
376	0.490474858589221\\
377	0.490443602290547\\
378	0.490416615095355\\
379	0.490392228058729\\
380	0.490369005982919\\
381	0.490347580097516\\
382	0.490326021078915\\
383	0.490306099536659\\
384	0.490287197037812\\
385	0.490269327854389\\
386	0.490254543020963\\
387	0.490243592510589\\
388	0.490236429074075\\
389	0.490232722051334\\
390	0.490230802564742\\
391	0.490226554288736\\
392	0.49021777208096\\
393	0.4902058045226\\
394	0.490192155899377\\
395	0.490178252059623\\
396	0.490167125779257\\
397	0.490160904747784\\
398	0.490157526541866\\
399	0.490153896415267\\
400	0.490147476534243\\
401	0.490136209101736\\
};
\addlegendentry{smooth \& steep}

\end{axis}
\end{tikzpicture}%

%% file: wheel_obj_scatter.tex
\definecolor{mycolor1}{rgb}{0.00000,0.44700,0.74100}%
\definecolor{mycolor2}{rgb}{0.85000,0.32500,0.09800}%

\begin{tikzpicture}

\begin{axis}[%
width=\textwidth,
height=.85\textwidth,
at={(1.517in,0.962in)},
%scale only axis,
xmin=0.5,
xmax=5.5,
ymin=0.615,
ymax=0.665,
axis background/.style={fill=white},
xmajorgrids,
ymajorgrids,
xtick={1,2,3,4,5},
xticklabels={4,8,16,32,64},
xlabel={Batch size},
legend style={legend cell align=left, align=left, draw=white!15!black},
ylabel={$\rvol$}
]
\addplot[only marks, mark=*, mark options={}, mark size=2.0pt, color=mycolor1, fill=mycolor1] table[row sep=crcr]{%
x	y\\
1	0.661602983620999\\
1	0.652942311227018\\
1	0.661228227170273\\
1	0.630049118746794\\
};
\addlegendentry{sMMA};

\addplot[only marks, mark=*, mark options={}, mark size=2.0pt, color=mycolor1, forget plot, fill=mycolor1] table[row sep=crcr]{%
x	y\\
2	0.632407245841023\\
2	0.642936134848562\\
2	0.633309618326712\\
2	0.636955230164887\\
};
\addplot[only marks, mark=*, mark options={}, mark size=2.0pt, color=mycolor1, forget plot, fill=mycolor1] table[row sep=crcr]{%
x	y\\
3	0.643697559129574\\
3	0.644541468424804\\
3	0.645754255627769\\
3	0.643675522959918\\
};
\addplot[only marks, mark=*, mark options={}, mark size=2.0pt, color=mycolor1, forget plot, fill=mycolor1] table[row sep=crcr]{%
x	y\\
4	0.636637722358925\\
4	0.64035958551182\\
4	0.632726691398583\\
4	0.642896895331209\\
};
\addplot[only marks, mark=*, mark options={}, mark size=2.0pt, color=mycolor1, forget plot, fill=mycolor1] table[row sep=crcr]{%
x	y\\
5	0.626088054319007\\
5	0.639481879389376\\
5	0.633452064098024\\
5	0.64152975966689\\
};
\addplot[only marks, mark=*, mark options={}, mark size=2.0pt, color=mycolor2, fill=mycolor2] table[row sep=crcr]{%
x	y\\
4	0.62572487625113\\
4	0.620347441406752\\
4	0.620990411579714\\
4	0.635079767172234\\
};
\addlegendentry{MMA};

\addplot[only marks, mark=*, mark options={}, mark size=2.0pt, color=mycolor2, forget plot, fill=mycolor2] table[row sep=crcr]{%
x	y\\
5	0.618849247004736\\
5	0.620025305876591\\
5	0.621154328035088\\
5	0.62899670306978\\
};
\end{axis}

\end{tikzpicture}%

%% file: wheel_phy_scatter.tex
\definecolor{mycolor1}{rgb}{0.00000,0.44700,0.74100}%
\definecolor{mycolor2}{rgb}{0.85000,0.32500,0.09800}%

\begin{tikzpicture}

\begin{axis}[%
width=\textwidth,
height=.85\textwidth,
at={(1.517in,0.962in)},
%scale only axis,
xmin=0.5,
xmax=5.5,
ymin=0.455,
ymax=0.5,
axis background/.style={fill=white},
xmajorgrids,
ymajorgrids,
xtick={1,2,3,4,5},
xticklabels={4,8,16,32,64},
xlabel={Batch size},
legend style={legend cell align=left, align=left, draw=white!15!black},
ylabel={$\pvol$}
]
\addplot[only marks, mark=*, mark options={}, mark size=2.0pt, color=mycolor1, fill=mycolor1] table[row sep=crcr]{%
x	y\\
1	0.473020674728141\\
1	0.464628110485212\\
1	0.465782079088995\\
1	0.478031668688297\\
};
\addlegendentry{sMMA};

\addplot[only marks, mark=*, mark options={}, mark size=2.0pt, color=mycolor1, forget plot, fill=mycolor1] table[row sep=crcr]{%
x	y\\
2	0.45731504469348\\
2	0.464073302598008\\
2	0.481055950346529\\
2	0.497055101561615\\
};
\addplot[only marks, mark=*, mark options={}, mark size=2.0pt, color=mycolor1, forget plot, fill=mycolor1] table[row sep=crcr]{%
x	y\\
3	0.474355995938482\\
3	0.477515773444111\\
3	0.474656965949866\\
3	0.469659288146445\\
};
\addplot[only marks, mark=*, mark options={}, mark size=2.0pt, color=mycolor1, forget plot, fill=mycolor1] table[row sep=crcr]{%
x	y\\
4	0.471709486307729\\
4	0.464321474083927\\
4	0.468792756947\\
4	0.47084898908096\\
};
\addplot[only marks, mark=*, mark options={}, mark size=2.0pt, color=mycolor1, forget plot, fill=mycolor1] table[row sep=crcr]{%
x	y\\
5	0.472211656813625\\
5	0.466773208489206\\
5	0.467762111571891\\
5	0.481381488794859\\
};
\addplot[only marks, mark=*, mark options={}, mark size=2.0pt, color=mycolor2, fill=mycolor2] table[row sep=crcr]{%
x	y\\
4	0.47824990594896\\
4	0.474127366173422\\
4	0.475463467807045\\
4	0.472745339668167\\
};
\addlegendentry{MMA};

\addplot[only marks, mark=*, mark options={}, mark size=2.0pt, color=mycolor2, forget plot, fill=mycolor2] table[row sep=crcr]{%
x	y\\
5	0.477878891251179\\
5	0.476729867755772\\
5	0.478775300297812\\
5	0.474901453032746\\
};
\end{axis}
\end{tikzpicture}%

%% file: 2dtable_setup.tex
\begin{tikzpicture}[scale=0.9\textwidth/40cm]
\filldraw[very thick,fill=black!5] (0,0) -- (36,0) -- (36,18) -- (0,18) -- (0,0);
\filldraw[thick,fill=black!55] (-2,0) -- (38,0) -- (38,-1) -- (-2,-1) -- (-2,0);
\draw[line width=3pt, magenta] (3.6,18)--(32.4,18);
\filldraw[thick,fill=blue!20, draw = blue] (6.5,18) -- (8.5,18) -- (8.5,18.5) -- (6.5,18.5) -- (6.5,18);
%\draw[<-,thick, blue] (7.8,15.6) -- (9.21,17.01);
\draw[very thick,dotted,black] (7.5,18.6) -- (7.5,21.6);
\draw[very thick,dotted,black] (7.5,18.6) -- (9.6621,20.721);
\draw[very thick,dotted,black] (7.5,18.6) -- (5.379,20.721);
\draw[thick, black] (8.58,19.665) arc (45:135:1.5);
%\filldraw[thick,dashed,magenta,fill=magenta!10] (4.5,2.25) -- (31.5,2.25) -- (31.5,15.75) -- (4.5,15.75) -- (4.5,2.25);
\draw[thick,dashed,teal] (4.5,2.25) -- (31.5,2.25) -- (31.5,15.75) -- (4.5,15.75) -- (4.5,2.25);
\filldraw[blue,fill=blue!25](19,11) circle (1);
\filldraw[blue] (19,11) circle(.15);
\draw[Latex-Latex] (-3,0)--node[left,midway]{$\ell$}(-3,18);
\draw[Latex-Latex] (0,-3)--node[below,midway]{$2\ell$}(36,-3);
\draw (2,5) node{\Large{$\mathscr{D}$}};
\draw (26,19.5) node{\color{magenta}{\Large{$\Omega$}}};
\draw (33,8) node{\color{teal}{\Large{$\Xi$}}};
\end{tikzpicture}

%% file: 3dSetupSpat.tex
\begin{tikzpicture}[scale=.75\textwidth/90cm]
\draw[thick,line cap=round](0,0,0)--(60,0,0)--(60,60,0)--(0,60,0)--cycle;
\draw[thick, line cap=round](0,0,0)--(0,0,60)--(0,60,60)--(0,60,0)--cycle;
\draw[thick, line cap=round](0,0,0)--(60,0,0)--(60,0,60)--(0,0,60)--cycle;
\draw[thick, line cap=round](0,60,0)--(60,60,0)--(60,60,60)--(0,60,60)--cycle;
\draw[thick, line cap=round](60,0,0)--(60,0,60)--(60,60,60)--(60,60,0)--cycle;
\filldraw[thick, fill=black!35] (3,0,3)--(3,0,6)--(6,0,6)--(6,0,3)--cycle;
\filldraw[thick, fill=black!35] (3,0,54)--(3,0,57)--(6,0,57)--(6,0,54)--cycle;
\filldraw[thick, fill=black!35] (54,0,3)--(57,0,3)--(57,0,6)--(54,0,6)--cycle;
\filldraw[thick,dashed,magenta,fill=magenta!10] (7.5,60,20)--(52.5,60,20)--(52.5,60,40)--(7.5,60,40)--cycle;
\filldraw[thick,fill=blue!20, draw = blue] (19,60,38)--(23.8,60,38)--(23.8,60,33.2)--(19,60,33.2)--cycle;
\filldraw[thick,fill=blue,draw=blue] (21.4,60,35.6) circle(.25);
\draw[Latex-Latex] (0,-1,61)--node[below,midway]{$\ell$}(60,-1,61);
\draw[Latex-Latex] (-1,0,61)--node[left,midway]{$\ell$}(-1,60,61);
\draw[Latex-Latex] (-1,61,0)--node[left,midway]{$\ell$}(-1,61,60);
\end{tikzpicture}

%% file: MCMSA_rel_comp_histo.tex
\definecolor{mycolor1}{rgb}{0.00000,0.44700,0.74100}%
\begin{tikzpicture}

\begin{axis}[%
width=\textwidth,
height=.75\textwidth,
at={(0.758in,0.481in)},
%scale only axis,
xmin=0,
xmax=4.1,
ymin=0,
ymax=0.6,
axis background/.style={fill=white},
xtick={0,1,2,3,4},
xticklabels={0,1,2,3,$\ge4$},
xmajorgrids,
ymajorgrids,
xlabel={Relative compliance},
ylabel={Probability}
]

\addplot[ybar interval, fill=mycolor1, fill opacity=0.6, draw=mycolor1, area legend] table[row sep=crcr] {%
x	y\\
0	0\\
0.1	0\\
0.2	0\\
0.3	0\\
0.4	0\\
0.5	0\\
0.6	0.0208\\
0.7	0.3512\\
0.8	0.5136\\
0.9	0.0692\\
1	0.0178\\
1.1	0.0104\\
1.2	0.0048\\
1.3	0.003\\
1.4	0.0026\\
1.5	0.0014\\
1.6	0.0014\\
1.7	0.0008\\
1.8	0.0002\\
1.9	0.0008\\
2	0.0002\\
2.1	0.0006\\
2.2	0\\
2.3	0.0002\\
2.4	0.0002\\
2.5	0.0002\\
2.6	0.0002\\
2.7	0\\
2.8	0\\
2.9	0\\
3	0.0002\\
3.1	0\\
3.2	0\\
3.3	0\\
3.4	0.0002\\
3.5	0\\
3.6	0\\
3.7	0\\
3.8	0\\
3.9	0\\
4	0\\
4.1	0\\
4.2	0\\
4.3	0\\
4.4	0\\
4.5	0\\
4.6	0\\
4.7	0\\
4.8	0\\
};
\addplot [color=black, dashed, forget plot,line width=1.0pt]
  table[row sep=crcr]{%
1	0\\
1	1\\
};
\end{axis}
\end{tikzpicture}%

%% file: MMA_rel_comp_histo.tex
\definecolor{mycolor2}{rgb}{0.85000,0.32500,0.09800}%
\begin{tikzpicture}

\begin{axis}[%
width=\textwidth,
height=.75\textwidth,
at={(0.758in,0.481in)},
%scale only axis,
xmin=0,
xmax=4.1,
ymin=0,
ymax=0.3,
xtick={0,1,2,3,4},
xticklabels={0,1,2,3,$\ge4$},
xmajorgrids,
ymajorgrids,
axis background/.style={fill=white},
xlabel={Relative compliance},
ylabel={Probability}
]

\addplot[ybar interval, fill=mycolor2, fill opacity=0.6, draw=mycolor2, area legend] table[row sep=crcr] {%
x	y\\
0	0\\
0.1	0\\
0.2	0\\
0.3	0\\
0.4	0\\
0.5	0\\
0.6	0\\
0.7	0.0048\\
0.8	0.0866\\
0.9	0.2538\\
1	0.2016\\
1.1	0.0876\\
1.2	0.0446\\
1.3	0.0262\\
1.4	0.0126\\
1.5	0.006\\
1.6	0.0048\\
1.7	0.0054\\
1.8	0.0042\\
1.9	0.004\\
2	0.0028\\
2.1	0.0022\\
2.2	0.0038\\
2.3	0.002\\
2.4	0.002\\
2.5	0.0014\\
2.6	0.0014\\
2.7	0.0026\\
2.8	0.002\\
2.9	0.0024\\
3	0.002\\
3.1	0.0024\\
3.2	0.002\\
3.3	0.002\\
3.4	0.0016\\
3.5	0.002\\
3.6	0.001\\
3.7	0.0016\\
3.8	0.0012\\
3.9	0.001\\
4 0.2184\\
4.1 0.2184 \\
% 4	0.001\\
% 4.1	0.001\\
% 4.2	0.001\\
% 4.3	0.0012\\
% 4.4	0.0012\\
% 4.5	0.0016\\
% 4.6	0.0008\\
% 4.7	0.0016\\
% 4.8	0.0014\\
% 4.9	0.0006\\
% 5	0.0014\\
% 5.1	0.0008\\
% 5.2	0.0014\\
% 5.3	0.002\\
% 5.4	0.0004\\
% 5.5	0.0014\\
% 5.6	0.0006\\
% 5.7	0.0014\\
% 5.8	0.001\\
% 5.9	0.0004\\
% 6	0.0014\\
% 6.1	0.0008\\
% 6.2	0.0014\\
% 6.3	0.001\\
% 6.4	0.0002\\
% 6.5	0.0018\\
% 6.6	0.001\\
% 6.7	0.001\\
% 6.8	0.0004\\
% 6.9	0.0008\\
% 7	0.0016\\
% 7.1	0.0004\\
% 7.2	0.0014\\
% 7.3	0.0008\\
% 7.4	0.0006\\
% 7.5	0.001\\
% 7.6	0.0002\\
% 7.7	0.0008\\
% 7.8	0.0008\\
% 7.9	0.0002\\
% 8	0.0004\\
% 8.1	0.0008\\
% 8.2	0.0002\\
% 8.3	0.0004\\
% 8.4	0.0006\\
% 8.5	0.001\\
% 8.6	0.0014\\
% 8.7	0.001\\
% 8.8	0.0002\\
% 8.9	0.0002\\
% 9	0.0006\\
% 9.1	0.0004\\
% 9.2	0.0006\\
% 9.3	0.001\\
% 9.4	0.001\\
% 9.5	0.0004\\
% 9.6	0.0012\\
% 9.7	0.0008\\
% 9.8	0.0002\\
% 9.9	0.0002\\
% 10	0.166\\
% 11	0.166\\
};
\addplot [color=black, dashed, forget plot,line width=1.0pt]
  table[row sep=crcr]{%
1	0\\
1	1\\
};
\end{axis}
\end{tikzpicture}%